\documentclass[11pt]{article}
\usepackage{paper,asb}

\usepackage{amsmath}
\usepackage{amsfonts}
\usepackage{amssymb}
\usepackage{amsthm}
\usepackage{hyperref}
\usepackage{bbm}
\newtheorem{theorem}{Theorem}
\newtheorem{lemma}[theorem]{Lemma}
\newtheorem{example}[theorem]{Example}
\newtheorem{termtest}[]{Termination Test}
\newtheorem{assumption}[theorem]{Assumption}
\newtheorem{remark}[theorem]{Remark}

\usepackage{marginnote}
\usepackage{subcaption}
\usepackage{booktabs}

\newcommand{\papertitle}{Optimistic Noise-Aware Sequential Quadratic Programming for Equality Constrained Optimization with Rank-Deficient Jacobians
}
\newcommand{\paperauthora}{Albert~S.~Berahas}
\newcommand{\paperauthoraaffiliation}{Dept. of Industrial and Operations Engineering, University of Michigan}
\newcommand{\paperauthoraemail}{albertberahas@gmail.com}
\newcommand{\paperauthorb}{Jiahao Shi}
\newcommand{\paperauthorbaffiliation}{Dept. of Industrial and Operations Engineering, University of Michigan}
\newcommand{\paperauthorbemail}{jiahaos@umich.edu}
\newcommand{\paperauthorc}{Baoyu~Zhou}
\newcommand{\paperauthorcaffiliation}{School of Computing and Augmented Intelligence, Arizona State University}
\newcommand{\paperauthorcemail}{baoyu.zhou@asu.edu}

\usepackage{booktabs}
\usepackage{caption}
\usepackage{float}
\usepackage{titlesec}
\usepackage{capt-of}

\usepackage{array}
\usepackage{arydshln}
\newcommand{\adasqp}{\texttt{AdaSQP}}
\newcommand{\lssqp}{\texttt{LSSQP}}
\newcommand{\adasqppes}{\texttt{AdaSQP-pes}}
\newcommand{\lssqppes}{\texttt{LSSQP-pes}}
\newcommand{\adasqpopt}{\texttt{AdaSQP-opt}}
\newcommand{\lssqpopt}{\texttt{LSSQP-opt}}
\newcommand{\adasqpinexact}{\texttt{AdaSQP-inexact}}
\newcommand{\lssqpinexact}{\texttt{LSSQP-inexact}}
\newcommand{\adasqpexact}{\texttt{AdaSQP-exact}}
\newcommand{\lssqpexact}{\texttt{LSSQP-exact}}
\newcommand{\ntsqp}{\texttt{NTSQP}}

\allowdisplaybreaks
                
\begin{document}
\title{\papertitle}

\author{\paperauthora\footnotemark[1]\ \footnotemark[2]
   \and \paperauthorb\footnotemark[3]
   \and \paperauthorc\footnotemark[4]}

\maketitle

\renewcommand{\thefootnote}{\fnsymbol{footnote}}
\footnotetext[1]{Corresponding author.}
\footnotetext[2]{\paperauthoraaffiliation. (\url{\paperauthoraemail})}
\footnotetext[3]{\paperauthorbaffiliation. (\url{\paperauthorbemail})}
\footnotetext[4]{\paperauthorcaffiliation. (\url{\paperauthorcemail})}
\renewcommand{\thefootnote}{\arabic{footnote}}

\begin{abstract}{
We propose and analyze a sequential quadratic programming algorithm for minimizing a noisy nonlinear smooth function subject to noisy nonlinear smooth equality constraints. The algorithm uses a step decomposition strategy and, as a result, is robust to potential rank-deficiency in the constraints, allows for two different step size strategies, and has an early stopping mechanism. Under the linear independence constraint qualification, convergence is established to a neighborhood of a first-order stationary point, where the radius of the neighborhood is proportional to the noise levels in the objective function and constraints. Moreover, in the rank-deficient setting, the merit parameter may converge to zero, and convergence to a neighborhood of an infeasible stationary point is established. Numerical experiments demonstrate the efficiency and robustness of the proposed method.}
\end{abstract}

\numberwithin{equation}{section}
\numberwithin{theorem}{section}

\section{Introduction}\label{sec.introduction}

Constrained optimization problems with noisy objective functions and noisy equality constraints are ubiquitous and arise in a plethora of applications, e.g., machine learning~\cite{marquez2017imposing}, statistics~\cite{chatterjee2016constrained}, engineering design~\cite{martins2021engineering}, and simulation optimization~\cite{amaran2016simulation}. In this paper, we consider noisy equality-constrained optimization problems of the form 
\begin{align}\label{prob.main}
  \min_{x\in\mathbb{R}^n}\ f(x)\ \text{ s.t. }\ c(x) = 0,
\end{align}
where $f: \mathbb{R}^n \rightarrow \mathbb{R}$ and $c: \mathbb{R}^n \rightarrow \mathbb{R}^m$ are smooth, possibly nonconvex functions, and where access to the exact values and derivatives of both the objective and constraint functions is not possible. Instead, noisy approximations, denoted by over-bars, are used that satisfy 
\begin{align} 
\label{eq.bounded_noise} 
  \|\bar{f}(x) - f(x) \| \leq \epsilon_f, \  \|\bar{g}(x) - \nabla f(x)\| \leq \epsilon_g, \  \|\bar{c}(x) - c(x) \| \leq  \epsilon_c,\text{ and } \|\bar{J}(x) - \nabla c(x)^T \| \leq  \epsilon_J,
\end{align} 
where $(\fbar(x),\gbar(x),\cbar(x),\bar{J}(x))$ denote objective function and gradient and constraint function and Jacobian approximations, respectively, at $x\in\mathbb{R}^n$, and $\{\epsilon_f,\epsilon_g,\epsilon_c,\epsilon_J\}\subset\R{}_{\geq 0}$ dictate the level of noise. Such problems pose significant challenges for optimization methods, particularly when noise impacts not only the objective function but also the constraints, potentially violating classical assumptions and making the notion of feasibility unclear.

Both the unconstrained and the constrained noisy settings have received much attention due to their wide applicability. In this paper, we focus on \textit{noise-aware} algorithms for solving such problems, i.e., algorithms that exploit information about the noise and that are adaptive. In the unconstrained and bounded noise setting, several noise-aware algorithms that leverage noise-level dependent constants (e.g., $\epsilon_f$ and $\epsilon_g$) to evaluate the acceptability of steps within line search \cite{berahas2019derivative,berahas2021global,jin2024high,xie2020analysis} or trust region \cite{bellavia2023impact,cao2023first,larson2024novel,sun2023trust} methods have been proposed. A natural extension of these algorithms to the constrained setting assumes bounded noise in the objective function and associated derivatives, and possibly in the constraint functions. In this context, and when noise is 
present only in the objective function, \cite{lou2024noise} proposed a line search gradient projection method with a noise level estimation approach, and \cite{na2023adaptive} proposed a line search sequential quadratic programming (SQP) method that uses a differentiable exact augmented Lagrangian as the merit function. 
Related to our problem setting, in the nonlinear equality-constrained setting with bounded noise in both the objective and constraint functions, \cite{oztoprak2023constrained} proposed a line search  SQP method for the setting in which the linear independence constraint qualification (LICQ) is satisfied, and recently a study proposed a trust-region SQP method for the setting in which the LICQ is violated~\cite{sun2024trust}. Under the same noise setting but for a different problem class, in \cite{dezfulian2024convergence} a line search interior-point method is proposed for solving bound-constrained problems with bounded noise in the objective function, and the methodology is extended in \cite{curtis2025interior} to tackle general inequality-constrained problems with bounded noise in the objective and constraint functions under the LICQ.

Another related but distinct problem class that has received significant attention in recent years involves problems with deterministic constraints and a stochastic objective function. In this setting, it is often assumed that an unbiased gradient approximation with bounded variance is available; 
see e.g., \cite{berahas2024modified,berahas2023stochastic,berahas2021sequential,berahas2023accelerating,curtis2023stochastic,curtis2024worst,curtis2024stochastic,fang2024fully,reddi2016stochastic}. This problem class presents a different set of challenges for algorithm design, analysis, and implementation, but shares some commonalities with our problem setting. Perhaps the most closely related works to ours are \cite{berahas2023stochastic,curtis2024stochastic}. In \cite{berahas2023stochastic}, the authors consider situations where the LICQ does not hold and propose an SQP method with a step decomposition strategy inspired by the Byrd-Omojokun approach \cite{omojokun1989trust}. In \cite{curtis2024stochastic}, an SQP method with step decomposition and termination tests for the subproblems is introduced to enhance computational efficiency. Our approach builds upon these two key ideas with the necessary modifications to solve noisy equality-constrained problems.

\subsection{Contributions}

The goal of this paper is to design, analyze, and implement an SQP algorithm for solving noisy equality-constrained optimization problems with possibly rank-deficient constraint Jacobians. Specifically, we focus on problems where both the objective and constraint functions, as well as their associated derivatives, are contaminated with bounded noise. The challenge lies in dealing with this noise in conjunction with the nonlinear constraints and possible rank deficiency. Our algorithm is noise-aware and requires knowledge or estimates of the noise levels. The main contributions can be summarized as follows. 
\begin{itemize}
	\item We propose an SQP algorithm to solve noisy equality-constrained problems \eqref{prob.main}. Our algorithm is built upon the classical line search SQP framework \cite[Chapter 18.4]{NoceWrig06} and stochastic SQP methods \cite{berahas2023stochastic,berahas2021sequential,curtis2024stochastic}, with the following adjustments to account for the noise, the possible rank deficiency, and for practicality.

    \textbf{Optimistic feasibility}: Given the inherent noise in the constraint functions, achieving strict feasibility is not practical. Therefore, we consider feasibility to be attained if the constraint violation falls below a specified threshold related to the noise level in the constraint function ($\epsilon_c$). We treat this as a feature of our algorithm, 
    and build upon the notion of optimistic feasibility. Based on this, we design an early termination condition. This approach not only improves efficiency by potentially reducing the number of iterations but also retains 
    the strong convergence guarantees associated with the SQP framework.

    \textbf{Step decomposition and inexact solution}: Since rank deficiency in the Jacobian of the constraints (violation of the LICQ) is possible due to noise or other factors, we employ a Byrd-Omojukun step decomposition strategy to compute the search direction \cite{omojokun1989trust}. This strategy involves computing normal and tangential steps, which requires solving two subproblems. For practical purposes and to account for the noise, we allow for inexact solutions to these subproblems. We carefully design termination conditions for both subproblems that ensure that the inexact solutions are sufficiently accurate to achieve global convergence while enhancing computational efficiency.

    \textbf{Flexible step size scheme}: We propose two efficient and theoretically sound schemes for selecting the step size parameter: 
    an adaptive scheme and a line search scheme. The adaptive scheme is motivated by \cite{berahas2023stochastic} and is objective-function-free, but requires 
    estimates of the Lipschitz constants of both the objective and constraint functions. 
    The line search scheme, motivated by \cite{berahas2019derivative,NoceWrig06}, requires function evaluations and knowledge of the noise levels (or estimates) in the objective and constraint functions.
    \item We prove that the iterates generated by the proposed SQP method, for both step size selection strategies, exhibit neighborhood convergence behavior, with the neighborhood size dictated by the noise 
    in the objective and constraint functions. Specifically, in the full-rank setting, the iterates converge to a neighborhood of a first-order stationary point of \eqref{prob.main}. When the noisy constraint Jacobian is rank-deficient, the merit parameter may converge to zero, and the iterates converge to a neighborhood of an infeasible stationary point. These results are consistent with those 
    for deterministic SQP variants~\cite{CurtNoceWach10}, and when all noise parameters are zero, we recover the deterministic results. Additionally, we establish a complexity analysis for both feasibility (or infeasible stationarity) and stationarity measures, 
highlighting the efficiency of our approach even under noisy conditions. Under similar noisy conditions, in \cite{sun2024trust} a trust-region SQP method is proposed, and, in addition to the guarantees in the full-rank setting, convergence to a neighborhood of an infeasible stationary point is established in the rank-deficient setting. To establish the result in the full-rank setting, an assumption is made on the discrepancy between the null spaces of the noisy and the true Jacobians \cite[Assumption 5]{sun2024trust}, which can be significant even when the noise levels are small. Our analysis does not require this assumption. Furthermore, contrary to both \cite{oztoprak2023constrained,sun2024trust} which establish only asymptotic convergence of the iterates, we derive non-asymptotic bounds for the average stationarity and feasibility measures across iterations and iteration complexity results. 
    \item We evaluate the performance of our proposed SQP method on a subset of equality-constrained problems from the CUTEst collection. By artificially adding different levels of noise, we conduct controlled experiments and test the robustness 
    of the method, and consider both the full-rank and rank-deficient constraint Jacobian settings. We compare the robustness of the two step size selection strategies to the noise 
    and also compare exact and inexact variants. Our experiments demonstrate that the method is robust to noise variations, with the early termination condition triggering appropriately, and that our method is competitive with existing noise-aware line search SQP methods \cite{oztoprak2023constrained}. 
\end{itemize}

\subsection{Organization}
We conclude this section by introducing the notation used throughout the paper. In Section~\ref{sec.algorithm}, we discuss the main assumptions and present our proposed algorithm. 
The convergence guarantees 
of the algorithm are given in Section~\ref{sec.convergence}. Empirical results are 
presented in Section~\ref{sec.experiments} and concluding remarks are provided in Section~\ref{sec.remarks}.

\subsection{Notation} 
Let $\Nmbb := \{0,1,2,\ldots\}$ denote the set of natural numbers, $\Rmbb$ denote the set of real numbers, $\R{n}$ denote the set of $n$-dimensional real vectors and $\R{m\times n}$ denote the set of $m$-by-$n$ dimensional real matrices. For any $r \in \R{}$, $\R{}_{\geq r}$ (resp., $\R{}_{>r}$) denotes the set of real numbers greater than or equal to $r$ (resp., strictly greater than $r$).  Unless specified, norms $\|\cdot\|$ are 2-norms.
Let $\Range(A)$ and $\Null(A)$ denote the range space and the null space of $A\in\R{m\times n}$, respectively. For any $k\in\Nmbb$, the subscript $k$ denotes the quantity 
evaluated at iterate $x_k$, e.g., $f_k:= f(x_k)$. Any quantity with an overbar denotes a noisy estimate, 
e.g., $\bar{f}_k$ is a noisy estimate of $f_k$.

\section{Noisy Sequential Quadratic Programming Method}\label{sec.algorithm}

In this section, we introduce our proposed noisy sequential quadratic programming (SQP) method. We begin with assumptions and preliminaries.

\subsection{Assumptions and Preliminaries}\label{sec.assumptions}

We make the following assumption about problem \eqref{prob.main} and the iterates generated by the algorithm.

\begin{assumption}\label{ass:prob}
We propose an iterative algorithm that generates a sequence of iterates $\{x_k\}\subset\R{n}$. We assume that there exists an open convex set $\mathcal{X} \subseteq \mathbb{R}^{n}$ containing all 
the iterates $\{x_k\}$ and the trial points $\{x_k+\dbar_k\}$ generated by our algorithm, where $\{\dbar_k\}$ is the sequence of search directions. 
The objective function $f : \R{n} \to \R{}$ is continuously differentiable and bounded below over $\mathcal{X}$, and its gradient $\nabla f : \mathbb{R}^{n} \to \mathbb{R}^{n}$ is Lipschitz continuous with constant $L\in\R{}_{>0}$ and bounded over $\mathcal{X}$. The constraint function $c : \mathbb{R}^{n} \to \mathbb{R}^{m}$ and its Jacobian $\nabla c^T : \mathbb{R}^{n} \to \mathbb{R}^{m \times n}$ are bounded over $\mathcal{X}$, and $\nabla c^T$ is Lipschitz continuous with constant $\Gamma \in\R{}_{\geq 0}$ over $\mathcal{X}$. Specifically, there exist constants $\left(f_{\text{inf }}, f_{\text {sup }}, \kappa_{g}, \kappa_c, \kappa_J\right) \in \mathbb{R} \times \mathbb{R} \times \mathbb{R}_{>0} \times \mathbb{R}_{>0} \times \mathbb{R}_{>0}$ such that for all $x\in\mathcal{X}$,
\begin{align*}
    f_{\text {inf }} \leq f(x) \leq f_{\text {sup }},  \quad\left\|\nabla f(x)\right\| \leq \kappa_{g}, \quad\left\|c(x)\right\| \leq \kappa_c, \qquad \text{and} \qquad \left\|\nabla c (x)\right\| \leq \kappa_J. 
\end{align*}
\end{assumption}

The next assumption pertains to the noisy evaluations of the objective and constraint functions, as well as their gradients.
\begin{assumption}\label{ass:error}
There exist constants $\{\epsilon_f,\epsilon_g,\epsilon_c,\epsilon_J\}\subset\R{}_{\geq 0}$ such that \eqref{eq.bounded_noise} holds for all $x\in\Xcal$. 
\end{assumption}

\begin{remark}  
The smoothness and boundedness assumptions are standard in the deterministic continuous constrained optimization literature~\cite{byrd1987trust,wachter2005line}, and the set $\mathcal{X}$ is assumed to be open and convex but not necessarily bounded. Throughout the paper, we assume that the noise in the objective and constraint functions, as well as their associated gradients, is bounded (Assumption~\ref{ass:error}). Under this noise model, we show that the components of our proposed algorithm (e.g., step size and merit parameter) are well defined. 
\end{remark}

Often, in nonlinear equality-constrained optimization, the linear independence constraint qualification (LICQ) is assumed to hold by imposing that the singular values of the Jacobian $\nabla c(x)^T$ are bounded away from zero; see, e.g., \cite[Assumption 18.1(a)]{NoceWrig06} and~\cite{powell2006fast}. We do not make such an assumption here; however, we consider this setting in our analysis (as a special case) in Section~\ref{sec.convergence}. Before proceeding, we make a few remarks about the singular values of the true and noisy Jacobians of the constraints, as they play a key role in equality-constrained optimization. The singular values of $\nabla c(x)^T$ along with the upper bound of the noise in the constraint Jacobians $\epsilon_J$ dictate 
the rank of $\bar{J}(x)$. Suppose that the singular values of $\nabla c(x)^T$ are lower bounded by 
$\kappa_{\sigma} \in \R{}_{>0}$ for all $x\in\Xcal$, i.e., the true constraint Jacobian is of full rank. If $\epsilon_J < \kappa_{\sigma}$ (as assumed in \cite{oztoprak2023constrained}), it then follows from Weyl's Theorem \cite{weyl1912asymptotische} that the singular values of $\bar J(x)$ are lower bounded by $\kappa_{\sigma}-\epsilon_J $ for all $x\in\Xcal$. However, if $\epsilon_J \ge \kappa_{\sigma}$, then the full-rank property of $\bar J(x)$ cannot be guaranteed. Now, suppose that $\nabla c(x)^T$ is rank-deficient. In this setting, under Assumption~\ref{ass:error} it is possible that $\bar J(x)$ is also rank-deficient. That said, there are also examples for which $\bar J(x)$ is of full rank with high probability; see~\cite{benaych2012singular,capitaine2007strong}. In this paper, we do not make explicit assumptions about the relationship between the singular values of $\nabla c(x)$ and $\epsilon_J$, but consider the full rank setting as a special case in our analysis.

Let $\mathcal{L}(x, y)=f(x)+c(x)^T y$ denote the Lagrangian corresponding to \eqref{prob.main}, where $y \in \mathbb{R}^m$ represents a vector of Lagrange multipliers. Under Assumption~\ref{ass:prob} and the LICQ, a first-order stationary point $x^*$ with respect to \eqref{prob.main} 
satisfies
\begin{align}
\label{eq.first_order_sta}
    \nabla f(x^*)+ \nabla c(x^*) y^*  = 0 \quad \text{and} \quad  c (x^*)  = 0  \;\ \text{for some } y^* \in \mathbb{R}^m.
\end{align}
Throughout this manuscript, we define the stationarity error at $(x,y) \in (\mathbb{R}^{n},\mathbb{R}^{m}) $ as $\| \nabla f(x)+\nabla c(x) y \|$ and the feasibility error at $x \in \mathbb{R}^n$ as $\|c(x)\|$. In the absence of constraint qualifications, it is possible that \eqref{prob.main} is infeasible, or \eqref{prob.main} can be degenerate, in which case \eqref{eq.first_order_sta} might not be a necessary condition. To address the difficulties that arise from the lack of constraint qualifications, it is essential to develop an algorithm that can transition 
its objective from finding a stationary point of \eqref{prob.main} to identifying a stationary point relative to a measure of infeasibility. To this end, we define the infeasibility measure $\varphi: \mathbb{R}^n \to \mathbb{R}$ as $\varphi(x) = \frac{1}{2}\|c(x)\|^2$ 
and a point $\tilde x \in \mathbb{R}^n$ is considered stationary with respect to $\varphi$ if one of the following conditions is satisfied: $(i)$ $c(\tilde x) = 0$, or $(ii)$ $c(\tilde x) \neq 0$ and $0=\nabla \varphi(\tilde x)=\nabla c(\tilde x) c(\tilde x)$, 
where case $(ii)$ is the definition of an infeasible stationary point. Henceforth, we define the infeasible stationarity error at $x \in \mathbb{R}^n$ as $\|\nabla c(x) c(x)\|$.

\subsection{A Noisy Sequential Quadratic Programming Method}\label{sec.noisy_SQP}

In this section, we describe our proposed SQP method. The algorithm has foundations in the line search SQP paradigm \cite[Chapter 18.4]{NoceWrig06}, with several modifications to account for the noise and possible violation of the LICQ. 

To address the challenges posed by the discrepancy between the observed noisy constraint functions $\bar{c}(x)$ and the true constraint functions $c(x)$, our algorithm has an early termination feature. The concept, termed \textit{optimistic feasibility}, is simple yet powerful. Due to the noise, even if $\bar{c}(x) = 0$, 
we cannot conclusively assert feasibility with respect to the true constraint functions. In fact, from a pessimistic standpoint, exact feasibility may never be guaranteed under our general noise model. To overcome this limitation, we adopt an optimistic perspective to feasibility. Specifically, if $\|\bar{c}(x)\| \leq \epsilon_o$ for some $\epsilon_o \in \mathbb{R}_{\geq 0}$, we consider feasibility to be approximately achieved. Consequently, the search direction at an approximately feasible point focuses on improving the objective function while maintaining approximate feasibility. This approach prevents the algorithm from becoming sensitive to the noise, which could result in the unnecessary minimization of constraint violation. 
Moreover, our algorithm has an early termination scheme that is triggered when $\|\bar{c}(x)\| \leq \epsilon_o$ and $\Delta \lbar(x_{k},\bar \tau_{k}, \bar d_{k})   \le  \epsilon_o$ and that returns an approximately first-order stationary iterate. 
In our algorithm, $\epsilon_o \in [0, \epsilon_c]$ to ensure that if the algorithm terminates early, it does so at an iterate that is sufficiently feasible with respect to the noise in the constraint functions. Moreover, $\epsilon_o = 0$ represents the pessimistic setting and $\epsilon_o = \epsilon_c$ represents the most optimistic setting.  In either setting, neighborhood convergence is established. In the optimistic setting, the algorithm may terminate early at an approximately first-order stationary point. 

At every iteration $k\in\N{}$ of our proposed algorithm, the search direction $\dbar_k$ is computed using a step decomposition technique~\cite{omojokun1989trust}, i.e., $\dbar_k = \vbar_k + \ubar_k$ where $\vbar_k \in \Range(\bar J_k^T)$ is the \textit{normal component} and $\ubar_k \in \Null(\bar J_k)$ is the \textit{tangential component}. Each component is computed by solving an optimization subproblem. In particular, when $\|\bar{J}_k^T\bar{c}_k\| \neq 0$, the algorithm first computes the normal component $\bar v_k$ as an approximate solution of the subproblem
\begin{align}
    \min _{\vbar \in \operatorname{Range}\left(\bar J_{k}^{T}\right)} \tfrac{1}{2}\left\|\bar c_{k}+\bar J_{k} \vbar\right\|^{2} \quad \text{ s.t. } \| \vbar  \| \le\sigma_{Jc} \|\bar J_{k}^T \bar c_{k}\|,
\label{eq.nolicq.tr.sub.p1}
\end{align}
where $\sigma_{Jc}\in\R{}_{>0}$,
that is required to satisfy the Cauchy-like decrease condition
\begin{align}
    \left\|\bar c_{k}\right\|-\left\|\bar c_{k}+\bar J_{k} \bar v_{k}\right\| \geq \bar \gamma_{c}\left(\left\|\bar c_{k}\right\|-\left\| \bar c_{k}+\bar \alpha_{k}^{c} \bar J_{k}\bar  v_{k}^{c}\right\|\right),
\label{eq.nolicq.tr.cauchy.decrease}
\end{align}
where $\bar \gamma_{c} \in(0,1]$ is a 
parameter, $\bar v_{k}^{c}:=-\bar J_{k}^{T} \bar c_{k}$ is the negative gradient direction of the objective function of \eqref{eq.nolicq.tr.sub.p1} at $\vbar=0$, and $\bar \alpha_{k}^{c} = \min \left\{ \sigma_{Jc}, \left\|\bar J_{k}^{T} \bar c_{k}\right\|^{2} /\left\|\bar J_{k} \bar J_{k}^{T} \bar c_{k}\right\|^{2} \right\} $ is the optimal step size along the direction of $\bar v_{k}^{c}$ that minimizes $\left\|\bar c_{k}+\alpha \bar J_{k}\bar v_{k}^{c}\right\|^2$ within the trust region constraint defined in \eqref{eq.nolicq.tr.sub.p1}. 
Subproblem \eqref{eq.nolicq.tr.sub.p1} minimizes the reduced linearized constraint violation subject to 
a trust-region constraint. In our optimistic approach, we skip subproblem \eqref{eq.nolicq.tr.sub.p1} 
when $\|\bar c_k\| \le \epsilon_o$ and set $\bar v_k = 0$, as reducing constraint violation is unnecessary. Efficient iterative techniques for solving \eqref{eq.nolicq.tr.sub.p1} include the Steihaug-Toint conjugate gradient method (CG-Steihaug) \cite{steihaug1983conjugate}  and GLTR \cite{gould1999solving}.

After computing 
the normal component $\vbar_k$, the algorithm computes the tangential component $\ubar_k$ as an approximate solution to the subproblem 
\begin{align}
\min_{\ubar \in \mathbb{R}^{n}} \left(\bar g_{k}+H_{k} \bar v_{k}\right)^{T} \ubar+\tfrac{1}{2} \ubar^{T} H_{k} \ubar  \quad \text { s.t. } \bar J_{k} \ubar=0,
\label{eq.nolicq.sub.p2}
\end{align}
whose optimal solution $\bar u_{k}^*$ and the associated dual multiplier $\bar y_{k}^*$ (under Assumption~\ref{ass:H}) can be computed via the following linear system of equations
\begin{align}
    \left[\begin{array}{cc}
H_{k} & \bar J_{k}^{T} \\
\bar J_{k} & 0
\end{array}\right]\left[\begin{array}{c}
\bar u_{k}^* \\
\bar y_{k}^*
\end{array}\right]=-\left[\begin{array}{c}
\bar g_{k}+H_{k} \bar v_{k}  \\
0
\end{array}\right].
\label{eq.nolicq.sub.p2.linear.system.star}
\end{align}
Note that even when $\bar J_k$ is rank-deficient, the solution component  $\bar u_{k}^*$ is unique while the dual multiplier component $\bar y_{k}^*$ is not; see Lemma \ref{lemma.y.exist}. Instead of computing an exact solution $(\bar u_{k}^*,\bar{y}_{k}^*)$, our algorithm relies on $(\bar{u}_k,\bar{y}_k)$, an inexact solution of~\eqref{eq.nolicq.sub.p2.linear.system.star}, which can be computed via iterative Krylov subspace methods (e.g., MINRES \cite{choi2011minres}) that only require matrix-vector multiplications. We further define the residual of inexact solution $(\bar{u}_k,\bar{y}_k)$ as
\begin{align}
\left[\begin{array}{c}
\bar\rho_{k} \\
\bar r_{k}
\end{array}\right]:=\left[\begin{array}{cc}
H_{k} & \bar J_{k}^{T} \\
\bar J_{k} & 0
\end{array}\right]\left[\begin{array}{c}
\bar u_{k} \\
\bar y_{k}
\end{array}\right]+\left[\begin{array}{c}
\bar g_{k}+H_{k}\bar  v_{k} \\
0
\end{array}\right].
\label{eq.nolicq.sub.p2.linear.system.res}
\end{align}
We make the following standard assumptions about the matrix $H_k$ in \eqref{eq.nolicq.sub.p2}--\eqref{eq.nolicq.sub.p2.linear.system.res}. 
\begin{assumption}\label{ass:H}
The sequence $\{H_k\}$ is bounded, e.g., there exists a constant $\kappa_H \in \mathbb{R}_{>0}$ such that $\|H_k\| \leq \kappa_H$ for all $k\in\N{}$. In addition, for all $k \in \mathbb{N}$, the matrix $H_k$ is sufficiently positive definite in the null space of both $\nabla c_k^T$ and $\bar{J}_k$, e.g., there exists a constant $\zeta \in \mathbb{R}_{>0}$ such that $u^TH_ku \geq \zeta \|u\|^2$ for all $u\in \Null(\nabla c_k^T) \cup \Null(\bar{J}_k)$.
\end{assumption}

\begin{remark}
Classical SQP methods assume that $H_k$ is sufficiently positive definite in the null space of $\nabla c_k^T$ \cite[Assumption 18.1(b)]{NoceWrig06}. Although $\nabla c_k$ is not computable in our setting, this assumption remains necessary for our analysis. Moreover, we require that $H_k$ is sufficiently positive definite in the null space of the noisy constraint Jacobian $\bar{J}_k$. The identity matrix satisfies Assumption~\ref{ass:H}. In~\cite{oztoprak2023constrained}, $H_k$ is a constant diagonal matrix that is strictly positive definite. We do not impose such structure on $H_k$ in this paper. 
\end{remark}

Next, our proposed algorithm relies on 
a merit function to balance the two possible competing goals of minimizing the objective function and minimizing the constraint violation. To this end, we use an exact $\ell_2$-merit function
\begin{align}\label{eq.merit}
    \bar\phi(x,{\bar\tau}) := \bar\tau \bar f(x) + \|\bar c(x)\| 
\end{align}
where $\bar \tau \in \mathbb{R}_{>0}$ is a merit parameter. Several classical SQP methods (e.g., \cite{Han77,powell2006fast}) use the $\ell_1$ norm in the merit function. We use the $\ell_2$ norm instead for ease of presentation and analysis, but note that the $\ell_1$ norm or other norms could be used. Our algorithm makes use of a local first-order model of the merit function $\bar l : \mathbb{R}^{n} \times \mathbb{R}_{>0} \times \mathbb{R}^{n} \to \mathbb{R}$, and more specifically, the reduction in the model of the merit function $\Delta \bar{l}:\R{n}\times\R{}_{>0}\times\R{n}\to\R{}$ 
\begin{equation}
\label{def.merit_model_reduction}
    \Delta \bar{l}(x,\bar\tau,\bar d) := \bar l(x,\bar\tau,0) - \bar l(x,\bar\tau,\bar d) =  -\bar \tau \gbar(x)^T \bar d + \|\cbar(x)\| - \|\cbar(x) + \bar{J}(x) \bar d\|. 
\end{equation}
The reduction in the model in the merit function is used in our early termination condition, i.e., when both $\|\cbar_k\| \le \epsilon_o$ and $\Delta \bar{l}(x_k,\bar\tau_k,\bar d_k) \le \epsilon_o$ the algorithm terminates. 

The algorithm also uses the reduction in the model of the merit function in the termination tests and in the merit parameter update scheme. Two different termination tests are used depending on whether $\|\cbar(x)\| \le \epsilon_o$ or $\|\cbar(x)\| > \epsilon_o$. Given inexact normal components $\bar{v}_k$ that satisfy \eqref{eq.nolicq.tr.cauchy.decrease}, inexact tangential components $\bar{u}_k$ that satisfy the following termination conditions are sufficient to guarantee convergence.

\begin{termtest}\label{tt1}
This termination test is only considered 
when $\|\cbar_k\| \le \epsilon_o$. Given $\lambda_{\rho r}\in (0,1), \kappa_{\rho r}  \in \mathbb{R}_{>0}, \lambda_u \in(0, \zeta)$ where $\zeta \in \mathbb{R}_{>0}$ is given in Assumption~\ref{ass:H}, 
and $\sigma_{u} \in(0,1)$, $\left( \bar  u_{k}, \bar  \rho_{k}, \bar  r_{k}\right)$ satisfy Termination Test 1 if
\begin{align}  
 \max\{ \left\|\bar\rho_{k}\right\|, \left\|\bar r_{k}\right\| \} & \leq  \lambda_{\rho r} \min \left\{  \max \left\{ \|\bar  u_{k}\|, \|\bar J_k^T \bar c_k\| \right\}, \kappa_{\rho r} \right\},  
\label{eq.tt1.cond1} \\
\bar u_{k}^{T} H_{k} \bar u_{k}  &\geq \lambda_u\left\|\bar u_{k}\right\|^{2} -  \epsilon_o \label{eq.tt1.cond2} \\ 
\bar g_{k}^{T}\bar u_{k}+\tfrac{1}{2} \bar u_{k}^{T} H_{k}\bar u_{k} &\leq 
\epsilon_o, \label{eq.tt1.cond3} \\ 
\text{and } \qquad  \Delta \bar  l\left(x_{k}, \bar\tau_{k-1},\bar u_k 
\right) &\geq \bar\tau_{k-1}\sigma_u  \max \left\{\bar u_{k}^{T} H_{k}\bar u_k, \lambda_u\left\|\bar u_{k}\right\|^{2}\right\} - \epsilon_o \label{eq.inexact.delta.l.desired.opt.case1}
\end{align}
hold. The algorithm then sets $\bar{\tau}_k = \bar{\tau}_{k-1}$.
\end{termtest}

\begin{termtest}\label{tt2}
This termination test is only considered when $\|\cbar_k\| > \epsilon_o$. Given $\lambda_{\rho r}\in (0,1), \kappa_{\rho r}  \in \mathbb{R}_{>0}$,  
$\lambda_{uv} \in \mathbb{R}_{>0},  \lambda_{u} \in(0, \zeta)$ where $\zeta \in \mathbb{R}_{>0}$ is given in Assumption~\ref{ass:H},  $\lambda_{v} \in \mathbb{R}_{>0}$, $\sigma_{u} \in(0,1)$, $\sigma_{c} \in(0,1)$, $\sigma_{r} \in(\sigma_{c},1)$, and $\bar v_{k} \in \operatorname{Range}\left(\bar J_{k}^{T}\right)$ satisfying~\eqref{eq.nolicq.tr.cauchy.decrease}, $\left( \bar  u_{k}, \bar  \rho_{k}, \bar  r_{k}\right)$ satisfy Termination Test 2 if condition \eqref{eq.tt1.cond1} holds and
\begin{align}  
    \left\|\bar u_{k}\right\| \leq \lambda_{uv} \left\|\bar v_{k}\right\| \text{ or }   \left\{\begin{aligned}\bar  u_{k}^{T} H_{k} \bar u_{k} & \geq \lambda_u\left\|\bar u_{k}\right\|^{2} \\\left(\bar g_{k}+H_{k} \bar  v_{k}\right)^{T}\bar u_{k}+ \max \left \{ \tfrac{1}{2}, 1 - \|\bar J_k^T \bar c_k\| \right \} \bar u_{k}^{T} H_{k}\bar u_{k} & \leq \lambda_{v} \left\|\bar v_{k}\right\|\end{aligned}\right\}\label{eq.tt2.cond2}
\end{align} 
together with at least one of the following two conditions,
\begin{align}
    \Delta  \bar l\left(x_{k}, \bar\tau_{k-1},  \bar d_k 
    \right) &\geq \bar\tau_{k-1}\sigma_u  \max \left\{\bar u_{k}^{T} H_{k}\bar u_k, \lambda_u\left\|\bar u_{k}\right\|^{2}\right\} \nonumber\\
    & \qquad + \sigma_{c}\left(\left\|\bar c_{k}\right\|-\left\|\bar c_{k}+\bar J_{k} \bar v_{k}\right\|\right) \label{eq.inexact.delta.l.desired.opt.case2} \\
    \text{or}\qquad \left\|\bar c_{k}\right\|-\|\bar c_{k}+\bar J_{k} \bar v_{k}+\bar r_{k}\| &\geq \sigma_r\left(\left\|\bar c_{k}\right\|-\left\|\bar c_{k}+\bar J_{k} \bar v_{k}\right\|\right)>0. \label{eq.tt2.cond1}
\end{align}
When~\eqref{eq.inexact.delta.l.desired.opt.case2} is satisfied, we set $\bar \tau_k = \bar \tau_{k-1}$. If~\eqref{eq.inexact.delta.l.desired.opt.case2} is not satisfied but~\eqref{eq.tt2.cond1} holds, given a user-defined parameter $\sigma_{\tau}\in(0,1)$, the algorithm sets
\begin{align} 
    \bar \tau_{k} \leftarrow \begin{cases} \bar\tau_{k-1} & \text { if } \bar\tau_{k-1} \leq \left(1-\sigma_{\tau}\right) \bar\tau_{k}^{\text {trial }} \\\left(1-\sigma_{\tau}\right)  \bar \tau_{k}^{\text {trial }}  & \text { otherwise, }\end{cases} \label{eq.tau.update.tt2} 
\end{align}
where
\begin{align}
    \bar\tau_{k}^{\text {trial }} \leftarrow \begin{cases}\infty & \text {if } \bar g_{k}^{T}\bar d_k
    +\max \left\{\bar u_{k}^{T} H_{k}\bar u_{k}, \lambda_u\left\|\bar u_{k}\right\|^{2}\right\} \leq 0 \\ \tfrac{\left(1-\tfrac{\sigma_{c}}{\sigma_{r}}\right)\left(\left\|\bar c_{k}\right\|-\|\bar c_{k}+\bar J_{k} \bar v_{k}+\bar r_{k}\|\right)}{\bar g_{k}^{T}\bar d_k
    +\max \left\{\bar u_{k}^{T} H_{k}\bar u_{k}, \lambda_u\left\|\bar u_{k}\right\|^{2}\right\}} & \text {otherwise. }\end{cases}\label{eq.tau.trial.tt2}
\end{align}
\end{termtest}

The conditions in Termination Test~\ref{tt2} are similar to those presented in~\cite{curtis2024stochastic}. The main motivation is to ensure that the reduction in the model of the merit function is sufficiently positive. This is achieved by either directly setting $\bar \tau_{k} = \bar \tau_{k-1}$ when \eqref{eq.inexact.delta.l.desired.opt.case2} is satisfied or via the merit parameter update rule \eqref{eq.tau.update.tt2}--\eqref{eq.tau.trial.tt2}. The conditions in Termination Test~\ref{tt1} are different, mainly due to the optimistic feasibility setting. 
In~\eqref{eq.tt1.cond2}--\eqref{eq.inexact.delta.l.desired.opt.case1}, we add a relaxation term ($\epsilon_o$) to ensure that Termination Test~\ref{tt1} can be satisfied when the residual (i.e., $\max\{\|\bar\rho_k\|,\|\bar{r}_k\|\}$) is sufficiently small. 
Minimizing the constraint violation is not a priority in this optimistic setting ($\|\bar v_{k}\| = 0$ since $\|\cbar_k\| \le \epsilon_o$), and we simply set $\bar \tau_k = \bar \tau_{k-1}$ and consider 
the model reduction condition given 
in \eqref{eq.inexact.delta.l.desired.opt.case1}. Note that the relaxation parameters in \eqref{eq.tt1.cond2}--\eqref{eq.inexact.delta.l.desired.opt.case1} are only required to be proportional to $\epsilon_o$. 
We use $\epsilon_o$ as the relaxation parameter as it simplifies the analysis and avoids the need for an additional parameter. 
The smaller the relaxation parameters, the more accurate $\bar u_k$ is required to be. When $\epsilon_o = 0$, conditions \eqref{eq.tt1.cond2}--\eqref{eq.inexact.delta.l.desired.opt.case1} reduce to those in Termination Test~\ref{tt2}. The $\bar \tau_k$ update rules in both termination tests ensure that $\{\bar \tau_k\}$ is monotonically non-increasing. It is important to note that \eqref{eq.inexact.delta.l.desired.opt.case1} does not necessarily ensure that the reduction in the model of the merit function is non-negative, i.e., $\Delta \bar l\left(x_{k}, \bar\tau_{k},\bar d_k \right) = \Delta \bar l\left(x_{k}, \bar\tau_{k},\bar u_k 
\right) \ge 0$, which poses a significant challenge as this quantity is no longer a reliable measure of optimality. Fortunately, when $\Delta \bar l\left(x_{k}, \bar\tau_{k},\bar d_k \right) < 0$, one can show that both the stationarity and feasibility errors are sufficiently small and comparable in magnitude to the noise levels in the problem. 
In particular, Theorem~\ref{theorem.early.termination}  
guarantees that if early termination is triggered due to $\|\cbar_k\| \le \epsilon_o$ and $\Delta \bar  l\left(x_{k}, \bar\tau_{k},\bar d_k \right) \le  \epsilon_o$ at some iteration, then both the stationarity and feasibility errors are sufficiently small.  Note that the termination condition for $\Delta \bar  l\left(x_{k}, \bar\tau_{k}, \bar d_k\right)$ remains robust even if the constant on the right-hand side is replaced by different constant that is proportional to $\epsilon_o$. 
Another termination condition arises when $\|\cbar_k\| > \epsilon_o$ and $\| \bar J_k^T \bar c_k \|  = 0$, corresponding to iterates that are approximately infeasible stationary. This termination condition aligns with SQP methods developed under the rank-deficient constraint Jacobian setting (see, e.g., \cite{CurtNoceWach10,curtis2024stochastic}). Unlike prior works, due to the presence of noise in the constraints, we can only conclude that the terminated iterate is approximately an infeasible stationary point, with the associated error characterized in Theorem~\ref{theorem.early.termination.infeasible}.
On the other hand, if the algorithm does not terminate early, 
then an averaged measure of 
stationarity and feasibility error and/or infeasible stationarity error is guaranteed to be small (Theorems~\ref{theorem.tau.lb} and \ref{theorem.tau.tozero}).

The final component of the algorithm is the mechanism for selecting the step sizes $\{\bar \alpha_k \} \subset \R{}_{>0}$ used to update the iterates, i.e., $x_{k+1} \gets x_k + \bar \alpha_k \bar d_k$. Our algorithm, see Algorithm~\ref{alg.dfo_sqp_LS}, allows for two 
step size schemes: an adaptive scheme inspired by \cite{berahas2023stochastic} and a line search scheme inspired by \cite{berahas2019derivative,NoceWrig06}. The goal of both step size procedures is to find a step size that ensures sufficient decrease as measured by the progress made on the merit function. A key difference between the two schemes is that the former is objective-function-free but requires knowledge (or estimates) of the Lipschitz constants of the objective and constraint gradients, while the latter requires function evaluations (that may be contaminated with noise) and knowledge (or estimate) of the noise level in the objective and constraint functions. We discuss both approaches in detail below.
 
The adaptive step size scheme (\textbf{Option I}) is very similar to that proposed in \cite{berahas2023stochastic}. For completeness and in order to introduce the notation used in the algorithm and analysis, we fully present the procedure. The key idea is to carefully select a step size based on the magnitude of the normal and tangential components, 
to safe-guard our algorithm from potential rank deficiency in the constraint Jacobians. The exact balance between 
these two components is problem-dependent and, to avoid the need for tuning, 
our algorithm generates two sequences $\{\bar \chi_k\}$ and $\{\bar \zeta_k\}$ to distinguish between tangentially and normally dominated steps. Specifically, at iteration $k \in \N{}$, after the two steps $(\bar{v}_k,\bar{u}_k)$ are computed, 
the conditions
\begin{align}
    \left\|\bar u_k\right\|^2 \geq \bar \chi_{k-1}\left\|\bar v_k\right\|^2 \qquad \text { and } \qquad \tfrac{1}{2} (\bar{v}_k + \bar{u}_k)^T H_k (\bar{v}_k + \bar{u}_k) < \tfrac{1}{4} \bar \zeta_{k-1}\left\|\bar u_k\right\|^2,
\label{eq.dominate.condition}
\end{align}
are employed to update the sequences via 
\begin{align}
    \left(\bar \chi_k, \bar \zeta_k\right) \leftarrow \begin{cases}\left(\left(1+\sigma_\chi\right) \bar \chi_{k-1},\left(1-\sigma_\zeta\right) \bar \zeta_{k-1}\right) & \text { if } \eqref{eq.dominate.condition} \text { holds } \\ \left(\bar \chi_{k-1}, \bar \zeta_{k-1}\right) & \text { otherwise}, \end{cases}
\label{eq.chi.update}
\end{align}
where $\sigma_\chi \in \mathbb{R}_{>0}$ and $\sigma_\zeta \in(0,1)$ are user-defined parameters. Under the conditions that $\bar\chi_{-1} > 0$ and $\bar\zeta_{-1} > 0$, by \eqref{eq.chi.update}, it follows that $\left\{\bar \chi_k\right\}$ is monotonically non-decreasing and $\left\{\bar \zeta_k\right\}$ is monotonically non-increasing. Moreover, our analysis will show that $\left\{\bar \chi_k\right\}$ (resp. $\left\{\bar \zeta_k\right\}$) is bounded above (resp.~below) uniformly by a positive real number, i.e., these sequences are bounded deterministically and $\left(\bar \chi_k, \bar \zeta_k\right)=\left(\bar \chi_{k-1}, \bar \zeta_{k-1}\right)$ for all sufficiently large $k \in \mathbb{N}$ in any run of the algorithm (see Lemma~\ref{lemma.chi.bound}). With $\bar d_k = \bar u_k + \bar v_k$ (note that $\bar d_k \ne 0$ by Lemma~\ref{lemma.d.nonzero.opt}), the algorithm then sets
\begin{align}
    \bar \xi_{k}:= \begin{cases}\bar \xi_{k-1} & \text { if }\bar  \xi_{k-1} \leq \bar \xi_{k}^{\text {trial}} \\ \min \left\{\left(1-\sigma_{\xi}\right) \bar \xi_{k-1}, \bar \xi_{k}^{\text {trial}}\right\} & \text { otherwise, }\end{cases}
\label{eq.inexact.xi.update}
\end{align}
where
\begin{align}
    \bar \xi_{k}^{\text{trial}}:= \begin{cases}\tfrac{\Delta \lbar(x_{k},\bar \tau_{k}, \bar d_{k}) }{\bar \tau_{k}\|\bar d_{k}\|^{2}} & \text { if }\left\|\bar u_k\right\|^2 \geq \bar \chi_{k}\left\|\bar v_k\right\|^2 \\\tfrac{\Delta \lbar(x_{k},\bar \tau_{k},\bar d_{k}) }{\|\bar d_{k}\|^{2}} & \text { otherwise, }\end{cases}
    \label{eq.inexact.xi.trial.update}
\end{align}
and $\sigma_\xi \in (0,1)$ is a user-defined parameter and $\bar{\chi}_k$ is given in \eqref{eq.chi.update}. By \eqref{eq.inexact.xi.update}
, $\bar\xi_{k} \leq \bar \xi_{k}^{\text{trial}}$ for all $k \in \mathbb{N}$, and $\{\bar\xi_{k}\}$ is monotonically non-increasing and bounded away from 0 (see Lemma~\ref{lemma.xi}). Then, the algorithm sets 
\begin{align}
    \bar\alpha_{k}^{\text {suff}} := \min \left\{\tfrac{2(1-\eta) \beta_{k} \Delta \lbar(x_{k},\bar \tau_{k}, \bar d_{k})   }{\left(\bar\tau_{k} L+\Gamma\right)\|\bar d_{k}\|^{2}}, 1\right\}, 
 \label{eq.inexact.alpha.suff}   
\end{align}
where $\eta\in(0,1)$ and $\left\{\beta_{k}\right\}$ are user-defined parameters satisfying $\tfrac{2(1-\eta) \beta_{k}\bar \xi_{-1} \max \{ \bar \tau_{-1}, 1 \}}{ \bar \tau_{-1}L + \Gamma} \in (0,1]$ for all $k \in \mathbb{N}$ with $(\bar\tau_{-1},\bar\xi_{-1})$ representing the initial values of sequences $\{\bar\tau_k\}$ and $\{\bar\xi_k\}$, respectively. We then employ a safeguarding scheme to ensure that the step size $\bar\alpha_k$ 
is not too large or too small. To this end, since $\bar\alpha_{k}^{\text {suff}}$ is based on stochastic quantities, we use the sequence $\{\bar\xi_k\}$, which is monotonically non-increasing and bounded from zero, as the proxy, and define
\begin{align}
    \bar\alpha_{k}^{\min }:= \begin{cases}\tfrac{2(1-\eta) \beta_{k}\bar \xi_{k} \bar\tau_{k}}{\bar\tau_{k} L+\Gamma} & \text { if }\left\|\bar u_k\right\|^2 \geq \bar  \chi_{k}\left\|\bar v_k\right\|^2 \\ \tfrac{2(1-\eta) \beta_{k}\bar \xi_{k} }{\bar\tau_{k} L+\Gamma} & \text { otherwise, }\end{cases}
 \label{eq.inexact.alpha.min}   
\end{align}
and $\bar\alpha_{k}^{\max }:= \bar \alpha_{k}^{\min }+\theta \beta_{k}$, where $\theta \in \mathbb{R}_{>0}$ is 
an input parameter. 
The algorithm then sets the step size $\bar\alpha_{k}$ by projecting $\bar\alpha_{k}^{\text {suff}}$ onto the interval $[\bar \alpha_{k}^{\min}, \bar \alpha_{k}^{\max}]$, i.e.,  
$\bar\alpha_{k} \gets \operatorname{Proj}_k(\bar\alpha_{k}^{\text {suff}})$, where
\begin{align}
    \label{eq.project}
    \operatorname{Proj}_k (\bar\alpha_{k}^{\text {suff}} )  =  \operatorname{Proj} \left(\bar\alpha_{k}^{\text {suff}}  \big| [\bar \alpha_{k}^{\min}, \bar \alpha_{k}^{\max}] \right) 
\end{align}
and $\operatorname{Proj}(\bar\alpha_{k}^{\text {suff}}| \mathcal{I})$ denotes the projection of $\bar\alpha_{k}^{\text {suff}}$ onto the interval $\mathcal{I} \subset \mathbb{R}$. Note that this projection step will never increase the step size (see Lemma~\ref{lemma.inexact.stepsize.opt}). This is necessary especially when $\bar{J}_k$ is rank-deficient and the merit parameter $\bar \tau_k$ is close to 0, which makes $\bar \alpha_{k}^{\min}$ close to 0 when $\|\bar u_k\|$ dominates  $\|\bar v_k\|$ . In this case, there could be a large gap between $\bar\alpha_{k}^{\text {suff}}$ and $\bar \alpha_{k}^{\min}$, and the projection scheme makes this gap controllable. 

The other option for the step size selection is to use a line search mechanism \textbf{Option II}. Of course, 
noisy evaluations of the objective and constraint functions 
are required by this option (while \textbf{Option I} is objective- and constraint-function-free). This approach is inspired by the relaxed Armijo condition proposed in \cite{berahas2019derivative} for the unconstrained setting with bounded noise and classical line search SQP methods \cite{NoceWrig06}. 
By Assumption \ref{ass:error} and $\phi(x_k,\bar\tau_k) := \bar\tau_kf(x_k) + \|c(x_k)\|$, it follows that 
\begin{align}
\label{eq.diff.merit.func}
    |\bar\phi(x_{k},\bar{\tau}_{k}) - \phi(x_{k},\bar{\tau}_{k}) | \leq \bar\tau_k|\bar f(x_k) - f(x_k)| + \|\bar c(x_k) - c(x_k)\| \le \bar{\tau}_{k} \epsilon_f + \epsilon_c.
\end{align}
We consider a modified version of the sufficient decrease Armijo condition, applied to the merit function \eqref{eq.merit},
\begin{equation}
   \bar\phi(x_{k} + \bar\alpha_k \bar{d}_k,\bar{\tau}_{k}) \le  \bar\phi(x_{k},\bar{\tau}_{k}) - \eta \bar\alpha_k \Delta \lbar(x_{k},\bar \tau_{k}, \bar d_{k})   + \epsilon_{A_k},
  \label{eq.modified.linesearch.1d}
\end{equation}
where $\eta \in (0,1)$ and $\epsilon_{A_k} \in \mathbb{R}_{\geq 0}$ is a user-defined parameter that depends on $(\epsilon_f,\epsilon_c,\epsilon_g,\epsilon_J)$ and accounts for the noise. The step size is selected via a standard backtracking mechanism, where the initial trial step size is a user-defined constant $ \alpha_u^L \in (0,1]$. 
If a given trial step size does not satisfy \eqref{eq.modified.linesearch.1d}, then the trial step size is reduced by a constant factor. This process is repeated until a step size that satisfies \eqref{eq.modified.linesearch.1d} is found (see Lemma \ref{lemma.step.size.ls}). The relaxation term $\epsilon_{A_k}$ is crucial for both theory (by guaranteeing a well-defined line search procedure) and practice (by potentially accepting longer step sizes). The explicit expression of the relaxation term $\epsilon_{A_k}$ is given in Lemma~\ref{lemma.step.size.ls}. 

The complete algorithm (with both step size options) is given in Algorithm~\ref{alg.dfo_sqp_LS}.
\begin{algorithm}[]
  \caption{Optimistic Noise-Aware SQP Algorithm}
  \label{alg.dfo_sqp_LS}
  \begin{algorithmic}[1]
    \Require $x_0 \in \mathbb{R}^{n}$ (initial iterate); $\bar \tau_{-1} \in \mathbb{R}_{>0}$ (initial merit parameter value); 
     $\epsilon_o \in [0,\epsilon_c]$ (optimistic feasibility parameter); 
     
     \hspace{-1cm}\textbf{Require (Option I: Adaptive):} $\{ \beta_k \} \subset (0,1]$ (adaptive step size parameter), $L \in \mathbb{R}_{>0}$, $\Gamma \in \mathbb{R}_{>0}$ (Lipschitz constant estimates), $\{\bar\chi_{-1}, \bar\xi_{-1}, \bar\zeta_{-1}\} \subset \mathbb{R}_{>0}$ (initial values for additional sequences) 
    
    \hspace{-1cm}\textbf{Require (Option II: Line Search):} $ \alpha_u^L \in (0,1]$ (initial step size),  $\nu \in (0,1)$ (search control parameters), $\epsilon_{A_k}$ (line search relaxation parameter)
    
    \For{\textbf{all} $k \in \mathbb{N}$}
     \State  Compute gradient  $\bar{g}_k$ and Jacobian  $\bar J_k$ estimates satisfying Assumption~\ref{ass:error}
     \If{$\|\bar c_k\| \leq \epsilon_o$}
      \State  Set $\bar v_k = 0$ and $\bar \tau_k = \bar \tau_{k-1}$ \label{line:vzero}
      \State Compute  $(\bar u_{k}, \bar\rho_k, \bar r_{k})$ via \eqref{eq.nolicq.sub.p2.linear.system.res} satisfying Termination Test~\ref{tt1} and set $\bar d_k = \bar u_k$
      \If{$\Delta \lbar(x_{k},\bar \tau_{k}, \bar d_{k})   \le  \epsilon_o  $}
 \State{\textbf{Terminate} and \textbf{return} $x_k$ \textit{(approximate first-order stationary point)}} \label{line:terminate} 
  \EndIf
       \Else
    \If{$\|\bar J_k^T \bar c_k\|= 0$} 
    \State{\textbf{Terminate} and \textbf{return} $x_k$ \textit{(approximate infeasible stationary point)}} \label{line.early.stopping}
    \EndIf
    \State Compute 
    $\bar v_{k}$ via \eqref{eq.nolicq.tr.sub.p1} satisfying \eqref{eq.nolicq.tr.cauchy.decrease}
    \State Compute  $(\bar u_{k}, \bar\rho_k, \bar r_{k})$ via \eqref{eq.nolicq.sub.p2.linear.system.res} satisfying Termination Test~\ref{tt2} and set $\bar d_k = \bar v_k + \bar u_k$
    \State Update $\bar\tau_k$ by \eqref{eq.tau.update.tt2} and \eqref{eq.tau.trial.tt2}
      \EndIf
        \State Choose step size $\bar\alpha_k$ following one of the two options:
          \State \hspace{0.75cm}
          \textbf{Option I:} Update $\bar \alpha_k$ via \eqref{eq.dominate.condition}--\eqref{eq.project}.
          \State \hspace{0.75cm}
          \textbf{Option II:} Set $\bar\alpha_k = \alpha_u^L$ and update $ \bar\alpha_k = \nu \bar\alpha_{k} $ until~\eqref{eq.modified.linesearch.1d} 
         is satisfied   
      \State Set $x_{k+1} \gets x_k + \bar\alpha_k \bar{d}_k$ \label{line.x.update}
    \EndFor
  \end{algorithmic}
\end{algorithm}

\begin{remark}
    We make a few remarks about Algorithm~\ref{alg.dfo_sqp_LS}.
    \begin{itemize}
        \item \textbf{Step size schemes}: Our proposed algorithm has two options for selecting the step size. Algorithm~\ref{alg.dfo_sqp_LS} with \textbf{Option I} (adaptive) is similar to \cite{berahas2023stochastic} and Algorithm~\ref{alg.dfo_sqp_LS} with \textbf{Option II} (line search) is similar to \cite{oztoprak2023constrained}. The main difference between these two schemes is that the adaptive scheme (\textbf{Option I}) requires estimates of the Lipschitz constants while the line search scheme (\textbf{Option II}) requires objective function and constraint evaluations and knowledge of $(\epsilon_f, \epsilon_g, \epsilon_c, \epsilon_J)$ (or at least estimates of these quantities).  
        \item \textbf{Comparison to \cite{berahas2023stochastic,curtis2024stochastic}}: 
        The algorithms proposed in \cite{berahas2023stochastic,curtis2024stochastic} are developed for solving optimization problems with stochastic objective functions and deterministic equality constaints. We consider general bounded noise in both the objective and constraint functions. Moreover, another differentiating factor is our optimistic early termination scheme. All these algorithms employ a step decomposition strategy to compute the step, but the normal and tangential component computations are not exactly the same. Due to potential rank-deficiency in the constraint Jacobians, similar to \cite{berahas2023stochastic} an inexact normal step is computed by solving a trust-region subproblem whose goal is to minimize linearized infeasibility. Then, similar to \cite{curtis2024stochastic} an inexact tangential step is computed with appropriately modified termination conditions to account for optimistic feasibility. 
        \item \textbf{Comparison to \cite{oztoprak2023constrained}}: 
        The algorithm and analysis in \cite{oztoprak2023constrained} is restricted to the setting in which the LICQ holds. We consider the setting in which the constraint Jacobians may be rank-deficient. Other differences include: $(i)$ the Hessian approximation in~\cite{oztoprak2023constrained} is assumed to be a scalar multiple of the identity matrix, while we allow for more general Hessian approximations; $(ii)$ we propose two step size selection schemes whereas \cite{oztoprak2023constrained} only proposes a line search strategy; $(iii)$ exact solutions to subproblems are considered in \cite{oztoprak2023constrained} whereas we allow for inexact solutions; and, $(iv)$ we adopt an optimistic approach to feasibility and develop early stopping mechanisms.
        \item \textbf{Comparison to \cite{sun2024trust}}: While our proposed algorithm addresses the same class of problems as \cite{sun2024trust}, it differs in the following aspects: $(i)$ \cite{sun2024trust} employs a trust-region approach whereas we use alternative step-size schemes; $(ii)$ our approach allows for inexact solutions and includes an early stopping mechanism, whereas such features are not present
        in \cite{sun2024trust}; and, $(iii)$ in the analysis presented in \cite{sun2024trust} a key assumption regarding the perturbation error of the orthogonal basis of the Jacobian matrix is made \cite[Assumption 5]{sun2024trust}, whereas our analysis does not impose such an assumption. 
    \end{itemize}
\end{remark}

\section{Convergence Analysis}
\label{sec.convergence}

In this section, we provide convergence guarantees for Algorithm~\ref{alg.dfo_sqp_LS}. Specifically, if the LICQ is satisfied, the algorithm generates a sequence of iterates whose averaged first-order stationary and feasibility measures converge to a sufficiently small value dictated by the noise. When the LICQ is not satisfied, there exists a subsequence of the iterates generated by the algorithm whose infeasible stationary measure converges to a sufficiently small value dictated by the noise. Furthermore, we show that the averaged stationary error is sufficiently small if the merit parameter is bounded away from zero. The convergence results are consistent with the results that can be established for the deterministic SQP method~\cite{CurtNoceWach10}, and when the noise levels are all zero we recover such results. To prove our main results (Section~\ref{sec:convergence.main}), we first establish the well-posedness of Algorithm~\ref{alg.dfo_sqp_LS} (Section~\ref{sec:convergence.wellposedness}). Next, we conduct perturbation (error) analysis for subproblems \eqref{eq.nolicq.tr.sub.p1} and \eqref{eq.nolicq.sub.p2.linear.system.res}, and demonstrate that the solutions to these subproblems remain stable under noise perturbations (Section~\ref{sec:convergence.Perturbation}), and also show the deviation of the lower bound of merit parameter values is proportional to the noise. We begin with three assumptions in Section~\ref{sec.assum}.

\subsection{Assumptions}\label{sec.assum}

In this subsection, we present three key assumptions made throughout Section~\ref{sec.convergence}. 
Assumption~\ref{ass.inexact.solver.accuracy} is related to the iterative solver used to compute an inexact tangential step, Assumption~\ref{ass.d.nonzero} is related to the iterates generated by Algorithm~\ref{alg.dfo_sqp_LS}, and Assumption~\ref{ass.noearlytermination} is related to the early termination feature of Algorithm~\ref{alg.dfo_sqp_LS}.

\begin{assumption}
\label{ass.inexact.solver.accuracy}
For all $k \in \mathbb{N}$, the iterative linear system solver employed in \eqref{eq.nolicq.sub.p2.linear.system.res}  generates a sequence of inexact solutions $\left\{\left(\bar u_{k, t}, \bar y_{k, t}, \bar\rho_{k, t},\bar r_{k, t}\right)\right\}_{t \in \mathbb{N}}$ satisfying
\begin{align*}
\left[\begin{array}{c}
\bar\rho_{k, t} \\
\bar r_{k, t}
\end{array}\right]=\left[\begin{array}{cc}
H_{k} &\bar J_{k}^{T} \\
\bar J_{k} & 0
\end{array}\right]\left[\begin{array}{c}
\bar u_{k, t} \\
\bar y_{k, t}
\end{array}\right]+\left[\begin{array}{c}
\bar g_{k}+H_{k}\bar v_{k} \\
0
\end{array}\right] \text { for all } t \in \mathbb{N}
\end{align*}
such that $\lim _{t \rightarrow \infty}\left\|\left(\bar u_{k, t}, \bar\rho_{k, t},\bar r_{k, t}\right)-\left(\bar u_{k}^*,  0,0\right)\right\|=0$, where 
$\bar u_{k}^*$ is the solution of \eqref{eq.nolicq.sub.p2.linear.system.star}. 
\end{assumption}

\begin{assumption}
For all $k \in \mathbb{N}$, $\bar{J}_k^T\bar c_k \neq 0$ or $\bar g_k  \notin \operatorname{Range}(\bar J_k^T)$.
 \label{ass.d.nonzero}
\end{assumption}

\begin{remark}
Assumption~\ref{ass.inexact.solver.accuracy} is satisfied if one employs a linear system solver, such as MINRES~\cite{choi2011minres}, given that $\bar u_{k}^*$ is unique (even though $\bar J_k$ could be rank-deficient).
Assumption~\ref{ass.d.nonzero} is necessary for the well-posedness of the proposed termination tests 
(see Lemma~\ref{lemma.algorithm.finite.termination}). Assumption~\ref{ass.d.nonzero} is violated only when $\bar{J}_k^T\bar c_k = 0$ and $\bar g_k$ lies precisely in $\operatorname{Range}(\bar J_k^T)$. When the algorithm reaches an iterate for which $\bar{J}_k^T\bar c_k = 0$, it terminates if $\|\bar c_k\| > \epsilon_o$ with an approximate infeasible stationary point.
\end{remark}

\begin{assumption}
\label{ass.noearlytermination}
Algorithm~\ref{alg.dfo_sqp_LS} does not terminate finitely, i.e., Algorithm~\ref{alg.dfo_sqp_LS} does not terminate on Lines~\ref{line:terminate} or~\ref{line.early.stopping} for any $k\in \mathbb{N}$.
\end{assumption}

\begin{remark}
Theorems~\ref{theorem.tau.lb},~\ref{theorem.tau.tozero} and~\ref{theorem.tau.lb.complexity} all rely on Assumption~\ref{ass.noearlytermination}. However, in general, Assumption~\ref{ass.noearlytermination} is not guaranteed to be satisfied, i.e., the algorithm may terminate finitely. To handle the case where Assumption~\ref{ass.noearlytermination} is violated, we will show in Theorem~\ref{theorem.early.termination} (resp., Theorem~\ref{theorem.early.termination.infeasible}) that if Algorithm~\ref{alg.dfo_sqp_LS} terminates on Line~\ref{line:terminate} (resp., Line~\ref{line.early.stopping}), a near stationary solution (resp., a near infeasible stationary solution) to \eqref{prob.main} is attained.
\end{remark}

\subsection{Well-posedness of Algorithm}
\label{sec:convergence.wellposedness}

In this section, we show that Algorithm~\ref{alg.dfo_sqp_LS} is well-posed, i.e., 
every iteration of the algorithm is well-defined and is guaranteed to terminate finitely. 
We first show that the termination tests for tangential components $\{ \bar u_k \}$ are well-posed and that at least one of them is satisfied in finite time. 
Then, we show that the sequences $\{\bar \chi_k\}$, $\{\bar \zeta_k\}$, and $\{\bar \xi_k\}$ defined in the adaptive step size scheme (\textbf{Option I}, \eqref{eq.dominate.condition}--\eqref{eq.inexact.xi.trial.update}) are bounded to ensure that the scheme is well-posed. Moreover, we show the boundedness of the merit parameter sequence $\{\bar \tau_k\}$ in the presence of the LICQ. Finally, we show that both the adaptive step size (\textbf{Option I})  and the line search (\textbf{Option II}) schemes are well-defined in that the step sizes are guaranteed to be upper bounded on all iterates 
and lower bounded under the condition that the merit parameter is lower bounded or that both the infeasibility stationarity and the feasibility errors are sufficiently large.

We present a lemma that pertains to the normal component and that is common in the constrained optimization literature; see e.g., \cite[Lemmas 3.5 \& 3.6]{CurtNoceWach10} and \cite[Lemma 1]{berahas2023stochastic}.  
\begin{lemma}
 \label{lemma.v.bound}
 Suppose that Assumptions~\ref{ass:prob} and \ref{ass.noearlytermination} hold. If $\|\bar c_k\|>
\epsilon_o$, there exist $\sigma_v \in \mathbb{R}_{>0}$ and $\underline{\sigma}_{Jc} \in \mathbb{R}_{>0}$ such that 
 \begin{align*}
     \left\|\bar c_k\right\|\left(\left\|\bar c_k\right\|-\left\|\bar c_k+\bar J_k\bar  v_k\right\|\right) \geq \sigma_v\left\|\bar J_k^T \bar c_k\right\|^2, \quad \text{ and } \quad
     \underline{\sigma}_{Jc} \left\|\bar J_k^T \bar c_k\right\|^2 \leq\left\|\bar v_k\right\| \leq \sigma_{Jc}\left\|\bar J_k^T\bar  c_k\right\|.  
 \end{align*}
\end{lemma}
\begin{proof}
The proof follows the same logic \cite[Lemmas 3.5 \& 3.6]{CurtNoceWach10} by replacing $(c_k, J_k, v_k)$ with their noisy evaluations $(\bar c_k,\bar J_k, \bar v_k)$. The result only holds when $\|\bar c_k\|> \epsilon_o$ since otherwise we do not solve subproblem \eqref{eq.nolicq.tr.sub.p1} and directly set $\bar v_k = 0$. 
\end{proof}

The following lemma guarantees the well-posedness of Termination Tests~\ref{tt1} and \ref{tt2} when computing the tangential step, i.e., a solver will compute $\bar u_k$  which satisfies our proposed termination tests within a finite number of iterations. Note that the MINRES method \cite{choi2011minres} guarantees an exact solution to the linear system within a uniformly bounded number of iterations.


\begin{lemma}
Suppose that Assumptions~\ref{ass:error}, \ref{ass:H}, \ref{ass.inexact.solver.accuracy}, \ref{ass.d.nonzero}, and \ref{ass.noearlytermination} hold. 
For all $k \in \mathbb{N}$, the iterative linear system solver computes $\left(\bar u_{k}, \bar\rho_k, \bar r_{k}\right) \leftarrow (\bar{u}_{k,t},\bar\rho_{k,t},\bar{r}_{k,t})$  satisfying either Termination Test~\ref{tt1} or Termination Test~\ref{tt2} in a finite number of iterations $t \in \N{}$. 
\label{lemma.algorithm.finite.termination}
\end{lemma}
\begin{proof}
By \eqref{eq.nolicq.tr.sub.p1} it follows that $\|\vbar_k\| \leq \sigma_{Jc}\|\bar{J}_k^T\bar{c}_k\|$ and by \eqref{eq.nolicq.sub.p2.linear.system.star} it follows that $H_k \bar u_k^* + \bar J_k^T \bar y_k^* + \bar g_k + H_k \bar v_k = 0$.  If $\|\bar u_k^*\| = \| \bar J_k^T \bar c_k\| = 0$, then $\vbar_k = \bar J_k^T \bar c_k = 0$ and $\bar J_k^T \bar y_k^* + \bar g_k = 0$, which violates Assumption~\ref{ass.d.nonzero}. Hence, $\max \{\|\bar u_k^*\|, \| \bar J_k^T \bar c_k\| \} > 0$. By Assumption~\ref{ass.inexact.solver.accuracy} it follows that \eqref{eq.tt1.cond1} is satisfied by $(\bar{u}_k,\bar\rho_k,\bar{r}_k) = (\bar{u}_{k,t},\bar\rho_{k,t},\bar{r}_{k,t})$ for sufficiently large $t \in \mathbb{N}$. We will first show that if $\|\bar{c}_k\| > \epsilon_o$ Termination Test~\ref{tt2} (\eqref{eq.tt2.cond2} and \eqref{eq.inexact.delta.l.desired.opt.case2} or \eqref{eq.tt2.cond1}) will eventually be satisfied, and then if $\|\bar{c}_k\|\leq\epsilon_o$ Termination Test~\ref{tt1} (\eqref{eq.tt1.cond1}--\eqref{eq.inexact.delta.l.desired.opt.case1}) will eventually be satisfied.

When $\left\|\bar c_k\right\|>\epsilon_o$, if $\|\bar J_k^T \bar c_k \|= 0$, Algorithm~\ref{alg.dfo_sqp_LS} terminates on Line~\ref{line.early.stopping}. Hence, we consider $\|\bar J_k^T \bar c_k \| \neq 0$ ($\|\bar v_k\| > 0$ by \eqref{eq.nolicq.tr.cauchy.decrease}) when computing $\bar u_k$. To satisfy condition \eqref{eq.tt2.cond2}, we consider two cases: $(i)$ if $\left\|\bar u_{k}^*\right\|=0$, then Assumption~\ref{ass.inexact.solver.accuracy} implies that $\left\{\left\|\bar u_{k, t}\right\|\right\}_{t\in\mathbb{N}} \rightarrow\left\|\bar u_{k}^*\right\|=0$  so that the former condition in  \eqref{eq.tt2.cond2} holds with  $\left(\bar{u}_k,\bar\rho_k, \bar r_k\right) \leftarrow \left(\bar{u}_{k,t},\bar\rho_{k, t},\bar r_{k, t}\right)$ for all sufficiently large $t \in \mathbb{N}$; $(ii)$ if $\left\|\bar u_{k}^*\right\|>0$, it follows by Assumption~\ref{ass:H} and \eqref{eq.nolicq.sub.p2.linear.system.star} that 
\begin{align*}
\bar u_{k}^{*T}\left(\bar g_k+H_k \bar v_k\right)+ \max \left \{ \tfrac{1}{2}, 1 - \|\bar J_k^T \bar c_k\| \right \} \bar u_{k}^{*T} H_k \bar u_{k}^* &<\bar u_{k}^{*T}\left(\bar g_k+H_k \bar v_k + H_k \bar u_{k}^*\right) \\
    &=-\bar u_{k}^{*T} \bar J_k^T \bar y_{k}^*=0.
\end{align*}
Since $\|\bar v_k\| > 0$, the latter condition in \eqref{eq.tt2.cond2} holds by Assumption~\ref{ass.inexact.solver.accuracy} ($(\bar{u}_k,\bar\rho_k,\bar{r}_k)\leftarrow (\bar{u}_{k,t},\bar\rho_{k,t},\bar{r}_{k,t})$ for all sufficiently large $t \in \mathbb{N}$). Moreover, by $\|\bar v_k\| > 0$ and \eqref{eq.nolicq.tr.cauchy.decrease} it follows that $\left\|\bar c_{k}\right\|-\left\|\bar c_{k}+\bar J_{k} \bar v_{k}\right\| > 0$, and by Assumption~\ref{ass.inexact.solver.accuracy} and $\sigma_r\in(\sigma_c,1)\subset
(0,1)$, 
\begin{align*}
    \lim_{t \to \infty} \left( \left\|\bar c_{k}\right\|-\left\|\bar c_{k}+\bar J_{k} \bar v_{k}+\bar r_{k,t}\right\| \right) = \left\|\bar c_{k}\right\|-\left\|\bar c_{k}+\bar J_{k} \bar v_{k}\right\| > \sigma_r  (\left\|\bar c_{k}\right\|-\left\|\bar c_{k}+\bar J_{k} \bar v_{k}\right\| ).
\end{align*}
Condition \eqref{eq.tt2.cond1} holds with $(\bar{u}_k,\bar\rho_k,\bar{r}_k)\leftarrow (\bar{u}_{k,t},\bar\rho_{k,t},\bar{r}_{k,t})$ for all sufficiently large $t \in \mathbb{N}$. Thus, Termination Test~\ref{tt2} is guaranteed to be satisfied in a finite number of iterations.

When $\left\|\bar c_k\right\|\le \epsilon_o$, then $\left\|\bar v_k\right\|=0$ (Line \ref{line:vzero}, Algorithm~\ref{alg.dfo_sqp_LS}). 
We consider two cases, $\epsilon_o > 0$ and $\epsilon_o = 0$. When $\epsilon_o > 0$, by Assumption~\ref{ass:H} and~\eqref{eq.nolicq.sub.p2.linear.system.star}, it follows that $\bar u_{k}^{*T} H_k \bar u_{k}^* \geq \zeta\|\bar u_{k}^*\|^2 \geq \lambda_u \|\bar u_{k}^*\|^2 > \lambda_u \|\bar u_{k}^*\|^2 - \epsilon_o$ and $\bar g_k^T\bar u_{k}^{*} +\tfrac{1}{2} \bar u_{k}^{*T} H_k \bar u_{k}^* = \bar u_{k}^{*T}\left(\bar g_k+H_k \bar v_k\right)+\tfrac{1}{2} \bar u_{k}^{*T} H_k \bar u_{k}^* = -\tfrac{1}{2} \bar u_{k}^{*T} H_k \bar u_{k}^* < \epsilon_o$. By Assumption~\ref{ass:H},~\eqref{eq.nolicq.sub.p2.linear.system.star},~\eqref{def.merit_model_reduction}, $\sigma_u\in(0,1)$ and $\bar{d}_k^* = \bar{u}_k^* + \bar{v}_k = \bar{u}_k^*$, 
\begin{align*}
    \Delta \lbar\left(x_k, \bar \tau_{k-1},  \bar d_{k}^*\right)  = \Delta \lbar\left(x_k, \bar \tau_{k-1},  \bar u_{k}^*\right)     & = -\bar\tau_{k-1}\bar g_k^T\bar u_{k}^* + \|\bar c_k\|-\|\bar c_k + \bar J_k \bar{u}_k^*\|   \\
    & =   -\bar\tau_{k-1}\left(-H_k \bar u_{k}^* -J_k^T \bar y_{k}^*\right)^T\bar u_{k}^* \\
    & = \bar\tau_{k-1}\bar u_{k}^{*T} H_k\bar u_{k}^* \\
    & > \bar \tau_{k-1} \sigma_u  \max \left\{\bar u_{k}^{*T} H_k \bar u_{k}^*, \lambda_u\left\|\bar u_{k}^*\right\|^2\right\} - \epsilon_o.
\end{align*}
By Assumption~\ref{ass.inexact.solver.accuracy}, 
for sufficiently large $t \in \N{}$, the inexact solution $\bar{u}_k \leftarrow \bar{u}_{k,t}$ satisfies~\eqref{eq.tt1.cond2}--\eqref{eq.inexact.delta.l.desired.opt.case1} when $\epsilon_o > 0$, and Termination Test~\ref{tt1} is guaranteed to be satisfied in finite iterations.

On the other hand, when $\|\bar{c}_k\| \leq \epsilon_o = 0$ by Assumption~\ref{ass.d.nonzero}, it follows that $\|\bar{c}_k\| = \|\bar{v}_k\| = 0$ and $\|\bar{u}_k^*\| > 0$. Therefore, by Assumption~\ref{ass:H},~\eqref{eq.nolicq.sub.p2.linear.system.star} and $\lambda_u\in(0,\zeta)$, it follows that $\bar u_{k}^{*T} H_k \bar u_{k}^* \geq \zeta\|\bar u_{k}^*\|^2 > \lambda_u \|\bar u_{k}^*\|^2$ and $\bar g_k^T\bar u_{k}^{*} +\tfrac{1}{2} \bar u_{k}^{*T} H_k \bar u_{k}^* = \bar u_{k}^{*T}\left(\bar g_k+H_k \bar v_k\right)+\tfrac{1}{2} \bar u_{k}^{*T} H_k \bar u_{k}^* = -\tfrac{1}{2} \bar u_{k}^{*T} H_k \bar u_{k}^* < 0 = \epsilon_o$. Additionally, by 
\eqref{def.merit_model_reduction}, $\sigma_u\in(0,1)$ and $\bar{d}_k^* = \bar{u}_k^* + \bar{v}_k = \bar{u}_k^*$, using similar logic as the former case, it follows that $\Delta \lbar\left(x_k, \bar \tau_{k-1},  \bar d_{k}^*\right) > \bar \tau_{k-1} \sigma_u  \max \left\{\bar u_{k}^{*T} H_k \bar u_{k}^*, \lambda_u\left\|\bar u_{k}^*\right\|^2\right\}$, since $\bar u_{k}^* \ne 0$. 
Thus, by Assumption~\ref{ass.inexact.solver.accuracy}, 
 for sufficiently large $t \in \N{}$, the inexact solution $\bar{u}_k \leftarrow \bar{u}_{k,t}$ satisfies~\eqref{eq.tt1.cond2}--\eqref{eq.inexact.delta.l.desired.opt.case1} when $\epsilon_o = 0$, and Termination Test~\ref{tt1} is guaranteed to be satisfied in finite iterations. 
\end{proof}

Next, we prove that the search direction $\bar d_k$ is non-zero at every iteration if the algorithm does not terminate early. This ensures that the step size schemes are 
well-posed.

\begin{lemma}
\label{lemma.d.nonzero.opt}
Suppose that Assumptions~\ref{ass:error}, \ref{ass:H}, \ref{ass.inexact.solver.accuracy}, \ref{ass.d.nonzero}, and \ref{ass.noearlytermination} hold. 
For all $k \in \mathbb{N}$, $\|\bar d_k \| > 0$.  
\end{lemma}
\begin{proof}
We consider two cases $\|\bar c_k\| \le \epsilon_o$ and $\|\bar c_k\| > \epsilon_o$ and prove the result in both cases using a contradiction argument. 

Suppose that $\|\bar d_k \| = 0$. If $\|\bar c_k\| \le \epsilon_o$ and $\epsilon_o = 0$, then it follows that $\|\bar c_k\| = \|\bar v_k\| = 0$, which implies that $\|\bar u_k\| = 0$. By \eqref{eq.tt1.cond1}, $ \|\bar\rho_k \| = \|\bar{r}_k \| = 0$ and consequently $\bar u_k^* =  \bar u_k = 0$, which violates Assumption~\ref{ass.d.nonzero}. If $\|\bar c_k\| \le \epsilon_o$ and $\epsilon_o > 0$, then it follows by \eqref{def.merit_model_reduction} that $\Delta \lbar\left(x_k, \bar \tau_{k},  \bar d_{k} \right) = 0 < \epsilon_o$. Algorithm~\ref{alg.dfo_sqp_LS} then terminates at line~\ref{line:terminate}, which causes contradiction.

Suppose that $\|\bar d_k \| = 0$. If $\|\bar c_k\| >  \epsilon_o$, by Assumption~\ref{ass.noearlytermination} the algorithm does not terminate on Line~\ref{line.early.stopping}, and thus $\|\bar c_k\| >  \epsilon_o$. This implies that $\|\bar v_k  \| > 0$ by~\eqref{eq.nolicq.tr.cauchy.decrease}, and so the right-hand-side of \eqref{eq.inexact.delta.l.desired.opt.case2} is positive. Since $\Delta \lbar\left(x_k, \bar \tau_{k},  \bar d_{k} \right) = 0$ by \eqref{def.merit_model_reduction},  \eqref{eq.inexact.delta.l.desired.opt.case2} cannot be satisfied. On the other hand, $\| \bar c_k \| - \| \bar c_k + \bar J_k \bar v_k + \bar r_k\| = \| \bar c_k \| - \| \bar c_k + \bar J_k \bar d_k\|  = 0$ and \eqref{eq.tt2.cond1}  cannot be satisfied. Thus, Termination Test~\ref{tt2} cannot be satisfied and we have reached a contradiction to Lemma~\ref{lemma.algorithm.finite.termination}.

In both cases, $\|\bar d_k \| = 0$ leads to a contradiction, and thus the desired result follows.    
\end{proof}
 
The following lemma characterizes the relationship between the reduction in the model of the merit function and the search direction under various conditions when early termination is not triggered, i.e., Algorithm~\ref{alg.dfo_sqp_LS} does not terminate on Lines~\ref{line:terminate} or \ref{line.early.stopping}. 
Similar results have been established in \cite[Lemma 3.4]{berahas2021sequential} and \cite[Lemma 8] {curtis2024stochastic}. It is worth noting that in our setting, a small value of $\Delta \lbar(x_{k}, \bar \tau_{k}, \bar d_{k})$ does not necessarily imply that the norm of the search direction is also small, which is different from~\cite{berahas2021sequential,curtis2024stochastic}. This distinction arises due to the potential rank deficiency in the 
noisy constraint Jacobians. 
We distinguish between two cases because $\bar v_k = 0$ is set when $\|\bar c_k\| \leq \epsilon_o$, resulting in no reduction in constraint violation in this scenario. Moreover, if the singular values of $\{\bar J_k\}$ are bounded below 
and $\|\bar c_k\| > \epsilon_o$, an alternative result can be derived by replacing the infeasible stationarity measure with the feasibility measure.

\begin{lemma}
\label{lemma.d.bar.bounded.by.Delta.l.opt}
Suppose that Assumptions~\ref{ass:prob}, \ref{ass:error}, \ref{ass:H}, \ref{ass.inexact.solver.accuracy}, \ref{ass.d.nonzero}, and~\ref{ass.noearlytermination} hold. 
If  
$\|\bar c_k\| \le \epsilon_o$, 
  \begin{align}
  \label{eq.Deltal.lb.LICQ.opt}
      \Delta \lbar(x_{k},\bar \tau_{k}, \bar d_{k})  \ge  \tfrac{\bar\tau_{k} \sigma_{u} \lambda_u}{2}  \|\bar u_k  \|^2= \tfrac{\bar\tau_{k} \sigma_{u} \lambda_u}{2}  \|\bar d_k  \|^2.   
  \end{align}
 If $\|\bar c_k\| > \epsilon_o$, it follows that
       \begin{align}
    \label{eq.Deltal.lb.NOLICQ}
      \Delta \lbar(x_{k},\bar \tau_{k}, \bar d_{k})  \ge  \bar\tau_{k}  \sigma_{u}\lambda_u\left\|\bar u_{k}\right\|^{2} + \tfrac{\sigma_{c}\sigma_v}{\kappa_c + \epsilon_c} \|\bar  J_k^T \bar  c_k \|^2 \ge \tfrac{\min\{\bar\tau_{k}  \sigma_{u}\lambda_u, \sigma_{c}\sigma_v/ (\kappa_c + \epsilon_c)\}}{\max \left\{2,2 \sigma_{Jc}^2 \right\}} \|\bar d_k  \|^2.
  \end{align}
Moreover, if  $\|\bar c_k\| > \epsilon_o$ and 
  the singular values of $\bar J_k$ are bounded below by $\kappa_{\sigma}$ for all $k \in \mathbb{N}$, 
    \begin{align}
    \label{eq.Deltal.lb.LICQ}
      \Delta \lbar(x_{k},\bar \tau_{k}, \bar d_{k})  \ge 
  \bar\tau_{k} \sigma_{u}  \lambda_u\left\|\bar u_{k}\right\|^{2} +  \kappa_{\sigma}^2 \sigma_{c}\sigma_v \| \bar  c_k \| \ge \tfrac{\min\{\bar\tau_{k} \sigma_{u}  \lambda_u, \kappa_{\sigma}^2 \sigma_{c}\sigma_v  \}}{\max \left\{2,2 \sigma_{Jc}^2 (\kappa_c+\epsilon_c)(\kappa_J+\epsilon_J)^2 \right\}} \|\bar d_k  \|^2 .  
  \end{align}
\end{lemma}
\begin{proof}
If $\|\bar c_k \| \le \epsilon_o$, then Termination Test~\ref{tt1} is satisfied, and
 \begin{equation}
 \label{eq.detal.lb1}
     \Delta \lbar\left(x_{k}, \bar\tau_{k},\bar d_k \right) \geq \bar\tau_{k} \sigma_u \max \left\{\bar u_{k}^{T} H_{k}\bar u_k, \lambda_u\left\|\bar u_{k}\right\|^{2}\right\}
     - \epsilon_o 
 \end{equation}
by~\eqref{eq.inexact.delta.l.desired.opt.case1} and $(\bar{d}_k,\bar\tau_k) = (\bar{u}_k,\bar\tau_{k-1})$. Since 
 $ \Delta \lbar(x_{k},\bar \tau_{k}, \bar d_{k}) > \epsilon_o$ by Assumption~\ref{ass.noearlytermination}, it follows that
\begin{align*}
    \Delta \lbar\left(x_{k}, \bar\tau_{k},\bar d_k \right) 
 &\geq  \bar\tau_{k} \sigma_{u} \lambda_u  \|\bar u_k  \|^2  - \epsilon_o >  \bar\tau_{k} \sigma_{u} \lambda_u  \|\bar u_k  \|^2  - \Delta \lbar\left(x_{k}, \bar\tau_{k},\bar d_k \right).  
\end{align*}
Re-arranging the above inequality and using the fact that 
$\bar v_k = 0$ we conclude~\eqref{eq.Deltal.lb.LICQ.opt}.

If $\|\bar c_k \| > \epsilon_o$, when~\eqref{eq.inexact.delta.l.desired.opt.case2} holds, we directly have $\bar\tau_k = \bar\tau_{k-1}$ and 
\begin{align*}
\Delta \lbar\left(x_{k}, \bar\tau_{k},\bar d_k \right) \geq \bar\tau_k \sigma_u \max \left\{\bar u_{k}^{T} H_{k}\bar u_{k}, \lambda_u\left\|\bar u_{k}\right\|^{2}\right\}   + \sigma_c ( \left\|\bar c_{k}\right\|-\left\|\bar c_{k}+\bar J_{k} \bar v_{k}\right\| ).
\end{align*}
Otherwise, when~\eqref{eq.inexact.delta.l.desired.opt.case2} does not hold, by Termination Test~\ref{tt2},~\eqref{eq.tt2.cond1} holds. By~\eqref{eq.tau.update.tt2} and~\eqref{eq.tau.trial.tt2}, it follows that $\bar\tau_k \leq (1-\sigma_{\tau})\bar\tau_{k}^{\text {trial}}$ and
\begin{align}\label{eq1}
    \bar\tau_{k}\left(\bar{g}_k^T\bar{d}_k + \max\{\bar{u}_k^TH_k\bar{u}_k,\lambda_u\|\bar{u}_k\|^2\}\right) \leq (1-\sigma_{\tau})(1-\tfrac{\sigma_c}{\sigma_r})(\|\bar{c}_k\| - \|\bar{c}_k + \bar{J}_k\bar{v}_k + \bar{r}_k\|).
\end{align}
By~\eqref{eq.nolicq.sub.p2.linear.system.res},~\eqref{def.merit_model_reduction},~\eqref{eq.tt2.cond1} and~\eqref{eq1}, it follows that
\begin{align}
 \Delta \lbar\left(x_{k}, \bar\tau_{k},\bar d_k \right) &= -\bar\tau_k \bar{g}_k^T\bar{d}_k + \|\bar{c}_k\| - \|\bar{c}_k + \bar{J}_k\bar{d}_k\| \notag \\ 
 &= -\bar\tau_k \bar{g}_k^T\bar{d}_k + \|\bar{c}_k\| - \|\bar{c}_k + \bar{J}_k\bar{v}_k + \bar{r}_k\| \notag \\
   &\ge  \bar\tau_{k} \max \left\{\bar u_{k}^{T} H_{k}\bar u_{k}, \lambda_u\left\|\bar u_{k}\right\|^{2}\right\}   + \tfrac{\sigma_c}{\sigma_r}(\|\bar{c}_k\| - \|\bar{c}_k + \bar{J}_k\bar{v}_k + \bar{r}_k\|) \notag \\
   &\ge \bar\tau_{k} \max \left\{\bar u_{k}^{T} H_{k}\bar u_{k}, \lambda_u\left\|\bar u_{k}\right\|^{2}\right\}   + \sigma_c( \left\|\bar c_{k}\right\|-\left\|\bar c_{k}+\bar J_{k} \bar v_{k}\right\| ). \label{eq.detal.lb2}
\end{align}
Combining with \eqref{eq.inexact.delta.l.desired.opt.case2}, we further know that when $\|\bar{c}_k\| > \epsilon_o$,
\begin{align}
\label{eq.detal.lb3}
 \Delta \lbar\left(x_{k}, \bar\tau_{k},\bar d_k \right) 
\ge\bar\tau_{k} \sigma_{u}   \max \left\{\bar u_{k}^{T} H_{k}\bar u_{k}, \lambda_u\left\|\bar u_{k}\right\|^{2}\right\}   + \sigma_c ( \left\|\bar c_{k}\right\|-\left\|\bar c_{k}+\bar J_{k} \bar v_{k}\right\| ). 
\end{align}
By the step update rule, Lemma~\ref{lemma.v.bound}, Assumptions~\ref{ass:prob} and~\ref{ass:error}, it follows that 
 \begin{align}
 \left\|\bar  d_k\right\|^2 
  \leq 2\left(\left\|\bar u_k\right\|^2+\left\|\bar v_k\right\|^2\right)  
 &\leq 2\left(\left\|\bar u_k\right\|^2+ \sigma_{Jc}^2\left\|\bar  J_k^T \bar  c_k \right\|^2\right) \label{eq.u^2.Jc^2}. 
 \end{align}
Furthermore, by \eqref{eq.detal.lb3}, \eqref{eq.u^2.Jc^2} and Lemma~
\ref{lemma.v.bound}, 
\begin{align*}
      \Delta \lbar\left(x_{k}, \bar\tau_{k},\bar d_k \right) 
 &\geq  \bar\tau_k \sigma_u \lambda_u\left\|\bar u_{k}\right\|^{2} + \sigma_{c}\sigma_v\tfrac{ \|\bar  J_k^T \bar  c_k \|^2 }{ \|\bar c_k  \|} \\
  &\geq  \bar\tau_k \sigma_u \lambda_u\left\|\bar u_{k}\right\|^{2} + \tfrac{\sigma_{c}\sigma_v 
 }{\kappa_c + \epsilon_c} \|\bar  J_k^T \bar  c_k \|^2,  
\end{align*}
we conclude 
\eqref{eq.Deltal.lb.NOLICQ}. 

If the singular values of $\{\bar J_k\}$ are bounded below by $\kappa_{\sigma} > 0$ for all $k \in \mathbb{N}$, it follows from the singular value decomposition of $\bar  J_k^T$ that 
\begin{align}
\label{eq.relation.c.Jc}
    \|\bar  J_k^T \bar  c_k \| \ge    \kappa_{\sigma}  \|\bar  c_k \|, 
\end{align}
and by Lemma~\ref{lemma.v.bound}, Assumptions~\ref{ass:prob} and~\ref{ass:error}
 \begin{align}
    \left\|\bar  d_k\right\|^2  
 &\leq 2\left(\left\|\bar u_k\right\|^2+ \sigma_{Jc}^2\left\|\bar  J_k^T \bar  c_k \right\|^2\right) 
  \leq 2 \left(\left\|\bar  u_k\right\|^2+ \sigma_{Jc}^2 (\kappa_c+\epsilon_c)(\kappa_J+\epsilon_J)^2\left\|\bar  c_k\right\|\right). \label{eq.u^2.c}
 \end{align}
If $\|\bar c_k\| > \epsilon_o$, by 
\eqref{eq.detal.lb3}, \eqref{eq.relation.c.Jc}, \eqref{eq.u^2.c} and Lemma~\ref{lemma.v.bound}, it follows that
\begin{align*}
 \Delta \lbar\left(x_{k}, \bar\tau_{k},\bar d_k \right) &\geq \bar\tau_k \sigma_u \max \left\{\bar u_{k}^{T} H_{k}\bar u_k, \lambda_u\left\|\bar u_{k}\right\|^{2}\right\}+\sigma_{c}\left(\left\|\bar c_{k}\right\|-\left\|\bar c_{k}+\bar J_{k} \bar v_{k}\right\|\right) \\
   &\geq  \bar\tau_k \sigma_u \lambda_u\left\|\bar u_{k}\right\|^{2} + \tfrac{\sigma_{c}\sigma_v \|\bar  J_k^T \bar  c_k \|^2 }{ \|\bar c_k  \|} \geq  \bar\tau_k \sigma_u \lambda_u\left\|\bar u_{k}\right\|^{2} +   \kappa_{\sigma}^2 \sigma_{c}\sigma_v \| \bar  c_k \|,
\end{align*}
from which we conclude 
\eqref{eq.Deltal.lb.LICQ}.
\end{proof}

Next, we introduce two lemmas related to the properties of the sequences $\{ \bar \chi_k,\bar \zeta_k \}$ 
(Lemma~\ref{lemma.chi.bound}) and $\{ \bar \xi_{k}\}$ (Lemma~\ref{lemma.xi}).


\begin{lemma}
\label{lemma.chi.bound}
Suppose that Assumptions~\ref{ass:H} and~\ref{ass.noearlytermination} hold. There exists $\left(\chi_{\max}, \zeta_{\min}\right) \in \mathbb{R}_{>0} \times \mathbb{R}_{>0}$ such that for some $K_{\chi,\zeta} \in \mathbb{N}$, $(\bar \chi_k, \bar \zeta_k) = (\bar \chi_{K_{\chi,\zeta}}, \bar\zeta_{K_{\chi,\zeta}})$ for all $k \in \mathbb{N}$
with $k \ge K_{\chi,\zeta}$, where $(\bar \chi_{K_{\chi,\zeta}}, \bar\zeta_{K_{\chi,\zeta}}) \in (0, \chi_{\max} ] \times [ \zeta_{\min}, \infty)$.
\end{lemma}
\begin{proof}
The proof proceeds in the same manner as 
\cite[Lemma 2]{berahas2023stochastic}. 
\end{proof}

\begin{lemma}
\label{lemma.xi}
Suppose Assumptions~\ref{ass:prob}, \ref{ass:error}, \ref{ass:H}, \ref{ass.inexact.solver.accuracy}, \ref{ass.d.nonzero}, and~\ref{ass.noearlytermination}  
hold. The sequence $\{ \bar \xi_{k}\}$ is non-increasing. In addition, there exists $\bar \xi_{\min} \in \mathbb{R}_{>0}$ such that $\bar \xi_{k} \ge \bar\xi_{\min}$ for all $k\in\mathbb{N}$. 
\end{lemma}
\begin{proof}
By 
\eqref{eq.inexact.xi.update} 
either $\bar \xi_{k}  =  \bar \xi_{k-1} $ or $\bar \xi_{k}  =  \min \{ (1- \sigma_{\xi})\bar \xi_{k-1} , \bar \xi_k^{\text{trial}} \} < \bar  \xi_{k-1}$. Therefore, $\{ \bar \xi_{k}\}$ is non-increasing. 

If $\|\bar c_k\| > \epsilon_o$, by \eqref{eq.Deltal.lb.NOLICQ}, 
\begin{align}
\label{eq.Delta.ge.max}
 \Delta \lbar\left(x_{k}, \bar\tau_{k},\bar d_k \right) &\geq  \max \left \{   \bar\tau_{k} \sigma_u\lambda_u\left\|\bar u_{k}\right\|^{2}  ,\sigma_c \sigma_v (\kappa_c+\epsilon_c)^{-1}\left\|\bar J_k^T \bar c_k\right\|^2 \right \}. 
\end{align}
We consider two cases. When $\left\|\bar u_k\right\|^2 \geq \bar \chi_k\left\|\bar v_k\right\|^2$, by \eqref{eq.inexact.xi.update} and \eqref{eq.inexact.xi.trial.update}, 
either $\bar\xi_k=\bar\xi_{k-1}$ or 
\begin{align*}
   \bar \xi_k \geq \left(1-\sigma_{\xi}\right)\tfrac{\Delta \lbar\left(x_k, \bar\tau_k, \bar d_k\right)}{\bar \tau_k\left\|\bar d_k\right\|^2} \geq \left(1 - 
 \sigma_{\xi}\right) \tfrac{ \bar\tau_k \sigma_u \lambda_u \left\|\bar u_k\right\|^2 }{2\bar\tau_k\left(1+\bar \chi_k^{-1}\right)\left\|\bar u_k\right\|^2}  = \tfrac{ \left(1-\sigma_{\xi}\right)\sigma_u \lambda_u}{2\left(1+\chi_{-1}^{-1}\right)}>0. 
\end{align*}
On the other hand, when $\left\|\bar u_k\right\|^2<\bar \chi_k\left\|\bar v_k\right\|^2$, by \eqref{eq.nolicq.tr.sub.p1},  \eqref{eq.inexact.xi.update}, \eqref{eq.inexact.xi.trial.update}, \eqref{eq.Deltal.lb.NOLICQ}, and Lemma~\ref{lemma.chi.bound} it follows that either $\bar \xi_k=\bar \xi_{k-1}$ or
\begin{align}
\label{eq.u_k.le.Chi.v_k}
    \bar \xi_k \geq 
    \left(1-\sigma_{\xi}\right)\tfrac{\Delta \lbar\left(x_k, \bar\tau_k, \bar d_k\right)}{\left\|\bar d_k\right\|^2}  \geq\left(1-\sigma_{\xi}\right) \tfrac{\sigma_c \sigma_v \left\|\bar J_k^T \bar c_k\right\|^2}{2\left(\bar \chi_k+1\right) (\kappa_c + \epsilon_c) \sigma_{Jc}^2\left\|\bar J_k^T \bar c_k\right\|^2} \geq \tfrac{\left(1-\sigma_{\xi}\right)\sigma_c \sigma_v }{2\left(\chi_{\max }+1\right) (\kappa_c + \epsilon_c) \sigma_{Jc}^2}>0.
\end{align}

If $\|\bar c_k\| \le \epsilon_o$, it follows that $\left\|\bar u_k\right\|^2 \geq \bar \chi_k\left\|\bar v_k\right\|^2$, and by \eqref{eq.inexact.xi.update}, \eqref{eq.inexact.xi.trial.update}, and \eqref{eq.Deltal.lb.LICQ.opt} 
\begin{align*}
   \bar \xi_k \geq \left(1-\sigma_{\xi}\right)\tfrac{\Delta \lbar\left(x_k, \bar\tau_k, \bar d_k\right)}{\bar \tau_k\left\|\bar d_k\right\|^2} \geq\left(1-\sigma_{\xi}\right) \tfrac{\bar\tau_k \sigma_u \lambda_u \left\|\bar d_k\right\|^2 }{2\bar\tau_k\left\|\bar d_k\right\|^2} 
   = \tfrac{\left(1-\sigma_{\xi}\right) \sigma_u \lambda_u}{2}>0.
\end{align*}

Setting $\bar \xi_{\min} = \min \left \{ \tfrac{\left(1-\sigma_{\xi}\right)\sigma_u\lambda_u}{ \max \left \{2\left(1+\chi_{-1}^{-1}\right) , 2 \right\} },  \tfrac{\left(1-\sigma_{\xi}\right)\sigma_c \sigma_v}{2\left(\chi_{\max }+1\right) (\kappa_c + \epsilon_c) \sigma_{Jc}^2} \right\}$ completes the proof. 
\end{proof}

The next lemma pertains to the boundedness of the tangential and full search directions.

\begin{lemma}
\label{lemma.dbar.ubar.bound}
Suppose that Assumptions~\ref{ass:prob}, \ref{ass:error}, \ref{ass:H}, \ref{ass.inexact.solver.accuracy}, \ref{ass.d.nonzero}, and  \ref{ass.noearlytermination} hold. Suppose that the minimum non-zero singular values of $\bar J_k$ are bounded below by $\kappa_{\sigma}\in \mathbb{R}_{>0} $ for all $k \in \mathbb{N}$. Then, for all $k\in \mathbb{N}$, there exists $\{ \kappa_{\bar u},  \kappa_{\bar d}\} \subset \mathbb{R}_{>0}$ such that $\| \bar u_{k}\| \le \kappa_{\bar u}$ and $\| \bar d_{k}\| \le \kappa_{\bar d}$. 
\end{lemma}
\begin{proof}
Suppose that $\text{rank}(\bar J_k ) = p$. First, let us discuss the case when $p = n \le m$. It follows by the rank-nullity theorem \cite{horn2012matrix} that $\text{nullity}(\bar J_k) = n - p= 0$, i.e., the null space of $\bar J_k$ is trivial (contains only the zero vector). Since $\bar J_k$ has full column rank, the matrix $\bar J_k^T \bar J_k$ is symmetric and positive definite, ensuring that it is invertible. It then follows by \eqref{eq.nolicq.sub.p2.linear.system.res}, \eqref{eq.tt1.cond1}, Assumptions~\ref{ass:prob}, \ref{ass:error}, and \ref{ass:H}, and the condition of the minimum non-zero singular values of $\bar J_k$ being bounded below by $\kappa_{\sigma} \in \R{}_{>0}$ that $\|\bar u_k\| = \|(\bar J_k^T \bar J_k)^{-1} \bar J_k^T \bar r_k\| \le \|(\bar J_k^T \bar J_k)^{-1}\| \|  \bar J_k \| \| \bar r_k \| \le  ( \kappa_J + \epsilon_J ) \kappa_{\sigma}^{-2} \lambda_{\rho r} \kappa_{\rho r} $.
Now, suppose that $ p = n \leq m $ does not hold. Since $ p \leq \min\{m, n\} $, we consider the following cases: $(i)$ $ p < n \leq m $ and $(ii)$ $ p \leq m < n $. In both cases, it follows that $ p < n $. The null space of $\bar J_k$ is therefore nontrivial. Suppose that $\bar \Sigma_{k}$ is a diagonal matrix containing the non-zero singular values of $\bar J_k $, where the full SVD decomposition of $\bar J_k $ is given as 
\begin{align}
\label{eq.J.bar.SVD}
 \bar J_k^T &= (\bar U_k \text{ } \bar Z_k) \left(\begin{array}{cc}\bar \Sigma_{k} & 0 \\ 0 & 0\end{array}\right) \left(\begin{array}{c} \bar V_k  \\ \bar  V_k^{\perp} \end{array}\right)   = \bar U_k \bar \Sigma_{k}\bar V_k, 
\end{align}
where $\bar U_k \in \mathbb{R}^{n\times p}$, $\bar Z_k \in \mathbb{R}^{n\times (n-p)}$, $\bar \Sigma_k \in \mathbb{R}^{p\times p} $, $\bar V_k \in \mathbb{R}^{p\times m}$,  $\bar V_k^{\perp} \in \mathbb{R}^{(m-p) \times m}$.  
Let $\tilde{\bar J}_k^T = \bar  U_k \bar \Sigma_{k} \in \mathbb{R}^{n\times p}$. 
We have 
$\text{Range}(\tilde{\bar J}_k^T) = \text{Range}({\bar J}_k^T)$, $\|\tilde{\bar J}_k \tilde{\bar J}_k^T \|=\| {\bar J}_k {\bar J}_k^T \|$. Since $\bar Z_k$ represents the orthonormal basis of null space of $\bar J_k $, by the linear system \eqref{eq.nolicq.sub.p2.linear.system.res} 
\begin{align} 
\label{eq.u.bar.expression}
\bar u_{k} &= 
\tilde{\bar J}_k^T \left(\tilde{\bar J}_k \tilde{\bar J}_k^{T}\right)^{-1}\bar V_k \bar r_{k} 
-\bar Z_{k}\left(\bar Z_{k}^{T} H_{k} \bar Z_{k}\right)^{-1} \bar Z_{k}^{T}\left(\bar g_{k}+H_{k}\bar v_{k}+ H_{k}\tilde{\bar J}_k^T \left(\tilde{\bar J}_k \tilde{\bar J}_k^{T}\right)^{-1} \bar V_k \bar r_{k} -\bar \rho_{k}\right).  
\end{align}
Since $\bar U_k$ and $\bar V_k$ are both submatrices of orthogonal matrices, $\|\bar U_k\| \le 1$, $\|\bar V_k\| \le 1$. By $\|(\tilde{\bar J}_k \tilde{\bar J}_k^{T})^{-1}\| =  \| \bar \Sigma_k^{-2}\|  \le \kappa_{\sigma}^{-2}$, Assumptions~\ref{ass:prob}, \ref{ass:error}, \ref{ass:H} and \eqref{eq.tt1.cond1} it follows that 
\begin{align*}
    &\   \| \bar u_{k}\| \\ \le &\ \left \| \tilde{\bar J}_k^T \left(\tilde{\bar J}_k \tilde{\bar J}_k^{T}\right)^{-1}\bar V_k \bar r_{k} \right \| +      \left \|  \bar Z_{k}\left(\bar Z_{k}^{T} H_{k} \bar Z_{k}\right)^{-1} \bar Z_{k}^{T}\left(\bar g_{k}+H_{k}\bar v_{k}+ H_{k}\tilde{\bar J}_k^T \left(\tilde{\bar J}_k \tilde{\bar J}_k^{T}\right)^{-1} \bar V_k \bar r_{k} -\bar \rho_{k}\right)   \right \| \\ 
    \le &\ ( \kappa_J + \epsilon_J ) \kappa_{\sigma}^{-2} \|\bar r_k\| + \zeta^{-1} \left( \kappa_g +\epsilon_g + \kappa_H\sigma_{Jc}    (\kappa_J + \epsilon_J)(\kappa_c+ \epsilon_c) + \kappa_H  (\kappa_J + \epsilon_J ) \kappa_{\sigma}^{-2}  \|\bar r_k\| + \|\bar \rho_k \|  \right) \\ 
    \le &\ \left( (1 + \zeta^{-1}  \kappa_H )  (\kappa_J + \epsilon_J ) \kappa_{\sigma}^{-2} 
     + \zeta^{-1}   \right) \lambda_{\rho r} \kappa_{\rho r}  + \zeta^{-1} \left( \kappa_g +\epsilon_g +  \kappa_H  \sigma_{Jc} (\kappa_J + \epsilon_J)(\kappa_c+ \epsilon_c) \right):= \kappa_{\bar u}.   
\end{align*}
Since $\kappa_{\bar u} > ( \kappa_J + \epsilon_J ) \kappa_{\sigma}^{-2} \lambda_{\rho r} \kappa_{\rho r} $, the bound derived when $\text{rank}(\bar J_k ) = p = n$, we have $\|\bar u_k\| \le  \kappa_{\bar u}$ for all $k \in \N{}$. 
Given that $\|\bar v_k\| \le \sigma_{Jc} \|\bar J_k^T \bar c_k\| \le \sigma_{Jc} (\kappa_J + \epsilon_J)(\kappa_c + \epsilon_c) $ and $\|\bar d_k\| \le \|\bar v_k\| + \|\bar u_k \|$, there exists $ \kappa_{\bar d}\in \mathbb{R}_{>0}$ such that $\|\bar d_k\|  \le \kappa_{\bar d}$. 
\end{proof}

Next, we prove that $\{\bar \tau_k\}$ is bounded away from zero when the singular values of $\bar J_k$ are bounded below by $\kappa_{\sigma} \in \mathbb{R}_{>0}$ for all $k \in \mathbb{N}$. The step sizes for both schemes (\textbf{Option I} and \textbf{Option II}) are also lower bounded in this case (Lemmas~\ref{lemma.inexact.stepsize.opt} and \ref{lemma.step.size.ls}).

\begin{lemma}
\label{lemma.tau.lb}
Suppose that Assumptions~\ref{ass:prob}, \ref{ass:error}, \ref{ass:H}, \ref{ass.inexact.solver.accuracy}, \ref{ass.d.nonzero}, and~\ref{ass.noearlytermination} hold. 
The merit parameter sequence  $\{\bar \tau_k\}$ is 
non-increasing. 
 If the singular values of $\bar J_k$ are bounded below by $\kappa_{\sigma} \in \mathbb{R}_{>0}$ for all $k \in \mathbb{N}$, 
 then there exists $\tilde \tau_{\min} \in \mathbb{R}_{>0}$ such that $\bar \tau_k \ge \tilde \tau_{\min}$ for all $k\in\mathbb{N}$. 
\end{lemma}
\begin{proof}
When Termination Test~\ref{tt1} is triggered, $\bar \tau_k = \bar \tau_{k-1}$. When Termination Test~\ref{tt2} is triggered, either $\bar \tau_k = \bar \tau_{k-1}$ or $\bar \tau_k = (1-\sigma_{\tau}) \bar \tau_k^{\text{trial}} < \bar \tau_{k-1}$. Therefore, $\{\bar \tau_k\}$ is a non-increasing sequence.  
By Algorithm~\ref{alg.dfo_sqp_LS}, 
\eqref{eq.tau.update.tt2}, and \eqref{eq.tau.trial.tt2}, it follows that $\bar \tau_k < \bar \tau_{k-1}$ only when  $\|\bar c_k \| > \epsilon_o$ (Termination Test~\ref{tt2}), $\bar g_{k}^{T}\bar d_{k}+\max \left\{\bar u_{k}^{T} H_{k}\bar u_{k}, \lambda_u\left\|\bar u_{k}\right\|^{2}\right\}>0$, and 
\begin{align*}
 \bar \tau_{k-1}  > 
     \tfrac{(1- \sigma_{\tau}) \left( 1-\tfrac{\sigma_{c}}{\sigma_r} \right) \left(\left\|\bar c_{k}\right\|-\|\bar c_{k}+\bar J_{k} \bar v_{k}+\bar r_{k}\|\right) }{\bar g_{k}^{T}\bar d_{k}+\max \left\{\bar u_{k}^{T} H_{k}\bar u_{k}, \lambda_u\left\|\bar u_{k}\right\|^{2}\right\} }. 
 \end{align*}
By Assumptions~\ref{ass:prob}, \ref{ass:error}, and~\ref{ass:H}, \eqref{eq.nolicq.tr.sub.p1}, and Lemma~\ref{lemma.dbar.ubar.bound}, considering the two cases in~\eqref{eq.tt2.cond2}, 
\begin{align*}
    \bar g_{k}^{T}\bar d_{k}+\max \left\{\bar u_{k}^{T} H_{k}\bar u_{k}, \lambda_u\left\|\bar u_{k}\right\|^{2}\right\} 
    = \ &\bar g_{k}^{T}\bar u_{k}+   \bar g_{k}^{T}\bar v_{k}+ \max \left\{\bar u_{k}^{T} H_{k}\bar u_{k}, \lambda_u\left\|\bar u_{k}\right\|^{2}\right\} \\  
    \le \ &(\kappa_g+ \epsilon_g) (\lambda_{uv} + 1) \|\bar v_k\| + \max\{\kappa_H, \lambda_u \} \|\bar u_k\|^2 \\  
    \le \ &\left( (\kappa_g+ \epsilon_g) (\lambda_{uv} + 1) + \kappa_{\bar u} \max\{\kappa_H, \lambda_u \} \lambda_{uv} \right) \sigma_{Jc} \|\bar J_k^T \bar c_k \|, 
\end{align*}
or 
\begin{align*}
    \bar g_{k}^{T}\bar d_{k}+\max \left\{\bar u_{k}^{T} H_{k}\bar u_{k}, \lambda_u\left\|\bar u_{k}\right\|^{2}\right\}  =  &\   \bar g_{k}^{T}\bar u_{k}+   \bar g_{k}^{T}\bar v_{k}+  \bar u_{k}^{T} H_{k}\bar u_{k} \\ 
    =&\  (\bar g_{k} +  H_k\bar v_k)^T \bar u_k + (1 - \|\bar J_k^T\bar c_k\|) \bar u_{k}^{T} H_{k}\bar u_{k} + \bar g_{k}^{T}\bar v_{k} \\
    & \qquad + \|\bar J_k^T\bar c_k\| \bar u_{k}^{T} H_{k}\bar u_{k} -   \bar v_k^T H_k \bar u_k \\ 
    \le &\  \lambda_u \|\bar v_{k}\| + (\kappa_g+ \epsilon_g) \|\bar v_{k}\| + \kappa_{\bar u}^2 \kappa_H \|\bar J_k^T\bar c_k\|  + \kappa_{\bar u} \kappa_H \|\bar v_{k}\|  \\  
    \le &\ \left( \left( \lambda_u +  \kappa_g+ \epsilon_g + \kappa_{\bar u} \kappa_H  \right) \sigma_{Jc} + \kappa_{\bar u}^2 \kappa_H  \right) \|\bar J_k^T \bar c_k \|. 
\end{align*}
Thus, there exists a constant $\kappa_D \in \mathbb{R}_{>0}$ such that $ \bar g_{k}^{T}\bar d_{k}+\max \left\{\bar u_{k}^{T} H_{k}\bar u_{k}, \lambda_u\left\|\bar u_{k}\right\|^{2}\right\}   \le \kappa_D  \|\bar J_k^T \bar c_k \|$. 
Hence, if the merit parameter $\bar\tau_k$ is ever updated, it follows by \eqref{eq.tt2.cond1},  \eqref{eq.relation.c.Jc}, and Lemma~\ref{lemma.v.bound} that 
\begin{equation}
\label{eq.tau.tilde}
\begin{split}
      \bar \tau_{k-1}  &>  \tfrac{ (1- \sigma_{\tau}) \left( 1-\tfrac{\sigma_{c}}{\sigma_r} \right)\left(\left\|\bar c_{k}\right\|-\|\bar c_{k}+\bar J_{k} \bar v_{k}+\bar r_{k}\| \right) }{\bar g_{k}^{T}\bar d_{k}+\max \left\{\bar u_{k}^{T} H_{k}\bar u_{k}, \lambda_u\left\|\bar u_{k}\right\|^{2}\right\} }\\
    &\ge \tfrac{(1- \sigma_{\tau})\left(\sigma_{r}-\sigma_{c}\right)\left(\left\|\bar c_{k}\right\|-\left\|\bar c_{k}+\bar J_{k} \bar v_{k}\right\|\right)}{\bar g_{k}^{T}\bar d_{k}+\max \left\{\bar u_{k}^{T} H_{k}\bar u_{k}, \lambda_u\left\|\bar u_{k}\right\|^{2}\right\} } \\
    &\ge \tfrac{(1- \sigma_{\tau})\left(\sigma_{r}-\sigma_{c}\right) \sigma_v \| \bar J_k^T \bar c_k\|^2  / \|\bar c_k\|}{ \kappa_D  \|\bar J_k^T \bar c_k \|}  \\
    &\ge  \tfrac{(1- \sigma_{\tau})\left(\sigma_{r}-\sigma_{c}\right) \sigma_v \kappa_{\sigma} }{ \kappa_D} := \tilde \tau_{\min} > 0.  
\end{split}    
\end{equation}
The desired result is therefore satisfied with the defined $\tilde \tau_{\min}$. 
\end{proof}

Lemma~\ref{lemma.tau.lb} shows that $\bar \tau_k$ is uniformly bounded away from zero when $\bar J_k$ is of full rank. 
In general, there are two possible outcomes for the sequence $\{\bar \tau_k\}$: 
$(i)$ there exists  $\bar \tau_{\min}\in \mathbb{R}_{>0}$ such that $\bar \tau_k \ge \bar \tau_{\min}$ for all $k\in\mathbb{N}$; and $(ii)$ $\bar\tau_k \rightarrow 0$. 
In Section~\ref{sec:convergence.main}, we derive convergence guarantees and complexity results for the different outcomes. The remainder of this subsection is about properties of step sizes.

Lemmas~\ref{lemma.inexact.stepsize.opt} and \ref{lemma.step.size.ls} provide lower and upper bounds for step sizes generated by both the adaptive and the line search schemes. The lower bounds of the step size are established when $\bar \tau_k \ge \bar \tau_{\min}$. However, when $\{\bar \tau_k\}$ diminishes to 0, the step size may converge to zero, which poses challenges for the convergence of the algorithm. We prove that the step size does not vanish under certain conditions. For this purpose, we define
\begin{align}\label{eq.gammas}
    \mathcal{X}_\gamma:=\left\{x \in \mathbb{R}^n:\left\|\bar J(x)^T 
\bar c(x)\right\| \geq \gamma, \|\bar c (x)\| >\epsilon_o \right\}
\end{align}
for some $\gamma \in \mathbb{R}_{>0}$ 
and show that the step sizes are uniformly bounded away from zero 
when $x_k \in \mathcal{X}_\gamma$ even when $\{\bar \tau_k\}$ diminishes to zero. We first show that the conclusions above hold for the adaptive step size scheme (\textbf{Option I}).

\begin{lemma}
Suppose that Assumptions~\ref{ass:prob}, \ref{ass:error}, \ref{ass:H}, \ref{ass.inexact.solver.accuracy}, \ref{ass.d.nonzero}, and  \ref{ass.noearlytermination} hold, 
and that \textbf{Option I} is used to select the step size in Algorithm~\ref{alg.dfo_sqp_LS}.
In addition, suppose that $\beta_k = \beta \in (0,1)$ satisfies $\tfrac{2(1-\eta) \beta \bar{\xi}_{-1} \max \left\{\bar{\tau}_{-1}, 1\right\}}{\bar{\tau}_{-1} \Gamma+L} \in (0,1]$ for all $k \in \mathbb{N}$. 
 For any $k\in\mathbb{N}$, it holds that $0<\bar \alpha_{k}^{\min} \leq \bar \alpha_{k}^{\text{suff}}$, $\bar \alpha_k \le \bar \alpha_k^{\text{suff}}\le 1$
  , and 
\begin{align}
\label{def.alpha.uA}
0<\bar\alpha_{k}^{\min } \leq \bar \alpha_{k}^{\max } \le \left(   \tfrac{2(1-\eta)\bar \xi_{-1} \max \{\bar\tau_{-1},1 \}}{\Gamma}  + \theta \right) \beta := \alpha_u^A.
\end{align}
\begin{itemize}
    \item 
    If $\bar \tau_k \ge \bar \tau_{\min}$ for all $k\in\mathbb{N}$, it follows that $\bar\alpha_{k}^{\min } \geq \tfrac{2(1-\eta) \beta\bar \xi_{\min} \min \{\bar\tau_{\min},1 \}}{\bar\tau_{\min} L+\Gamma} := \alpha_{l\tau}^{A}$.
    \item 
    If $\{\bar \tau_k\}$ diminishes to 0, for $\gamma \in \mathbb{R}_{>0}$, there exists $\alpha_{l0}^{A} \in \mathbb{R}_{>0}$ such that if $x_k \in \mathcal{X}_\gamma$, then $\bar \alpha_k \geq \alpha_{l0}^{A}$.
\end{itemize}
\label{lemma.inexact.stepsize.opt}
\end{lemma}
\begin{proof}
Since $\bar \xi_{k} \le \bar \xi_{k}^{\text{trial}}$ by 
\eqref{eq.inexact.xi.update}, by the update rule of $ \bar\alpha_{k}^{\min}$~\eqref{eq.inexact.alpha.min} it follows that 
\begin{align*}
    \bar\alpha_{k}^{\min } \le \tfrac{2(1-\eta) \beta_{k} \Delta \lbar(x_{k},\bar \tau_{k}, \bar d_{k})   }{\left(\bar\tau_{k} L+\Gamma\right)\|\bar d_{k}\|^{2}}. 
\end{align*}
In addition, by Lemmas~\ref{lemma.xi} and~\ref{lemma.tau.lb}, \eqref{eq.inexact.alpha.min} and the prescribed condition for $\{\beta_k\}$ 
\begin{align*}
 \bar\alpha_{k}^{\min } \le \tfrac{2(1-\eta) \beta_k \bar{\xi}_{-1} \max \left\{\bar{\tau}_{-1}, 1\right\}}{\bar{\tau}_{-1} \Gamma+L} \le 1.
\end{align*}
By~\eqref{eq.inexact.alpha.suff}, $ \bar\alpha_{k}^{\min } \le \min \left\{ \tfrac{2(1-\eta) \beta_{k} \Delta \lbar(x_{k},\bar \tau_{k}, \bar d_{k})   }{\left(\bar\tau_{k} L+\Gamma\right)\|\bar d_{k}\|^{2}}, 1\right\} = \bar\alpha_{k}^{\text {suff}}$, and hence, the projection operator will never increase the step size, i.e., 
$\bar \alpha_k \le \bar \alpha_k^{\text{suff}} \le 1$. Moreover, by $\beta_k = \beta \in(0,1)$ and the definitions of $\bar\alpha_k^{\min}$ and $\bar\alpha_k^{\max}$, it follows that $\bar  \alpha_{k}^{\max } =  \bar \alpha_{k}^{\min } + \theta \beta \le  \left(   \tfrac{2(1-\eta)\bar \xi_{-1} \max \{\bar\tau_{-1},1 \}}{\Gamma}  + \theta \right) \beta$.

If $\bar \tau_k \ge \bar \tau_{\min}$, it follows by \eqref{eq.inexact.alpha.min} and Lemma~\ref{lemma.tau.lb} that 
\begin{align*}
 \tfrac{2(1-\eta) \beta\bar \xi_{\min} \min \{\bar\tau_{\min},1 \}}{\bar\tau_{\min} L+\Gamma} \le \tfrac{2(1-\eta) \beta\bar \xi_k \min \{\bar\tau_k,1 \}}{\bar\tau_k L+\Gamma} \le  \bar \alpha_{k}^{\min}. 
\end{align*}
The desired conclusion is therefore satisfied. 

If $\{\bar \tau_k\}$ is vanishes to zero, 
define $\mathcal{K}_\gamma:=\left\{k \in \mathbb{N}: x_k \in \mathcal{X}_\gamma\right\}$. It follows by Lemma \ref{lemma.v.bound} that $\left\|\bar v_k\right\| \geq \underline{\sigma}_{Jc}\left\|\bar J_k^T \bar c_k\right\|^2 \geq \underline{\sigma}_{Jc} \gamma^2$ for all $k \in \mathcal{K}_\gamma$. Furthermore, $\left\|\bar u_k\right\| \leq \tfrac{\kappa_{\bar u}}{\underline{\sigma}_{Jc} \gamma^2}\left\|\bar v_k\right\| $ for all $k \in \mathcal{K}_\gamma$, by Lemma~\ref{lemma.dbar.ubar.bound}. 
If $\left\|\bar u_k\right\|^2<\bar \chi_k\left\|\bar v_k\right\|^2$, then $\bar \alpha_k \geq \tfrac{2(1-\eta) \beta \bar \xi_{\min} }{\left(\bar\tau_{-1} L+\Gamma\right)}$ by~\eqref{eq.inexact.alpha.min} and Lemma~\ref{lemma.xi}. If $\left\|\bar u_k\right\|^2 \geq \bar \chi_k\left\|\bar v_k\right\|^2$, by \eqref{eq.project}, \eqref{eq.Deltal.lb.NOLICQ},  Lemma~\ref{lemma.v.bound}, 
and $\bar\alpha_{k}^{\min } \le \bar\alpha_k^{\text{suff}}$,
\begin{align}
    \bar\alpha_k& = \min \left \{ \bar\alpha_k^{\text{suff}},  \bar\alpha_k^{\text{max}}   \right \} \notag \\
    & =\min \left\{\tfrac{2(1-\eta) \beta \Delta \lbar\left(x_k, \bar\tau_k, \bar d_k\right)}{\left(\bar\tau_k L+\Gamma\right)\left\|\bar d_k\right\|^2}, 1, \tfrac{2(1-\eta) \beta\bar \xi_k \bar\tau_k}{\bar\tau_k L+\Gamma}+\theta \beta \right\} \notag \\
    & \geq \min \left\{\tfrac{(1-\eta)\beta \sigma_c \sigma_v\left\|\bar J_k^T\bar  c_k\right\|^2}{\left(\bar \tau_k L+\Gamma\right) (\kappa_c + \epsilon_c) \left(\left\|\bar u_k\right\|^2+\left\|\bar v_k\right\|^2\right)}, 1, \tfrac{2(1-\eta) \beta\bar  \xi_k \bar \tau_k}{\bar \tau_k L+\Gamma}+\theta \beta \right\} \notag \\& \geq \min \left\{\tfrac{(1-\eta)\beta \sigma_c \sigma_v\left\|\bar J_k^T\bar  c_k\right\|^2}{\left(\bar \tau_k L+\Gamma\right) (\kappa_c + \epsilon_c) \left(\kappa_{\bar u}^2/(\underline{\sigma}_{Jc}^2 \gamma^4)+1\right) \sigma_{Jc}^2\left\|\bar J_k^T \bar c_k\right\|^2}, 1, \tfrac{2(1-\eta) \beta \bar \xi_k \bar \tau_k}{\bar \tau_k L+\Gamma}+\theta \beta \right\} \notag \\
    & \geq \min \left\{\tfrac{(1-\eta) \beta \sigma_c \sigma_v \underline{\sigma}_{Jc}^2 \gamma^4}{(\bar \tau_{-1} L+\Gamma) (\kappa_c + \epsilon_c) \sigma_{Jc}^2\left(\kappa_{\bar u}^2+\underline{\sigma}_{Jc}^2 \gamma^4\right)}, 1, \theta \beta \right\} \label{eq.alpha.min.NOLICQ.u_k.ge.v_k}. 
\end{align}
Let $\alpha_{l0}^{A} := \min \left\{\tfrac{(1-\eta) \beta \sigma_c \sigma_v \underline{\sigma}_{Jc}^2 \gamma^4}{(\bar \tau_{-1} L+\Gamma) (\kappa_c + \epsilon_c) \sigma_{Jc}^2\left(\kappa_{\bar u}^2+\underline{\sigma}_{Jc}^2 \gamma^4\right)}, 1, \theta \beta, \tfrac{2(1-\eta) \beta \bar \xi_{\min}}{\bar\tau_{-1} L+\Gamma}  \right\}$, and the result holds. 
\end{proof}




The next lemma proves that similar results apply to the line search scheme (\textbf{Option II}), albeit with different lower and upper bounds. 

\begin{lemma}
\label{lemma.step.size.ls}
Suppose that Assumptions~\ref{ass:prob}, \ref{ass:error}, \ref{ass:H}, \ref{ass.inexact.solver.accuracy}, \ref{ass.d.nonzero}, and  \ref{ass.noearlytermination} hold, 
and that \textbf{Option II} is used to select the step size in Algorithm~\ref{alg.dfo_sqp_LS} with 
\begin{align}
\label{def.epsilon_Ak}
    \epsilon_{A_k} = 2\bar \tau_k \epsilon_f + 4\epsilon_c +  \bar\tau_{k}  \alpha_u^L\epsilon_g \|\bar d_k\| +   \alpha_u^L\epsilon_J \|\bar d_k\|.
\end{align}
\begin{itemize}
    \item If $\bar \tau_k \ge \bar{\tau}_{\min}$ for all $k
    \in\mathbb{N}$,  there exists $\kappa_{\Delta l} \in \mathbb{R}_{>0}$ (that depends on $ \bar{\tau}_{\min}$) such that 
   $\bar \alpha_k \ge \min \left\{  \alpha_u^L, \tfrac{\nu (1-\eta) \kappa_{\Delta l} }{\bar\tau_{-1} L+\Gamma }    \right \} := \alpha_{l\tau}^{L}$. 
    \item If $\{\bar \tau_k\}$ diminishes to 0, for $\gamma \in \mathbb{R}_{>0}$, there exists $\alpha_{l0}^{L} \in \mathbb{R}_{>0}$ such that if $x_k \in \mathcal{X}_\gamma$, then $\bar \alpha_k \geq \alpha_{l0}^{L}$.

\end{itemize}
\end{lemma}
\begin{proof}
If $\bar \tau_k \ge \bar{\tau}_{\min}$, by Lemma~\ref{lemma.d.bar.bounded.by.Delta.l.opt} 
there exists $\kappa_{\Delta l} \in \mathbb{R}_{>0}$ (that depends on $ \bar{\tau}_{\min}$) 
such that $\Delta \lbar\left(x_k, \bar\tau_k, \bar d_k\right) \ge \kappa_{\Delta l} \|\bar d_k\|^2$. By 
Assumptions~\ref{ass:prob} and \ref{ass:error}, \eqref{eq.merit}, \eqref{def.merit_model_reduction}, 
and $\bar\alpha_k \le \alpha_u^L\leq 1$ 
    \begin{equation}
    \label{eq.model.reduction.c.large}
        \begin{split}
    &\ \phi\left(x_{k}+ \bar \alpha_k \bar d_{k}, \bar \tau_{k}\right)-\phi\left(x_{k}, \bar \tau_{k}\right)\\  
    =&\
    \bar  \tau_{k} f_{k+1}-\bar\tau_{k} f_{k}+\left\| c_{k+1}\right\|-\left\| c_{k}\right\|\\
    \leq  &\  \bar \alpha_k \bar \tau_{k} g_{k}^{T} \bar d_{k}+\left\| c_{k}+\bar \alpha_k  J_{k} \bar d_{k}\right\|-\left\| c_{k}\right\|+\tfrac{1}{2}\left(\bar\tau_{k} L+\Gamma\right) \bar \alpha_k ^{2}\left\|\bar d_{k}\right\|^{2}  \\
    \leq  &\  \bar \alpha_k \bar \tau_{k} g_{k}^{T} \bar d_{k}+\left\|\bar c_{k}+\bar \alpha_k  \bar J_{k} \bar d_{k}\right\|-\left\|\bar c_{k}\right\|+\tfrac{1}{2}\left(\bar\tau_{k} L+\Gamma\right) \bar \alpha_k ^{2}\left\|\bar d_{k}\right\|^{2}   + \| c_k - \bar c_k  \| + \| c_k  - \bar c_k  + \bar\alpha_k (J_k - \bar J_k) \bar d_k \|
     \\
     \leq &\ \bar \alpha_k \bar \tau_{k} \bar  g_{k}^{T} \bar d_{k} +  \bar \alpha_k \bar \tau_{k} ( g_{k} -\bar g_{k} )^{T}  \bar d_{k} - \bar \alpha_k \left\|\bar c_{k}\right\| +  \bar \alpha_k \left\|\bar c_{k}+  \bar J_{k} \bar d_{k}\right\| \\  &\ + \left\|\bar c_{k}+\bar \alpha_k  \bar J_{k} \bar d_{k}\right\|-\left\|\bar c_{k}\right\| +  \bar \alpha_k \left\|\bar c_{k}\right\| -  \bar \alpha_k \left\|\bar c_{k}+  \bar J_{k} \bar d_{k}\right\| + \tfrac{1}{2}\left(\bar\tau_{k} L+\Gamma\right) \bar \alpha_k ^{2}\left\|\bar d_{k}\right\|^{2}  + 2\epsilon_c +  \bar \alpha_k\epsilon_J \|\bar d_k\|  \\ 
     \leq  &\ - \bar \alpha_k \Delta \lbar\left(x_k, \bar\tau_k, \bar d_k\right)  +  \bar \alpha_k \bar \tau_{k} ( g_{k} -\bar g_{k} )^{T}  \bar d_{k} \\  &\ + \left\|\bar c_{k}+\bar \alpha_k  \bar J_{k} \bar d_{k}\right\|-\left\|\bar c_{k}\right\| +  \bar \alpha_k \left\|\bar c_{k}\right\| -  \bar \alpha_k \left\|\bar c_{k}+  \bar J_{k} \bar d_{k}\right\| + \tfrac{1}{2}\left(\bar\tau_{k} L+\Gamma\right) \bar \alpha_k ^{2}\left\|\bar d_{k}\right\|^{2}  + 2\epsilon_c + \bar \alpha_k \epsilon_J \|\bar d_k\|  
     \\
     \le  &  - \bar \alpha_k \Delta \lbar\left(x_k, \bar\tau_k, \bar d_k\right) 
     + \bar \alpha_k \bar \tau_{k} ( g_{k} -\bar g_{k} )^{T}  \bar d_{k}
     + (1-\bar \alpha_k) \|\bar c_k \| + (\bar \alpha_k-1) \|\bar c_k \| \\  &\  + 
     \tfrac{1}{2}\left(\bar\tau_{k} L+\Gamma\right) \bar \alpha_k ^{2}\left\|\bar d_{k}\right\|^{2} + 
     2\epsilon_c +  \bar \alpha_k\epsilon_J \|\bar d_k\|\\ \le  &  - \bar \alpha_k \Delta \lbar\left(x_k, \bar\tau_k, \bar d_k\right) 
     + \bar \alpha_k \bar \tau_{k} ( g_{k} -\bar g_{k} )^{T}  \bar d_{k} +  
     \tfrac{\bar\tau_{-1} L+\Gamma}{2\kappa_{\Delta l} } \bar \alpha_k^{2} 
     \Delta \lbar\left(x_k, \bar\tau_k, \bar d_k\right) + 2\epsilon_c +  \bar \alpha_k\epsilon_J \|\bar d_k\|
     \\ \le  &  - \bar \alpha_k \left( 1-  \tfrac{\bar\tau_{-1} L+\Gamma}{2\kappa_{\Delta l} } \bar \alpha_k 
  \right) \Delta \lbar\left(x_k, \bar\tau_k, \bar d_k\right)+ \bar\tau_{k}  \alpha_u^L\epsilon_g \|\bar d_k\| + 2\epsilon_c + \alpha_u^L \epsilon_J \|\bar d_k\|.
        \end{split}
    \end{equation}
Moreover, \eqref{eq.diff.merit.func} implies that when $\bar \alpha_k \le \tfrac{(1-\eta) \kappa_{\Delta l} }{\bar\tau_{-1} L+\Gamma }$, 
\eqref{eq.modified.linesearch.1d} is satisfied by the definition of $\epsilon_{A_k}$.  By setting $\alpha_{l\tau}^{L} := \min \left\{ \alpha_u^L, \tfrac{\nu (1-\eta) \kappa_{\Delta l} }{\bar\tau_{-1} L+\Gamma }   \right \}$, the desired conclusion then follows.

If $\{\bar \tau_k\}$ vanishes to zero, 
define $\mathcal{K}_\gamma:=\left\{k \in \mathbb{N}: x_k \in \mathcal{X}_\gamma\right\}$. It follows by the same arguments as in the proof of 
Lemma~\ref{lemma.inexact.stepsize.opt} that $\left\|\bar u_k\right\| \leq \tfrac{\kappa_{\bar u}}{\underline{\sigma}_{Jc} \gamma^2}\left\|\bar v_k\right\| $ for all  $ k \in \mathcal{K}_\gamma $ and 
\begin{align*}
    \tfrac{\Delta \lbar\left(x_k, \bar\tau_k, \bar d_k\right)}{\left\|\bar d_k\right\|^2} \ge \tfrac{\sigma_c \sigma_v\left\|\bar J_k^T\bar  c_k\right\|^2}{2(\kappa_c + \epsilon_c) \left(\left\|\bar u_k\right\|^2+\left\|\bar v_k\right\|^2\right)}  \ge \tfrac{ \sigma_c \sigma_v\left\|\bar J_k^T\bar  c_k\right\|^2}{ 2(\kappa_c + \epsilon_c) \left(\kappa_{\bar u}^2/(\underline{\sigma}_{Jc}^2 \gamma^4)+1\right) \sigma_{Jc}^2\left\|\bar J_k^T \bar c_k\right\|^2} \ge \tfrac{\sigma_c \sigma_v \underline{\sigma}_{Jc}^2 \gamma^4}{2(\kappa_c + \epsilon_c) \left(\kappa_{\bar u}^2+\underline{\sigma}_{Jc}^2 \gamma^4\right)\sigma_{Jc}^2}. 
\end{align*}
The desired conclusion then follows by the definition of $\kappa_{\Delta l} (\gamma)$, $\alpha_{l0}^{L}$ and a similar rationale as in the previous case. 
\end{proof}

\begin{remark} Lemmas~\ref{lemma.inexact.stepsize.opt} and \ref{lemma.step.size.ls} provide useful lower and upper bounds of the step size for the adaptive and line search schemes, respectively.  The results in Lemma~\ref{lemma.inexact.stepsize.opt} are similar to those in~\cite[Theorem 2 \& Lemma 12]{berahas2023stochastic}. The relaxation parameter $\epsilon_{A_k} \in \mathbb{R}_{\ge0}$ (defined in \eqref{def.epsilon_Ak}) ensures that the step size that satisfies the line search condition \eqref{eq.modified.linesearch.1d} is bounded away from zero under reasonable conditions, ensuring that the backtracking line search procedure is well defined. Relaxation parameters are common in adaptive algorithms in the noisy setting. 
It is important to note that the relaxation parameters used in \cite{berahas2019derivative,berahas2021global,berahas2025sequential,dezfulian2024convergence,jin2024high,lou2024noise} depend on $\epsilon_f$ for problems with noisy objective functions, and the relaxation parameter in \cite{oztoprak2023constrained} depends on both $\epsilon_f$ and $\epsilon_c$ for problems with both noisy objective and constraint functions. 
Similar relaxation strategies have been developed in the context of trust-region methods; see e.g., \cite{cao2023first,larson2024novel,sun2023trust} for problems with noisy objective functions and \cite{sun2024trust} for problems with both noisy objective and constraint functions. 
The fact that all relaxation parameters in previous works do not depend on $\epsilon_g$ and $\epsilon_J$ can be attributed to the following reasons: $(i)$ stronger assumptions are made on the gradient approximations \cite{berahas2019derivative,dezfulian2024convergence} or $(ii)$ the desired 
stationarity measure is assumed to be sufficiently large compared to $\epsilon_g$ and/or $\epsilon_J$ in the convergence analysis \cite{berahas2021global,berahas2025sequential,jin2024high,lou2024noise,oztoprak2023constrained,cao2023first,larson2024novel,sun2023trust}.   
The additional terms in the relaxation parameter \eqref{def.epsilon_Ak}, $ \bar\tau_{k}  \alpha_u^L\epsilon_g \|\bar d_k\| +   \alpha_u^L\epsilon_J \|\bar d_k\|$, are necessary in order to establish non-asymptotic guarantees  (Theorem~\ref{theorem.tau.lb}) without the need for additional assumptions or restrictions. It is important to note that the additional terms can be replaced by quantities proportional to $ (\bar\tau_{-1} \epsilon_g + \epsilon_J )^2$. Finally, when $\{\bar \tau_k\}$ diminishes to zero, the noise levels in the constraint functions, $\epsilon_c$ and $\epsilon_J$, dominate the relaxation parameter $\epsilon_{A_k}$.
\end{remark}

Table~\ref{table.result.bound.alpha} summarizes the results from Lemmas \ref{lemma.inexact.stepsize.opt} and \ref{lemma.step.size.ls}.

\begin{table}[h]
\centering
\caption{Lower and upper bounds for $\bar \alpha_k$ in Algorithm~\ref{alg.dfo_sqp_LS} (\textbf{Option I} and \textbf{Option II}).}
\label{table.result.bound.alpha}
\resizebox{\columnwidth}{!}{
\begin{tabular}{lccc}
\toprule
 & Upper Bound & Lower Bound ($\bar \tau_k \ge \bar \tau_{\min}$) & Lower Bound ($\bar \tau_k \to 0$ and $x_k \in \mathcal{X}_\gamma$) \\
\midrule
\textbf{Option I} &  $ \left(   \tfrac{2(1-\eta)\bar \xi_{-1} \max \{\bar\tau_{-1},1 \}}{\Gamma}  + \theta \right) \beta$  & $\tfrac{2(1-\eta) \bar \xi_{\min} \min \{\bar\tau_{\min},1 \}}{\bar\tau_{\min} L+\Gamma}\beta$ & $\min \left\{\tfrac{(1-\eta)  \sigma_c \sigma_v \underline{\sigma}_{Jc}^2 \gamma^4}{(\bar \tau_{-1} L+\Gamma) (\kappa_c + \epsilon_c) \sigma_{Jc}^2\left(\kappa_{\bar u}^2+\underline{\sigma}_{Jc}^2 \gamma^4\right)}, \tfrac{1}{\beta}, \theta, \tfrac{2(1-\eta) \bar \xi_{\min}}{\bar\tau_{-1} L+\Gamma}  \right\} \beta $ \\
\textbf{Option II} &  $ \alpha_u^L \in (0,1]$ & $\min \left\{ \alpha_u^L, \tfrac{\nu (1-\eta) \kappa_{\Delta l} }{\bar\tau_{-1} L+\Gamma }   \right \}$ & $\min \left\{ \alpha_u^L, \tfrac{\nu (1-\eta) \kappa_{\Delta l} (\gamma) }{\bar\tau_{-1} L+\Gamma }   \right \}$ \\
\bottomrule
\end{tabular}}
\end{table}

\subsection{Perturbation Analysis}
\label{sec:convergence.Perturbation}

In this section, we present theoretical guarantees that characterize the perturbations 
of the (possibly inexact) solutions to subproblems \eqref{eq.nolicq.tr.sub.p1} and \eqref{eq.nolicq.sub.p2}. Specifically, we establish how the solutions of the noisy subproblems deviate from those of the fictitious (never solved) true subproblems. Moreover, we derive useful upper bounds on the 
deviation of the lower bound of the merit parameter value as compared to the fictitious (never computed) true merit parameter value. The former results provide insight into the stationarity, feasibility, or infeasible stationarity measures of the true problem, and the latter results provide insights into the stability of the algorithm. 

To simplify the notation, 
let $g_k = \nabla f(x_k)$ and $J_k = \nabla c(x_k)^T$, and let $v_k$, 
$u_k^*$ and $\tau_k$ denote the true normal and tangential 
components and true merit parameter (defined explicitly below in \eqref{eq.nolicq.tr.sub.p1.det}, \eqref{eq.nolicq.sub.p2.linear.system.star.det}, and \eqref{def.tau.det}--\eqref{def.tau.trial.det}), respectively. As shown in \cite{berahas2024modified}, when the constraint evaluations are exact (i.e., $\epsilon_c = \epsilon_J = 0$) and the singular values of $J_k$ are bounded away from zero, the difference between the noisy and true search directions ($\bar d_k$ and $d_k$) is upper bounded by a quantity proportional to $\epsilon_g$. Deriving such a result in our noisy setting is not as straightforward. The main reason for this is that $J_k$ or $\bar J_k$ can be rank-deficient and the fact that our algorithm makes use of inexact subproblem solutions. As such, we do not directly prove that $\|\bar{d}_k - d_k\|$ is sufficiently small. Instead, we establish the convergence of our algorithm by deriving upper bounds for $\|\bar{v}_k - v_k\|$ when the LICQ holds and either exact solutions or Cauchy solutions are employed (see Lemmas~\ref{lemma.diff.v} and \ref{lemma.diff.v.exact}) and 
$\|\bar{u}_k - u_k^*\|$ regardless of whether $J_k$ or $\bar{J}_k$ are rank-deficient (see Lemma~\ref{lemma.diff.u}). 
When the LICQ is violated, $\bar{v}_k$ and $v_k$ can be non-unique and the norm difference $\|\bar{v}_k - v_k\|$ can be large. 
The trust region constraint $\|\bar{v}_k\| \leq \sigma_{Jc} \|\bar{J}_k^T \bar{c}_k\|$ is sufficient for algorithmic convergence, as we will demonstrate that $\|\bar{J}_k^T \bar{c}_k\|$ eventually becomes small (see Theorem~\ref{theorem.tau.lb}).

The first lemma bounds the difference between $\|\bar{J}_k^T \bar{c}_k\|$ and $\|{J}_k^T {c}_k\|$.

\begin{lemma}
\label{lemma.mathcalE.jc}
Suppose that Assumptions~\ref{ass:prob} and \ref{ass:error} hold. 
There exists $\mathcal{E}_{Jc}$, $\mathcal{E}_{Jc^2} \in \mathbb{R}_{\ge 0}$
such that, for all $k \in \N{}$
\begin{align*}
\| J_k^T c_k - \bar  J_k^T \bar  c_k\| & \le  \mathcal{E}_{Jc} \quad 
     \text{ and } \quad \big| \| J_k^T c_k \|^2 - \|\bar  J_k^T \bar  c_k \|^2 \mathbbm{1}\{ \|\bar c_k\| >\epsilon_o \} \big| \le \mathcal{E}_{Jc^2}. 
\end{align*}
\end{lemma}
\begin{proof}
By 
Assumptions~\ref{ass:prob} and \ref{ass:error} it follows that 
\begin{align}
\label{eq.diff.Jc}
\| J_k^T c_k - \bar  J_k^T \bar  c_k\| \le \|\bar{J}_k^T(c_k - \bar{c}_k)\| + \|(J_k - \bar{J}_k)^Tc_k\| \le ( (\kappa_J + \epsilon_J) \epsilon_c + \kappa_c \epsilon_J )   :=\mathcal{E}_{Jc}, 
\end{align}
and 
\begin{equation}
\begin{aligned}
\label{eq.diff.Jc.sto.det}
    & \   \big| \| J_k^T c_k \|^2 - \|\bar  J_k^T \bar  c_k \|^2 \mathbbm{1}\{ \|\bar c_k\| >\epsilon_o \} \big|   \\   \le &\
    \max \left \{  \big| \|\bar  J_k^T \bar  c_k \|^2  - \| J_k^T c_k \|^2 \big|,  \| J_k^T c_k \|^2 \mathbbm{1}\{ \|\bar c_k\| \le \epsilon_o \}    \right \}   \\   
    \le &\ \max \left \{ \| J_k^T c_k - \bar  J_k^T \bar  c_k\| \left(  (\kappa_J + \epsilon_J)  (\kappa_c + \epsilon_c) + \kappa_J \kappa_c \right) ,  \kappa_J^2 \epsilon_o^2\right \}   \\   
    \le &\ \max \left \{  ( (\kappa_J + \epsilon_J) \epsilon_c + \kappa_c \epsilon_J ) \left(  (\kappa_J + \epsilon_J)  (\kappa_c + \epsilon_c) + \kappa_J \kappa_c \right),  \kappa_J^2 \epsilon_c^2\right \}:=\mathcal{E}_{Jc^2},
\end{aligned}
\end{equation}
which completes the proof.
\end{proof}
Lemma~\ref{lemma.mathcalE.jc} is useful for establishing upper bounds of the perturbation errors of inexact solutions to normal \eqref{eq.nolicq.tr.sub.p1} and tangential \eqref{eq.nolicq.sub.p2.linear.system.res} subproblems, as shown in Lemmas~\ref{lemma.diff.v} and~\ref{lemma.diff.u}, respectively. 
By \eqref{eq.diff.Jc} and \eqref{eq.diff.Jc.sto.det}, when $\epsilon_c = \epsilon_J = 0$, it follows $\mathcal{E}_{Jc} = \mathcal{E}_{Jc^2} = 0$. 

Let $v_k$ denote the inexact solution to the deterministic trust region subproblem 
\begin{align}
    \min _{v \in \operatorname{Range}\left( J_{k}^{T}\right)} \tfrac{1}{2}\left\| c_{k}+ J_{k} v \right\|^{2} \quad \text{ s.t. } \| v  \| \le\sigma_{Jc} \| J_{k}^T  c_{k}\|
\label{eq.nolicq.tr.sub.p1.det}
\end{align}
that is required to satisfy the Cauchy-like decrease condition
\begin{align}
    \left\| c_{k}\right\|-\left\| c_{k}+ J_{k}  v_{k}\right\| \geq  \gamma_{c}\left(\left\| c_{k}\right\|-\left\|  c_{k}+ \alpha_{k}^{c}  J_{k}  v_{k}^{c}\right\|\right),
\label{eq.nolicq.tr.cauchy.decrease.det}
\end{align}
where $ \gamma_{c} \in(0,1]$ is a user-defined parameter, $ v_{k}^{c}:=- J_{k}^{T}  c_{k}$ is the negative gradient direction for the objective of \eqref{eq.nolicq.tr.sub.p1.det} at $v=0$, and $ \alpha_{k}^{c} = \min \left\{ \sigma_{Jc}, \left\| J_{k}^{T}  c_{k}\right\|^{2} /\left\| J_{k}  J_{k}^{T}  c_{k}\right\|^{2} \right \} $ is the optimal step length along the direction of $ v_{k}^{c}$ that minimizes $\left\| c_{k}+\alpha  J_{k} v_{k}^{c}\right\|^2$ under the trust region constraint defined in \eqref{eq.nolicq.tr.sub.p1.det} over $\alpha \in \mathbb{R}$. The next lemma bounds the error between the true and noisy Cauchy steps. In the ideal setting, our goal is to show that $\|\bar v_k - v_k\|$ is sufficiently small for any inexact solutions  $\bar v_k$ and $v_k$ that satisfy Cauchy-like decrease. However, the fact that the difference between Cauchy step and exact step for either true or noisy problem cannot be bounded by $\epsilon$ terms under all conditions, prohibits us from being able to establish such a result. 
As a compromise, we provide bounds on the error between the true and noisy normal components in two extreme cases, where either both steps are Cauchy steps or both steps are the exact steps.

\begin{lemma}
\label{lemma.diff.v}
Suppose that Assumptions~\ref{ass:prob}, \ref{ass:error}, and \ref{ass.noearlytermination} hold. In addition, suppose that the singular values of both $J_k$ and $\bar J_k$ are bounded below by $\kappa_{\sigma} \in \R{}_{>0}$ for all $k \in \mathbb{N}$. Furthermore, suppose that $\bar v_k  = \bar  \alpha_{k}^{c}  \bar v_{k}^{c}$ is the Cauchy solution to \eqref{eq.nolicq.tr.sub.p1} and $v_k  = \alpha_{k}^{c} v_{k}^{c}$ is the Cauchy solution to \eqref{eq.nolicq.tr.sub.p1.det}. Then, there exists $\mathcal{E}_{v} \in \R{}_{\geq 0}$ 
such that $\|\bar v_k - v_k\| \le   \mathcal{E}_{v}$. 
\end{lemma}
\begin{proof}
\allowdisplaybreaks{
When $\|\bar c_k\|>\epsilon_o$, we have $ \bar v_k^c = - \bar J_k^T \bar c_k$. 
Let $\bar \alpha_{k}^{*} = \left\|\bar J_{k}^{T} \bar c_{k}\right\|^{2} /\left\|\bar J_{k} \bar J_{k}^{T} \bar c_{k}\right\|^{2}$ and $\alpha_{k}^{*} =  \left\| J_{k}^{T}  c_{k}\right\|^{2} /\left\| J_{k}  J_{k}^{T}  c_{k}\right\|^{2}$. It follows from Lemma \ref{lemma.mathcalE.jc} that 
\begin{align}
            \left\|J_{k} J_{k}^{T} c_{k}\right\| -  \left\|\bar J_{k} \bar J_{k}^{T} \bar  c_{k}\right\| &\le \left\|J_{k} J_{k}^{T} c_{k} - \bar J_{k} \bar J_{k}^{T} \bar  c_{k}\right\| \notag\\  
            &\le \left\|J_{k} J_{k}^{T} c_{k} -\bar J_{k} J_{k}^{T} c_{k}\right\| + \left\|\bar J_{k} J_{k}^{T} c_{k}  - \bar J_{k} \bar J_{k}^{T} \bar  c_{k}\right\| \notag \\  
            &\le   \kappa_J \kappa_c \epsilon_J + (\kappa_J + \epsilon_J) ( (\kappa_J + \epsilon_J) \epsilon_c + \kappa_c \epsilon_J ):= \mathcal{E}_{JJc}. \label{def.epsilon_JJc}  
    \end{align}

Let us discuss the following four cases:
\begin{itemize}
    \item[(a)] $\bar \alpha_{k}^{*} \le \sigma_{Jc}$ and $\alpha_{k}^{*} \le \sigma_{Jc} $: $ |\bar \alpha_{k}^{c}  - \alpha_{k}^{c} | = | \bar \alpha_{k}^{*}  - \alpha_{k}^{*} |$,
    \item[(b)]  $\bar \alpha_{k}^{*} \le \sigma_{Jc}$ and $\alpha_{k}^{*} > \sigma_{Jc} $: $ |\bar \alpha_{k}^{c}  - \alpha_{k}^{c} | = | \bar \alpha_{k}^{*}  - \sigma_{Jc} | = \sigma_{Jc} - \bar \alpha_{k}^{*} \le   \alpha_{k}^{*} - \bar \alpha_{k}^{*} = |\alpha_{k}^{*} - \bar \alpha_{k}^{*}|$,
    \item[(c)] $\bar \alpha_{k}^{*} > \sigma_{Jc}$ and $\alpha_{k}^{*} \le \sigma_{Jc} $: $ |\bar \alpha_{k}^{c}  - \alpha_{k}^{c} | = |  \alpha_{k}^{*}  - \sigma_{Jc} | = \sigma_{Jc} - \alpha_{k}^{*} \le   \bar  \alpha_{k}^{*} - \alpha_{k}^{*} = |\alpha_{k}^{*} - \bar \alpha_{k}^{*}|$, 
    \item[(d)] $\bar \alpha_{k}^{*} > \sigma_{Jc}$ and $\alpha_{k}^{*}> \sigma_{Jc} $: $ |\bar \alpha_{k}^{c}  - \alpha_{k}^{c} | = |  \sigma_{Jc} - \sigma_{Jc} | = 0 \le |\alpha_{k}^{*} - \bar \alpha_{k}^{*}|$, 
\end{itemize}
where $\sigma_{Jc}$ is the constant in the trust region constraint in \eqref{eq.nolicq.tr.sub.p1} and \eqref{eq.nolicq.tr.sub.p1.det}. In all cases, $ |\bar \alpha_{k}^{c}  - \alpha_{k}^{c} | \le  |\alpha_{k}^{*} - \bar \alpha_{k}^{*}|$, from which it follows that 
\begin{align*}
  \ &   \|\alpha_k^{c}  v_k^c - \bar \alpha_k^{c} \bar v_k^c\| 
\\  \le \ &   \|\alpha_k^{c}  v_k^c - \bar \alpha_k^{c}   v_k^c   \| + \|\bar  \alpha_k^{c}     v_k^c - \bar \alpha_k^{c}  \bar v_k^c   \| \\  \le \ &  |\alpha_k^{*} -\bar \alpha_k^{*}  |  \|     v_k^c   \| + \|\bar  \alpha_k^{c}    v_k^c - \bar \alpha_k^{c}  \bar v_k^c   \| \\  \le \ &  \left|\tfrac{\left\|J_{k}^{T} c_{k}\right\|^{2} }{\left\|J_{k} J_{k}^{T} c_{k}\right\|^{2} } - \tfrac{\left\|\bar J_{k}^{T}\bar  c_{k}\right\|^{2} }{\left\|\bar J_{k} \bar J_{k}^{T} \bar  c_{k}\right\|^{2} }  \right|   \|   J_{k}^{T} c_{k}   \| + \sigma_{Jc} \|  J_{k}^{T}  c_{k} -  \bar J_{k}^{T} \bar  c_{k}   \|   \\  = \ &  \left|\tfrac{\left \|J_{k}^{T} c_{k}\right\|^{2} \left\|\bar J_{k} \bar J_{k}^{T} \bar  c_{k}\right\|^{2} - \left \|\bar  J_{k}^{T}\bar  c_{k}\right\|^{2}\left\|\bar  J_{k}\bar  J_{k}^{T}\bar  c_{k}\right\|^{2}  + \left \|\bar J_{k}^{T} \bar c_{k}\right\|^{2} \left\|\bar J_{k}\bar  J_{k}^{T}\bar  c_{k}\right\|^{2} -  \left\|\bar J_{k}^{T}\bar  c_{k}\right\|^{2}  \left\|J_{k} J_{k}^{T} c_{k}\right\|^{2}  }{\left\|J_{k} J_{k}^{T} c_{k}\right\|^{2} \left\|\bar J_{k} \bar J_{k}^{T} \bar  c_{k}\right\|^{2} } \right| \left\|J_{k}^{T} c_{k}\right\| \\ 
& + \sigma_{Jc} \|  J_{k}^{T}  c_{k} -  \bar J_{k}^{T} \bar  c_{k}   \|  \\  
\leq \ &  \tfrac{\left\|  J_{k}^{T}    c_{k}\right\| }{\left\|J_{k} J_{k}^{T}  c_{k}\right\|^{2}} ( \| \bar J_{k}^{T}\bar  c_{k} \| + \| J_{k}^{T}  c_{k} \| ) (  \left|\| \bar J_{k}^{T}\bar  c_{k} \| - \| J_{k}^{T}  c_{k} \|\right| ) \\ 
& +\tfrac{\left\|  J_{k}^{T}    c_{k}\right\| }{\left\|J_{k} J_{k}^{T}  c_{k}\right\|^{2}}   \tfrac{\left\|  \bar J_{k}^{T}  \bar  c_{k}\right\|^{2} }{\left\|\bar J_{k} \bar J_{k}^{T} \bar  c_{k}\right\|^{2}}  \left( \left\|J_{k} J_{k}^{T} c_{k}\right\| + \left\|\bar J_{k} \bar J_{k}^{T} \bar  c_{k}\right\|  \right)\left( \left| \left\|J_{k} J_{k}^{T} c_{k}\right\| - \left\|\bar J_{k} \bar J_{k}^{T} \bar  c_{k}\right\| \right| \right) \\ 
& + \sigma_{Jc} \|  J_{k}^{T}  c_{k} -  \bar J_{k}^{T} \bar  c_{k}   \|  \\  
\le \ &  \tfrac{\left\|  J_{k}^{T}    c_{k}\right\| }{\left\|J_{k} J_{k}^{T}  c_{k}\right\|^{2}} ( \| \bar J_{k}^{T}\bar  c_{k} \| - \| J_{k}^{T}  c_{k} \|  + 2\| J_{k}^{T}  c_{k} \|  ) \| \bar J_{k}^{T}\bar  c_{k} - J_{k}^{T}  c_{k} \|  \\ 
& +\tfrac{\left\|  J_{k}^{T}    c_{k}\right\| }{ \left\|J_{k} J_{k}^{T}  c_{k}\right\|^{2}}  \kappa_{\sigma}^{-2} \left( \left\|\bar J_{k} \bar J_{k}^{T} \bar  c_{k}\right\| -  \left\|J_{k} J_{k}^{T} c_{k}\right\|    + 2\left\| J_{k}  J_{k}^{T}   c_{k}\right\|  \right) 
\left\|J_{k} J_{k}^{T} c_{k} - \bar J_{k} \bar J_{k}^{T} \bar  c_{k}\right\|  \\ 
& + \sigma_{Jc} \|  J_{k}^{T}  c_{k} -  \bar J_{k}^{T} \bar  c_{k}   \|   \\  
\le \ &  \tfrac{2 \left\|  J_{k}^{T}    c_{k}\right\|^2 }{\left\|J_{k} J_{k}^{T}  c_{k}\right\|^{2}} \mathcal{E}_{Jc} + \tfrac{ \left\|  J_{k}^{T}    c_{k}\right\|}{\left\|J_{k} J_{k}^{T}  c_{k}\right\|^{2}} \mathcal{E}_{Jc}^2 +  \tfrac{2 \left\|  J_{k}^{T}    c_{k}\right\|^2 }{\left\|J_{k} J_{k}^{T}  c_{k}\right\|^{2}} \kappa_J \kappa_{\sigma}^{-2} \mathcal{E}_{JJc} + \tfrac{\left\|  J_{k}^{T}    c_{k}\right\| }{\left\|J_{k} J_{k}^{T}  c_{k}\right\|^{2}} \kappa_{\sigma}^{-2} \mathcal{E}_{JJc}^2  + \sigma_{Jc}  \mathcal{E}_{Jc}, 
\end{align*}
where the second to last inequality follows from the imposed assumption about singular values of $\bar J_k$, and the last inequality follows from Lemma~\ref{lemma.mathcalE.jc} and \eqref{def.epsilon_JJc}.  
When $\left\|J_{k} J_{k}^{T}  c_{k}\right\| > \max \{ \mathcal{E}_{Jc}, \mathcal{E}_{JJc} \} $, 
\begin{align*}
      \ &   \|\alpha_k^{c}  v_k^c - \bar \alpha_k^{c} \bar v_k^c\|  \\  \le \ &  \tfrac{2 \left\|  J_{k}^{T}    c_{k}\right\|^2 }{\left\|J_{k} J_{k}^{T}  c_{k}\right\|^{2}} \mathcal{E}_{Jc} + \tfrac{ \left\|  J_{k}^{T}    c_{k}\right\|}{\left\|J_{k} J_{k}^{T}  c_{k}\right\|} \mathcal{E}_{Jc} +  \tfrac{2 \left\|  J_{k}^{T}    c_{k}\right\|^2 }{\left\|J_{k} J_{k}^{T}  c_{k}\right\|^{2}} \kappa_J \kappa_{\sigma}^{-2} \mathcal{E}_{JJc} + \tfrac{\left\|  J_{k}^{T}    c_{k}\right\| }{\left\|J_{k} J_{k}^{T}  c_{k}\right\|} \kappa_{\sigma}^{-2}  \mathcal{E}_{JJc}  + \sigma_{Jc}  \mathcal{E}_{Jc}  \\  \le \ &  2 \kappa_{\sigma}^{-2}  \mathcal{E}_{Jc} + \kappa_{\sigma}^{-1}  \mathcal{E}_{Jc} + 2  \kappa_{\sigma}^{-4} \kappa_J  \mathcal{E}_{JJc} +  \kappa_{\sigma}^{-3}  \mathcal{E}_{JJc} + \sigma_{Jc}  \mathcal{E}_{Jc}. 
\end{align*}
Otherwise, $\left\| J_{k}^{T}  c_{k}\right\|  \le \kappa_{\sigma}^{-1}  \left\|J_{k} J_{k}^{T}  c_{k}\right\| \le \kappa_{\sigma}^{-1} \max  \{ \mathcal{E}_{Jc}, \mathcal{E}_{JJc} \} $, and it follows by Lemma~\ref{lemma.mathcalE.jc} 
\begin{align*}
      \|\alpha_k^{c}  v_k^c - \bar \alpha_k^{c} \bar v_k^c\| &\le \max \{\alpha_k^{*},  \bar \alpha_k^{*} \}   \|     v_k^c   \| + \|\bar  \alpha_k^{c}    v_k^c - \bar \alpha_k^{c}  \bar v_k^c   \| \\ 
       &\le  \sigma_{Jc} \kappa_{\sigma}^{-1} \max  \{ \mathcal{E}_{Jc}, \mathcal{E}_{JJc} \}  + \sigma_{Jc}  \mathcal{E}_{Jc}. 
\end{align*}
If $\|\bar c_k \| \le \epsilon_o$, then $ \bar v_k^c = 0$ and $\| c_k\| \le \|\bar c_k \| + \| c_k  -\bar c_k \| \le 2 \epsilon_c$, from which it follows that 
\begin{align*}
      \|\alpha_k^{c}  v_k^c - \bar \alpha_k^{c} \bar v_k^c\| = 
\| \alpha_k^{c}  J_k^T  c_k \| \le  
    2  \sigma_{Jc}  \kappa_J  \epsilon_c. 
\end{align*}
The desired conclusion follows by defining 

\noindent
$\mathcal{E}_{v} := \max \left\{(2 \kappa_{\sigma}^{-2} + \kappa_{\sigma}^{-1} +  \sigma_{Jc}  ) \mathcal{E}_{Jc}  + (  2  \kappa_{\sigma}^{-4} \kappa_J + \kappa_{\sigma}^{-3}   ) \mathcal{E}_{JJc}   ,\sigma_{Jc} \kappa_{\sigma}^{-1} \max  \left\{ \mathcal{E}_{Jc}, \mathcal{E}_{JJc} \right\}  + \sigma_{Jc}  \mathcal{E}_{Jc}, 2  \sigma_{Jc}  \kappa_J  \epsilon_c  \right\}$. 
}
\end{proof}

The next lemma provides an upper bound of the perturbation errors of the exact (optimal) solutions to \eqref{eq.nolicq.tr.sub.p1} and \eqref{eq.nolicq.tr.sub.p1.det}.
\begin{lemma}
\label{lemma.diff.v.exact}
Suppose that Assumptions~\ref{ass:prob}, \ref{ass:error}, and  \ref{ass.noearlytermination}  hold. In addition, suppose that the singular values of both $J_k$ and $\bar J_k$ are bounded below by $\kappa_{\sigma} \in \R{}_{>0}$ for all $k \in \mathbb{N}$. 
Furthermore, suppose that  ${\bar v_k^*}$ and  ${v_k^*}$ are the exact solutions to \eqref{eq.nolicq.tr.sub.p1} and \eqref{eq.nolicq.tr.sub.p1.det}, respectively.
Then, there exists $\mathcal{E}_{v^*} 
\in \mathbb{R}_{\geq 0}$ 
such that $\|{\bar v_k^*}  - {v_k^*}\| \le   \mathcal{E}_{v^*}$. 
\end{lemma}
\begin{proof}
\allowdisplaybreaks{
When $\|\bar c_k\| \le \epsilon_o \le \epsilon_c$, by \eqref{eq.nolicq.tr.sub.p1} and Assumption~\ref{ass:prob}, it follows that $\| \bar v_k^*  - v_k^*\| \le \| v_k^*\| \le \sigma_{Jc} \| J_k^T c_k\| \le  \sigma_{Jc} \kappa_J \epsilon_c$.

For the remainder of the proof we consider the case where $\|\bar c_k\| > \epsilon_o$. 
Let $\bar v_k = \bar J_k^T \bar s_k $ and $ v_k =  J_k^T s_k $ with vectors $\{\bar s_k, s_k\} \subset \mathbb{R}^m$. We have $\| \bar{J}_k^T \bar{c}_k \| \ne 0$ since otherwise Algorithm~\ref{alg.dfo_sqp_LS} terminates on Line~\ref{line.early.stopping} and violates Assumption~\ref{ass.noearlytermination}. When $\| J_k^T c_k \| = 0$, we have $\|v_k\| = 0$ and $\|\bar J_k^T \bar c_k \| \le \mathcal{E}_{Jc}$ by Lemma~\ref{lemma.mathcalE.jc}, which implies that $\|\bar v_k - v_k\|\le \|\bar v_k\| \le \sigma_{Jc} \mathcal{E}_{Jc}$ and the desired conclusion holds. We therefore assume that $\| J_k^T c_k \| \ne 0$ in the remainder of the analysis. 
Using the vectors 
$\{\bar s_k, s_k\} \subset \mathbb{R}^m$ we reformulate problems \eqref{eq.nolicq.tr.sub.p1} and \eqref{eq.nolicq.tr.sub.p1.det} 
\begin{align}
 &\  \min_{\bar s \in\mathbb{R}^m} \;\ \tfrac{1}{2}\left\|\bar c_{k}+\bar J_{k}\bar J_{k}^T \bar  s \right\|^{2} \quad \text{ s.t. } \|\bar J_{k}^T \bar  s  \| \le  \sigma_{Jc} \|\bar J_{k}^T \bar c_{k}\|,  
\label{eq.nolicq.tr.sub.p1.s} \\ 
 \text{and} \qquad  &\  \min_{s \in\mathbb{R}^m}  \;\ \tfrac{1}{2}\left\| c_{k}+ J_{k}J_{k}^T s\right\|^{2} \quad \text{ s.t. } \|J_{k}^T s\| \le  \sigma_{Jc} \| J_{k}^T  c_{k}\|. 
\label{eq.nolicq.tr.sub.p1.det.s}
\end{align} 
Let $\mathcal{C} = \{s \in\mathbb{R}^m | \|\bar J_{k}^T  s  \| \le  \sigma_{Jc} \|\bar J_{k}^T \bar c_{k}\| \}$, where $\mathcal{C}$ is a nonempty, closed, convex set, and 
\begin{align}
\label{eq.hat.s}
    \hat  s  =  \tfrac{\|\bar J_k^T \bar c_k \|}{\| J_k^T c_k \| } \left( \bar J_k  \bar J_k^T  \right)^{-1/2} \left(  J_k  J_k^T  \right)^{1/2} s.  
\end{align}
Note that $\bar J_k  \bar J_k^T$ and $ J_k  J_k^T$ are both invertible since $\bar J_k $ and $J_k $ are of full rank, as well as $\| \bar{J}_k^T \bar{c}_k \| \ne 0$ and $\| {J}_k^T {c}_k \| \ne 0$. 
Then $\hat  s $ in~\eqref{eq.hat.s} satisfies
\begin{align*}
    \|\bar J_k^T \hat s \|^2  = 
    \tfrac{\|\bar J_k^T \bar c_k \|^2}{\| J_k^T c_k \|^2 } s^T \left(  J_k  J_k^T  \right)^{1/2}  \left( \bar J_k  \bar J_k^T  \right)^{-1/2}  \bar J_k \bar J_k^T  \left( \bar J_k  \bar J_k^T  \right)^{-1/2} \left(  J_k  J_k^T  \right)^{1/2} s  = \tfrac{\|\bar J_k^T \bar c_k \|^2}{\| J_k^T c_k \|^2 }   \| J_k^T  s \|^2. 
\end{align*}
Thus, we can reformulate \eqref{eq.nolicq.tr.sub.p1.det.s} by changing the variable   from $s$ to $\hat s$ as 
\begin{align}
  \min_{\hat  s \in \mathbb{R}^m} \;\ \tfrac{1}{2}\left\| c_{k}+ \tfrac{\| J_k^T c_k \| }{\|\bar J_k^T \bar c_k \|} \left( J_{k} J_k^T  \right)^{1/2}\left( \bar J_{k}\bar J_k^T  \right)^{1/2}  \hat  s \right\|^{2} \quad \text{ s.t. } \|\bar J_{k}^T \hat s  \| \le  \sigma_{Jc} \|\bar J_{k}^T \bar c_{k}\|. 
\label{eq.nolicq.tr.sub.p1.det.s.hat}
\end{align} 
Let us further reformulate \eqref{eq.nolicq.tr.sub.p1.s} and \eqref{eq.nolicq.tr.sub.p1.det.s.hat} to \eqref{eq.nolicq.tr.sub.p1.s.proj} and \eqref{eq.nolicq.tr.sub.p1.det.s.proj}, respectively, by rewriting the objective functions
\begin{align}
 &\  \min_{\bar s \in\mathbb{R}^m} \tfrac{1}{2}\left\|\bar s   + (\bar J_{k}\bar J_{k}^T)^{-1} \bar c_{k} \right\|_{\bar J_{k} \bar J_{k}^T \bar J_{k} \bar J_{k}^T}^{2} \;\ \text{ s.t. } \|\bar J_{k}^T \bar s  \| \le  \sigma_{Jc} \|\bar J_{k}^T \bar c_{k}\|, \qquad \text{and}
  \label{eq.nolicq.tr.sub.p1.s.proj} \\
  &\   \min_{\hat s \in\mathbb{R}^m} \tfrac{\|J_k^Tc_k\|^2}{2\|\bar{J}_k^T\bar{c}_k\|^2}\left\|\hat s  + \tfrac{\|\bar J_k^T \bar c_k \|}{\| J_k^T c_k \| } (\bar J_{k}\bar J_k^T)^{-1/2} (J_{k} J_k^T)^{-1/2}   c_k  \right\|_{(\bar J_{k} \bar J_{k}^T)^{1/2} J_{k}  J_{k}^T (\bar J_{k} \bar J_{k}^T)^{1/2} }^{2} \;\ \text{ s.t. } \|\bar J_{k}^T \hat s\| \le  \sigma_{Jc} \|\bar J_{k}^T \bar c_{k}\|. 
\label{eq.nolicq.tr.sub.p1.det.s.proj}
\end{align} 
Note that \eqref{eq.nolicq.tr.sub.p1.s.proj} and \eqref{eq.nolicq.tr.sub.p1.det.s.proj} share the same constraint. 
Let $ \bar s^*$ and $\hat s^*$ be the optimal solutions of \eqref{eq.nolicq.tr.sub.p1.s.proj} and \eqref{eq.nolicq.tr.sub.p1.det.s.proj}, respectively. It then follows from the projection theorem \cite[Theorem 4.3-1]{ciarlet2013linear} that $ \bar s^*$ and $\hat s^*$ are both unique solutions to \eqref{eq.nolicq.tr.sub.p1.s.proj} and \eqref{eq.nolicq.tr.sub.p1.det.s.proj}, and
\begin{align}
    \forall t \in \mathcal{C}, \quad  \left  \langle -(\bar J_{k} \bar J_{k}^T)^{-1} \bar c_k - \bar s^*, t - \bar s^*   \right \rangle_{\bar J_{k} \bar J_{k}^T \bar J_{k} \bar J_{k}^T} &\le 0, \label{eq.projection1} \\
     \forall t \in \mathcal{C}, \quad \left \langle \hat s^* + \tfrac{\|\bar J_k^T \bar c_k \|}{\| J_k^T c_k \| } (\bar J_{k}\bar J_k^T)^{-1/2} (J_{k} J_k^T)^{-1/2}   c_k  , \hat s^* - t  \right \rangle_{ (\bar J_{k} \bar J_{k}^T)^{1/2} J_{k}  J_{k}^T (\bar J_{k} \bar J_{k}^T)^{1/2} } &\le 0, \label{eq.projection2}
\end{align}
where $\langle x,y\rangle_A = x^TAy$, and $x,y \in \mathbb{R}^n$ and $A \in \mathbb{R}^{n \times n}$ is a positive definite matrix. 
Plugging $\hat s^*$ and $\bar s^*$ into \eqref{eq.projection1} and \eqref{eq.projection2}, respectively, 
\begin{align}
    \left \langle - (\bar J_{k} \bar J_{k}^T)^{-1} \bar c_k  - \bar s^* , \hat s^* - \bar s^* \right  \rangle_{\bar J_{k} \bar J_{k}^T \bar J_{k} \bar J_{k}^T} &\le 0,  \label{eq.projection3} \\
   \text{and} \quad \left   \langle \hat s^* + \tfrac{\|\bar J_k^T \bar c_k \|}{\| J_k^T c_k \| } (\bar J_{k}\bar J_k^T)^{-1/2} (J_{k} J_k^T)^{-1/2}  c_k , \hat s^* - \bar s^*  \right \rangle_{(\bar J_{k} \bar J_{k}^T)^{1/2} J_{k}  J_{k}^T (\bar J_{k} \bar J_{k}^T)^{1/2}} &\le 0.  \label{eq.projection4}  
\end{align}
Summing \eqref{eq.projection3} and \eqref{eq.projection4}, we have
\begin{align*}
    & \left \langle - (\bar J_{k} \bar J_{k}^T)^{-1} \bar c_k  - \bar s^* +  \hat s^* + \tfrac{\|\bar J_k^T \bar c_k \|}{\| J_k^T c_k \| } (\bar J_{k}\bar J_k^T)^{-1/2} (J_{k} J_k^T)^{-1/2} c_k , \hat s^* - \bar s^*  \right \rangle_{\bar J_{k} \bar J_{k}^T \bar J_{k} \bar J_{k}^T} \\
    & \quad + \left( \hat s^* + \tfrac{\|\bar J_k^T \bar c_k \|}{\| J_k^T c_k \| }(\bar J_{k}\bar J_k^T)^{-1/2} (J_{k} J_k^T)^{-1/2}  c_k\right)^T \left( (\bar J_{k} \bar J_{k}^T)^{1/2} J_{k}  J_{k}^T (\bar J_{k} \bar J_{k}^T)^{1/2} -  \bar J_{k} \bar J_{k}^T \bar J_{k} \bar J_{k}^T \right) (\hat s^* - \bar s^* )  \le 0.     
\end{align*}
Since the minimum singular value of $\bar{J}_k$ is lower bounded by $\kappa_{\sigma}$, it follows that the minimum eigenvalue of $(\bar J_{k} \bar J_{k}^T)^2$ is bounded below by $\kappa_{\sigma}^4$. Thus, 
\begin{align}
  &\  \|\hat s^* - \bar s^* \|^2 \notag \\   \le &\ \kappa_{\sigma}^{-4} 
   \|\hat s^* - \bar s^* \|_{\bar J_{k} \bar J_{k}^T \bar J_{k} \bar J_{k}^T}^2 \notag \\
    \le  &\  \kappa_{\sigma}^{-4} \left \langle  (\bar J_{k} \bar J_{k}^T)^{-1} \bar c_k  -  \tfrac{\|\bar J_k^T \bar c_k \|}{\| J_k^T c_k \| } (\bar J_{k}\bar J_k^T)^{-1/2} (J_{k} J_k^T)^{-1/2}  c_k , \hat s^* - \bar s^*\right    \rangle_{\bar J_{k} \bar J_{k}^T \bar J_{k} \bar J_{k}^T}  \notag \\
    \quad &\ + \kappa_{\sigma}^{-4} \left( \hat s^* + \tfrac{\|\bar J_k^T \bar c_k \|}{\| J_k^T c_k \| } (\bar J_{k}\bar J_k^T)^{-1/2} (J_{k} J_k^T)^{-1/2}  c_k \right)^T \left(  \bar J_{k} \bar J_{k}^T \bar J_{k} \bar J_{k}^T - (\bar J_{k} \bar J_{k}^T)^{1/2} J_{k}  J_{k}^T (\bar J_{k} \bar J_{k}^T)^{1/2} \right) (\hat s^* - \bar s^* ). \label{eq.bound.s2.s1}
\end{align}
Since the minimum singular value of $\bar J_{k}$ and $J_{k}$ are both bounded below by $\kappa_{\sigma}$, 
$\| (\bar J_{k}\bar J_k^T)^{-1/2} \| \le \kappa_{\sigma}^{-1}$ and  $\| (J_{k} J_k^T)^{-1/2} \| \le \kappa_{\sigma}^{-1}$.   
It follows by \cite[Lemma 2.2]{schmitt1992perturbation} (perturbation bound for matrix square root) and \cite[Theorem 2.5]{stewart1990matrix} (perturbation bound for matrix inverse) that
\begin{equation}
    \begin{split}
           \left  \|(\bar J_{k}\bar J_k^T)^{-1/2} - (J_{k} J_k^T)^{-1/2} \right \| & \le \left \| (\bar J_{k}\bar J_k^T)^{-1/2} \right \| \left  \| (J_{k} J_k^T)^{-1/2} \right\|   \left \| (\bar J_{k}\bar J_k^T)^{1/2} - (J_{k} J_k^T)^{1/2}  \right \| \\ & \le \kappa_{\sigma}^{-2} \tfrac{1}{\kappa_{\sigma}+ \kappa_{\sigma}} \|\bar J_{k}\bar J_k^T - J_{k} J_k^T \|  \le \tfrac{ (2 \kappa_J + \epsilon_J ) \epsilon_J }{ 2 \kappa_{\sigma}^{3} }.  
    \end{split}
         \label{eq.diff.JJ.half}
\end{equation}
Moreover, by 
\eqref{eq.diff.JJ.half}, Lemma~\ref{lemma.mathcalE.jc}, and Assumptions~\ref{ass:prob} and \ref{ass:error}, it follows that
\begin{align}
   &  \quad \left\| (\bar J_{k} \bar J_{k}^T)^{-1} \bar c_k  -  \tfrac{\|\bar J_k^T \bar c_k \|}{\| J_k^T c_k \| } (\bar J_{k}\bar J_k^T)^{-1/2} (J_{k} J_k^T)^{-1/2}  c_k  \right \| \notag\\ 
   &\le \left\| (\bar J_{k} \bar J_{k}^T)^{-1} \bar c_k  - (\bar J_{k}\bar J_k^T)^{-1/2} (J_{k} J_k^T)^{-1/2} \bar c_k  \right \| + \left\|  (\bar J_{k}\bar J_k^T)^{-1/2} (J_{k} J_k^T)^{-1/2} (\bar c_k - c_k)  \right \| \notag\\ 
   &\quad +\left\| \tfrac{\| J_k^T c_k \| - \|\bar J_k^T \bar c_k \|}{\| J_k^T c_k \| }(\bar J_{k}\bar J_k^T)^{-1/2} (J_{k} J_k^T)^{-1/2} c_k \right \|  \notag\\ 
   &\le \left\| (\bar J_{k}\bar J_k^T)^{-1/2} \left( (\bar J_{k}\bar J_k^T)^{-1/2}  - (J_{k} J_k^T)^{-1/2} \right) \bar c_k  \right \| + \left\|  (\bar J_{k}\bar J_k^T)^{-1/2} (J_{k} J_k^T)^{-1/2} (\bar c_k - c_k)  \right \| \notag\\ 
   &\quad + \tfrac{|\| J_k^T c_k \| - \|\bar J_k^T \bar c_k \||}{\| J_k^T c_k \| } \left\| (\bar J_{k}\bar J_k^T)^{-1/2} (J_{k} J_k^T)^{-1/2} c_k \right \| \notag\\ 
      &\le \left\| (\bar J_{k}\bar J_k^T)^{-1/2} \right \|   \left\|  (\bar J_{k}\bar J_k^T)^{-1/2}  - (J_{k} J_k^T)^{-1/2} \right \|   \left\| \bar c_k  \right \| + \left\|  (\bar J_{k}\bar J_k^T)^{-1/2} \right \|  
     \left\|  (J_{k} J_k^T)^{-1/2} \right \|    \left\|  \bar c_k - c_k  \right \| 
   \notag\\ 
   &\quad + \| J_k^T c_k  - \bar J_k^T \bar c_k \| \tfrac{\|c_k\|}{\| J_k^T c_k \| } \left\| (\bar J_{k}\bar J_k^T)^{-1/2}  \right \|   \left\|  (J_{k} J_k^T)^{-1/2} \right \| 
   \notag\\ 
   &\le 
   \tfrac{ \left( (2 \kappa_J + \epsilon_J ) (\kappa_c + \epsilon_c)  + 2\kappa_{\sigma} \kappa_c  \right) \epsilon_J + ( 2\kappa_{\sigma}^{2} + 2\kappa_{\sigma}(\kappa_J + \epsilon_J) )\epsilon_c  }{ 2 \kappa_{\sigma}^{4} } .
  \label{eq.first.term}
\end{align}
The first term of inequality \eqref{eq.bound.s2.s1} (excluding constant) can 
be bounded by 
\begin{align*}
 & \quad  \left \langle  (\bar J_{k} \bar J_{k}^T)^{-1} \bar c_k  -  \tfrac{\|\bar J_k^T \bar c_k \|}{\| J_k^T c_k \| } (\bar J_{k}\bar J_k^T)^{-1/2} (J_{k} J_k^T)^{-1/2}  c_k , \hat s^* - \bar s^*\right    \rangle_{\bar J_{k} \bar J_{k}^T \bar J_{k} \bar J_{k}^T} \\ &\le   \tfrac{ \left( (2 \kappa_J + \epsilon_J ) (\kappa_c + \epsilon_c)  + 2\kappa_{\sigma} \kappa_c  \right) \epsilon_J + ( 2\kappa_{\sigma}^{2} + 2\kappa_{\sigma}(\kappa_J + \epsilon_J) )\epsilon_c  }{ 2 \kappa_{\sigma}^{4} }    (\kappa_J + \epsilon_J)^4  \| \hat s^* - \bar s^* \|.
\end{align*}
Additionally, it follows by Assumptions~\ref{ass:prob} and \ref{ass:error} that 
\begin{align*}
 \| \bar J_{k} \bar J_{k}^T \bar J_{k} \bar J_{k}^T -    (\bar J_{k} \bar J_{k}^T)^{1/2} J_{k}  J_{k}^T (\bar J_{k} \bar J_{k}^T)^{1/2} \|   &\le 
   \|   (\bar J_{k} \bar J_{k}^T)^{1/2} (   \bar J_{k} \bar J_{k}^T -   J_{k}  J_{k}^T )  (\bar J_{k} \bar J_{k}^T)^{1/2}  \| \\&\le (\kappa_J + \epsilon_J)^2 (2\kappa_J + \epsilon_J  ) \epsilon_J. 
\end{align*}
The second term of inequality \eqref{eq.bound.s2.s1} (excluding constant) can 
be bounded by 
\begin{equation}
    \begin{split}
   & \quad    \left ( \hat s^* + \tfrac{\|\bar J_k^T \bar c_k \|}{\| J_k^T c_k \| } (\bar J_{k}\bar J_k^T)^{-\tfrac12} (J_{k} J_k^T)^{-\tfrac12}  c_k \right )^T \left(  \bar J_{k} \bar J_{k}^T \bar J_{k} \bar J_{k}^T - (\bar J_{k} \bar J_{k}^T)^{\tfrac12} J_{k}  J_{k}^T (\bar J_{k} \bar J_{k}^T)^{\tfrac12} \right) (\hat s^* - \bar s^* )  \\  &\le \left( \|  \hat s^* \|  +  \tfrac{\|\bar J_k^T \bar c_k \|}{\| J_k^T c_k \| }  \left\| (\bar J_{k}\bar J_k^T)^{-\tfrac12} (J_{k} J_k^T)^{-\tfrac12}  c_k  \right\| \right)  (\kappa_J + \epsilon_J)^2 (2\kappa_J + \epsilon_J  ) \epsilon_J  \| \hat s^* - \bar s^* \| \\  &\le \left( \tfrac{\|\hat s^*\|}{\| \bar J_k^T \hat s^* \|}  \| \bar J_k^T \hat s^* \|  +  \|\bar J_k^T \bar c_k \| \tfrac{\| c_k \|}{\| J_k^T c_k \| }  \left\| (\bar J_{k}\bar J_k^T)^{-\tfrac12} \right\|  \left\|(J_{k} J_k^T)^{-\tfrac12}  \right\| \right)  (\kappa_J + \epsilon_J)^2 (2\kappa_J + \epsilon_J  ) \epsilon_J  \| \hat s^* - \bar s^* \| \\
   &\le \left( \kappa_{\sigma}^{-1}  \sigma_{Jc}  + \kappa_{\sigma}^{-3}    \right) (\kappa_c + \epsilon_c) (\kappa_J + \epsilon_J)^3 (2\kappa_J + \epsilon_J)  \epsilon_J  \| \hat s^* - \bar s^* \|. 
    \end{split}
       \label{eq.second.term}
\end{equation}
As a result, by plugging \eqref{eq.first.term} and \eqref{eq.second.term} into \eqref{eq.bound.s2.s1}, it follows that
\begin{equation}
    \begin{split}
\label{eq.def.diff.s.bar}
    \| \hat s^* - \bar s^* \| &\le (\kappa_J + \epsilon_J)^4\tfrac{ \left( (2 \kappa_J + \epsilon_J ) (\kappa_c + \epsilon_c)  + 2\kappa_{\sigma} \kappa_c  \right) \epsilon_J + ( 2\kappa_{\sigma}^{2} + 2\kappa_{\sigma}(\kappa_J + \epsilon_J) )\epsilon_c  }{ 2 \kappa_{\sigma}^{8} }     \\  & \quad + \kappa_{\sigma}^{-4} \left( \kappa_{\sigma}^{-1}  \sigma_{Jc}  + \kappa_{\sigma}^{-3}    \right) (\kappa_c + \epsilon_c) (\kappa_J + \epsilon_J)^3 (2\kappa_J + \epsilon_J)  \epsilon_J:= \mathcal{E}_{\bar s}.  
    \end{split}
\end{equation}
By Assumption~\ref{ass:prob}, \eqref{eq.nolicq.tr.sub.p1.det.s}, \eqref{eq.hat.s}, \eqref{eq.diff.JJ.half}, and the fact that all singular values of $\bar J_k$ are bounded below by $\kappa_{\sigma}$, 
\begin{equation}
    \begin{split}
\label{eq.def.diff.s}
    \| \hat s^* -  s^* \| &= \left\| \tfrac{\|\bar J_k^T \bar c_k \|}{\| J_k^T c_k \| }(\bar J_{k}\bar J_k^T)^{-1/2} (J_{k} J_k^T)^{1/2} s^* -  s^* \right \| 
    \\ &\le \tfrac{\|J_{k} J_k^T\|^{1/2}}{ \| J_k^T c_k \|  } \left\| \|\bar J_k^T \bar c_k \|  (\bar J_{k}\bar J_k^T)^{-1/2} - \| J_k^T c_k \| (J_{k} J_k^T)^{-1/2} \right\| \tfrac{\|s^*\|}{ \| J_k^T s^*\|}  \| J_k^T s^*\| 
       \\ &\le \kappa_J \sigma_{Jc} \kappa_{\sigma}^{-1} \left\| \|\bar J_k^T \bar c_k  -  J_k^T c_k \| (\bar J_{k}\bar J_k^T)^{-1/2}  + \| J_k^T c_k \|  \left( (\bar J_{k} \bar J_k^T)^{-1/2}  - (J_{k} J_k^T)^{-1/2} \right) \right\| 
       \\ &\le \kappa_J \sigma_{Jc} \kappa_{\sigma}^{-1} \left( \|\bar J_k^T \bar c_k  -  J_k^T c_k \| \left\|(\bar J_{k}\bar J_k^T)^{-1/2} \right\|  + \| J_k^T c_k \|  \left\| (\bar J_{k} \bar J_k^T)^{-1/2}  - (J_{k} J_k^T)^{-1/2} \right\| \right)
    \\ &\le \kappa_J \sigma_{Jc} \kappa_{\sigma}^{-1} \left(  \kappa_{\sigma}^{-1} (  \kappa_c \epsilon_J + (\kappa_J + \epsilon_J) \epsilon_c ) +  \tfrac{\kappa_J \kappa_c  (2 \kappa_J + \epsilon_J ) \epsilon_J }{ 2 \kappa_{\sigma}^{3} }  \right): = \mathcal{E}_s. 
    \end{split}
\end{equation}
It then follows 
that $\| \bar s^* -  s^* \| \le  \| \hat s^* -  s^* \| +  \| \hat s^* - \bar   s^* \| \le \mathcal{E}_s+ \mathcal{E}_{\bar s}$. 
By Assumption~\ref{ass:prob}, \ref{ass:error}, \eqref{eq.def.diff.s}, and the fact that 
$\bar v_k^* = \bar J_k^T \bar s_k $ and $ v_k^* =  J_k^T s_k $, 
\begin{align*}
    \| \bar v_k^* -  v_k^* \|  
    \le \|\bar J_k^T  \bar s^* - \bar   J_k^T  s^* \|  +  \|\bar J_k^T   s^* -   J_k^T  s^* \| \le (\kappa_J + \epsilon_J) ( \mathcal{E}_s+ \mathcal{E}_{\bar s}) + \sigma_{Jc} \kappa_{\sigma}^{-1} \kappa_J \kappa_c \epsilon_J. 
\end{align*}

Defining $\mathcal{E}_{v^*} := \max \left \{ (\kappa_J + \epsilon_J) ( \mathcal{E}_s+ \mathcal{E}_{\bar s}) + \sigma_{Jc} \kappa_{\sigma}^{-1} \kappa_J \kappa_c \epsilon_J, \sigma_{Jc} \kappa_J \epsilon_c   \right\}$ 
concludes the proof.}
\end{proof}

\begin{remark}
Lemma~\ref{lemma.diff.v} shows that when the LICQ holds for both true and noisy constraint Jacobians, 
Cauchy solutions to problems \eqref{eq.nolicq.tr.sub.p1} and \eqref{eq.nolicq.tr.sub.p1.det} are sufficiently close. Lemma~\ref{lemma.diff.v.exact} shows that a similar result holds for the optimal solutions of problems \eqref{eq.nolicq.tr.sub.p1} and \eqref{eq.nolicq.tr.sub.p1.det}. It can be seen from the definitions of $\mathcal{E}_v$ and $\mathcal{E}_{v^*}$ that both perturbation errors depend on $\epsilon_c$ and $\epsilon_J$ linearly. Moreover, when $\epsilon_c = \epsilon_J = 0$, both errors are zero, i.e., $\mathcal{E}_v =\mathcal{E}_{v^*} =  0$. The analysis in Lemmas~\ref{lemma.diff.v} and \ref{lemma.diff.v.exact} relies critically on the LICQ. Without it, one cannot establish bounds for the ratios $\left\|\bar J_{k}^{T} \bar c_{k}\right\| /\left\|\bar J_{k} \bar J_{k}^{T} \bar c_{k}\right\|$ and $\left\| J_{k}^{T}  c_{k}\right\|/\left\| J_{k}  J_{k}^{T}  c_{k}\right\|$ 
(proof of Lemma~\ref{lemma.diff.v}), and the non-uniqueness of optimal solutions to \eqref{eq.nolicq.tr.sub.p1} makes it challenging to derive the results in Lemma~\ref{lemma.diff.v.exact}. We note that Lemmas~\ref{lemma.diff.v} and \ref{lemma.diff.v.exact} are not used to derive the main convergence guarantees (Section~\ref{sec:convergence.main}).
\end{remark}

As we will see in Section~\ref{sec:convergence.main}, our main interest is the distance between the true and the noisy tangential 
components. Let $u_{k}^*$ denote the exact solution to the deterministic version of the linear system \eqref{eq.nolicq.sub.p2.linear.system.star}, i.e.,  
\begin{align}
    \left[\begin{array}{cc}
H_{k} &  J_{k}^{T} \\
J_{k} & 0
\end{array}\right]\left[\begin{array}{c}
  {u_{k}^*} \\
 {y_{k}^*}
\end{array}\right]=-\left[\begin{array}{c}
 g_{k}+H_{k}  v_{k} \\
0
\end{array}\right],
\label{eq.nolicq.sub.p2.linear.system.star.det}
\end{align}
where $y_{k}^* \in \mathbb{R}^m$ is the vector of dual multipliers.  We aim to derive an upper bound for 
$\|u_k^* - \bar u_{k}\|$, the difference between the exact solution to the deterministic linear system \eqref{eq.nolicq.sub.p2.linear.system.star.det} and an inexact solution to the noisy linear system \eqref{eq.nolicq.sub.p2.linear.system.star}. The following results hold regardless of whether the LICQ holds for either the noisy or deterministic problems. First present a technical lemma that bounds the true and noisy tangential steps, the true search direction, and the dual multipliers.

\begin{lemma}
\label{lemma.y.exist}
Suppose that Assumptions~\ref{ass:prob}, \ref{ass:error}, \ref{ass:H}, 
and \ref{ass.noearlytermination} hold. In addition, suppose that the minimum non-zero singular values of both $J_k$ and $\bar J_k$ are bounded below by $\kappa_{\sigma} \in \R{}_{>0}$ for all $k \in \mathbb{N}$. 
Then, for all $k\in \mathbb{N}$, there exist $\{ \kappa_{\bar u^*}, \kappa_{u^*},  \kappa_{d^*}
\} \subset \mathbb{R}_{>0}$ 
such that $\| \bar u_{k}^*\| \le \kappa_{\bar u^*}$, $\|   {u_{k}^*}\| \le \kappa_{u^*}$ and $\|   {d_{k}^*}\| \le \kappa_{d^*}$, where $d_{k}^* = u_{k}^* + v_{k}$, where $v_k$ is the inexact solution to \eqref{eq.nolicq.tr.sub.p1.det} satisfying \eqref{eq.nolicq.tr.cauchy.decrease.det}. 
    Moreover, for all $k\in \mathbb{N}$, there exist $\{\kappa_{\bar y^*}, \kappa_{y^*} \} \subset \mathbb{R}_{>0}$ and solutions $\{\bar y_{k}^*, y_{k}^*\}\subset \mathbb{R}^m$ from linear systems \eqref{eq.nolicq.sub.p2.linear.system.res} and \eqref{eq.nolicq.sub.p2.linear.system.star.det} 
such that $\|  \bar y_{k}^* \|\le \kappa_{\bar y^*}$ and $\| {y_k^*} \| \le  \kappa_{y^*}$.
\end{lemma}
\begin{proof}
Suppose that $\text{rank}(\bar J_k ) = p$ and $\text{rank}( J_k ) = q$. First, let 
$p = n \le m$, in which case the null space of $\bar J_k$ is trivial (contains only the zero vector). 
By \eqref{eq.nolicq.sub.p2.linear.system.star}, 
it follows that $\bar u_k^* = 0$ and $\bar y_k^* = (\bar J_k^T \bar J_k)^{-1} \bar J_k^T (\bar g_k + \bar H_k \bar v_k )$, where the latter is bounded by Assumptions~\ref{ass:prob} and~\ref{ass:error}.  Similarly, when $q = n \le m$, it follows that $u_k^* = 0$ and $y_k^* = (J_k^T J_k)^{-1} J_k^T (g_k + H_k v_k )$, where the latter is bounded by Assumptions~\ref{ass:prob} and~\ref{ass:error}. The desired conclusion therefore holds for $p = n \le m$ or $q = n \le m$.

Hence, for the remainder of the proof, we consider 
$\text{rank}(\bar J_k ) = p < n$  and $\text{rank}(J_k ) = q < n$. 
Suppose that   $\Sigma_{k}$ are diagonal matrices containing the non-zero singular values of $J_k$. Similar to \eqref{eq.J.bar.SVD}, it follows that 
    \begin{align*}
  J_k^T &= (U_k \text{ } Z_k) \left(\begin{array}{cc}\Sigma_{k} & 0 \\  0 & 0\end{array}\right) \left(\begin{array}{c}  V_k  \\ V_k^{\perp} \end{array}\right) = U_k \Sigma_{k} V_k, 
\end{align*}
where $Z_k$ is the orthonormal basis of null space of $J_k $. 
By the linear systems \eqref{eq.nolicq.sub.p2.linear.system.star} and \eqref{eq.nolicq.sub.p2.linear.system.star.det}, 
\begin{align} 
 \bar u_{k}^* &= - \bar Z_{k}\left(\bar Z_{k}^{T} H_{k} \bar Z_{k}\right)^{-1}  \bar Z_{k}^{T}\left(\bar g_{k}+H_{k}\bar v_{k}\right) \label{eq.u.star.expression} \\ 
\text{and}\qquad u_{k}^* &=-Z_{k}\left(Z_{k}^{T} H_{k} Z_{k}\right)^{-1} Z_{k}^{T}\left(g_{k}+H_{k} v_{k}\right),  \label{eq.u.expression} 
\end{align}
where $\bar Z_{k}$ is defined in \eqref{eq.J.bar.SVD}. 
In addition, it follows by Assumption~\ref{ass:prob} that 
\begin{align}
    \| \bar u_{k}^*\| &\le \zeta^{-1} \left( \kappa_g + \epsilon_g +  \kappa_H   \| \bar v_{k}\|  \right) \le \zeta^{-1} \left( \kappa_g + \epsilon_g + \kappa_H \sigma_{Jc}  \| \bar J_{k}^T \bar c_k \|  \right)   \notag \\ &\le
     \zeta^{-1} \left( \kappa_g +\epsilon_g + \kappa_H \sigma_{Jc}   (\kappa_J + \epsilon_J)(\kappa_c+ \epsilon_c)  \right):= \kappa_{\bar u^*} \label{eq.u.sto.star} \\
    \text{and} \qquad \| u_{k}^*\|  &\le \zeta^{-1} \left( \kappa_g + \kappa_H \sigma_{Jc}  \| J_{k}^T  c_k \|  \right)   \le \zeta^{-1} \left( \kappa_g + \kappa_H  \sigma_{Jc}   \kappa_J \kappa_c \right): = \kappa_{u^*}. \label{eq.u}
\end{align}
Given that $\| v_k\| \le \sigma_{Jc} \|J_k^T c_k\| \le \sigma_{Jc} \kappa_J \kappa_c$ (\eqref{eq.nolicq.tr.sub.p1.det} and Assumption~\ref{ass:prob}) and \eqref{eq.u}, 
there exists $\kappa_{d^*} \in \mathbb{R}_{>0}$ such that  $\|d_k^*\|  \le \kappa_{d^*}$. 
Let $\bar V_k^{\dagger} = \bar V_k^T (\bar V_k \bar V_k^T)^{-1} $ denote the Moore-Penrose pseudo-inverse of $\bar V_k$ (defined in \eqref{eq.J.bar.SVD}) and $ V_k^{\dagger} = V_k^T ( V_k V_k^T)^{-1}$. Although the Lagrange multipliers are not unique for the linear systems \eqref{eq.nolicq.sub.p2.linear.system.star} 
and \eqref{eq.nolicq.sub.p2.linear.system.star.det}, 
\begin{align} 
    \bar y_{k}^* &= -\bar V_k^{\dagger} \left(\tilde{\bar J}_k \tilde{\bar J}_k^{T}\right)^{-1}  \tilde{\bar J}_k^T \left(\bar  g_k + H_k (\bar u_{k}^* + \bar v_k) \right),  
 \label{eq.delta.sto.star} \\ 
 \text{and}\qquad y_{k}^* &= - V_k^{\dagger} \left(\tilde{ J}_k \tilde{J}_k^{T}\right)^{-1}  \tilde{J}_k^T \left( g_k + H_k ( u_k^* + v_k) \right)  \label{eq.delta}
\end{align}
are one of the solutions 
to the linear systems \eqref{eq.nolicq.sub.p2.linear.system.star}
and \eqref{eq.nolicq.sub.p2.linear.system.star.det}, respectively. Suppose that $\|\bar V_k^{\dagger}\| \le \kappa_{\bar V^{\dagger}}$ and $\| V_k^{\dagger}\| \le \kappa_{ V^{\dagger}}$, it then follows by Assumptions~\ref{ass:prob} and~\ref{ass:error}, \eqref{eq.u.sto.star}, \eqref{eq.u}, \eqref{eq.delta.sto.star}, and \eqref{eq.delta} that 
\begin{align*}
  \|  \bar y_{k}^* \|&\le \kappa_{\bar V^{\dagger}}   \kappa_{\sigma}^{-2} (\kappa_{J} + \epsilon_J ) ( \kappa_g + \epsilon_g + \kappa_H ( \kappa_{\bar u^*} + \sigma_{Jc} (\kappa_J + \epsilon_J)(\kappa_c+ \epsilon_c) ))  \\
  \text{and}\qquad \|  y_k^* \|&\le  \kappa_{ V^{\dagger}}  \kappa_{\sigma}^{-2} \kappa_{J}  ( \kappa_g + \kappa_H ( \kappa_{u^*} + \sigma_{Jc} \kappa_J \kappa_c )), 
\end{align*}
which concludes the proof.
\end{proof}

Lemma~\ref{lemma.y.exist} provides upper bounds for the tangential 
components, the search directions, and the dual multipliers required in the analysis. The next lemma (Lemma~\ref{lemma.diff.u}) shows that the distance between inexact noisy solution $\bar u_{k}$ and optimal noisy solution $\bar u_{k}^*$ depends on the norm of residuals ($\|\bar \rho_k\|$ and $\|\bar r_k\|$) and the errors in the approximations. Moreover, it provides an upper bound on the distance between the noisy tangential 
component $\bar u_{k}$ in Algorithm~\ref{alg.dfo_sqp_LS} and the true exact tangential  
component 
$u_k^*$.  
{Before we proceed we make an additional technical assumption about the null spaces of the true and noisy constraint Jacobians. Under a similar noise setting to that considered in this paper, in~\cite{sun2024trust}, the authors show that the true stationarity error measure is sufficiently small under the assumption that $\|Z_k - \bar Z_k\|$ is sufficiently small.} 
Instead, we directly analyze the error in the tangential 
components under the following assumption regarding the orthogonal basis of null spaces. 

\begin{assumption}
\label{ass:Z.exist}
    The null spaces of $J_k$ and $\bar J_k$ are 
    nontrivial, i.e., both $Z_k$ and $\bar Z_k$ exist.  
\end{assumption} 
\begin{remark}
The null space of $J_k$ (resp., $\bar J_k$) is trivial if and only if  $n \le m$ and  $\text{rank}( J_k) = n$ (resp.,  $\text{rank}( \bar J_k) = n$).   
    Assumption~\ref{ass:Z.exist} is therefore equivalent to requiring that when $n \le m$, $\text{rank}(\bar{J}_k) < n$ and $\text{rank}(J_k) < n$. This assumption is critical since our stationarity measures involve $Z_k$ (resp., $\bar Z_k$), the null space of $J_k$ (resp., $\bar J_k$). Assumption~\ref{ass:Z.exist} is common in the literature; see, e.g., \cite{CurtNoceWach10,sun2024trust}. Note that \cite{sun2024trust} additionally assumes that $Z_k$ and $\bar{Z}_k$ have the same dimension (otherwise $\|Z_k - \bar Z_k\|$ is not computable), whereas we do not impose this restriction.
\end{remark}

When $\text{rank}(\bar J_k) = \text{rank}(J_k) = n \le m$, the null spaces of $J_k$  and $\bar J_k$ are both trivial. In this case, $u_k^* =0 $ and $\bar u_k^*  = 0$.
We use the following counterexample to show that when the null space of $\bar J_k$ is trivial, the difference between $u_k^*$ and $\bar u_k^*$ can be significant.
\begin{example}
Let $n = m = 2$, and set $H_k = I$. Consider the matrices
\begin{align*}
J_k = \begin{bmatrix} 1.0 & 2.0 \\ 0.5 & 1.0 \end{bmatrix}, \quad  \bar J_k= \begin{bmatrix} 1.0 & 2.0 \\ 0.5 & 1.0 + \epsilon_J \end{bmatrix}
\end{align*}
within the context of equations \eqref{eq.nolicq.sub.p2.linear.system.star} and \eqref{eq.nolicq.sub.p2.linear.system.star.det}. By construction, $J_k$ is rank-deficient, whereas $\bar J_k$ has full rank, with the perturbation satisfying $\|J_k - \bar J_k\| = \epsilon_J$. 
Next, let $g_k + H_kv_k =g_k + v_k = [z_1, z_2]^T$ for some $(z_1, z_2) \in \mathbb{R}^2$. By \eqref{eq.u.expression}, it follows that 
\begin{align*}
u_k^* = \begin{bmatrix} -0.8 z_1 + 0.4 z_2 \\ z_1 + 0.4 z_2 \end{bmatrix}.
\end{align*}
Conversely, since the null space of $\bar J_k$ contains only the zero vector and is therefore independent of $\bar g_k + H_k \bar v_k$, we obtain $\bar u_k^* = 0$. 
Since $u_k^*$ can attain arbitrarily large values, it is evident that the discrepancy between $u_k^*$ and $\bar u_k^*$ can be significant, and cannot be controlled by $\epsilon_J$. 
\end{example}

\begin{lemma}
\label{lemma.diff.u}
Suppose that Assumptions~\ref{ass:prob}, \ref{ass:error}, \ref{ass:H}, \ref{ass.inexact.solver.accuracy}, \ref{ass.d.nonzero}, 
 \ref{ass.noearlytermination}, and \ref{ass:Z.exist} hold.  In addition, suppose that the minimum non-zero singular values of both $J_k$ and $\bar J_k$ are bounded below by $\kappa_{\sigma} \in \R{}_{>0}$ for all $k \in \mathbb{N}$. 
Then, for all $k \in \mathbb{N}$
\begin{align*}
    \|{u_k^*} - \bar u_k \|^2  &\le  6\zeta^{-2} \left(  \epsilon_g^2+ \kappa_{\bar y^*}^2 \epsilon_J^2  + \kappa_{y^*}^2 \epsilon_J^2  \right)  +  2 \lambda_{\rho r}^2  \left(\zeta^{-1}  + (\kappa_{J} + \epsilon_J )   \kappa_{\sigma}^{-2}  \right)^2 \max \{ \|\bar u_k \|^2, \| \bar J_k^T \bar c_k \|^2 \}
  \\ &\qquad +  12\zeta^{-2} \kappa_H^2 \sigma_{Jc}^2 (\|\bar J_k^T \bar c_k\|^2 \mathbbm{1}\{ \|\bar c_k\| >\epsilon_o \} + \|J_k^T c_k\|^2).  
\end{align*}
\end{lemma}
\begin{proof}
\allowdisplaybreaks{
By Assumption~\ref{ass:Z.exist}, both $Z_k$ and $\bar Z_k$ exist. 
It follows from \cite[equation (37)]{curtis2024stochastic} that $\| I -  \bar Z_{k} (\bar Z_{k}^{T} H_{k} \bar Z_{k})^{-1}  \bar Z_{k}^{T} H_k \| \le 1$. 
Since $\bar U_k$, $\bar V_k$ (defined in \eqref{eq.J.bar.SVD}) are both submatrices of orthogonal matrices, $\|\bar U_k\| \le 1$, $\|\bar V_k\| \le 1$.  
Additionally, it follows by $\|(\tilde{\bar J}_k \tilde{\bar J}_k^{T})^{-1}\| =  \| \bar \Sigma_k^{-2}\|  \le \kappa_{\sigma}^{-2}$,  Assumptions~\ref{ass:error} and~\ref{ass:H}, and \eqref{eq.u.bar.expression}, \eqref{eq.u.expression}  
that 
\begin{align*}
    \|\bar u_{k}^* -\bar u_{k}^*  \| &\le  \|\bar Z_{k}(\bar Z_{k}^{T} H_{k} \bar Z_{k})^{-1} \bar Z_{k}^{T} \bar \rho_k \| +  \| (I -  \bar Z_{k}(\bar Z_{k}^{T} H_{k} \bar Z_{k})^{-1}  \bar Z_{k}^{T} H_k) \bar U_k  \bar \Sigma_k (\tilde{\bar J}_k \tilde{\bar J}_k^{T})^{-1}\bar V_k  \bar r_{k}  \|  \\ 
    &\le \zeta^{-1} \|\bar \rho_k \| + (\kappa_{J} + \epsilon_J )  \kappa_{\sigma}^{-2} \|\bar r_k \|.  
\end{align*} 
Let $C_k$ denote the orthonormal basis of $\operatorname{Null}(\bar J_k) \cap \operatorname{Null}( J_k)$. Then, one can construct the orthonormal basis of $\operatorname{Null}(\bar J_k)$ and denote it as $\bar Z_k = [C_k ,\bar Z_k^c] \in \mathbb{R}^{n \times (n-p)}$. Similarly, one can construct the orthonormal basis of $\operatorname{Null}(J_k)$ and denote it as $Z_k = [Z_k^c, C_k ]\in \mathbb{R}^{n \times (n-q)}$. Suppose that $\bar u_{k}^*= \bar Z_k \bar w_k$ and $ u_{k}^* = Z_k  w_k$ where $\bar w_k \in \mathbb{R}^{n-p}$ and $w_k \in \mathbb{R}^{n-q}$, respectively. By \eqref{eq.nolicq.sub.p2.linear.system.star} and \eqref{eq.nolicq.sub.p2.linear.system.star.det}
\begin{align*}
 H_k [Z_k^c  \text{ } C_k   \text{ } \bar Z_k^c ] \left[\begin{array}{c}
0 \\
\bar w_k
\end{array}\right] = H_k [Z_k^c  \text{ } \bar Z_k ] \left[\begin{array}{c}
0 \\
\bar w_k
\end{array}\right]  &=    H_k \bar Z_k \bar w_k   \\
&=  H_k \bar u_{k}^* = -(\bar g_k + H_k \bar v_k + \bar J_k^T \bar y_{k}^*) \\
\text{and} \qquad
 H_k [Z_k^c  \text{ } C_k   \text{ } \bar Z_k^c ] \left[\begin{array}{c}
 w_k \\
0
\end{array}\right]= H_k [Z_k  \text{ } \bar Z_k^c ] \left[\begin{array}{c}
 w_k \\
0
\end{array}\right]  & =  H_k  Z_k  w_k  \\
    & =   H_k   u_{k}^*  = -( g_k + H_k v_k + J_k^T y_k^*). 
\end{align*}
Multiplying $ [Z_k^c  \text{ } C_k   \text{ } \bar Z_k^c]^T$ to both sides, it follows from $J_k Z_k = 0$ and $\bar J_k  \bar Z_k = 0$ that 
\begin{align}
     [Z_k^c  \text{ } C_k   \text{ } \bar Z_k^c ]^T H_k [Z_k^c  \text{ } C_k   \text{ } \bar Z_k^c ] \left[\begin{array}{c}
0 \\
\bar w_k
\end{array}\right]  &=  -[Z_k^c  \text{ } C_k   \text{ } \bar Z_k^c ]^T (\bar g_k + H_k \bar v_k + \bar J_k^T \bar y_{k}^*) \notag \\ 
&=  -\left[\begin{array}{c}
Z_k^T(\bar g_k + H_k \bar v_k + ( \bar J_k -  J_k )^T   \bar y_{k}^*  )   \\
\bar{Z_k^c}^T (\bar g_k + H_k \bar v_k)
\end{array}\right],
\label{eq.multiplying.long1}
\\
\text{and} \qquad  [Z_k^c  \text{ } C_k   \text{ } \bar Z_k^c ]^T  H_k [Z_k^c  \text{ } C_k   \text{ } \bar Z_k^c ] \left[\begin{array}{c}
 w_k \\
0
\end{array}\right] &=  -[Z_k^c  \text{ } C_k   \text{ } \bar Z_k^c ]^T ( g_k + H_k v_k + J_k^T y_k^* ) \notag \\ 
&= -\left[\begin{array}{c}
Z_k^T( g_k + H_k  v_k) \\
\bar{Z_k^c}^T (g_k + H_k v_k + (\bar J_k - J_k)^T  y_k^*)
\end{array}\right]. \label{eq.multiplying.long2}
\end{align}
Since $[Z_k^c  \text{ } C_k   \text{ } \bar Z_k^c ]$ is an orthonormal basis of $\operatorname{Null}(J_k) \cup \operatorname{Null}(\bar J_k)$, it follows by Assumption~\ref{ass:H} that $\|(  [Z_k^c  \text{ } C_k   \text{ } \bar Z_k^c ]^T  H_k [Z_k^c  \text{ } C_k   \text{ } \bar Z_k^c ]  )^{-1} \| \le \zeta^{-1}$. Thus, 
 \begin{align}
   &\ \|u_k^* - \bar u_{k}^* \|^2   = \left \| Z_k w_k - \bar Z_k \bar  w_k \right\|^2 \notag \\  = &\ \left \| [Z_k^c  \text{ } C_k   \text{ } \bar Z_k^c ] \left[\begin{array}{c} 
 w_k \\
0 
\end{array}\right] -  [Z_k^c  \text{ } C_k   \text{ } \bar Z_k^c ] \left[\begin{array}{c}
 0 \\
\bar w_k
\end{array}\right]  \right\|^2 \notag \\
\le  &\ \left \| [Z_k^c  \text{ } C_k   \text{ } \bar Z_k^c ] \left (  [Z_k^c  \text{ } C_k   \text{ } \bar Z_k^c ]^T  H_k [Z_k^c  \text{ } C_k   \text{ } \bar Z_k^c ]    \right)^{-1} \left[\begin{array}{c}
Z_k^T(\bar g_k + H_k \bar v_k + ( \bar J_k -  J_k )^T   \bar y_{k}^*  )   \\
\bar{Z_k^c}^T (\bar g_k + H_k \bar v_k)
\end{array}\right]
\right. \notag \\ &\ - 
\left.  [Z_k^c  \text{ } C_k   \text{ } \bar Z_k^c ] \left (  [Z_k^c  \text{ } C_k   \text{ } \bar Z_k^c ]^T  H_k [Z_k^c  \text{ } C_k   \text{ } \bar Z_k^c ]    \right)^{-1}  \left[\begin{array}{c}
Z_k^T( g_k + H_k  v_k) \\
\bar{Z_k^c}^T (g_k + H_k v_k + (\bar J_k - J_k)^T  y_k^*)
\end{array}\right] 
\right\|^2\notag \\
=  &\ \left \| [Z_k^c  \text{ } C_k   \text{ } \bar Z_k^c ] \left (  [Z_k^c  \text{ } C_k   \text{ } \bar Z_k^c ]^T  H_k [Z_k^c  \text{ } C_k   \text{ } \bar Z_k^c ]    \right)^{-1} 
\left[\begin{array}{c}
Z_k^T  \\
\bar{Z_k^c}^T 
\end{array}\right] \left( (\bar g_k - g_k)  +  H_k (\bar v_k - v_k )  \right) \right.  \notag \\ 
&\ + 
\left.  [Z_k^c  \text{ } C_k   \text{ } \bar Z_k^c ]  \left(  [Z_k^c  \text{ } C_k  \text{ } \bar Z_k^c ]^T  H_k [Z_k^c  \text{ } C_k   \text{ } \bar Z_k^c ]    \right)^{-1}  \left[\begin{array}{c}
Z_k^T \left( ( \bar J_k -  J_k )^T   \bar y_{k}^* \right) \\ 
-\bar{Z_k^c}^T \left((\bar J_k - J_k)^T  y_k^*  \right)
\end{array}\right] 
\right\|^2 \notag \\
\le  &\ 3\zeta^{-2} \left(  \| \bar g_k - g_k \|^2 + \| H_k (\bar v_k - v_k ) \|^2 + \left \| \left[\begin{array}{c}
Z_k^T \left( ( \bar J_k -  J_k )^T   \bar y_{k}^* \right) \\ 
-\bar{Z_k^c}^T \left((\bar J_k - J_k)^T  y_k^*  \right)
\end{array}\right] 
\right\|^2  \right)    \notag   \\ 
\le  &\  3\zeta^{-2} \left(\epsilon_g^2 + \kappa_H^2 \|\bar v_k - v_k\|^2  + \kappa_{\bar y^*}^2 \epsilon_J^2  + \kappa_{y^*}^2 \epsilon_J^2   \right)  \notag   \\ 
\le  &\  3\zeta^{-2} \left(  \epsilon_g^2+ \kappa_{\bar y^*}^2 \epsilon_J^2  + \kappa_{y^*}^2 \epsilon_J^2  \right)  
 +  6\zeta^{-2} \kappa_H^2 \sigma_{Jc}^2 (\|\bar J_k^T \bar c_k\|^2\mathbbm{1}\{ \|\bar c_k\| >\epsilon_o \} + \|J_k^T c_k\|^2),  
 \label{eq.diff.u.u.star}
\end{align}
where the first inequality follows multiplying $(  [Z_k^c  \text{ } C_k   \text{ } \bar Z_k^c ]^T  H_k [Z_k^c  \text{ } C_k   \text{ } \bar Z_k^c ]  )^{-1}$ to both sides of  \eqref{eq.multiplying.long1} and \eqref{eq.multiplying.long2}, the second  inequality follows by the 
orthonormality of $[Z_k^c  \text{ } C_k   \text{ } \bar Z_k^c ]$, Assumption~\ref{ass:H}, and the fact that $\|a_1+a_2+a_3\|^2 \le 3(\|a_1\|^2 + \|a_2\|^2 + \|a_3\|^2)$, 
 the third inequality follows by 
 Assumptions~\ref{ass:error} and~\ref{ass:H}, 
Lemma~\ref{lemma.y.exist}, and the last inequality follows by \eqref{eq.nolicq.tr.sub.p1} and \eqref{eq.nolicq.tr.sub.p1.det}. It then follows by Young's inequality that 
\begin{align*}
    \|u_k^* - \bar u_k \|^2 &\le 2\|u_k^* - \bar u_{k}^* \|^2 + 
2\|\bar u_{k}^* - \bar u_k \|^2  
\\ &\le  6\zeta^{-2} \left(  \epsilon_g^2+ \kappa_{\bar y^*}^2 \epsilon_J^2  + \kappa_{y^*}^2 \epsilon_J^2  \right) 
 +  12\zeta^{-2} \kappa_H^2 \sigma_{Jc}^2 (\|\bar J_k^T \bar c_k\|^2\mathbbm{1}\{ \|\bar c_k\| >\epsilon_o \} + \|J_k^T c_k\|^2) 
 \\ &\qquad + 2\left(\zeta^{-1} \|\bar \rho_k \| + (\kappa_{J} + \epsilon_J )   \kappa_{\sigma}^{-2} \|\bar r_k \| \right)^2.  
\end{align*}
The desired conclusion (second statement) then follows by \eqref{eq.tt1.cond1}. 
}
\end{proof}

\begin{remark}
    Lemma~\ref{lemma.diff.u} shows that the upper bound of $\|u_k^* - \bar u_k \|^2$ depends on noise levels $\epsilon_g, \epsilon_J$, residual norm of problem \eqref{eq.nolicq.sub.p2.linear.system.res} (which is bounded by $\max \{ \|\bar u_k \|^2, \| \bar J_k^T \bar c_k \|^2 \}$), and the infeasible stationarity errors $\|\bar J_k^T \bar c_k\|^2$  and  $\|J_k^T c_k\|^2 $. Note that by 
    \eqref{eq.diff.u.u.star} 
    the dependence on $\|\bar J_k^T \bar c_k\|^2 + \|J_k^T c_k\|^2 $ can be replaced by $\| \bar v_k - v_k \|^2$. When $\epsilon_c = \epsilon_J = \epsilon_g = 0$, $\bar v_k$ and $v_k$ are either both Cauchy or both exact solutions to subproblems \eqref{eq.nolicq.tr.sub.p1} and \eqref{eq.nolicq.tr.sub.p1.det}, respectively, then $\bar u_k=u_k^*$. 
    However, our proof does not require that $\|u_k^* - \bar u_k \|^2$ is sufficiently small at each iteration. Instead, $\|u_k^* - \bar u_k \|^2$ is sufficiently small once both the stationarity ($\| \bar u_k\|$) and infeasible stationarity  ($\|\bar J_k^T \bar c_k\|$) errors are sufficiently small. 
\end{remark}


Given the fictitious (never computed in Algorithm~\ref{alg.dfo_sqp_LS}) normal solution $v_k$ and tangential solution $u_k^*$ for the true SQP subproblems \eqref{eq.nolicq.tr.sub.p1.det} and \eqref{eq.nolicq.sub.p2.linear.system.star.det}, respectively, we define a fictitious (never computed in Algorithm~\ref{alg.dfo_sqp_LS}) merit parameter $\tau_k$ update rule for the deterministic counterpart of our proposed algorithm, 
\begin{equation}
    \label{def.tau.det}
    \tau_k \leftarrow \begin{cases}{\bar \tau}_{k-1} & \text { if } \bar {\tau}_{k-1} \leq\left(1-\sigma_\tau\right) {\tau}_k^{\text {trial }} 
\\ \left(1-\sigma_\tau\right) {\tau}_k^{\text {trial }} & \text { otherwise},\end{cases}
\end{equation}
where
\begin{equation}
 \label{def.tau.trial.det}
{\tau}_k^{\text {trial }} \leftarrow \begin{cases}\infty & \text { if } {g}_k^T {d}_{k}^*+\max \left\{u_k^{*T} H_k u_k^*, \lambda_u\left\|{u}_k^*\right\|^2\right\} \leq 0  \\ \tfrac{\left(1-\frac{\sigma_{c}}{\sigma_{r}} \right)(\left\|{c}_k\right\| - \left\|{c}_k + J_k v_k \right\|)}{{g}_k^T {d}_k^*+\max \left\{u_k^{*T} H_k u_k^*, \lambda_u\left\|{u}_k^*\right\|^2\right\}} & \text { otherwise. }\end{cases}
\end{equation}
The next lemma characterizes the relationship between noisy and true merit parameters. Note that it requires $\bar u_k$ to be the exact solution to subproblem  \eqref{eq.nolicq.sub.p2.linear.system.star} to maintain consistency with the exactness level of $u_k^*$, which restricts $\lambda_{\rho r} = 0 $ in \eqref{eq.tt1.cond1}. 




\begin{lemma}
\label{lemma.diff.tau}
Suppose that Assumptions~\ref{ass:prob}, \ref{ass:error}, \ref{ass:H}, 
\ref{ass.inexact.solver.accuracy}, \ref{ass.d.nonzero}, 
 \ref{ass.noearlytermination}, and \ref{ass:Z.exist}  hold. 
In addition, suppose that $\bar v_k$ and $v_k$ are either both Cauchy solutions or both exact solutions to subproblems \eqref{eq.nolicq.tr.sub.p1} and \eqref{eq.nolicq.tr.sub.p1.det}, respectively, and that  
$\bar u_k$ and $u_k^*$ are the exact solutions to subproblems \eqref{eq.nolicq.sub.p2.linear.system.star} and \eqref{eq.nolicq.sub.p2.linear.system.star.det}, respectively, i.e., $\bar u_k = \bar u_k^*$. 
Then, there exists $ \mathcal{E}_{\tau} \in \mathbb{R}_{\ge 0}$ such that for all $k \in \N{}$, $| (\bar{\tau}_k-\tau_k) (g_k^T d_k^* + u_k^{*T} H_k u_k^*) | \le \mathcal{E}_{\tau}$.   
\end{lemma}

\begin{proof}
We first consider the case when $\|\bar c_k\| \le  \epsilon_o$, it follows by Assumptions~\ref{ass:prob}, \ref{ass:error}, and~\ref{ass:H}, Lemma~\ref{lemma.y.exist}, \eqref{eq.nolicq.tr.sub.p1.det}, and \eqref{eq.nolicq.sub.p2.linear.system.star.det} that 
\begin{align*}
   \left| g_k^T d_k^* + u_k^{*T} H_k u_k^*  \right| &=  \left| g_k^T u_k^*  + u_k^{*T} H_k u_k^*   +  g_k^T v_{k} \right| = \left|  -{u}_k^{*T}H_k  v_{k}  +  g_k^T v_{k} \right| \\
   &\le ( \kappa_{ u^*} \kappa_H    +  \kappa_g ) \sigma_{Jc} \kappa_J  \| c_k\|
    \le 2( \kappa_{ u^*} \kappa_H    +  \kappa_g ) \sigma_{Jc} \kappa_J  \epsilon_c: = \mathcal{E}_c,
\end{align*}
which implies that  $| (\bar{\tau}_k-\tau_k) (g_k^T d_k^* + u_k^{*T} H_k u_k^*) | \le 2 
\bar\tau_{-1} \mathcal{E}_c$ by Lemma~\ref{lemma.tau.lb}.

Next, we consider the case when $\|\bar c_k\| > \epsilon_o$. 
Since we allow $\bar v_k$ and $v_k$ to be either both Cauchy or both exact solutions to subproblems \eqref{eq.nolicq.tr.sub.p1} and \eqref{eq.nolicq.tr.sub.p1.det}, respectively, it follows that $\|v_k - \bar v_k\| \le \max \{ \mathcal{E}_v, \mathcal{E}_{v^*} \}$, where $\mathcal{E}_v$ and $\mathcal{E}_{v^*}$ are defined in 
Lemmas~\ref{lemma.diff.v} and~\ref{lemma.diff.v.exact}, respectively. 
The difference between the exact noisy tangential 
component $\bar u_k = \bar u_k^*$ and the true tangential 
component $u_k^*$ can be derived by \eqref{eq.diff.u.u.star},
 \begin{align*}
   &\ \|u_k^* - \bar u_{k} \|   \\
   = &\left \| Z_k w_k - \bar Z_k \bar  w_k \right\| \\
   \le  &\ \left \| [Z_k^c  \text{ } C_k   \text{ } \bar Z_k^c ] \left (  [Z_k^c  \text{ } C_k   \text{ } \bar Z_k^c ]^T  H_k [Z_k^c  \text{ } C_k   \text{ } \bar Z_k^c ]    \right)^{-1} 
\left[\begin{array}{c}
Z_k^T  \\
\bar{Z_k^c}^T 
\end{array}\right] \left( (\bar g_k - g_k)  +  H_k (\bar v_k - v_k )  \right) \right.  \notag \\ 
&\ + 
\left.  [Z_k^c  \text{ } C_k   \text{ } \bar Z_k^c ]  \left(  [Z_k^c  \text{ } C_k  \text{ } \bar Z_k^c ]^T  H_k [Z_k^c  \text{ } C_k   \text{ } \bar Z_k^c ]    \right)^{-1}  \left[\begin{array}{c}
Z_k^T \left( ( \bar J_k -  J_k )^T   \bar y_{k}^* \right) \\ 
-\bar{Z_k^c}^T \left((\bar J_k - J_k)^T  y_k^*  \right)
\end{array}\right] 
\right\| \notag \\
\le  &\  \zeta^{-1}  \left( (\kappa_{y^*} + \kappa_{\bar y^*}) \epsilon_J+ \epsilon_g + \kappa_H \max \{ \mathcal{E}_v, \mathcal{E}_{v^*} \} \right)  \notag :=   \mathcal{E}_u.  
 \end{align*}

We divide the proof into four cases that depend on whether $\left(\tau_k, \bar{\tau}_k\right)$ is updated 
conditioned on $\bar \tau_{k-1}$. 
For conciseness, we define $h_k = g_k^T d_k^*+\max \left\{u_k^{*T} H_k u_k^*, \lambda_u\left\|{u}_k^*\right\|^2\right\}$ and $\bar h_k =\bar{g}_k^T  \bar d_k+\max \left\{\bar{u}_k^T H_k \bar{u}_k, \lambda_u\left\|\bar{u}_k\right\|^2\right\}$. We begin with a general bound that holds in all cases. By Assumptions~\ref{ass:prob}, \ref{ass:error} and \ref{ass:H}, and the fact that $d_k^* = u_k^* + v_k$ and $\bar d_k = \bar u_k + \bar v_k$, 
\begin{align*}
    |h_k - \bar h_k| &= \left|{g}_k^T d_k^*+\max \left\{u_k^{*T} H_k u_k^*, \lambda_u\left\|{u}_k^*\right\|^2\right\} - \bar{g}_k^T  \bar d_k -  \max \left\{\bar{u}_k^T H_k \bar{u}_k, \lambda_u\left\|\bar{u}_k\right\|^2\right\} \right| \\ 
    &\le \left|{g}_k^T d_k^*+u_k^{*T} H_k u_k^*  - \bar{g}_k^T  \bar d_k -  \bar{u}_k^T H_k \bar{u}_k  \right| \\ 
    &\le \left|{g}_k^T {d}_k^* - {g}_k^T\bar {d}_k  + {g}_k^T\bar d_k  - \bar{g}_k^T  \bar {d}_k  +u_k^{*T} H_k u_k^*  -  \bar u_k^{*T} H_k u_k^* +\bar u_k^{*T} H_k u_k^*   -  \bar{u}_k^T H_k \bar{u}_k  \right| \\ 
     &\le \left|{g}_k^T {d}_k^* - {g}_k^T\bar {d}_k\right|  +\left|{g}_k^T\bar {d}_k  - \bar{g}_k^T  \bar {d}_k \right| +\left|u_k^{*T} H_k u_k^*  -  \bar u_k^{*T} H_k u_k^* \right| +\left| \bar u_k^{*T} H_k u_k^*   -  \bar{u}_k^T H_k \bar{u}_k  \right| \\ 
     &\le \kappa_g (\mathcal{E}_v +\mathcal{E}_u)  + \kappa_{\bar d} \epsilon_g + \kappa_H (\kappa_u + \kappa_{\bar u^*}) \mathcal{E}_u : = \mathcal{E}_h.  
\end{align*}
Additionally, when $\bar v_k$ and $v_k$ are either both Cauchy solutions or both exact solutions
to subproblems \eqref{eq.nolicq.tr.sub.p1} and \eqref{eq.nolicq.tr.sub.p1.det}, respectively, it follows by Lemmas~\ref{lemma.v.bound} and~\ref{lemma.diff.v}, and Assumption~\ref{ass:error}
\begin{align}
\label{eq.diff.c+Jv}
        \left\|\bar 
c_{k}+\bar J_{k} \bar v_{k} \right\| - \left\|{c}_k + J_k v_k \right\|  &\le \left\|{c}_k + J_k v_k  - \bar  c_{k} - \bar J_{k} \bar v_{k} \right\|  \notag \\ 
    &\le \left\|{c}_k  -  \bar  c_{k}  \right\| + \left\|   J_k v_k  - J_k  \bar v_{k} \right\|  + \left\| J_k  \bar v_{k}  -  \bar J_{k} \bar v_{k} \right\| \notag\\ 
    &\le \epsilon_c + \kappa_J \max \{ \mathcal{E}_v, \mathcal{E}_{v^*} \}  +  \sigma_{Jc} (\kappa_J + \epsilon_J ) (\kappa_c + \epsilon_c)\epsilon_J := \mathcal{E}_{cJv}.
\end{align}

Case $(i)$: Neither $\bar{\tau}_k$ nor $\tau_k$ is updated, i.e., $\tau_k=\bar{\tau}_k=\bar{\tau}_{k-1}$. 
The result holds in case $(i)$. 

Case $(ii)$: Both $\bar{\tau}_k$ and $\tau_k$ are updated. 
 By the update rule for $\tau_k$ and $\bar\tau_k$, we have $h_k > 0$, $\bar h_k > 0$, and  $(1- \sigma_{\tau}) \left(1-\sigma_{c}/\sigma_{r} \right)\left\| \bar {c}_k\right\| < \bar \tau_{k-1} \bar  h_k \le \bar \tau_{-1} \bar h_k$. 
 It then follows by $\bar{u}_k = \bar{u}_k^*$, 
 \eqref{eq.tau.trial.tt2}, \eqref{def.tau.trial.det}, and \eqref{eq.diff.c+Jv} that 
\begin{align*}
    \left|\bar{\tau}_k-\tau_k\right| & = (1- \sigma_{\tau})\left|\bar{\tau}_k^{trial}-\tau_k^{trial}\right| 
    \\ & = (1- \sigma_{\tau}) \left| \tfrac{\left(1-\sigma_{c}/\sigma_{r} \right)(\left\|{c}_k\right\| - \left\|{c}_k + J_k v_k \right\|)}{h_k} - \tfrac{\left(1-\sigma_{c}/\sigma_{r}\right)(\left\|\bar c_{k}\right\|-\|\bar c_{k}+\bar J_{k} \bar v_{k}\|)}{\bar h_k} \right| \\ 
    & = (1- \sigma_{\tau})\left(1-\tfrac{\sigma_{c}}{\sigma_{r}} \right)  \left| \tfrac{(\left\|{c}_k\right\|- \left\|{c}_k + J_k v_k \right\| - \left\|\bar c_{k}\right\| + \|\bar c_{k}+\bar J_{k} \bar v_{k}\|  ) \bar h_k + (\left\|\bar c_{k}\right\| - \|\bar c_{k}+\bar J_{k} \bar v_{k}\|) (\bar h_k - h_k)}{h_k \bar h_k} \right|   \\ 
        & \leq (1- \sigma_{\tau})\left(1-\tfrac{\sigma_{c}}{\sigma_{r}} \right)  \left| \tfrac{\left\|\bar{c}_k\right\| \mathcal{E}_h + \bar h_k (\epsilon_c + \mathcal{E}_{cJv})}{h_k \bar h_k} \right| \\ 
  & \leq \tfrac{\bar \tau_{-1} \mathcal{E}_h }{ h_k}  +  (1- \sigma_{\tau})\left(1-\tfrac{\sigma_{c}}{\sigma_{r}} \right) \tfrac{\epsilon_c + \mathcal{E}_{cJv}}{ h_k}. 
\end{align*}
Thus, in case $(ii)$, $| (\bar{\tau}_k-\tau_k) h_k | \le \bar \tau_{-1} \mathcal{E}_h  + (1- \sigma_{\tau})\left(1-\tfrac{\sigma_{c}}{\sigma_{r}} \right) (\epsilon_c + \mathcal{E}_{cJv})$.


Case $(iii)$: $\bar{\tau}_k$ is updated and $\tau_k$ is not, i.e., $\bar{\tau}_k 
 = \bar{\tau}_k^{\text {trial }} <  \bar{\tau}_{k-1}$ and $ \tau_k=\bar{\tau}_{k-1} \le (1-\sigma_{\tau}) \bar{\tau}_k^{\text {trial }}$.
 By the update rule of $\bar\tau_k$, we have $\bar h_k > 0$. Let us first discuss the case where  $\bar h_k > \mathcal{E}_h$, which implies that $h_k > 0$ and ${\tau}_k^{\text {trial }} < \infty$. 
 It then follows by~\eqref{def.tau.trial.det}
\begin{align*}
    \left|\bar{\tau}_k-\tau_k\right| & = \tau_{k} - (1- \sigma_{\tau}) \left(1 - \tfrac{\sigma_{c}}{\sigma_{r}}\right) \tfrac{\left\|\bar c_{k}\right\|-\|\bar c_{k}+\bar J_{k} \bar v_{k}+\bar r_{k}\|}{\bar h_k} \\ 
    & \le (1- \sigma_{\tau})\left(1-\tfrac{\sigma_{c}}{\sigma_{r}} \right)  \left( \tfrac{\left\|{c}_k\right\|- \left\|{c}_k + J_k v_k \right\|}{h_k} - \tfrac{\|\bar c_{k}\| - \|\bar 
c_{k}+\bar J_{k} \bar v_{k} \| }{\bar h_k}  \right), 
\end{align*}
from which we conclude that the result from 
case $(ii)$ holds. When $\bar h_k  \le \mathcal{E}_h$, it follows that 
 $ h_k  \le 2\mathcal{E}_h$ and therefore $| (\bar{\tau}_k-\tau_k) h_k | \le 4\bar\tau_{-1} \mathcal{E}_h$.

Case $(iv)$: $\tau_k$ is updated and $\bar{\tau}_k$ is not, i.e., $\tau_k \leq \bar{\tau}_k=\bar{\tau}_{k-1}$ and $\bar\tau_k = \bar\tau_{k=1}$. 
By the update rule of $\tau_k$, we have $h_k > 0$. Let us first discuss the case when $ h_k > \mathcal{E}_h$, which implies that $\bar h_k > 0$ and $ \bar {\tau}_k^{\text {trial }} < \infty$. 
 It then follows by~\eqref{def.tau.trial.det}  
\begin{align*}
    \left|\bar{\tau}_k-\tau_k\right| & = \bar \tau_{k-1} - (1- \sigma_{\tau}) \left(1 - \tfrac{\sigma_{c}}{\sigma_{r}}\right) \tfrac{\left\|{c}_k\right\| - \left\|{c}_k + J_k v_k \right\|}{h_k}\\ 
    & \le (1- \sigma_{\tau})\left(1-\tfrac{\sigma_{c}}{\sigma_{r}} \right)  \left(\tfrac{\left\|\bar c_{k}\right\| - \left\|\bar 
c_{k}+\bar J_{k} \bar v_{k} \right\| }{\bar h_k}  - \tfrac{\left\|{c}_k\right\|- \left\|{c}_k + J_k v_k \right\|}{h_k}   \right), 
\end{align*}
from which we conclude that the result from 
case $(ii)$ holds.
When $h_k  \le \mathcal{E}_h$, we have $| (\bar{\tau}_k-\tau_k) h_k | \le 2\bar\tau_{-1} \mathcal{E}_h$.

Finally, defining $\mathcal{E}_\tau : = \max\left\{ \bar \tau_{-1} \mathcal{E}_h  + (1- \sigma_{\tau})\left(1-\tfrac{\sigma_{c}}{\sigma_{r}} \right) (\epsilon_c + \mathcal{E}_{cJv}),   4\bar\tau_{-1} \mathcal{E}_h,   2\bar\tau_{-1} \mathcal{E}_c \right\} $ completes the proof.
 \end{proof}
 
\begin{remark}
    In Lemma~\ref{lemma.diff.tau}, we bound $| (\bar{\tau}_k-\tau_k) (g_k^T d_k^* + u_k^{*T} H_k u_k^*) |$ by quantities related to the noise in the objective and constraint functions. When $\epsilon_g = \epsilon_c = \epsilon_J = 0$, we have $\mathcal{E}_v = \mathcal{E}_u = \mathcal{E}_{h} = \mathcal{E}_{cJv} =\mathcal{E}_{c} = 0$ and the upper bound is zero ($\bar{\tau}_k = \tau_k$). When $g_k^T d_k^* + u_k^{*T} H_k u_k^*$ is sufficiently small, i.e., the algorithm has reached an approximate first-order stationary point, the bound becomes useless. However, when $g_k^T d_k^* + u_k^{*T} H_k u_k^*$ is sufficiently large, the difference between $\tau_k$ and $\bar{\tau}_k$ is sufficiently small and proportional to the noise levels in the problem. 
    Lemma~\ref{lemma.diff.tau} imposes conditions on the normal and tangential components, i.e., the fact that the exactness of both the noisy normal component and the noisy tangential component must be of the same level as the true components. It is important to note that we do not require such conditions in the main convergence analysis (Section~\ref{sec:convergence.main}).  
\end{remark}

Careful consideration of the behavior of the merit parameter $\bar \tau_k$ is essential as it plays a pivotal role in both the formulation and analysis of the proposed noisy SQP method. It is shown in \cite[Lemma 3.16]{CurtNoceWach10} that when the problem has no noise and the LICQ holds, the merit parameter sequence is guaranteed to be bounded away from zero. Next, we discuss the deviation in the lower bound of the merit parameter sequence ($\{\bar \tau_k\}$) from a fictitious lower bound of a true merit parameter sequence ($\{\tau_k\}$, never computed) under the LICQ.

\begin{lemma}
Suppose that Assumptions~\ref{ass:prob}, \ref{ass:error}, \ref{ass:H}, 
\ref{ass.inexact.solver.accuracy}, \ref{ass.d.nonzero}, 
and \ref{ass.noearlytermination} hold. In addition, suppose that the singular values of both $\bar J_k$ and $J_k$  are bounded below by $\kappa_{\sigma} \in \R{}_{>0}$ for all $k \in \mathbb{N}$. Furthermore, suppose that $\bar v_k$ and $v_k$ are the exact solutions of \eqref{eq.nolicq.tr.sub.p1} and \eqref{eq.nolicq.tr.sub.p1.det}, respectively, and that $\bar u_k$ and $u_k^*$ are exact solutions to subproblems \eqref{eq.nolicq.sub.p2.linear.system.star} and \eqref{eq.nolicq.sub.p2.linear.system.star.det}, respectively, i.e., $\left\|\bar c_{k}+\bar J_{k} \bar v_{k}\right\| = 0$, $\|\bar r_{k}\| = \|\bar \rho_{k}\| = 0$, and $\bar u_k = \bar u_k^*$. Then, there exist lower bounds $\{\tau_{\min}$, $\tilde \tau_{\min}\}\subset\mathbb{R}_{>0}$ such that 
$\tau_k \ge  \tau_{\min}$, $\bar \tau_k \ge \tilde \tau_{\min}$, and $\mathcal{E}_{\tau,\min} \in \mathbb{R}_{\ge 0}$ such that $| \tau_{\min} - \tilde \tau_{\min}| \le \mathcal{E}_{\tau,\min}$. 
\label{lemma.diff.tau.min}
\end{lemma}
\begin{proof}
Since $\|\bar r_{k}\| = \|\bar \rho_{k}\| = 0$, we have $\bar u_k  = \bar u_k^*$. Note that $\bar \tau_k < \bar \tau_{k-1}$ only when $\|\bar c_k \| > \epsilon_o$, $\bar g_{k}^{T}\bar d_{k}+\max \left\{\bar u_{k}^{T} H_{k}\bar u_{k}, \lambda_u\left\|\bar u_{k}\right\|^{2}\right\}>0$, and $\bar \tau_{k-1}  \ge  \tfrac{ (1- \sigma_{\tau}) \left( 1 - \sigma_{c}/\sigma_r \right)\left\|\bar c_{k}\right\| }{\bar g_{k}^{T}\bar d_{k}+\max \left\{\bar u_{k}^{*T} H_{k}\bar u_{k}^*, \lambda_u\left\|\bar u_{k}^*\right\|^{2}\right\} }$. 
Hence, if the merit parameter is ever updated, it follows by \eqref{eq.tau.tilde} that 
\begin{align*}
    \bar \tau_{k-1}  
    &\ge \tfrac{(1- \sigma_{\tau})\left(\sigma_{r}-\sigma_{c}\right) \sigma_v \| \bar J_k^T \bar c_k\|^2  / \|\bar c_k\|}{\bar g_{k}^{T}\bar u_{k}^* + \bar g_{k}^{T}\bar v_{k} +\max \left\{\bar u_{k}^{*T} H_{k}\bar u_{k}^*, \lambda_u\left\|\bar u_{k}^*\right\|^{2}\right\}} \\
    &\ge \tfrac{(1- \sigma_{\tau})\left( 1-\sigma_{c}/\sigma_r \right)\sigma_v \| \bar J_k^T \bar c_k\|^2  / \|\bar c_k\| }{ -\bar u_{k}^{*T} H_k \bar v_{k} + \bar g_{k}^{T}\bar v_{k} } \\
        &\ge \tfrac{(1- \sigma_{\tau})\left( 1-\sigma_{c}/\sigma_r \right) \sigma_v \| \bar J_k^T \bar c_k\|^2  / \|\bar c_k\| }{ \kappa_{\bar u^*} \kappa_H \sigma_{Jc} (\kappa_J  + \epsilon_J) \|\bar c_k\| +  (\kappa_g + \epsilon_g)  \sigma_{Jc} (\kappa_J  + \epsilon_J) \|\bar c_k\| } \\
    &= \tfrac{(1- \sigma_{\tau})\left( 1-\sigma_{c}/\sigma_r \right) \sigma_v \kappa_{\sigma}^{2}}{  \kappa_{\bar u^*} \kappa_H \sigma_{Jc} (\kappa_J  + \epsilon_J) + (\kappa_g + \epsilon_g)  \sigma_{Jc} (\kappa_J  + \epsilon_J)}:= \tilde \tau_{\min}. 
\end{align*} 

As for deterministic SQP method (no noise), it follows from \cite[Lemmas 3.5 \& 3.6]{CurtNoceWach10} that there exists $\sigma_v \in \mathbb{R}_{>0}$ such that 
\begin{align}
\label{eq.constraint.reduction.det}
     \left\|c_k\right\|\left(\left\|c_k\right\|-\left\|c_k+J_k v_k\right\|\right) &\geq \sigma_v\left\|J_k^T c_k\right\|^2. 
\end{align}
Although $\sigma_v$ is not necessarily the same constant as in Lemma \ref{lemma.v.bound}, this discrepancy does not result in a loss of generality. Specifically, $\sigma_v$ can be defined as the minimum of the two constants appearing in \eqref{eq.constraint.reduction.det} and its noisy counterpart. 
Whenever $\tau_k$ is updated, we have 
$ g_{k}^{T} d_{k}^*+\max \left\{ u_{k}^{*T} H_{k} u_{k}^*, \lambda_u\left\| u_{k}^*\right\|^{2}\right\}>0$, from which it follows that 
\begin{align*}
     \tau_{k-1}  
    &\ge \tfrac{(1- \sigma_{\tau})\left( 1-\sigma_{c}/\sigma_r \right) \left(\left\| c_{k}\right\|-\left\| c_{k}+ J_{k}  v_{k} \right\| \right) }{ g_{k}^{T} u_{k}^* +  g_{k}^{T} v_{k} +\max \left\{u_{k}^{*T} H_{k} u_{k}^*, \lambda_u\left\| u_{k}^* \right\|^{2}\right\}} \\
    &\ge \tfrac{(1- \sigma_{\tau})\left( 1-\sigma_{c}/\sigma_r \right) \sigma_v \|  J_k^T  c_k\|^2  / \| c_k\|  }{ {u}_k^{*T}H_k  v_{k} +  g_{k}^{T} v_{k} } \\
        &\ge \tfrac{(1- \sigma_{\tau})\left( 1-\sigma_{c}/\sigma_r \right) \sigma_v \kappa_{\sigma}^2 \left\| c_{k}\right\| }{ \kappa_{ u^*} \kappa_H \sigma_{Jc} \kappa_J  \| c_k\| +  \kappa_g  \sigma_{Jc} \kappa_J  \| c_k\| } \\
    &= \tfrac{(1- \sigma_{\tau})\left( 1-\sigma_{c}/\sigma_r \right) \sigma_v \kappa_{\sigma}^2  }{ \kappa_{ u^*} \kappa_H \sigma_{Jc} \kappa_J  +  \kappa_g  \sigma_{Jc} \kappa_J  }:=  \tau_{\min}.
\end{align*} 
Given \eqref{eq.u.sto.star} and \eqref{eq.u}, it follows by Assumption~\ref{ass:error} that 
\begin{align*}
 \kappa_{\bar u^*}  -  \kappa_{u^*} = \zeta^{-1} \left( \epsilon_g + \kappa_H \sigma_{Jc} (\kappa_J \epsilon_c + \kappa_c \epsilon_J + \epsilon_c  \epsilon_J  ) \right):= \epsilon_{u^*}.   
\end{align*}
By Assumption \ref{ass:error}, Lemma \ref{lemma.y.exist},  the definitions of $\tau_{\min}$ and $\tilde \tau_{\min}$, it follows that
\begin{align*}
    | \tau_{\min} - \tilde \tau_{\min}| &\le \left |  \tfrac{(1- \sigma_{\tau})\left( 1-\sigma_{c}/\sigma_r \right)  \sigma_v \kappa_{\sigma}^2}{ \kappa_{ u^*} \kappa_H \sigma_{Jc} \kappa_J  +  \kappa_g  \sigma_{Jc} \kappa_J  } -  \tfrac{(1- \sigma_{\tau})\left( 1-\sigma_{c}/\sigma_r \right)  \sigma_v \kappa_{\sigma}^2}{ \kappa_{\bar u^*} \kappa_H \sigma_{Jc} (\kappa_J  + \epsilon_J) +  (\kappa_g + \epsilon_g)  \sigma_{Jc} (\kappa_J  + \epsilon_J) } \right|  
    \\ &\le (1- \sigma_{\tau})\left( 1-\tfrac{\sigma_{c}}{\sigma_r} \right) \sigma_v \kappa_{\sigma}^2  \left | 
    \tfrac{\kappa_{\bar u^*} \kappa_H \sigma_{Jc} (\kappa_J  + \epsilon_J) +  (\kappa_g + \epsilon_g)  \sigma_{Jc} (\kappa_J  + \epsilon_J)  - \kappa_{ u^*} \kappa_H \sigma_{Jc} \kappa_J  +  \kappa_g  \sigma_{Jc} \kappa_J  }{( \kappa_{ u^*} \kappa_H \sigma_{Jc} \kappa_J  +  \kappa_g  \sigma_{Jc} \kappa_J )(\kappa_{\bar u^*} \kappa_H \sigma_{Jc} (\kappa_J  + \epsilon_J) +  (\kappa_g + \epsilon_g)  \sigma_{Jc} (\kappa_J  + \epsilon_J) )}  \right | 
    \\ &\le (1- \sigma_{\tau})\left( 1-\tfrac{\sigma_{c}}{\sigma_r} \right)  \sigma_v \kappa_{\sigma}^2   \left | 
    \tfrac{ \kappa_H \sigma_{Jc} \kappa_J  \epsilon_{u^*} + \kappa_{\bar u^*} \kappa_H \sigma_{Jc} \epsilon_J + \kappa_g \sigma_{Jc} \epsilon_J + \kappa_J \sigma_{Jc} \epsilon_g+ \epsilon_J \sigma_{Jc} \epsilon_g    }{( \kappa_{ u^*} \kappa_H \sigma_{Jc} \kappa_J  +  \kappa_g  \sigma_{Jc} \kappa_J )(\kappa_{\bar u^*} \kappa_H \sigma_{Jc} (\kappa_J  + \epsilon_J) +  (\kappa_g + \epsilon_g)  \sigma_{Jc} (\kappa_J  + \epsilon_J) )}  \right |  \\ &:=  \mathcal{E}_{\tau,\min}, 
\end{align*}
which completes the proof.
\end{proof} 

\begin{remark}
    Lemma~\ref{lemma.diff.tau.min} bounds the deviation in the lower bound of the merit parameter sequence from a fictious true (never computed) lower bound when the LICQ holds. 
    By Lemma~\ref{lemma.tau.lb}, this 
     deviation 
    is sufficiently small and proportional to the noise levels in the problem. When $\epsilon_g = \epsilon_c = \epsilon_J = 0$, we have $\mathcal{E}_{\tau,\min} = 0$ and $(\bar{\tau}_k$, $\tau_k)$ share the same positive lower bound.     
\end{remark}







\subsection{Main Convergence Theorems}
\label{sec:convergence.main}

In this section, we establish the main convergence theorems for Algorithm \ref{alg.dfo_sqp_LS} under both step size strategies. Before presenting the main results, we note that our results are analogous to results that can be established in the determistic setting, 
see e.g., \cite[Theorem 3.3]{CurtNoceWach10}, and are summarized in Table~\ref{tab:measure_table}.
In the most general case, convergence can only be established in terms of the \textit{infeasible stationarity error} $\|J_k^T c_k\|$. 
When $J_k$ is of full rank, we can prove a stronger convergence result in terms of the \textit{feasibility error} $\|c_k\|$. Furthermore, as we will show later in this section, when $\bar J_k$ is of full rank, we can also establish the convergence of $\|c_k\|$. Whether or not convergence can be established in terms of the \textit{stationarity error} 
\begin{align}
\label{eq.stationarity.error}
    \| Z_{k}^T(\nabla f(x_{k})+H_{k} v_{k})\|
\end{align} 
heavily depends on the behavior of the merit parameter sequence. When $\{\bar \tau_k\}$ diminishes to zero, no meaningful conclusions can be drawn about the \textit{stationarity error}. When $\{\bar \tau_k\}$ is bounded below  by some $\bar \tau_{\min}\in\mathbb{R}_{>0}$, we prove that the averaged \textit{stationarity error} is sufficiently small. It is worth noting that the reduced gradient $\left\|Z_{k}^T \nabla f\left(x_{k}\right) \right\| = \left\|Z_{k}^T ( \nabla f\left(x_{k}  \right) + J_k^T y_k^*) \right\|$ 
is used as a stationarity measure in \cite{sun2024trust}. 
We argue that \eqref{eq.stationarity.error} 
is a reasonable stationarity measure for problem \eqref{prob.main} provided that $\|J_k^T c_k \|$ vanishes or is sufficiently small. Since $ \| Z_{k}^T H_{k} v_{k}\| \leq \| H_{k}\|\|v_{k}\| \le \kappa_H \sigma_{Jc} \| J_k^T c_k \|$ by Assumption~\ref{ass:H} and~\eqref{eq.nolicq.tr.sub.p1.det}, as long as $\| J_k^T c_k \|$ is sufficiently small,  $\| Z_{k}^T(\nabla f(x_{k})+H_{k} v_{k})\|$ closely approximates the norm of the reduced gradient 
$\left\|Z_{k}^T \nabla f\left(x_{k}\right) \right\|$ and is thus a reasonable stationarity measure.  The stationarity error is related to the norm of tangential 
component.


\begin{table}[]
\centering
\resizebox{\columnwidth}{!}{
\begin{tabular}{c c ccccc}
\toprule
\multirow{2}{*}{} & \multirow{2}{*}{$\sigma_{\min} (\bar J_k) \ge \kappa_{\sigma}$, $\bar \tau_k \ge \bar \tau_{\min} $} & \multicolumn{2}{c}{$\sigma_{\min} (\bar J_k) =0$, $\bar \tau_k \ge \bar \tau_{\min} $} & \multicolumn{2}{c}{$\sigma_{\min} (\bar J_k) =0$, $\bar \tau_k \to 0$} \\ 
\cmidrule(l){3-7}
   &  & $\sigma_{\min} ( J_k) \ge \kappa_{\sigma}$ & $\sigma_{\min} ( J_k) =0$ & $\sigma_{\min} ( J_k) \ge \kappa_{\sigma}$& $\sigma_{\min} ( J_k) =0$ \\  
\cmidrule(r){1-7}
(in)fea. & $\|c_k\|$ & $\|c_k\|^2$ & $\|J_k^T c_k\|^2$  & $\|c_k\|^2$ & $\|J_k^T c_k\|^2$ \\
sta. & $\left\| Z_k^T\left(\nabla f\left(x_k\right)+H_k v_k\right)\right\|^2$    & \multicolumn{2}{c}{$\left\| Z_k^T\left(\nabla f\left(x_k\right)+H_k v_k\right)\right\|^2$}  & \multicolumn{2}{c}{--} \\
\bottomrule
\end{tabular}}
\caption{ Summary of feasibility or infeasible stationarity measure and  stationarity measure used in different cases.} 
\label{tab:measure_table}
\end{table}

We first derive upper bounds of the model reduction for both step size schemes when $\bar \tau_k \ge \bar{\tau}_{\min}$. However, when 
$\{\bar\tau_k\}$ diminishes to zero, the step size may vanish, and a sufficient reduction in the merit function can no longer be ensured. In this case, we show that the algorithm converges to a nearly infeasible stationary solution. To address this scenario, we provide a result on the reduction of constraint violation, which supports the convergence analysis when the merit parameter vanishes. We emphasize the importance of demonstrating that the true merit function achieves sufficient reduction at each iteration under the proposed step size. To this end, we define 
\begin{align}\label{def.merit}
   \phi(x,\tau) := \tau  f(x) + \| c(x)\|.  
\end{align}

\begin{lemma}
\label{lemma.model.reduction.combined}
 Suppose that Assumptions~\ref{ass:prob}, \ref{ass:error}, \ref{ass:H}, 
\ref{ass.inexact.solver.accuracy}, \ref{ass.d.nonzero}, 
and \ref{ass.noearlytermination} hold. 
 Let  
 \begin{align}
   \mathcal{E}^A  &:=\alpha_u^A \bar \tau_{-1} \kappa_{\bar d} \epsilon_g  + 2 \epsilon_c +   \alpha_u^A  \kappa_{\bar d} \epsilon_J \label{def.epsilon_A} \\
    \text{and} \qquad \mathcal{E}^L  &:=   4\bar{\tau}_{-1} \epsilon_f  +  \bar\tau_{-1}  \alpha_u^L \kappa_{\bar d} \epsilon_g + 
        6 \epsilon_c 
         +  \alpha_u^L \kappa_{\bar d} \epsilon_J.     \label{eq.epsilon.L}
\end{align}
If $\bar \tau_k \ge \bar{\tau}_{\min}$, for all $k \in \N{}$
\begin{align}
\label{eq.lemma.model.reduction.ada.LICQ}
    \phi\left(x_{k}+ \bar \alpha_k \bar d_{k}, \bar \tau_{k}\right)-\phi\left(x_{k}, \bar \tau_{k}\right) \le & - \eta \alpha_{l\tau}^i \Delta \lbar\left(x_k, \bar\tau_k, \bar d_k\right)  + \mathcal{E}^i, 
\end{align}
where $i \in \{A, L\}$ denotes the adaptive (\textbf{Option I}) and line search (\textbf{Option II}) scheme, respectively, in Algorithm~\ref{alg.dfo_sqp_LS}. 
If $\{\bar \tau_k\}$ diminishes to 0 
and $x_k \in \mathcal{X}_\gamma$, for all $k \in \N{}$
\begin{align}
\label{eq.lemma.model.reduction.ada.NOLICQ}
 \left\|c_{k+1}\right\| -\left\| c_{k}\right\| \le - \eta \alpha_{l0}^i \Delta \lbar\left(x_k, \bar\tau_k, \bar d_k\right) - \bar  \tau_{k} \left( f_{k+1}  - f_{k} \right)    + \mathcal{E}^i.   
\end{align}
\end{lemma}
\begin{proof}
We first prove that \eqref{eq.lemma.model.reduction.ada.LICQ} and  \eqref{eq.lemma.model.reduction.ada.NOLICQ} hold when $i = A$. By Lemma~\ref{lemma.inexact.stepsize.opt} it follows that $\bar\alpha_{k} \le 1$. 
When $\bar \tau_k \ge \bar{\tau}_{\min}$, 
it follows that  
\begin{align}
    &\ \phi\left(x_{k}+ \bar \alpha_k \bar d_{k}, \bar \tau_{k}\right)-\phi\left(x_{k}, \bar \tau_{k}\right) \notag 
    \\ \le &\ - \bar \alpha_k \Delta \lbar\left(x_k, \bar\tau_k, \bar d_k\right)  +  \bar \alpha_k \bar \tau_{k} ( g_{k} -\bar g_{k} )^{T}  \bar d_{k}  + \tfrac{1}{2}\left(\bar\tau_{k} L+\Gamma\right) \bar \alpha_k ^{2}\left\|\bar d_{k}\right\|^{2} + 
     2\epsilon_c +  \bar \alpha_k\epsilon_J \|\bar d_k\| \notag \\ 
    \le  &  - \bar \alpha_k \Delta \lbar\left(x_k, \bar\tau_k, \bar d_k\right)  +  \bar \alpha_k \bar \tau_{k} ( g_{k} -\bar g_{k} )^{T}  \bar d_{k}  + (1-\eta) \beta \bar \alpha_k \Delta \lbar\left(x_k, \bar\tau_k, \bar d_k\right)   + 2\epsilon_c + \bar \alpha_k\epsilon_J \|\bar d_k\|
    \notag \\
        \le  &  - \eta \bar \alpha_k  \Delta \lbar\left(x_k, \bar\tau_k, \bar d_k\right)  +  \bar \alpha_k  \bar \tau_{k} \kappa_{\bar d} \epsilon_g  + 2\epsilon_c +  \bar \alpha_k  \kappa_{\bar d} \epsilon_J
    \label{eq.model.red.ada.second}
    \\ \le & - \eta \alpha_{l\tau}^A \Delta \lbar\left(x_k, \bar\tau_k, \bar d_k\right)  +  \alpha_u^A \bar \tau_{-1} \kappa_{\bar d} \epsilon_g  + 2\epsilon_c +   \alpha_u^A  \kappa_{\bar d} \epsilon_J,   \label{eq.model.red.ada}
\end{align}
where the first inequality follows by the third last  inequality in \eqref{eq.model.reduction.c.large}, the second inequality follows by the fact that $\bar \alpha_k \le \bar \alpha_k^{\text{suff}} \le \tfrac{2(1-\eta) \beta_{k} \Delta \lbar(x_{k},\bar \tau_{k}, \bar d_{k})   }{\left(\bar\tau_{k} L+\Gamma\right)\|\bar d_{k}\|^{2}}$ (see Lemma~\ref{lemma.inexact.stepsize.opt} and \eqref{eq.inexact.alpha.suff}), the third 
 inequality follows by Lemma~\ref{lemma.dbar.ubar.bound}, 
 Assumption~\ref{ass:error}, 
 and the fact that $\beta \le 1$, and the last inequality follows by Lemmas~\ref{lemma.tau.lb} and~\ref{lemma.inexact.stepsize.opt}. 
 \eqref{eq.lemma.model.reduction.ada.LICQ} is therefore satisfied.  When $\{\bar \tau_k\}$ diminishes to 0 and $x_k \in \mathcal{X}_\gamma$, it follows from Lemma~\ref{lemma.inexact.stepsize.opt} that $\bar \alpha_k \le \alpha_{l0}^A$. Finally, \eqref{eq.lemma.model.reduction.ada.NOLICQ} follows by \eqref{def.merit} and \eqref{eq.model.red.ada.second}.

Next, we prove that \eqref{eq.lemma.model.reduction.ada.LICQ} and~\eqref{eq.lemma.model.reduction.ada.NOLICQ} hold when $i = L$. When $\bar \tau_k \ge \bar{\tau}_{\min}$, it follows by Assumption~\ref{ass:error}, \eqref{eq.diff.merit.func}, \eqref{eq.modified.linesearch.1d}, Lemma~\ref{lemma.dbar.ubar.bound} and Lemma~\ref{lemma.step.size.ls} that 
\begin{align}
      \phi(x_{k} + \bar\alpha_k \bar{d}_k,\bar{\tau}_{k}) -  \phi(x_{k},\bar{\tau}_{k})  &\le 
       \bar\phi(x_{k} + \bar\alpha_k \bar{d}_k,\bar{\tau}_{k}) - \bar\phi(x_{k},\bar{\tau}_{k}) + 2 \bar \tau_k \epsilon_f + 2\epsilon_c \notag \\ 
        &\le  - \eta \bar\alpha_k \Delta \lbar(x_{k},\bar \tau_{k}, \bar d_{k}) + \epsilon_{A_k} + 2 \bar \tau_k \epsilon_f + 2\epsilon_c \label{eq.model.red.ls.second} \\
        &\le  - \eta  \alpha_{l\tau}^L  \Delta \lbar(x_{k},\bar \tau_{k}, \bar d_{k}) + 4\bar{\tau}_{-1} \epsilon_f  +  \bar\tau_{-1}  \alpha_u^L \kappa_{\bar d} \epsilon_g  +       6 \epsilon_c 
         +  \alpha_u^L \kappa_{\bar d} \epsilon_J,   \label{eq.model.red.ls}
\end{align}
and \eqref{eq.lemma.model.reduction.ada.LICQ} is therefore satisfied.
When $\{\bar \tau_k\}$ diminishes to 0 and $x_k \in \mathcal{X}_\gamma$, 
\eqref{eq.lemma.model.reduction.ada.NOLICQ}
follows by \eqref{def.merit}, \eqref{eq.model.red.ls}, and Lemma~\ref{lemma.step.size.ls}. 
\end{proof}

Table~\ref{tbl.alpha} summarizes the results of Lemma~\ref{lemma.model.reduction.combined}.

\begin{table}[]
\centering
\caption{Table summarizing the notation (upper and lower bounds and merit function increase) from Lemma~\ref{lemma.model.reduction.combined} for the two step size schemes $i \in \{A,L\}$ and different scenarios for the merit parameter ($\bar \tau_k \ge \bar \tau_{\min}$ and $\bar \tau_k \to 0$, $x_k \in \mathcal{X}_\gamma$). 
\label{tbl.alpha}}
\resizebox{\columnwidth}{!}{
\begin{tabular}{lcccc}
\toprule
& \multicolumn{1}{c}{Upper Bound} & \multicolumn{1}{c}{Lower Bound} & \multicolumn{1}{c}{Lower Bound} & \multicolumn{1}{c}{Merit Function} \\
& & ($\bar \tau_k \ge \bar \tau_{\min}$) & ($\bar \tau_k \to 0$ and $x_k \in \mathcal{X}_\gamma$) & Increase \\
\midrule
\textbf{Option I} & $ \alpha_u^A $ & $ \alpha_{l\tau}^A $ & $\alpha_{l0}^A $ & $\mathcal{E}^A= \mathcal{O}(\epsilon_g+\epsilon_c+\epsilon_J)$ \\
\textbf{Option II} & $ \alpha_u^L $ & $ \alpha_{l\tau}^L $ & $\alpha_{l0}^L $ & $\mathcal{E}^L = \mathcal{O}(\epsilon_f+\epsilon_g+\epsilon_c+\epsilon_J)$ \\
\bottomrule
\end{tabular}
}
\end{table}

The following lemma shows a relationship between the true tangential step $ u_k^*$ and the stationarity error \eqref{eq.stationarity.error}.

\begin{lemma}
\label{lemma.u.lower.bound}
Suppose that Assumptions~\ref{ass:H},  \ref{ass.noearlytermination}, and \ref{ass:Z.exist} hold. 
For all $k \in \mathbb{N}$, 

\noindent
$\left\| u_k^*\right\| \geq \kappa_H^{-1}\left\| Z_k^T\left(\nabla f\left(x_k\right)+H_k v_k\right)\right\|$. 
\end{lemma}
\begin{proof}
The proof follows the same logic as \cite[Lemma 9]{berahas2023stochastic}.  
\end{proof}

The following lemma establishes bounds for feasibility and stationarity errors in the event that Algorithm \ref{alg.dfo_sqp_LS} terminates finitely on Line \ref{line:terminate}.


\begin{theorem}
\label{theorem.early.termination}
Suppose that Assumptions~\ref{ass:prob}, \ref{ass:error}, \ref{ass:H}, \ref{ass.inexact.solver.accuracy}, \ref{ass.d.nonzero}, and \ref{ass:Z.exist} hold. If Algorithm~\ref{alg.dfo_sqp_LS} terminates finitely at iteration $k = K_T \in \mathbb{N}$ on Line~\ref{line:terminate} (i.e., $\|\bar c_{K_T}\| \le \epsilon_o$ and $\Delta \bar l(x_{K_T},\bar \tau_{K_T}, \bar d_{K_T} ) \le \epsilon_o$), where $\epsilon_o \in [0,\epsilon_c]$, then $\| c_{K_T}\| \le 2\epsilon_c$ and
\begin{equation}\label{eq.opt.early.termination}
\begin{aligned}
    &\left\| Z_{K_T}^T\left(\nabla f\left(x_{K_T}\right)+H_{K_T} v_{K_T}\right)\right\|^2 \\
    \le & \tfrac{ 4 \kappa_H^2 \epsilon_c}{ \bar \tau_{K_T} \sigma_u  \lambda_u } + 12 \kappa_H^2\zeta^{-2} \left(  \epsilon_g^2+ \kappa_{\bar y^*}^2 \epsilon_J^2  + \kappa_{y^*}^2 \epsilon_J^2  \right) + 96 \kappa_H^4 \zeta^{-2}\sigma_{Jc}^2 \kappa_J^2 \epsilon_c^2\\
     &\quad +  4  \kappa_H^2\lambda_{\rho r}^2  \left(\zeta^{-1}  + (\kappa_{J} + \epsilon_J )   \kappa_{\sigma}^{-2}  \right)^2 \max \left \{ \tfrac{ 2\epsilon_c }{ \bar \tau_{K_T} \sigma_u  \lambda_u }, (\kappa_J + \epsilon_J)^2 \epsilon_c^2  \right \} 
\end{aligned}
\end{equation}
\end{theorem}
\begin{proof}
    Since Algorithm~\ref{alg.dfo_sqp_LS} terminates at line~\ref{line:terminate} only when $\|\bar c_k\| \le \epsilon_o$ and $\Delta \bar l(x_k,\bar \tau_k, \bar d_k ) \le \epsilon_o$, it follows by Assumption~\ref{ass:error} that $\|c_k\| \le \epsilon_c + \epsilon_o \le 2\epsilon_c$. Additionally, by \eqref{eq.detal.lb1}
\begin{align*}
 \bar \tau_{K_T} \sigma_u  \lambda_u \|\bar u_{K_T}\|^2 - \epsilon_o  \le \Delta \lbar(x_{K_T},\bar\tau_{K_T},\bar{d}_{K_T}) \le \epsilon_o, 
\end{align*}
which implies that $\| \bar u_{K_T}\|^2 \le \tfrac{ 2\epsilon_o }{ \bar \tau_{K_T} \sigma_u  \lambda_u }\le \tfrac{ 2\epsilon_c }{ \bar \tau_{K_T} \sigma_u  \lambda_u }$. Moreover, by Assumptions~\ref{ass:prob} and \ref{ass:error} it follows that $\|\bar J_{K_T}^T \bar c_{K_T}\| \le (\kappa_J + \epsilon_J) \epsilon_o \le   (\kappa_J + \epsilon_J) \epsilon_c $ and $\|J_{K_T}^T c_{K_T}\| \le 2\kappa_J \epsilon_c$. Thus, the stationarity result \eqref{eq.opt.early.termination} holds 
by Lemmas~\ref{lemma.diff.u} and~\ref{lemma.u.lower.bound}, and the fact that $\| u_{K_T}^*\|^2 \le 2 (\|\bar u_{K_T}-  u_{K_T}^*\|^2 + \|\bar u_{K_T}\|^2 )$.
\end{proof}

Next, we show that if Algorithm~\ref{alg.dfo_sqp_LS} terminates finitely on Line~\ref{line.early.stopping}, the \textit{infeasible stationarity error} is sufficiently small.
\begin{theorem}
\label{theorem.early.termination.infeasible}
Suppose that Assumptions~\ref{ass:prob}, \ref{ass:error}, \ref{ass:H}, \ref{ass.inexact.solver.accuracy}, and \ref{ass.d.nonzero} hold. If Algorithm~\ref{alg.dfo_sqp_LS} terminates finitely at iteration $k = K_T \in \mathbb{N}$ on Line~\ref{line.early.stopping} (i.e., $ \| \bar c_{K_T} \| > \epsilon_o$ and  $\|\bar J_{K_T}^T \bar c_{K_T}\| = 0$), where $\epsilon_o \in (0,\epsilon_c]$, then $\|J_{K_T}^T c_{K_T}\| \le \kappa_c \epsilon_J +  ( \kappa_J + \epsilon_J) \epsilon_c$.
\end{theorem}
\begin{proof}
    By 
    Lemma~\ref{lemma.mathcalE.jc}, 
   \begin{align*}
          \|J_{K_T}^T c_{K_T}\| \le \|\bar  J_{K_T}^T \bar c_{K_T}  \|  + \|J_{K_T}^T c_{K_T} - \bar  J_{K_T}^T \bar c_{K_T}  \|  \le \kappa_c \epsilon_J +  ( \kappa_J + \epsilon_J) \epsilon_c, 
    \end{align*}
    which concludes the statement.
\end{proof}

\begin{remark}
By Theorem~\ref{theorem.early.termination}, we note that if Algorithm \ref{alg.dfo_sqp_LS} ever terminates on Line~\ref{line:terminate}, the 
terminated iterate is sufficiently feasible even if the constraint Jacobians are rank-deficient. Moreover, we provide an upper bound on the first-order stationarity error that depends on the noise level and merit parameter value at the terminated iterate ($\bar \tau_{K_T}$). If the singular values of $\bar J_k$ are bounded away from 0, then we have $\bar \tau_{K_T} \ge \bar \tau_{\min}$. However, if $\bar J_k$ is rank-deficient, then $\bar \tau_{K_T}$ can be arbitrarily close to zero, in which case the stationarity error can be large. By Theorem~\ref{theorem.early.termination.infeasible}, we note that if Algorithm \ref{alg.dfo_sqp_LS} ever terminates on Line~\ref{line.early.stopping}, the infeasible stationarity error at terminated iterate is sufficiently small.
\end{remark}

We set up some general conditions necessary for the step size strategies required for Theorem~\ref{theorem.tau.lb}, \ref{theorem.tau.tozero}, and~\ref{theorem.tau.lb.complexity}. If Algorithm~\ref{alg.dfo_sqp_LS} employs \textbf{Option I} for the step size scheme, $\beta$ must satisfy the conditions outlined in Lemma \ref{lemma.inexact.stepsize.opt}. Alternatively, if Algorithm~\ref{alg.dfo_sqp_LS} employs \textbf{Option II}, the relaxation parameter $\epsilon_{A_k}$ is defined in \eqref{def.epsilon_Ak}. Under these conditions, we show that when $\bar \tau_k$ is bounded away from zero, the averaged stationarity and infeasible stationarity errors converge to a sufficiently small value, 
the size of which depends on the noise levels.

\begin{theorem}
\label{theorem.tau.lb}
Suppose Assumptions~\ref{ass:prob}, \ref{ass:error}, \ref{ass:H}, \ref{ass.inexact.solver.accuracy}, \ref{ass.d.nonzero}, \ref{ass.noearlytermination}, and \ref{ass:Z.exist} hold. 
In addition, suppose that there exists $\bar \tau_{\min} \in \mathbb{R}_{>0}$ such that $\bar \tau_k \ge \bar \tau_{\min}$ for all $k\in \mathbb{N}$. Let $\phi_{\min} = \min \{0,\bar\tau_{-1} f_{\inf}\}$
, and let $i \in \{A, L\}$ denote the adaptive (\textbf{Option I}) and line search (\textbf{Option II}) scheme, respectively, in Algorithm~\ref{alg.dfo_sqp_LS}. Then, for all $k \in \N{}$, 
\begin{align}
\label{eq.convergence.theorem.ub.Jc}
    \tfrac{1}{k+1}  \sum_{j=0}^{k}    \| J_j^T  c_j \|^2  &\le \tfrac{ ( \bar \tau_{-1} |f_{\inf}| + \phi(x_ 0 , \bar{\tau}_{0})  - \phi_{\min})(\kappa_c + \epsilon_c)  }{(k+1)\sigma_{c}\sigma_v\eta\alpha_{l\tau}^i} + \tfrac{\kappa_c + \epsilon_c}{\sigma_{c}\sigma_v\eta\alpha_{l\tau}^i}  \mathcal{E}^i + \mathcal{E}_{Jc^2},  
\end{align}
and
\begin{equation}
     \begin{split}
      &\   \tfrac{1}{k+1} \sum_{j=0}^{k} \tfrac{\left\| Z_j^T\left(\nabla f\left(x_j\right)+H_j v_j\right)\right\|^2}{\kappa_H^2}   \\ 
      \le &\ \left( \tfrac{4+ 2 \kappa_1}{\bar\tau_{\min} \sigma_{u} \lambda_u \eta\alpha_{l\tau}^i} + \tfrac{\left(  \kappa_1 + 2 \kappa_2  \right)(\kappa_c + \epsilon_c)}{\sigma_{c}\sigma_v\eta\alpha_{l\tau}^i} \right) \tfrac{\bar \tau_{-1} |f_{\inf}| + \phi(x_ 0 , \bar{\tau}_{0})  - \phi_{\min}}{k+1}\\
   &\ + \left(  \tfrac{4 +  2 \kappa_1 }{\bar\tau_{\min} \sigma_{u} \lambda_u \eta\alpha_{l\tau}^i}  + \tfrac{\left(  \kappa_1 + 2 \kappa_2  \right)(\kappa_c + \epsilon_c)}{\sigma_{c}\sigma_v\eta\alpha_{l\tau}^i}   \right) \mathcal{E}^i  \\
   &\ + 2\kappa_2 \mathcal{E}_{Jc^2} +  \kappa_1  (\kappa_J+\epsilon_J)^2 \epsilon_c^2 + 12 \zeta^{-2} \left(  \epsilon_g^2+ \kappa_{\bar y^*}^2 \epsilon_J^2  + \kappa_{y^*}^2 \epsilon_J^2  \right),  
     \end{split} 
  \label{eq.convergence.theorem.ub.u}
 \end{equation}
 where $\kappa_1 = 4 \lambda_{\rho r}^2  \left(\zeta^{-1}  + (\kappa_{J} + \epsilon_J )   \kappa_{\sigma}^{-2}  \right)^2$ and  $\kappa_2 = 24\zeta^{-2} \kappa_H^2 \sigma_{Jc}^2 $.

If the singular values of $ J_k$ are bounded below by $\kappa_{\sigma} \in \R{}_{>0}$ for all $k \in \mathbb{N}$, 
 \begin{align}
\label{eq.convergence.theorem.ub.c2}
    \tfrac{1}{k+1}  \sum_{j=0}^{k}    \|c_j \|^2  &\le \tfrac{ ( \bar \tau_{-1} |f_{\inf}| + \phi(x_ 0 , \bar{\tau}_{0})  - \phi_{\min})(\kappa_c + \epsilon_c)  }{(k+1)\sigma_{c}\sigma_v\eta\alpha_{l\tau}^i \kappa_{\sigma}^2} + \tfrac{\kappa_c + \epsilon_c}{\sigma_{c}\sigma_v\eta\alpha_{l\tau}^i \kappa_{\sigma}^2}  \mathcal{E}^i + 
 \tfrac{1}{ \kappa_{\sigma}^2} \mathcal{E}_{Jc^2}. 
\end{align}

 If the singular values of $\bar J_k$ are bounded below by $\kappa_{\sigma}  \in \R{}_{>0}$ for all $k \in \mathbb{N}$,
 \begin{align}
\label{eq.convergence.theorem.ub.c}
    \tfrac{1}{k+1}  \sum_{j=0}^{k}    \| c_j \|   &\le  \tfrac{\bar \tau_{-1} |f_{\inf}| + \phi(x_ 0 , \bar{\tau}_{0})  - \phi_{\min}  }{  (k+1)   \kappa_{\sigma}^{2} \sigma_{c}\sigma_v  \eta\alpha_{l\tau}^i } + \tfrac{1}{\kappa_{\sigma}^{2} \sigma_{c}\sigma_v  \eta\alpha_{l\tau}^i } \mathcal{E}^i  + 2\epsilon_c .  
\end{align}
\end{theorem}
\begin{proof}
 \allowdisplaybreaks{
It follows by Assumption~\ref{ass:prob} and Lemmas~\ref{lemma.d.bar.bounded.by.Delta.l.opt},~\ref{lemma.tau.lb} and~\ref{lemma.model.reduction.combined} that  
   \begin{align}\label{eq.112}
        &\ \phi\left(x_{k}+ \bar \alpha_k \bar d_{k},  \bar\tau_{k+1}\right)-\phi\left(x_{k}, \bar\tau_{k}\right)  \notag \\ =&\    \phi\left(x_{k}+ \bar \alpha_k \bar d_{k},  \bar\tau_{k}\right)-\phi\left(x_{k}, \bar\tau_{k}\right) + \phi\left(x_{k}+ \bar \alpha_k \bar d_{k},  \bar\tau_{k+1}\right)- \phi\left(x_{k}+ \bar \alpha_k \bar d_{k},  \bar\tau_{k}\right) \notag  \\ \le &   - \eta \alpha_{l\tau}^i \Delta \lbar\left(x_k, \bar\tau_k, \bar d_k\right)  +  \mathcal{E}^i + (   \bar\tau_{k+1} -  \bar\tau_{k} ) f_{\inf} 
       \notag  \\ \le&\   - \eta\alpha_{l\tau}^i \left( \tfrac{\bar\tau_{k} \sigma_{u} \lambda_u}{2} \left\|\bar u_{k}\right\|^{2}   + \tfrac{\sigma_{c}\sigma_v}{\kappa_c + \epsilon_c}   \|\bar  J_k^T \bar  c_k \|^2   \mathbbm{1}\{ \|\bar c_k\| >\epsilon_o \}  \right)  +  \mathcal{E}^i  + (\bar\tau_{k+1} -  \bar\tau_{k}) f_{\inf}.  
\end{align}
When $f_{\inf} \ge 0$ it follows that $\phi\left(x_{k}, \bar\tau_{k}\right) = \bar \tau_k f_k + \|c_k\| \ge 0$ and when $f_{\inf} < 0$ it follows that $\phi\left(x_{k}, \bar\tau_{k}\right) \ge \bar \tau_{-1} f_{\inf}$. Hence, for all $k \in \mathbb{N}$, it follows that 
\begin{align}
\label{eq.phi.min}
    \phi(x_k, \bar \tau_k) \ge \min \{0, \bar \tau_{-1} f_{\inf} \}:= \phi_{\min}. 
\end{align}
Summing \eqref{eq.112} for $j \in \{0, 1, \dots , k\}$, it follows by \eqref{eq.phi.min} that
\begin{equation}\label{eq.key_cond_tele}
\begin{aligned}
       &\ \phi_{\min}  -  \phi(x_ 0 , \bar{\tau}_{0}) \\ \le &\ \phi(x_{ k + 1}, \bar{\tau}_{k+1})-  \phi(x_ 0 , \bar{\tau}_{0})  
       \\ \leq &\ -  \tfrac{\bar\tau_{k} \sigma_{u} \lambda_u \eta\alpha_{l\tau}^i}{2} \sum_{j=0}^{k}   \left\|\bar u_{j}\right\|^{2} -  \tfrac{\sigma_{c}\sigma_v\eta\alpha_{l\tau}^i}{\kappa_c + \epsilon_c}      \sum_{j=0}^{k}  \left(  \|\bar  J_j^T \bar  c_j \|^2  \mathbbm{1}\{ \|\bar c_j\| >\epsilon_o \} \right)+ \bar\tau_{-1} |f_{\inf}| + (k+1)  \mathcal{E}^i. 
\end{aligned}
\end{equation}
It then follows by \eqref{eq.key_cond_tele} that 
\begin{align}
\label{eq.convergence.ub.Jc.bar}
\tfrac{\sigma_{c}\sigma_v\eta\alpha_{l\tau}^i}{\kappa_c + \epsilon_c}      \sum_{j=0}^{k}    \|\bar  J_j^T \bar  c_j \|^2  \mathbbm{1}\{ \|\bar c_j\| >\epsilon_o \} \le (k+1) \mathcal{E}^i + \bar \tau_{-1} |f_{\inf}| + \phi(x_ 0 , \bar{\tau}_{0})  - \phi_{\min}.  
\end{align}
Since $\left | \|\bar  J_j^T \bar  c_j \|^2  \mathbbm{1}\{ \|\bar c_j\| >\epsilon_o \}  - \|\bar  J_j^T \bar  c_j \|^2 \right | \le  \|\bar  J_j^T \bar  c_j \|^2 \mathbbm{1}\{ \|\bar c_j\| \le \epsilon_o \} \le (\kappa_J+\epsilon_J)^2 \epsilon_c^2 $, we have
\begin{align}
\label{eq.convergence.ub.Jc.bar.wo}
\tfrac{\sigma_{c}\sigma_v\eta\alpha_{l\tau}^i}{\kappa_c + \epsilon_c}     \sum_{j=0}^{k}    \|\bar  J_j^T \bar  c_j \|^2 & \le   \tfrac{\sigma_{c}\sigma_v\eta\alpha_{l\tau}^i}{\kappa_c + \epsilon_c}     \sum_{j=0}^{k}    \left( \|\bar  J_j^T \bar  c_j \|^2     \mathbbm{1}\{ \|\bar c_j\| >\epsilon_o \} + \|\bar  J_j^T \bar  c_j \|^2  \mathbbm{1}\{ \|\bar c_j\| \le \epsilon_o \} \right) \notag \\ &\le (k+1) \left( \mathcal{E}^i + \tfrac{\sigma_{c}\sigma_v\eta\alpha_{l\tau}^i}{\kappa_c + \epsilon_c}    (\kappa_J+\epsilon_J)^2 \epsilon_c^2  \right) + \bar \tau_{-1} |f_{\inf}| + \phi(x_ 0 , \bar{\tau}_{0})  - \phi_{\min}.  
\end{align}
By \eqref{eq.diff.Jc.sto.det} and \eqref{eq.convergence.ub.Jc.bar.wo}, it follows that 
\begin{align}
\label{eq.convergence.ub.Jc}
\tfrac{\sigma_{c}\sigma_v\eta\alpha_{l\tau}^i}{\kappa_c + \epsilon_c}      \sum_{j=0}^{k}    \| J_j^T  c_j \|^2  &\le \tfrac{\sigma_{c}\sigma_v\eta\alpha_{l\tau}^i}{\kappa_c + \epsilon_c}      \sum_{j=0}^{k}      \left(\|\bar  J_j^T \bar  c_j \|^2  \mathbbm{1}\{ \|\bar c_j\| >\epsilon_o \} + \left|  \| J_j^T  c_j \|^2  -  \|\bar  J_j^T \bar  c_j \|^2  \mathbbm{1}\{ \|\bar c_j\| >\epsilon_o \}   \right|\right)  \notag 
\\  &\le
 (k+1) \left( \mathcal{E}^i + \tfrac{\sigma_{c}\sigma_v\eta\alpha_{l\tau}^i}{\kappa_c + \epsilon_c}   
 \mathcal{E}_{Jc^2} \right) + \bar \tau_{-1} |f_{\inf}| + \phi(x_ 0 , \bar{\tau}_{0})  - \phi_{\min}. 
\end{align}
By re-arranging the above, it follows that \eqref{eq.convergence.theorem.ub.Jc} is satisfied. 

Additionally, it follows by \eqref{eq.key_cond_tele} that 
\begin{align}
\label{eq.convergence.ub.u2}
     \tfrac{\bar\tau_{\min} \sigma_{u} \lambda_u \eta\alpha_{l\tau}^i}{2} \sum_{j=0}^{k}   \left\|\bar u_{j}\right\|^{2} &\le  \tfrac{\bar\tau_{k} \sigma_{u} \lambda_u \eta\alpha_{l\tau}^i}{2} \sum_{j=0}^{k}   \left\|\bar u_{j}\right\|^{2} \notag\\
     & \le  \bar \tau_{-1} |f_{\inf}| + \phi(x_ 0 , \bar{\tau}_{0})  - \phi_{\min} + (k+1) \mathcal{E}^i.  
\end{align}
By Lemmas~\ref{lemma.diff.u}  and \ref{lemma.u.lower.bound}, it follows that
 \begin{align*}
 &\ \sum_{j=0}^{k} \left\| Z_j^T\left(\nabla f\left(x_j\right)+H_j v_j\right)\right\|^2  \\
 \le & \   \sum_{j=0}^{k} \kappa_H^2  \left\|u_{j}^*\right\|^{2} \\
 \le  & \ \sum_{j=0}^{k} \kappa_H^2  \left( 2\left\|\bar u_{j}\right\|^{2} + 2 \left\|u_{j}^* - \bar u_{j}\right\|^{2} \right)  \\  \le &\ 2\kappa_H^2 \sum_{j=0}^{k}  \left\|\bar u_{j}\right\|^{2}  +  4 \kappa_H^2 \lambda_{\rho r}^2  \left(\zeta^{-1}  + (\kappa_{J} + \epsilon_J )   \kappa_{\sigma}^{-2}  \right)^2 \sum_{j=0}^{k}  ( \|\bar u_j \|^2 +  \| \bar J_j^T \bar c_j\|^2 )  
 \\ &\ + 24\zeta^{-2} \kappa_H^4 \sigma_{Jc}^2 \sum_{j=0}^{k}  (\|\bar J_j^T \bar c_j\|^2 \mathbbm{1}\{ \|\bar c_j\| >\epsilon_o \} + \|J_j^T c_j\|^2)    +  12(k+1) \kappa_H^2  \zeta^{-2} \left(  \epsilon_g^2+ \kappa_{\bar y^*}^2 \epsilon_J^2  + \kappa_{y^*}^2 \epsilon_J^2  \right) 
 \\  \le &\ (k+1)  \tfrac{4\kappa_H^2+  8 \kappa_H^2 \lambda_{\rho r}^2  \left(\zeta^{-1}  + (\kappa_{J} + \epsilon_J )   \kappa_{\sigma}^{-2}  \right)^2}{\bar\tau_{\min} \sigma_{u} \lambda_u \eta\alpha_{l\tau}^i} \mathcal{E}^i + \tfrac{4\kappa_H^2+  8 \kappa_H^2 \lambda_{\rho r}^2  \left(\zeta^{-1}  + (\kappa_{J} + \epsilon_J )   \kappa_{\sigma}^{-2}  \right)^2}{\bar\tau_{\min} \sigma_{u} \lambda_u \eta\alpha_{l\tau}^i} \left( \bar \tau_{-1} |f_{\inf}| + \phi(x_ 0 , \bar{\tau}_{0})  - \phi_{\min}  \right)  
 \\  &\ + 4 \kappa_H^2 \lambda_{\rho r}^2  \left(\zeta^{-1} \!+\!  (\kappa_{J}  \!+\! \epsilon_J )   \kappa_{\sigma}^{-2}  \right)^2 \tfrac{\kappa_c \!+\! \epsilon_c}{\sigma_{c}\sigma_v\eta\alpha_{l\tau}^i}  \left( (k\!+\!1) \left( \mathcal{E}^i \!+\!\tfrac{\sigma_{c}\sigma_v\eta\alpha_{l\tau}^i}{\kappa_c \!+\! \epsilon_c}   
(\kappa_J\!+\!\epsilon_J)^2 \epsilon_c^2  \right) \!+\!  \bar \tau_{-1} |f_{\inf}| \!+\! \phi(x_ 0 , \bar{\tau}_{0})  \!-\! \phi_{\min}    \right) 
 \\  &\ +  24\zeta^{-2} \kappa_H^4 \sigma_{Jc}^2  \tfrac{\kappa_c + \epsilon_c}{\sigma_{c}\sigma_v\eta\alpha_{l\tau}^i}   \left( (k+1) \left( 2\mathcal{E}^i + \tfrac{\sigma_{c}\sigma_v\eta\alpha_{l\tau}^i}{\kappa_c + \epsilon_c}   
 \mathcal{E}_{Jc^2} \right) + 2\bar \tau_{-1} |f_{\inf}| +2 \phi(x_ 0 , \bar{\tau}_{0})  - 2\phi_{\min}    \right) 
 \\  &\   +  12(k+1) \kappa_H^2  \zeta^{-2} \left(  \epsilon_g^2+ \kappa_{\bar y^*}^2 \epsilon_J^2  + \kappa_{y^*}^2 \epsilon_J^2  \right) 
  \\  \le &\ (k+1) \left(  \tfrac{4\kappa_H^2 +  8 \kappa_H^2 \lambda_{\rho r}^2  \left(\zeta^{-1}  + (\kappa_{J} + \epsilon_J )   \kappa_{\sigma}^{-2}  \right)^2 }{\bar\tau_{\min} \sigma_{u} \lambda_u \eta\alpha_{l\tau}^i}  + \left( 48\zeta^{-2} \kappa_H^4 \sigma_{Jc}^2 + 4 \kappa_H^2 \lambda_{\rho r}^2  \left(\zeta^{-1}  + (\kappa_{J} + \epsilon_J )   \kappa_{\sigma}^{-2}  \right)^2   \right) \tfrac{\kappa_c + \epsilon_c}{\sigma_{c}\sigma_v\eta\alpha_{l\tau}^i}   \right) \mathcal{E}^i 
  \\  &\ + 24 (k+1)  \zeta^{-2} \kappa_H^4 \sigma_{Jc}^2 \mathcal{E}_{Jc^2} +  4  (k+1)  \kappa_H^2 \lambda_{\rho r}^2  \left(\zeta^{-1}  + (\kappa_{J} + \epsilon_J )   \kappa_{\sigma}^{-2}  \right)^2  (\kappa_J+\epsilon_J)^2 \epsilon_c^2  \\  &\   +  12(k+1) \kappa_H^2  \zeta^{-2} \left(  \epsilon_g^2+ \kappa_{\bar y^*}^2 \epsilon_J^2  + \kappa_{y^*}^2 \epsilon_J^2  \right) 
    \\  &\ +  \tfrac{4\kappa_H^2+  8 \kappa_H^2 \lambda_{\rho r}^2  \left(\zeta^{-1}  + (\kappa_{J} + \epsilon_J )   \kappa_{\sigma}^{-2}  \right)^2 }{\bar\tau_{\min} \sigma_{u} \lambda_u \eta\alpha_{l\tau}^i} \left( \bar \tau_{-1} |f_{\inf}| + \phi(x_ 0 , \bar{\tau}_{0})  - \phi_{\min}  \right)   
    \\  &\ + \left( 48\zeta^{-2} \kappa_H^4 \sigma_{Jc}^2 + 4 \kappa_H^2 \lambda_{\rho r}^2  \left(\zeta^{-1}  + (\kappa_{J} + \epsilon_J )   \kappa_{\sigma}^{-2}  \right)^2   \right)  \tfrac{\kappa_c + \epsilon_c}{\sigma_{c}\sigma_v\eta\alpha_{l\tau}^i}    \left( \bar \tau_{-1} |f_{\inf}| + \phi(x_ 0 , \bar{\tau}_{0})  - \phi_{\min}  \right). 
 \end{align*}
Thus, \eqref{eq.convergence.theorem.ub.u} is satisfied.

Starting with \eqref{eq.convergence.ub.Jc},  if the singular values of $J_k$ are bounded below by $\kappa_{\sigma}$ for all $k \in \mathbb{N}$,  \eqref{eq.convergence.theorem.ub.c2} is  then satisfied by \eqref{eq.relation.c.Jc}.

Finally, if the singular values of $\bar J_k$ are bounded below by $\kappa_{\sigma} \in \R{}_{>0}$ for all $k \in \mathbb{N}$, 
it follows by Assumption~\ref{ass:prob}, \eqref{eq.merit} and Lemmas \ref{lemma.d.bar.bounded.by.Delta.l.opt} and \ref{lemma.model.reduction.combined} that 
   \begin{align}\label{eq.113}
        &\ \phi\left(x_{k}+ \bar \alpha_k \bar d_{k},  \bar\tau_{k+1}\right)-\phi\left(x_{k}, \bar\tau_{k}\right)  \notag 
        \\ \le &   - \eta \alpha_{l\tau}^i \Delta \lbar\left(x_k, \bar\tau_k, \bar d_k\right)  +  \mathcal{E}^i + (   \bar\tau_{k+1} -  \bar\tau_{k} ) f_{\inf} 
       \notag  \\ \le&\   - \eta\alpha_{l\tau}^i \left( \tfrac{\bar\tau_{k} \sigma_{u} \lambda_u}{2} \left\|\bar u_{k}\right\|^{2}   + \kappa_{\sigma}^{2} \sigma_{c}\sigma_v \| \bar  c_k \| \mathbbm{1}\{ \|\bar c_j\| >\epsilon_o \} \right)  +  \mathcal{E}^i  + (   \bar\tau_{k+1} -  \bar\tau_{k} ) f_{\inf}.  
   \end{align}
Summing \eqref{eq.113} for $j \in \{0, 1, \dots , k\}$ and re-arranging, it follows that 
\begin{align*}
    \eta\alpha_{l\tau}^i  \kappa_{\sigma}^{2} \sigma_{c}\sigma_v     \sum_{j=0}^{k}    \| \bar  c_j \| \mathbbm{1}\{ \|\bar c_j\| >\epsilon_o \} \le (k+1) \mathcal{E}^i + \bar \tau_{-1} |f_{\inf}| + \phi(x_ 0 , \bar{\tau}_{0})  - \phi_{\min}.   
\end{align*}
Since
\begin{align*}
    \  \left| \| \bar  c_k \| \mathbbm{1}\{ \|\bar c_k\| >\epsilon_o \}  - \| c_k \| \right|
    &\le 
    \max \left \{  \left| \| \bar  c_k \|  - \| c_k \| \right|, \left| \| c_k \| \mathbbm{1}\{ \|\bar c_k\| \le \epsilon_o \}   \right| \right \} \notag\\
    & \le \max \{\epsilon_c, \epsilon_c + \epsilon_o \} \le 2\epsilon_c, 
\end{align*} 
it follows that
\begin{align*}
    \eta\alpha_{l\tau}^i  \kappa_{\sigma}^{2} \sigma_{c}\sigma_v    \sum_{j=0}^{k}    \| c_j \|  \le (k+1) \left( \mathcal{E}^i + 2\eta\alpha_{l\tau}^i  \kappa_{\sigma}^{2} \sigma_{c}\sigma_v \epsilon_c \right) + \bar \tau_{-1} |f_{\inf}| + \phi(x_ 0 , \bar{\tau}_{0})  - \phi_{\min},   
\end{align*}
which completes the proof.
}
\end{proof}

\begin{remark}
    Theorem \ref{theorem.tau.lb} provides the convergence guarantees for Algorithm \ref{alg.dfo_sqp_LS} with either step size strategy when the merit parameter is bounded away from zero. In the most general setting, we can only prove that the averaged square of both the stationarity error and the infeasibility stationarity error converges to a neighborhood of zero as $k \to \infty$. By Lemmas~\ref{lemma.mathcalE.jc} and \ref{lemma.model.reduction.combined}, it follows that $\mathcal{E}^A = \mathcal{O}(\epsilon_g + \epsilon_c+ \epsilon_J)$, $\mathcal{E}^L = \mathcal{O}(\epsilon_f + \epsilon_g + \epsilon_c+ \epsilon_J)$, and $\mathcal{E}_{Jc^2} = \mathcal{O}(\epsilon_c +\epsilon_J) $ when $\epsilon_c \le 1$. Thus, the neighborhood size is therefore $\mathcal{O}(\epsilon_f + \epsilon_g + \epsilon_c+ \epsilon_J)$. 
If the singular values of $J_k$ are bounded below by $\kappa_{\sigma}$, we prove that the averaged square of the feasibility error \eqref{eq.convergence.theorem.ub.c2} converges to a neighborhood of zero as $k \to \infty$, the size of which is $\mathcal{O}(\epsilon_f + \epsilon_g + \epsilon_c+ \epsilon_J)$.
Furthermore, if the singular values of $\bar{J}_k$ are bounded below by $\kappa_{\sigma}$, we can prove a stronger result; the averaged feasibility error \eqref{eq.convergence.theorem.ub.c} converges to 
a neighborhood of zero as 
$k \to \infty$, the size of which is still $\mathcal{O}(\epsilon_f + \epsilon_g + \epsilon_c+ \epsilon_J)$.
When $\epsilon_f = \epsilon_g = \epsilon_c = \epsilon_J = 0$, Theorem \ref{theorem.tau.lb} recovers (3.9) and (3.11) in \cite[Theorem 3.3]{CurtNoceWach10} for a noise-free SQP method with rank-deficient Jabobians.
\end{remark}

\begin{theorem}
\label{theorem.tau.tozero}
Suppose Assumptions~\ref{ass:prob}, \ref{ass:error}, \ref{ass:H}, \ref{ass.inexact.solver.accuracy}, \ref{ass.d.nonzero}, and \ref{ass.noearlytermination} hold. In addition, suppose that the merit parameter sequence $\{\bar \tau_k\}$ diminishes to zero. 
Let $\mathcal{E}^A_{c}  :=  2 \epsilon_c +   \alpha_u^A  \kappa_{\bar d} \epsilon_J$ and $\mathcal{E}^L_{c}  :=   6  \epsilon_c   +  \alpha_u^L \kappa_{\bar d} \epsilon_J$, respectively, and let $i \in \{A, L\}$ denote the adaptive (\textbf{Option I}) and line search (\textbf{Option II}) scheme, respectively, in Algorithm~\ref{alg.dfo_sqp_LS}. 
Then, 
for some $\sigma_{\gamma} \in (0,1)$, for all $k \in \mathbb{N}$,
   \begin{align*}
       \liminf_{k \to \infty}  \|J_k^T c_k\| \le  \max \left\{ \epsilon_c (\epsilon_J + \kappa_J) + \kappa_c \epsilon_J + \sqrt{ \tfrac{ (1+\sigma_{\gamma})(\kappa_c + \epsilon_c) \mathcal{E}^i_{c} }{\eta \alpha_{l0}^i \sigma_c \sigma_v}}, 2\kappa_J \epsilon_c \right\}. 
   \end{align*}
    
    In addition, if the singular values of $J_k$ are bounded below by $\kappa_{\sigma} \in \R{}_{>0}$, for all $k \in \mathbb{N}$, 
  \begin{align*}
       \liminf_{k \to \infty}  \|c_k\| \le  \kappa_{\sigma}^{-1} \max \left\{  \epsilon_c (\epsilon_J + \kappa_J) + \kappa_c \epsilon_J + \sqrt{ \tfrac{ (1+\sigma_{\gamma})(\kappa_c + \epsilon_c) \mathcal{E}^i_{c} }{\eta \alpha_{l0}^i \sigma_c \sigma_v}}, 2\kappa_J \epsilon_c  \right\}. 
   \end{align*}
\end{theorem}
\begin{proof}
Our goal is to prove that for any $\delta \in \mathbb{R}_{>0}$, there exist infinitely many iterates $x_k$ for which 
\begin{align}
\label{eq.to.satisfy}
\|\bar J_k^T \bar c_k\| < \sqrt{ \tfrac{ (1+\sigma_{\gamma})(\kappa_c + \epsilon_c) \mathcal{E}^i_{c} }{\eta \alpha_{l0}^i \sigma_c \sigma_v}}  + \delta \text{ or }  \|\bar c_k \| < \epsilon_o  + \delta. 
\end{align}
To derive a contradiction, we 
assume that there exists some $\delta\in  \mathbb{R}_{>0}$ and some $k_\gamma \in \mathbb{N}$ such that for all $k \geq k_\gamma$,
\begin{align*}
\|\bar J_k^T \bar c_k\| \ge \sqrt{ \tfrac{ (1+\sigma_{\gamma})(\kappa_c + \epsilon_c) \mathcal{E}^i_{c} }{\eta \alpha_{l0}^i \sigma_c \sigma_v}}  + \delta \text{ and }  \|\bar c_k \| \ge \epsilon_o  + \delta.
\end{align*} 
Let $\gamma: = \sqrt{ \tfrac{ (1+\sigma_{\gamma})(\kappa_c + \epsilon_c) \mathcal{E}^i_{c} }{\eta \alpha_{l0}^i \sigma_c \sigma_v}}$ for some $\sigma_{\gamma} \in (0,1)$. By the definition of $\mathcal{X}_\gamma$ (see \eqref{eq.gammas}), we have $x_k \in \mathcal{X}_\gamma$ for all $k \geq k_\gamma$. 

Let $ \mathcal{E}^A_{f} :=  \alpha_u^A  \kappa_{\bar d} \epsilon_g  $ and $\mathcal{E}^L_{f}  :=  4\epsilon_f +  \alpha_u^L \kappa_{\bar d} \epsilon_g$, respectively. It then follows that for all $k \geq k_\gamma$, 
    \begin{equation}
    \label{eq.c.k.minus}
        \begin{split}
             \left\| c_{k}\right\| - \left\|c_{k+1}\right\| 
 \ge&\  \bar  \tau_{k} \left( f_{k+1}  - f_{k} \right)  + \eta \alpha_{l0}^i \Delta \lbar\left(x_k, \bar\tau_k, \bar d_k\right)  -   \mathcal{E}^i_{f} \bar \tau_k -   \mathcal{E}^i_{c} \\ 
 \ge &\   \eta \alpha_{l0}^i \sigma_c \sigma_v (\kappa_c + \epsilon_c)^{-1}  \|\bar J_k^T \bar c_k\|^2 +   \bar  \tau_{k} \left(  f_{\inf} - f_{\sup}  \right)  -   \mathcal{E}^i_{f} \bar \tau_k -   \mathcal{E}^i_{c} \\ \ge &\ (1 + \sigma_{\gamma}) \mathcal{E}^i_{c}  + \bar  \tau_{k} \left(  f_{\inf} - f_{\sup}  \right)  -   \mathcal{E}^i_{f} \bar \tau_k -   \mathcal{E}^i_{c}  \\ = &\  \sigma_{\gamma} \mathcal{E}^i_{c}   +    \bar  \tau_{k} \left(  f_{\inf} - f_{\sup}  \right) -  \mathcal{E}^i_{f} \bar \tau_k,
        \end{split}
    \end{equation}
where the first inequality follows by \eqref{eq.merit} and Lemmas~\ref{lemma.inexact.stepsize.opt} and~\ref{lemma.step.size.ls}, \eqref{eq.model.red.ada.second}, 
\eqref{eq.model.red.ls.second}, the second inequality follows by Assumption~\ref{ass:prob} and  Lemma~\ref{lemma.d.bar.bounded.by.Delta.l.opt}, and the third inequality follows by the fact that $x_k \in \mathcal{X}_\gamma$. Since $\bar \tau_k \to 0$, there exists sufficiently large $ K_\gamma \geq k_\gamma$ such that when  $k \ge K_\gamma$, $\bar  \tau_{k} \left(  f_{\sup} - f_{\inf}   \right) + \mathcal{E}^i_{f} \bar \tau_k \le \tfrac12 \sigma_{\gamma}  \mathcal{E}^i_{c}$, 
which implies that $\left\|c_{k}\right\| - \left\|c_{k+1}\right\|   \ge \tfrac12 \sigma_{\gamma} \mathcal{E}^i_{c} > 0$. 

As a result, for some iteration $k_{c} $ it follows that $\left\| c_{k_{c}}\right\| < \epsilon_o + \delta$, which leads to a contradiction. Hence, there exists a subsequence $ \mathcal{K} \subset \{k_{\gamma}, k_{\gamma}+ 1, \dots \}$ such that \eqref{eq.to.satisfy} is satisfied for all $k \in \mathcal{K}$ and any $\delta \in \mathbb{R}_{>0}$.  
If the former case in \eqref{eq.to.satisfy} is satisfied, it follows by Lemma~\ref{lemma.mathcalE.jc} that  
\begin{align*}
   \| J_k^T c_k \| \le  \| J_k^T  c_k  - \bar J_k^T \bar c_k\| +   \|\bar J_k^T \bar c_k\|  \le 
    (\kappa_J + \epsilon_J) \epsilon_c + \kappa_c \epsilon_J + \sqrt{ \tfrac{ (1+\sigma_{\gamma})(\kappa_c + \epsilon_c) \mathcal{E}^i_{c} }{\eta \alpha_{l0}^i \sigma_c \sigma_v}} +  \delta
\end{align*}
and if the latter case in \eqref{eq.to.satisfy} is satisfied
it follows by Assumption~\ref{ass:prob} and~\ref{ass:error} that 
\begin{align*}
    \| J_k^T c_k \| \le \| J_k^T \bar c_k \|  + \| J_k^T( \bar c_k - c_k) \|  \le \kappa_J (\epsilon_o + \delta + \epsilon_c) \le 2\kappa_J \epsilon_c + \kappa_J\delta.
\end{align*}
Hence, for any $\delta \in \mathbb{R}_{>0}$ and all $k \in \mathcal{K}$, we have
   \begin{align*}
    \|J_k^T c_k\| <  \max \left\{ \epsilon_c (\epsilon_J + \kappa_J) + \kappa_c \epsilon_J + \sqrt{ \tfrac{ (1+\sigma_{\gamma})(\kappa_c + \epsilon_c) \mathcal{E}^i_{c} }{\eta \alpha_{l0}^i \sigma_c \sigma_v}}, 2\kappa_J \epsilon_c \right\} + \max\{\kappa_J,1\}\delta. 
   \end{align*}
By choosing an arbitrarily small value for  $\delta>0$, the first conclusion is therefore satisfied. 

If in addition we assume that the singular values of $J_k$ are bounded below by $\kappa_{\sigma}$ for all $k \in \mathbb{N}$, 
the second conclusion holds as 
$\|J_k^Tc_k\| \geq \kappa_{\sigma}\|c_k\|$.
\end{proof}

\begin{remark}
Theorem \ref{theorem.tau.tozero} provides convergence guarantees for Algorithm \ref{alg.dfo_sqp_LS} with either step size strategy when the merit parameter sequence vanishes to zero. In the most general setting, when $\{\epsilon_c, \epsilon_J\} \subset (0,1]$, 
we prove that there exists a subsequence of ${\|J_k^T c_k\|}$ (infeasible stationarity measure) that converges to a neighborhood of zero whose size is $\mathcal{O}( (\epsilon_c + \epsilon_J)^{1/2})$. Contrary to Theorem \ref{theorem.tau.lb}, we cannot prove convergence for the entire sequence of iterates since the step size is not necessarily lower bounded if $x_k \notin \mathcal{X}_\gamma$. This result is consistent with \cite[Proposition 1]{sun2024trust}, which shows that the iterates visit a critical region related to ${\|J_k^T c_k\|}$ infinitely often. If we additionally assume that the singular values of $J_k$ are bounded below by $\kappa_{\sigma}$, we can prove a similar $\liminf$-type result for ${\|c_k\|}$. 
Finally, under the condition that $\epsilon_c = \epsilon_J = 0$, Theorem~\ref{theorem.tau.tozero} implies that $\liminf_{k\to \infty} \| J_k^T c_k\| = 0$. This result differs from equation (3.10) in \cite[Theorem 3.3]{CurtNoceWach10}, which demonstrates the stronger conclusion that $\lim_{k\to \infty} \| J_k^T c_k\| = 0$. The discrepancy arises due to 
the step size scheme in the presence of noise, making our result the strongest possible under such conditions. Our $\liminf$-type result recovers the corresponding conclusion in \cite[Theorem 3]{berahas2023stochastic} when $\epsilon_c = \epsilon_J = 0$.
\end{remark}

We derive worst-case complexity bounds for the iterates generated by Algorithm~\ref{alg.dfo_sqp_LS} in different settings. We use  $K_{inf}$ and $K_{fea}$ to denote the number of iterations required until 
\begin{align}
    \|J_k^T c_k \| &\le \epsilon \label{eq.complexity.inf}\\ \text{ and } \qquad \qquad \| c_k \| &\le \epsilon \label{eq.complexity.fea}
\end{align}
are guaranteed to be satisfied, respectively. Moreover, we use  $K_{si}$ to denote the number of iterations required until 
\begin{align}\label{eq.complexity.si}
    \left\| Z_k^T\left(\nabla f\left(x_k\right)+H_k v_k\right)\right\| \le \epsilon \qquad \text{ and } \qquad \|J_k^T c_k \| \le \epsilon
\end{align}
are guaranteed to both be satisfied, we use $K_{sf}$  to denote the number of iterations required until 
\begin{align}\label{eq.complexity.sf}
    \left\| Z_k^T\left(\nabla f\left(x_k\right)+H_k v_k\right)\right\| \le \epsilon \qquad \text{ and } \qquad \|c_k \| \le \epsilon
\end{align}
are guaranteed to both be satisfied, and, finally, we use $K_{\overline{sf}}$ to denote the number of iterations required until 
\begin{align}\label{eq.complexity.sf.bar}
    \left\| Z_k^T\left(\nabla f\left(x_k\right)+H_k v_k\right)\right\| \le \epsilon \qquad \text{ and } \qquad  \sqrt{\|c_k \|} \le \epsilon
\end{align}
are both satisfied. Theorems \ref{theorem.tau.lb} and \ref{theorem.tau.tozero} present neighborhood convergence guarantees for Algorithm \ref{alg.dfo_sqp_LS} in different settings. As a result, the worst case complexity results presented below  explicitly define this neighborhood, i.e., for the different convergence measures stated above the target accuracy levels $\epsilon \in \mathbb{R}_{>0}$ are bounded below. Such conditions are analogous to those for noise-aware methods for noisy unconstrained optimization problems \cite{berahas2021global,cao2023first,jin2024high}.

\begin{theorem}
\label{theorem.tau.lb.complexity}
Suppose Assumptions~\ref{ass:prob}, \ref{ass:error}, \ref{ass:H}, \ref{ass.inexact.solver.accuracy}, \ref{ass.d.nonzero}, \ref{ass.noearlytermination}, and \ref{ass:Z.exist} hold. Let $\phi_{\min} = \min \{0,\bar\tau_{-1} f_{\inf}\}$, and let $i \in \{A, L\}$ denote the adaptive (\textbf{Option I}) and line search (\textbf{Option II}) scheme, respectively, in Algorithm~\ref{alg.dfo_sqp_LS}. In addition, let $\sigma_{\epsilon} \in \mathbb{R}_{>1}$ and suppose that \begin{align}
    \epsilon^2 &>  \max \left \{ 4\mathcal{E}_{Jc}^2, 
\sigma_{\epsilon} \left(  \mathcal{E}_{Jc^1} +  \tfrac{ \kappa_c + \epsilon_c}{\eta  \alpha_{l0}^i   \sigma_c \sigma_v } \mathcal{E}^i \right), \sigma_{\epsilon} \left(  \mathcal{E}_{Jc^2} + \tfrac{\kappa_c + \epsilon_c}{\eta\alpha_{l\tau}^i  \sigma_{c}\sigma_v}  \mathcal{E}^i \right) \right \}, 
   \label{eq.epsilon.4}
\end{align}
where $\mathcal{E}_{Jc^1}: =  ( (\kappa_J + \epsilon_J) \epsilon_c + \kappa_c \epsilon_J ) \left(  (\kappa_J + \epsilon_J)  (\kappa_c + \epsilon_c) + \kappa_J \kappa_c \right)$. Then, \eqref{eq.complexity.inf} is guaranteed to be satisfied after 
\begin{align*}
    K_{inf} = \left\lceil \tfrac{ (\kappa_c + \epsilon_c)  ( \bar \tau_{-1} |f_{\inf}| + \phi(x_ 0 , \bar{\tau}_{0})  - \phi_{\min} )}{\eta \min \{ \alpha_{l0}^i, \alpha_{l\tau}^i \} \sigma_c \sigma_v (1- {\sigma_{\epsilon}^{-1}} ) \epsilon^2} \right\rceil
\end{align*}
iterations. Further, if the singular values of $J_k$ are bounded below by $\kappa_{\sigma} \in \mathbb{R}_{>0}$ for all $k \in \mathbb{N}$ and 
\begin{align}
    \epsilon^2 &>  \max \left \{ \tfrac{ 4\mathcal{E}_{Jc}^2}{\kappa_{\sigma}^2}, 
\sigma_{\epsilon} \left( \tfrac{1}{\kappa_{\sigma}^2} \mathcal{E}_{Jc^1} +  \tfrac{ \kappa_c + \epsilon_c}{\eta \min \{ \alpha_{l0}^i, \alpha_{l\tau}^i \} \sigma_c \sigma_v \kappa_{\sigma}^2} \mathcal{E}^i \right) , \sigma_{\epsilon} \left( \tfrac{1}{\kappa_{\sigma}^2} \mathcal{E}_{Jc^2} +  \tfrac{ \kappa_c + \epsilon_c}{\eta \alpha_{l\tau}^i \sigma_c \sigma_v \kappa_{\sigma}^2} \mathcal{E}^i \right)    \right \}, 
\label{eq.epsilon.5}
\end{align}
then, \eqref{eq.complexity.fea} is guaranteed to be satisfied after 
\begin{align*}
    K_{fea} = \left\lceil \tfrac{ (\kappa_c + \epsilon_c)  ( \bar \tau_{-1} |f_{\inf}| + \phi(x_ 0 , \bar{\tau}_{0})  - \phi_{\min} )}{\eta \min \{ \alpha_{l0}^i, \alpha_{l\tau}^i \} \sigma_c \sigma_v \kappa_{\sigma}^2 (1- {\sigma_{\epsilon}^{-1}} ) \epsilon^2} \right\rceil
\end{align*}
iterations. 

Moreover, suppose that there exists $\bar \tau_{\min} \in \mathbb{R}_{>0}$ such that $\bar \tau_k \ge \bar \tau_{\min}$ for all $k\in \mathbb{N}$ and 
\begin{align}
    \epsilon^2 &> \sigma_{\epsilon} \left(\kappa_{si} 
 \mathcal{E}^i + (2\kappa_H^2 \kappa_2 + 1) \mathcal{E}_{Jc^2} + \kappa_H^2 \left(  \kappa_1 (\kappa_J+\epsilon_J)^2 \epsilon_c^2 +  12   \zeta^{-2} \left(  \epsilon_g^2+ \kappa_{\bar y^*}^2 \epsilon_J^2  + \kappa_{y^*}^2 \epsilon_J^2  \right) \right) \right), 
   \label{eq.epsilon.1}
\end{align}
then, \eqref{eq.complexity.si} is guaranteed to be satisfied after 
\begin{align*}
    K_{si} = \left\lceil \tfrac{\kappa_{si} ( \bar \tau_{-1} |f_{\inf}| + \phi(x_ 0 , \bar{\tau}_{0})  - \phi_{\min} )}{(1- {\sigma_{\epsilon}^{-1}} ) \epsilon^2} \right\rceil
\end{align*}
iterations. Further, if the singular values of $J_k$ are bounded below by $\kappa_{\sigma} \in \mathbb{R}_{>0}$ for all $k \in \mathbb{N}$ and
  \begin{align}
    \epsilon^2 &> \sigma_{\epsilon} \left(\kappa_{sf} 
 \mathcal{E}^i + (2\kappa_H^2 \kappa_2 + \kappa_{\sigma}^{-2}) \mathcal{E}_{Jc^2} + \kappa_H^2 \left(  \kappa_1 (\kappa_J+\epsilon_J)^2 \epsilon_c^2 +  12   \zeta^{-2} \left(  \epsilon_g^2+ \kappa_{\bar y^*}^2 \epsilon_J^2  + \kappa_{y^*}^2 \epsilon_J^2  \right) \right) \right),  
   \label{eq.epsilon.2}
\end{align}
then, \eqref{eq.complexity.sf} is guaranteed to be satisfied after 
\begin{align*}
    K_{sf} = \left\lceil \tfrac{\kappa_{sf} ( \bar \tau_{-1} |f_{\inf}| + \phi(x_ 0 , \bar{\tau}_{0})  - \phi_{\min} )}{(1- {\sigma_{\epsilon}^{-1}} ) \epsilon^2} \right\rceil
\end{align*}
iterations. Finally, if the singular values of $\bar J_k$ are bounded below by $\kappa_{\sigma} \in \mathbb{R}_{>0}$ for all $k \in \mathbb{N}$, and 
\begin{align}
    \epsilon^2 &> \sigma_{\epsilon} \left(\kappa_{\overline{sf}} 
 \mathcal{E}^i + 2\kappa_H^2 \kappa_2 \mathcal{E}_{Jc^2} + 
 2\epsilon_c +\kappa_H^2 \left(  \kappa_1 (\kappa_J+\epsilon_J)^2 \epsilon_c^2 +  12   \zeta^{-2} \left(  \epsilon_g^2+ \kappa_{\bar y^*}^2 \epsilon_J^2  + \kappa_{y^*}^2 \epsilon_J^2  \right) \right) \right),  
   \label{eq.epsilon.3}
\end{align}
then, \eqref{eq.complexity.sf.bar} is guaranteed to be satisfied after
\begin{align*}
    K_{\overline{sf}} = \left\lceil \tfrac{\kappa_{\overline{sf}} ( \bar \tau_{-1} |f_{\inf}| + \phi(x_ 0 , \bar{\tau}_{0})  - \phi_{\min} )}{(1- {\sigma_{\epsilon}^{-1}} ) \epsilon^2} \right\rceil
\end{align*}
iterations. Note that $\kappa_{si}$, $\kappa_{sf}$, and $\kappa_{\overline{sf}}$ are all problem-dependent constants and are explicitly defined in the proof.  
\end{theorem}
\begin{proof}
We prove the first two general results with regard to $K_{inf}$ and $K_{fea}$ by discussing two cases: $(i)$ $\lim_{k\to\infty}\bar\tau_k = 0$ and $(ii)$ $\bar \tau_k \ge \bar \tau_{\min}$ for all $k\in\mathbb{N}$.  

When $\lim_{k\to\infty}\bar\tau_k = 0$, 
we use  $K_{inf}$ to denote the number of iterations required until \eqref{eq.complexity.inf} is  satisfied at an iteration. 
Since for all $k < K_{inf}$, $\|J_k^T  c_k\| > \epsilon$, it follows from Lemma~\ref{lemma.mathcalE.jc} that  $\|\bar J_k^T  \bar c_k\| \geq \|J_k^T  c_k\| - \mathcal{E}_{Jc} >  \epsilon- \mathcal{E}_{Jc} > \mathcal{E}_{Jc}$. Additionally, it follows from \eqref{eq.diff.Jc} and Cauchy-Schwarz inequality that  $\|\bar c_k\| \ge \|c_k\| - \epsilon_c \geq \tfrac{\|J_k^T  c_k\|}{\kappa_J} - \epsilon_c > \tfrac{2\mathcal{E}_{Jc}}{\kappa_J} - \epsilon_c \geq \epsilon_c \ge \epsilon_o$. The lower bounds $\alpha_{l0}^A$ and $\alpha_{l0}^L$ provided in  Lemmas~\ref{lemma.inexact.stepsize.opt} and~\ref{lemma.step.size.ls} are therefore valid with $\gamma =  \mathcal{E}_{Jc}$. It then follows from  \eqref{eq.Deltal.lb.NOLICQ}, \eqref{eq.diff.Jc.sto.det}, \eqref{def.merit}, and \eqref{eq.lemma.model.reduction.ada.NOLICQ}  that 
\begin{align*}
      \phi\left(x_{k}+ \bar \alpha_k \bar d_{k}, \bar \tau_{k}\right)-\phi\left(x_{k}, \bar \tau_{k}\right)  & \le - \eta \alpha_{l0}^i \Delta \lbar\left(x_k, \bar\tau_k, \bar d_k\right)  + \mathcal{E}^i \\ 
     &  \le - \eta \alpha_{l0}^i \sigma_c \sigma_v (\kappa_c + \epsilon_c)^{-1}  \|\bar J_k^T \bar c_k\|^2 + \mathcal{E}^i \\  &  \le - \eta \alpha_{l0}^i \sigma_c \sigma_v (\kappa_c + \epsilon_c)^{-1}  \|J_k^T c_k\|^2 + \eta \alpha_{l0}^i \sigma_c \sigma_v (\kappa_c + \epsilon_c)^{-1} \mathcal{E}_{Jc^1} + \mathcal{E}^i,    
\end{align*}
where $\mathcal{E}_{Jc^1}: =  ( (\kappa_J + \epsilon_J) \epsilon_c + \kappa_c \epsilon_J ) \left(  (\kappa_J + \epsilon_J)  (\kappa_c + \epsilon_c) + \kappa_J \kappa_c \right)$. The last inequality is satisfied since $| \|J_k^T c_k\|^2 -  \|\bar J_k^T \bar  c_k\|^2  | \le \mathcal{E}_{Jc^1}$ is proved in \eqref{eq.diff.Jc.sto.det}. It then follows that
\begin{align*}
&\  \phi\left(x_{k}+ \bar \alpha_k \bar d_{k}, \bar \tau_{k+1}\right)-\phi\left(x_{k}, \bar \tau_{k}\right) \\  \le &\ - \eta \alpha_{l0}^i \sigma_c \sigma_v (\kappa_c + \epsilon_c)^{-1}  \|J_k^T c_k\|^2 + \eta \alpha_{l0}^i \sigma_c \sigma_v (\kappa_c + \epsilon_c)^{-1} \mathcal{E}_{Jc^1} + \mathcal{E}^i  + (   \bar\tau_{k+1} -  \bar\tau_{k} ) f_{\inf}.  
\end{align*}
Summing the inequality above for all iterations within $\{0, 1, \dots , K_{inf} - 1\}$ and doing rearrangement, we have from~\eqref{eq.epsilon.4} and   \eqref{eq.phi.min} that  
   \begin{equation}\label{eq.final_thm_key_cond.4}
     \begin{split}
 &\  K_{inf} \eta \alpha_{l0}^i \sigma_c \sigma_v (\kappa_c + \epsilon_c)^{-1} \min_{j \in\{0,1,\cdots,K_{inf} - 1\}}   \|J_j^T c_j\|^2   \le \sum_{j=0}^{K_{inf} - 1} \eta \alpha_{l0}^i \sigma_c \sigma_v (\kappa_c + \epsilon_c)^{-1}  \|J_j^T c_j\|^2  \\  \le &\  K_{inf}  \left( \eta \alpha_{l0}^i \sigma_c \sigma_v (\kappa_c + \epsilon_c)^{-1} \mathcal{E}_{Jc^1} + \mathcal{E}^i \right) + \bar \tau_{-1} |f_{\inf}| + \phi(x_ 0 , \bar{\tau}_{0})  - \phi_{\min} 
 \\  \le &\  K_{inf} \eta \alpha_{l0}^i \sigma_c \sigma_v (\kappa_c + \epsilon_c)^{-1} \sigma_{\epsilon}^{-1} \epsilon^2 + \bar \tau_{-1} |f_{\inf}| + \phi(x_ 0 , \bar{\tau}_{0})  - \phi_{\min}. 
 \end{split}
 \end{equation}
   Since for all $k < K_{inf}$, we always have $\|J_k^T  c_k\| > \epsilon$. It then follows by \eqref{eq.final_thm_key_cond.4} that  
\begin{align*}
\tfrac{\kappa_c + \epsilon_c}{\eta \alpha_{l0}^i \sigma_c \sigma_v } >  \tfrac{K_{inf}} {\bar \tau_{-1} |f_{\inf}| + \phi(x_0,  \bar{\tau}_{0})  - \phi_{\min}  } (1- {\sigma_{\epsilon}^{-1}} ) \epsilon^2, 
\end{align*}
and we require at most  
$ K_{inf} =  \left\lceil \tfrac{ (\kappa_c + \epsilon_c)  ( \bar \tau_{-1} |f_{\inf}| + \phi(x_ 0 , \bar{\tau}_{0})  - \phi_{\min} )}{\eta \alpha_{l0}^i \sigma_c \sigma_v (1- {\sigma_{\epsilon}^{-1}} ) \epsilon^2} \right\rceil$ 
many iterations such that \eqref{eq.complexity.inf} is satisfied. 

When $\bar \tau_k \ge \bar \tau_{\min}$ for all $k\in\mathbb{N}$, 
it follows from \eqref{eq.convergence.theorem.ub.Jc} that 
\begin{align*}
     \min_{j \in \{0,\ldots,K_{inf}-1 \}}    \| J_j^T  c_j \|^2  &\le \tfrac{\kappa_c + \epsilon_c}{\sigma_{c}\sigma_v\eta\alpha_{l\tau}^i}  \mathcal{E}^i + 
 \mathcal{E}_{Jc^2}  + \tfrac{ ( \bar \tau_{-1} |f_{\inf}| + \phi(x_ 0 , \bar{\tau}_{0})  - \phi_{\min})(\kappa_c + \epsilon_c)  }{K_{inf}\sigma_{c}\sigma_v\eta\alpha_{l\tau}^i},  
\end{align*}
which implies by \eqref{eq.epsilon.4} that 
\begin{align*}
 \tfrac{ ( \bar \tau_{-1} |f_{\inf}| + \phi(x_ 0 , \bar{\tau}_{0})  - \phi_{\min})(\kappa_c + \epsilon_c)  }{K_{inf}\sigma_{c}\sigma_v\eta\alpha_{l\tau}^i}  >   {\epsilon^2} - \tfrac{\kappa_c + \epsilon_c}{\sigma_{c}\sigma_v\eta\alpha_{l\tau}^i}  \mathcal{E}^i - 
 \mathcal{E}_{Jc^2} > (1- {\sigma_{\epsilon}^{-1}} ) {\epsilon^2}, 
\end{align*}
and we require at most 
$ K_{inf} =  \left\lceil \tfrac{ (\kappa_c + \epsilon_c)  ( \bar \tau_{-1} |f_{\inf}| + \phi(x_ 0 , \bar{\tau}_{0})  - \phi_{\min} )}{\eta \alpha_{l\tau}^i \sigma_c \sigma_v (1- {\sigma_{\epsilon}^{-1}} ) \epsilon^2} \right\rceil$ 
iterations such that \eqref{eq.complexity.inf} is satisfied. The first result is therefore satisfied by combining the two cases.

If the singular values of $J_k$ are bounded below by $\kappa_{\sigma}$ for all $k\in \mathbb{N}$, the second desired conclusion (see $K_{fea}$ and \eqref{eq.complexity.fea}) is directly satisfied by $\|J_k^T  c_k\| \ge  \kappa_{\sigma} \|c_k\|$.

When $\bar \tau_k \ge \bar \tau_{\min}$ for all $k\in\mathbb{N}$, combining \eqref{eq.convergence.theorem.ub.Jc} with a rearrangement of 
\eqref{eq.convergence.theorem.ub.u} yields 
  \begin{equation}\label{eq.final_thm_key_cond}
     \begin{split}
      &\       \min_{{j\in\{0,\ldots,K_{si}-1\}}}  \left\| Z_j^T\left(\nabla f\left(x_j\right)+H_j v_j\right)\right\|^2  +  \| J_j^T  c_j \|^2 
  \\  \le &\ \left( \tfrac{\kappa_H^2 ( 4+ 2 \kappa_1) }{\bar\tau_{\min} \sigma_{u} \lambda_u \eta\alpha_{l\tau}^i} + \tfrac{ \kappa_H^2  \left(  \kappa_1 + 2 \kappa_2  \right)(\kappa_c + \epsilon_c)}{\sigma_{c}\sigma_v\eta\alpha_{l\tau}^i} + \tfrac{ \kappa_c + \epsilon_c}{\sigma_{c}\sigma_v\eta\alpha_{l\tau}^i} \right) \tfrac{\bar \tau_{-1} |f_{\inf}| + \phi(x_ 0 , \bar{\tau}_{0})  - \phi_{\min}}{K_{si}}\\
   &\ + \left(  \tfrac{\kappa_H^2 ( 4+ 2 \kappa_1)}{\bar\tau_{\min} \sigma_{u} \lambda_u \eta\alpha_{l\tau}^i}  + \tfrac{\kappa_H^2 \left(  \kappa_1 + 2 \kappa_2  \right)(\kappa_c + \epsilon_c)}{\sigma_{c}\sigma_v\eta\alpha_{l\tau}^i} + \tfrac{\kappa_c + \epsilon_c}{\sigma_{c}\sigma_v\eta\alpha_{l\tau}^i}   \right) \mathcal{E}^i  \\
   &\ + (2\kappa_H^2 \kappa_2 +1 )\mathcal{E}_{Jc^2} +  \kappa_H^2 \kappa_1  (\kappa_J+\epsilon_J)^2 \epsilon_c^2 + 12  \kappa_H^2  \zeta^{-2} \left(  \epsilon_g^2+ \kappa_{\bar y^*}^2 \epsilon_J^2  + \kappa_{y^*}^2 \epsilon_J^2  \right).   
     \end{split} 
 \end{equation}
Let $ \kappa_{si} =  \tfrac{\kappa_H^2 ( 4+ 2 \kappa_1) }{\bar\tau_{\min} \sigma_{u} \lambda_u \eta\alpha_{l\tau}^i} + \tfrac{ 
     \left(\kappa_H^2  \left(  \kappa_1 + 2 \kappa_2  \right) + 1 \right)(\kappa_c + \epsilon_c)}{\sigma_{c}\sigma_v\eta\alpha_{l\tau}^i}$. 
     We use  $K_{si}$ to denote the number of iterations required until \eqref{eq.complexity.si} is satisfied at an iteration.  
  Since for all $k < K_{si}$, $\left\| Z_k^T\left(\nabla f\left(x_k\right)+H_k v_k\right)\right\|^2  +  \| J_k^T  c_k \|^2 > \epsilon^2$, it follows by 
  \eqref{eq.epsilon.1} and \eqref{eq.final_thm_key_cond} that 
\begin{align*}
\kappa_{si} > \tfrac{K_{si} } {\bar \tau_{-1} |f_{\inf}| + \phi(x_ 0 , \bar{\tau}_{0})  - \phi_{\min}  } (1- {\sigma_{\epsilon}^{-1}} ) \epsilon^2, 
\end{align*}
and the third desired conclusion (see $K_{si}$ and \eqref{eq.complexity.si}) is satisfied. 
 If the singular values of $J_k$ are bounded below by $\kappa_{\sigma}$   for all $k\in \mathbb{N}$, combining~\eqref{eq.convergence.theorem.ub.u} with a rearrangement of~\eqref{eq.convergence.theorem.ub.c2} yields \begin{equation}\label{eq.final_thm_key_cond.2}
     \begin{split}
      &\       \min_{{j\in\{0,\ldots,K_{sf}-1\}}}  \left\| Z_j^T\left(\nabla f\left(x_j\right)+H_j v_j\right)\right\|^2  +  \| c_j \|^2 
  \\  \le &\ \left( \tfrac{\kappa_H^2 ( 4+ 2 \kappa_1) }{\bar\tau_{\min} \sigma_{u} \lambda_u \eta\alpha_{l\tau}^i} + \tfrac{ \kappa_H^2  \left(  \kappa_1 + 2 \kappa_2  \right)(\kappa_c + \epsilon_c)}{\sigma_{c}\sigma_v\eta\alpha_{l\tau}^i } + \tfrac{ \kappa_c + \epsilon_c}{\sigma_{c}\sigma_v\eta\alpha_{l\tau}^i\kappa_{\sigma}^2} \right) \tfrac{\bar \tau_{-1} |f_{\inf}| + \phi(x_ 0 , \bar{\tau}_{0})  - \phi_{\min}}{K_{sf}}\\
   &\ + \left(  \tfrac{\kappa_H^2 ( 4+ 2 \kappa_1)}{\bar\tau_{\min} \sigma_{u} \lambda_u \eta\alpha_{l\tau}^i}  + \tfrac{\kappa_H^2 \left(  \kappa_1 + 2 \kappa_2  \right)(\kappa_c + \epsilon_c)}{\sigma_{c}\sigma_v\eta\alpha_{l\tau}^i} + \tfrac{\kappa_c + \epsilon_c}{\sigma_{c}\sigma_v\eta\alpha_{l\tau}^i\kappa_{\sigma}^2}   \right) \mathcal{E}^i  \\
   &\ + (2\kappa_H^2 \kappa_2 +\kappa_{\sigma}^{-2} )\mathcal{E}_{Jc^2} +  \kappa_H^2 \kappa_1  (\kappa_J+\epsilon_J)^2 \epsilon_c^2 + 12  \kappa_H^2  \zeta^{-2} \left(  \epsilon_g^2+ \kappa_{\bar y^*}^2 \epsilon_J^2  + \kappa_{y^*}^2 \epsilon_J^2  \right).   
     \end{split} 
 \end{equation}
 Let $ \kappa_{sf} =  \tfrac{\kappa_H^2 ( 4+ 2 \kappa_1) }{\bar\tau_{\min} \sigma_{u} \lambda_u \eta\alpha_{l\tau}^i} + \tfrac{ 
     \left(\kappa_{\sigma}^2\kappa_H^2  \left(  \kappa_1 + 2 \kappa_2  \right) + 1 \right)(\kappa_c + \epsilon_c)}{\sigma_{c}\sigma_v\eta\alpha_{l\tau}^i\kappa_{\sigma}^2}$. We use  $K_{sf}$ to denote the number of iterations required until \eqref{eq.complexity.sf} is  satisfied at an iteration. 
   Since for all $k < K_{sf}$, $\left\| Z_k^T\left(\nabla f\left(x_k\right)+H_k v_k\right)\right\|^2  +  \| c_k \|^2 > \epsilon^2$, it follows by 
   \eqref{eq.epsilon.2} and \eqref{eq.final_thm_key_cond.2} that 
\begin{align*}
\kappa_{sf} >  \tfrac{K_{sf} } {\bar \tau_{-1} |f_{\inf}| + \phi(x_ 0 , \bar{\tau}_{0})  - \phi_{\min}  } (1- {\sigma_{\epsilon}^{-1}} ) \epsilon^2, 
\end{align*}
and the fourth desired conclusion (see $K_{sf}$ and \eqref{eq.complexity.sf}) is satisfied. Finally, if the singular values of $\bar J_k$ are bounded below by $\kappa_{\sigma}$ for all $k\in \mathbb{N}$, by~\eqref{eq.convergence.theorem.ub.u} and \eqref{eq.convergence.theorem.ub.c}, we have
   \begin{equation}\label{eq.final_thm_key_cond.3}
     \begin{split}
      &\       \min_{{j\in\{0,\ldots,K_{\overline{sf}}-1\}}}  \left\| Z_j^T\left(\nabla f\left(x_j\right)+H_j v_j\right)\right\|^2  +  \| c_j \|  
  \\  \le &\ \left( \tfrac{\kappa_H^2 ( 4+ 2 \kappa_1) }{\bar\tau_{\min} \sigma_{u} \lambda_u \eta\alpha_{l\tau}^i} + \tfrac{ \kappa_H^2  \left(  \kappa_1 + 2 \kappa_2  \right)(\kappa_c + \epsilon_c)}{\sigma_{c}\sigma_v\eta\alpha_{l\tau}^i } + \tfrac{1}{\sigma_{c}\sigma_v\eta\alpha_{l\tau}^i\kappa_{\sigma}^2} \right) \tfrac{\bar \tau_{-1} |f_{\inf}| + \phi(x_ 0 , \bar{\tau}_{0})  - \phi_{\min}}{K_{\overline{sf}}}\\
   &\ + \left(  \tfrac{\kappa_H^2 ( 4+ 2 \kappa_1)}{\bar\tau_{\min} \sigma_{u} \lambda_u \eta\alpha_{l\tau}^i}  + \tfrac{\kappa_H^2 \left(  \kappa_1 + 2 \kappa_2  \right)(\kappa_c + \epsilon_c)}{\sigma_{c}\sigma_v\eta\alpha_{l\tau}^i} + \tfrac{1}{\sigma_{c}\sigma_v\eta\alpha_{l\tau}^i\kappa_{\sigma}^2}   \right) \mathcal{E}^i  \\
   &\ + 2\kappa_H^2 \kappa_2\mathcal{E}_{Jc^2} + 2\epsilon_c +  \kappa_H^2 \kappa_1  (\kappa_J+\epsilon_J)^2 \epsilon_c^2 + 12  \kappa_H^2  \zeta^{-2} \left(  \epsilon_g^2+ \kappa_{\bar y^*}^2 \epsilon_J^2  + \kappa_{y^*}^2 \epsilon_J^2  \right).   
     \end{split} 
 \end{equation}
 Let $\kappa_{\overline{sf}} =   \tfrac{\kappa_H^2 ( 4+ 2 \kappa_1) }{\bar\tau_{\min} \sigma_{u} \lambda_u \eta\alpha_{l\tau}^i} + \tfrac{ 
     \kappa_{\sigma}^2 \kappa_H^2  \left( \kappa_1 + 2 \kappa_2  \right) (\kappa_c + \epsilon_c) + 1}{\sigma_{c}\sigma_v\eta\alpha_{l\tau}^i\kappa_{\sigma}^2}$.  We use  $K_{\overline{sf}}$ to denote the number of iterations required until \eqref{eq.complexity.sf.bar} is  satisfied at an iteration. 
    Since for all $k < K_{\overline{sf}}$, we always have $\left\| Z_k^T\left(\nabla f\left(x_k\right)+H_k v_k\right)\right\|^2  +  \| c_k \| > \epsilon^2$, it follows by 
    \eqref{eq.epsilon.3} and \eqref{eq.final_thm_key_cond.3} that 
\begin{align*}
\kappa_{\overline{sf}} >  \tfrac{K_{\overline{sf}} } {\bar \tau_{-1} |f_{\inf}| + \phi(x_ 0 , \bar{\tau}_{0})  - \phi_{\min}  } (1- {\sigma_{\epsilon}^{-1}} ) \epsilon^2, 
\end{align*}
and the final desired conclusion (see $K_{\overline{sf}}$ and \eqref{eq.complexity.sf.bar}) is satisfied. 
\end{proof}

\begin{remark}
Theorem~\ref{theorem.tau.lb.complexity} provides worst-case iteration complexity guarantees for Algorithm \ref{alg.dfo_sqp_LS} under the conditions that $\epsilon$ is sufficiently large compared to $\mathcal{O}(( \epsilon_g + \epsilon_c + \epsilon_J)^{1/2})$ where $\{\epsilon_g,\epsilon_c,\epsilon_J\} \subset (0,1]$. The results of the theorem can be summarized as follows: $(i)$ 
 In the most general case,  Algorithm \ref{alg.dfo_sqp_LS}  requires $\mathcal{O}(\epsilon^{-2})$ iterations to achieve the desired accuracy in terms of the infeasible  stationarity  error; $(ii)$  If additionally $\{J_k\}$ is with full rank, then Algorithm \ref{alg.dfo_sqp_LS} can also achieve the desired feasibility error in at most $\mathcal{O}(\epsilon^{-2})$ iterations;  
$(iii)$ When $\bar \tau_k \ge \bar \tau_{\min}$ for all $k\in\mathbb{N}$, Algorithm \ref{alg.dfo_sqp_LS}  requires $\mathcal{O}(\epsilon^{-2})$  iterations to achieve both the desired accuracy in terms of stationarity error and infeasible  stationarity  error.   
$(iv)$ When $\bar \tau_k \ge \bar \tau_{\min}$ and $\{J_k\}$ or $\{\bar J_k\}$ are of full rank, then Algorithm \ref{alg.dfo_sqp_LS} can also achieve the desired accuracy in terms of stationarity error and feasibility error in at most $\mathcal{O}(\epsilon^{-2})$ iterations. The only difference is that when  $\{\bar J_k\}$ is of full rank, the desired accuracy in terms of feasibility error can be improved (see \eqref{eq.complexity.sf.bar}).    
The target accuracies for each case are consistent with Table~\ref{tab:measure_table}. 
In the setting with no noise, i.e., $\epsilon_f = \epsilon_g = \epsilon_c = \epsilon_J = 0$, the conditions on $\epsilon$ only require that these accuracy parameters are positive, which could be arbitrarily close to zero. Under the LICQ, 
the algorithm recovers the complexity results of the deterministic SQP method \cite{curtis2024worst}. 
\end{remark}

\section{Numerical Experiments}\label{sec.experiments}

In this section, we present numerical results that showcase the performance of the proposed algorithm applied to standard equality-constrained nonlinear optimization problems from the CUTEst collection \cite{gratton2024s2mpj}, with different levels of manually introduced noise. We selected a subset of equality-constrained optimization problems from the CUTEst collection where the objective function is non-constant and the LICQ holds at all iterations for the deterministic line search SQP method \cite{NoceWrig06}. This resulted in 88 problems. The performance of our algorithm is evaluated under two scenarios: $(i)$ when the LICQ is satisfied and $(ii)$ when the LICQ is violated. The latter scenario is carefully induced to ensure fair experiments.

To simulate noise, we add uniformly distributed random perturbations to the exact function and gradient evaluations of both the objective and constraints for each problem. This results in the noisy quantities  $(\bar f_k, \bar g_k, \bar c_k, \bar J_k)$, which are used in the algorithms. Specifically, 
for $e_f \sim \mathcal{U}\left(-\epsilon_f, \epsilon_f\right)$, $e_g \sim \mathcal{U}\left(-\tfrac{\epsilon_g}{\sqrt{n}} \mathbf{1}_n, \tfrac{\epsilon_g}{\sqrt{n}} \mathbf{1}_n \right)$, $e_c \sim \mathcal{U}\left(-\tfrac{\epsilon_c}{\sqrt{m}} \mathbf{1}_m, \tfrac{\epsilon_c}{\sqrt{m}}\mathbf{1}_m \right)$,  and $e_J \sim \mathcal{U}\left(- \tfrac{\epsilon_J}{\sqrt{mn}} \mathbf{1}_{m \times n} , \tfrac{\epsilon_J}{\sqrt{mn}} \mathbf{1}_{m \times n}  \right)$,  we set
\begin{align*}
    \bar{f}(x) =f(x)+e_f, \quad \bar{c}(x)=c(x)+e_c, \quad \bar{g}(x) =g(x)+e_g, \quad \text{and} \quad \bar{J}(x)=J(x)+e_J, 
\end{align*} 
where $\{\epsilon_f,\epsilon_g,\epsilon_c,\epsilon_J\}\subset\R{}_{>0}$ dictate the level of the noise. We consider four different
noise levels in the objective function and constraint functions, $\epsilon_f \in \{10^{-1}, 10^{-2}, 10^{-4},  10^{-8}\}$ and $\epsilon_c \in \{10^{-1}, 10^{-2}, 10^{-4},  10^{-8}\}$, and all combinations, and set $\epsilon_g = \sqrt{\epsilon_f}$ and $\epsilon_J = \sqrt{\epsilon_c}$ to align the gradient noise levels with the upper bounds established for common DFO methods \cite{larson2019derivative}.




To implement the adaptive step size scheme (\adasqp{}), we set $\eta = 0.5$ and $\beta = 1$ and estimate the Lipschitz constants of the objective gradient and the constraint Jacobian using gradient differences near the initial point and keep these estimates constant. As for the line search scheme (\lssqp{}), we set $ \alpha_u^L = 1$, $\eta = 10^{-3}$, and $\nu = 0.5$. The other parameters were set as $\bar \tau_{-1} = 1$, $\lambda_{u} = 5 \times 10^{-9}$, $\sigma_{Jc} = 10^{2}$, $\sigma_{u} = 0.99$, $\sigma_{c} = 0.1$, $\sigma_{r} = 0.9999$, $\sigma_{\tau} = 10^{-2}$, $\xi_{-1} = 1$, $\chi_{-1} = 10^{-3}$, $\zeta_{-1} = 10^3$, and $\theta = 10^4$, to ensure that the parameters interfered as minimally as possible with the algorithm. We set $\epsilon_o = 0$ for the pessimistic algorithms (\adasqppes{} and \lssqppes{}) and set $\epsilon_o = \epsilon_c$ for the optimistic algorithms (\adasqpopt{} and \lssqpopt{}). 
We employ the CG-Steihaug method \cite{NoceWrig06} to solve \eqref{eq.nolicq.tr.sub.p1} potentially inexactly; the method is terminated when 
\eqref{eq.nolicq.tr.cauchy.decrease} holds and the residual vectors, defined as $\bar  R_{k,t} =\bar J_k^T \bar J_k \bar v_{k,t}  + \bar J_k^T \bar c_k$ where $ \bar v_{k,t}$ denotes the $t$th iteration of the CG-Steihaug algorithm, satisfy 
\begin{equation}
\label{eq.inexact.terminate.v}
    \left\|\bar R_{k,t} \right\|_{\infty} \le \kappa_v \min \{ \epsilon_c , \epsilon_f\} \max \{1, \|\bar J_k^T \bar c_k \|_{\infty} \}.
\end{equation}
We employ the MINRES method \cite{choi2011minres} to potentially solve \eqref{eq.nolicq.sub.p2} inexactly; the method 
is terminated when at least one of Termination Test~\ref{tt1} or \ref{tt2} holds and the 
residual vectors $\bar \rho_{k,t}$ and $ \bar r_{k,t}$ (defined in \eqref{eq.nolicq.sub.p2.linear.system.res} where the second subscript defines the iterations of the MINRES procedure) satisfy 
 \begin{equation}
 \label{eq.inexact.terminate.u_old}
    \left\|\left[\begin{array}{c}
\bar\rho_{k,t} \\
\bar r_{k,t}
\end{array}\right]\right\|_{\infty} \leq \kappa_u \min \{ \epsilon_c , \epsilon_f\}  \max \left\{ \min \left\{  \max \left\{ \|\bar  u_{k}\|_{\infty}, \|\bar J_k^T \bar c_k\|_{\infty} \right\}, 10^2 \right\} , 10^{-2} \right\} ,   
\end{equation}
Compared to \eqref{eq.tt1.cond1}, we add a $\max \{ \ \cdot \ , 10^{-2} \}$ term to ensure that the right-hand-side is not too small. 
For the inexact algorithms (\adasqpinexact{} and \lssqpinexact{}), we use $\kappa_u = \kappa_v = 10^{-2}$ in \eqref{eq.inexact.terminate.v} and \eqref{eq.inexact.terminate.u_old}, and for the exact variants (\adasqpexact{} and \lssqpexact{}) we replace $\kappa_v \min \{ \epsilon_c , \epsilon_f\}$  in  \eqref{eq.inexact.terminate.v}  and  $\kappa_u \min \{ \epsilon_c , \epsilon_f\}$ in \eqref{eq.inexact.terminate.u_old} with $10^{-10}$.  

\begin{figure}[htbp]
\centering    
\includegraphics[width=0.35\textwidth,clip=true,trim=10 180 50 150]{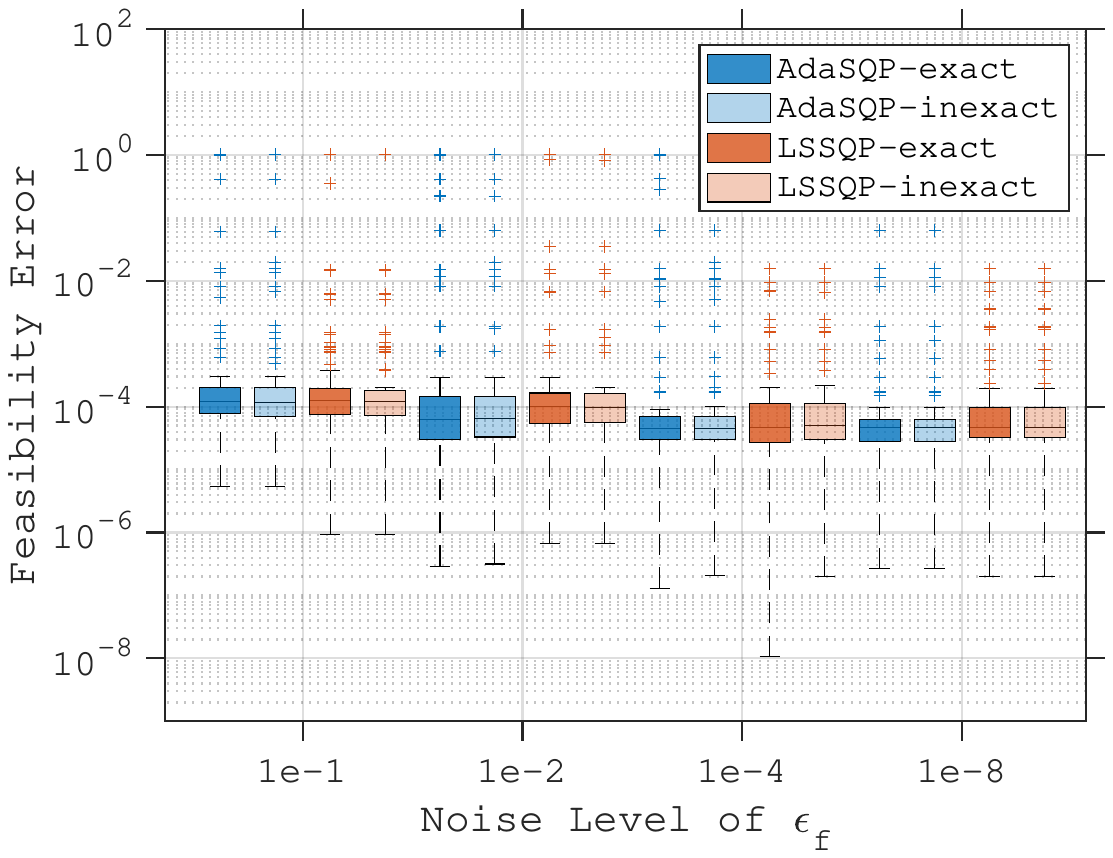}
\includegraphics[width=0.35\textwidth,clip=true,trim=10 180 50 150]{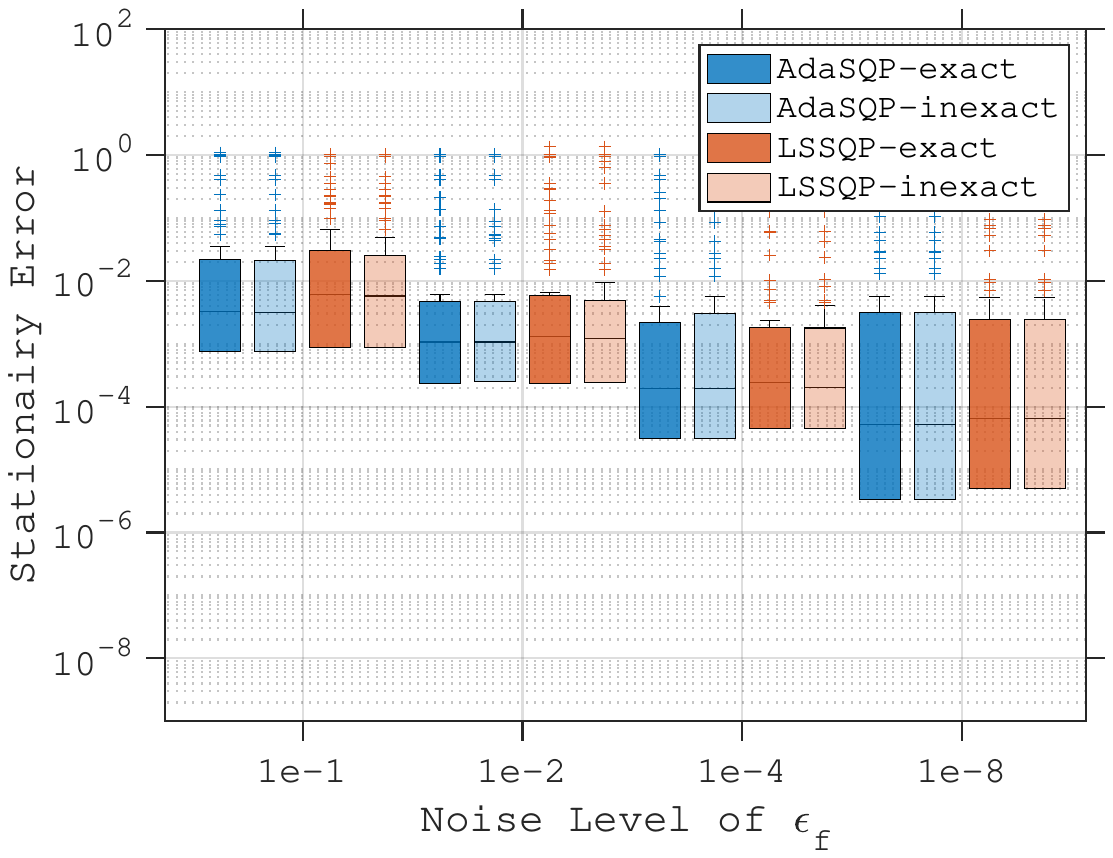}
\vspace{-0.25cm}

\includegraphics[width=0.35\textwidth,clip=true,trim=10 180 50 150]{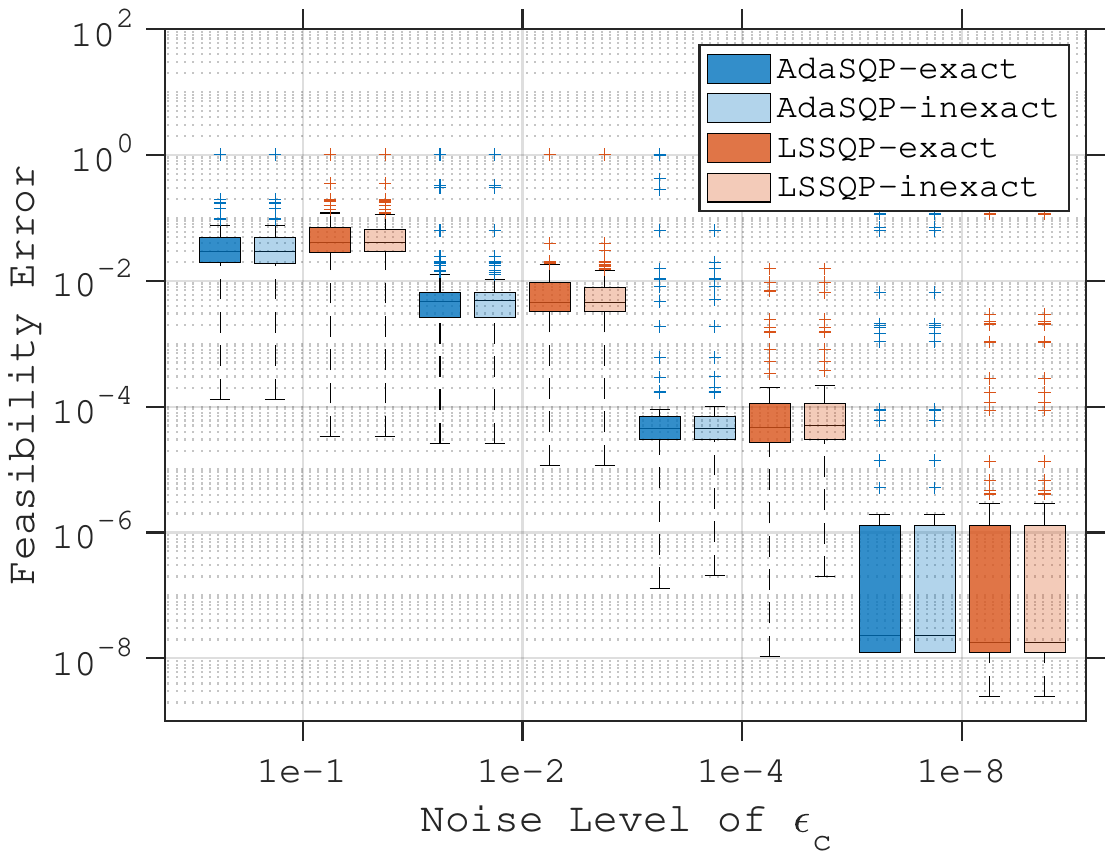}
\includegraphics[width=0.35\textwidth,clip=true,trim=10 180 50 150]{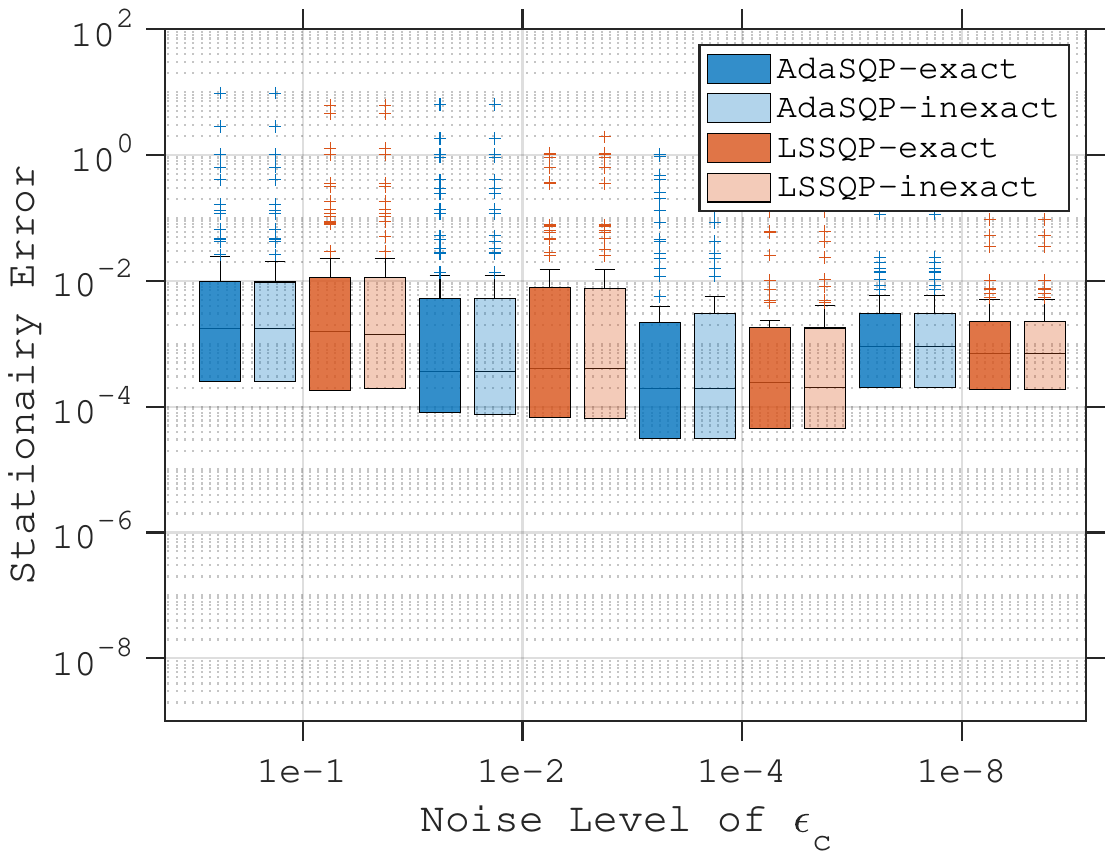}
\vspace{-0.15cm}

\caption{Box plots illustrating feasibility errors (left) and stationarity errors (right) for exact and inexact \adasqpopt{} and \lssqpopt{} methods on CUTEst problems that satisfy the LICQ. First row: $\epsilon_c = 10^{-4}$ and $\epsilon_f \in \{10^{-1}, 10^{-2}, 10^{-4},  10^{-8}\}$. Second row: $\epsilon_f = 10^{-4}$ and $\epsilon_c \in \{10^{-1}, 10^{-2}, 10^{-4},  10^{-8}\}$.}
\label{fig.boxplot.LICQ}
\end{figure}

We first test the performance of \adasqpopt{} and \lssqpopt{}
under different noise levels and exactness levels. We start by investigating the numerical performance of the proposed algorithms on the original CUTEst problem instances where the LICQ is satisfied. The feasibility and stationarity errors, defined as $\|c_k\|_{\infty}$ and $\| \nabla f(x_k) + J_k^T y_k \|_{\infty}$, respectively, where $y_k$ is a least-squares multiplier obtained using the true gradient $\nabla f(x_k)$ and the true Jacobian $J_k$, are reported based on the \emph{best iterate} from a maximum of 1000 iterations. The best iterate was selected as follows. For each problem instance, we selected the iterate with the lowest stationarity error, from the subset of iterates that satisfied the feasibility condition $\|c_k\|_{\infty} \leq 2 \max \{\epsilon_c, \epsilon_f\}$. If no iterate satisfied the feasibility condition, we instead reported the iterate with the smallest feasibility error. If the algorithm terminates early because of $\|\bar c_k\| \le \epsilon_o$ and $ \Delta \lbar(x_{k},\bar \tau_{k}, \bar d_{k}) \le  \epsilon_o$, the \emph{best iterate} is selected similarly, but only from the iterates up to the termination. It can be observed from Figure \ref{fig.boxplot.LICQ} that the performance of the inexact method is comparable to that of the exact method for both step size schemes. When $\epsilon_c$ is fixed, the feasibility error remains at a consistent level, while smaller values of $\epsilon_f$ lead to lower stationarity errors. Similarly, when $\epsilon_f$ is fixed, the stationarity error remains stable and the feasibility error scales with the noise level $\epsilon_c$.

\begin{figure}[htbp]
\centering    
\includegraphics[width=0.32\textwidth,clip=true,trim=10 180 50 150]{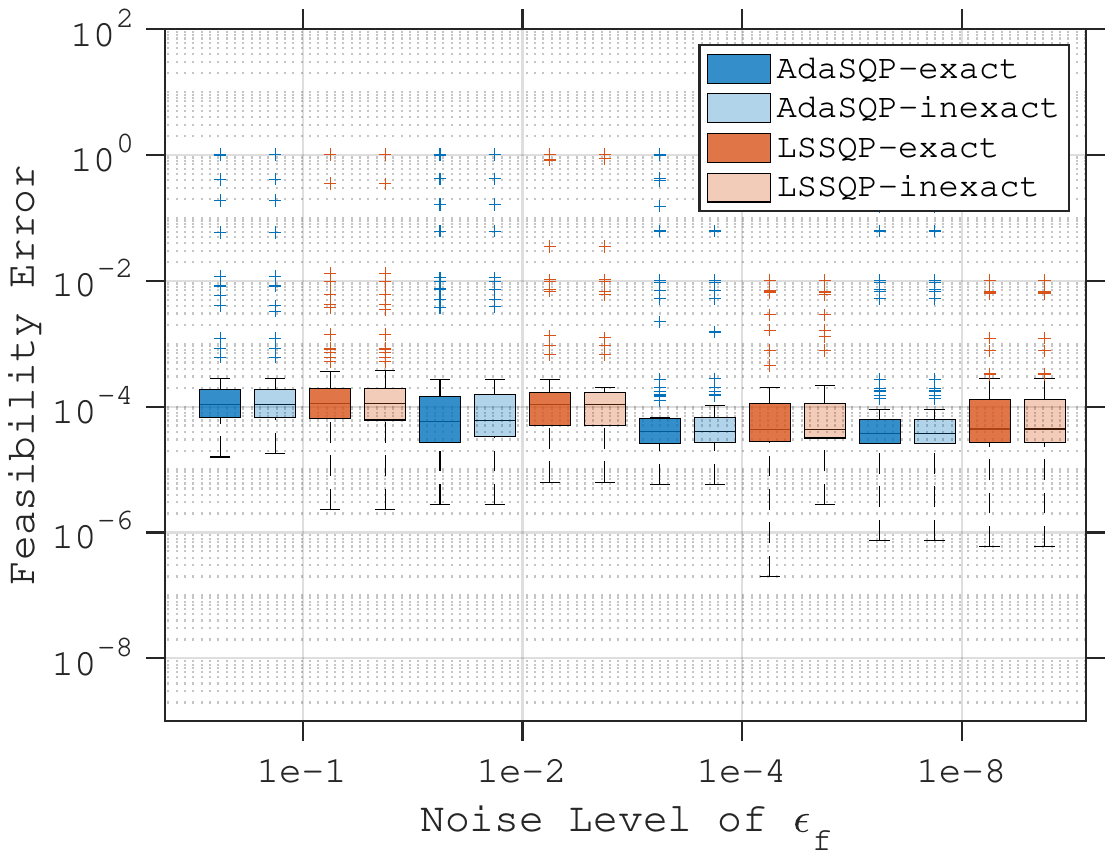}
\includegraphics[width=0.32\textwidth,clip=true,trim=10 180 50 150]{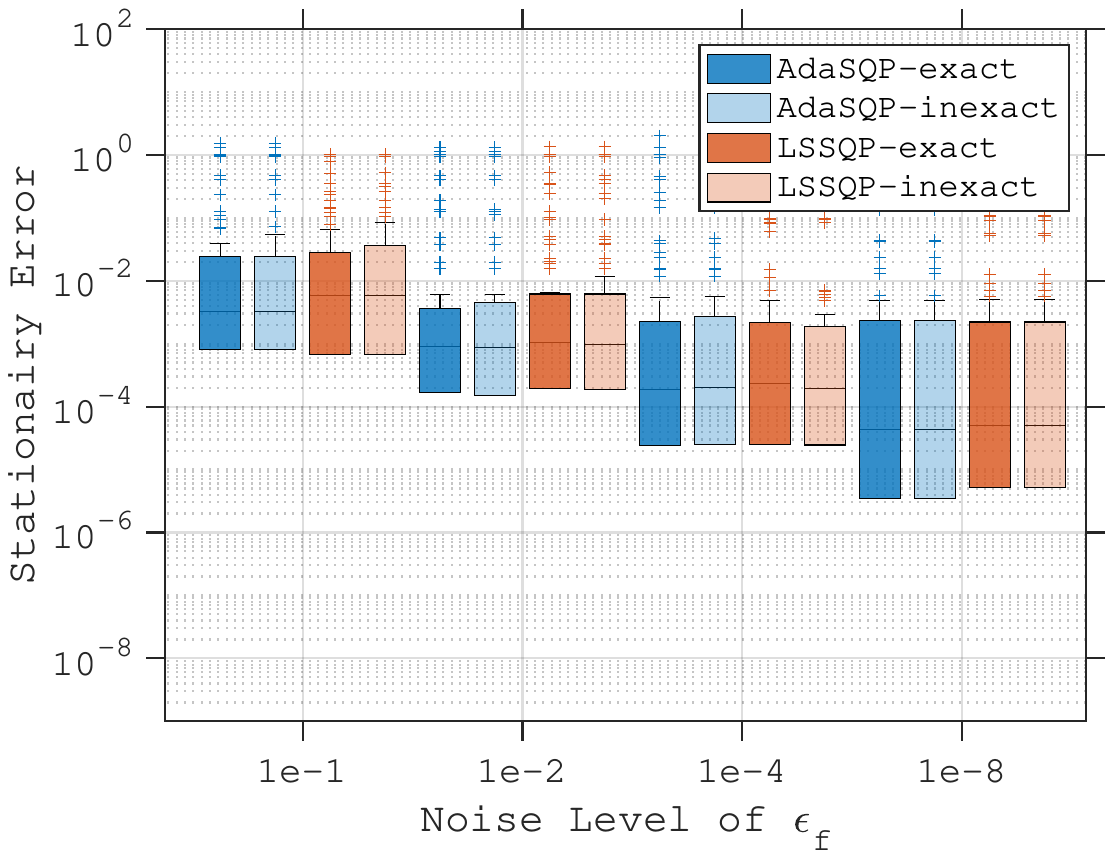}
\includegraphics[width=0.32\textwidth,clip=true,trim=10 180 50 150]{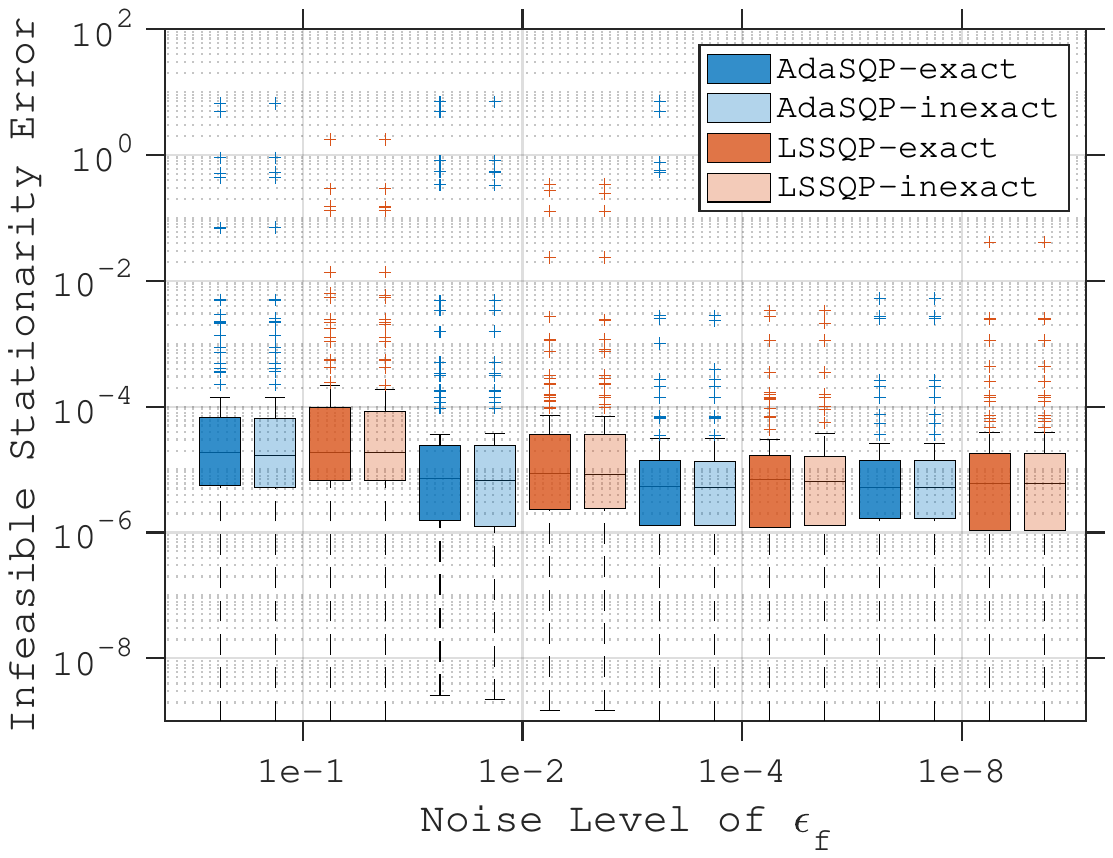}
\vspace{-0.25cm}

\includegraphics[width=0.32\textwidth,clip=true,trim=10 180 50 150]{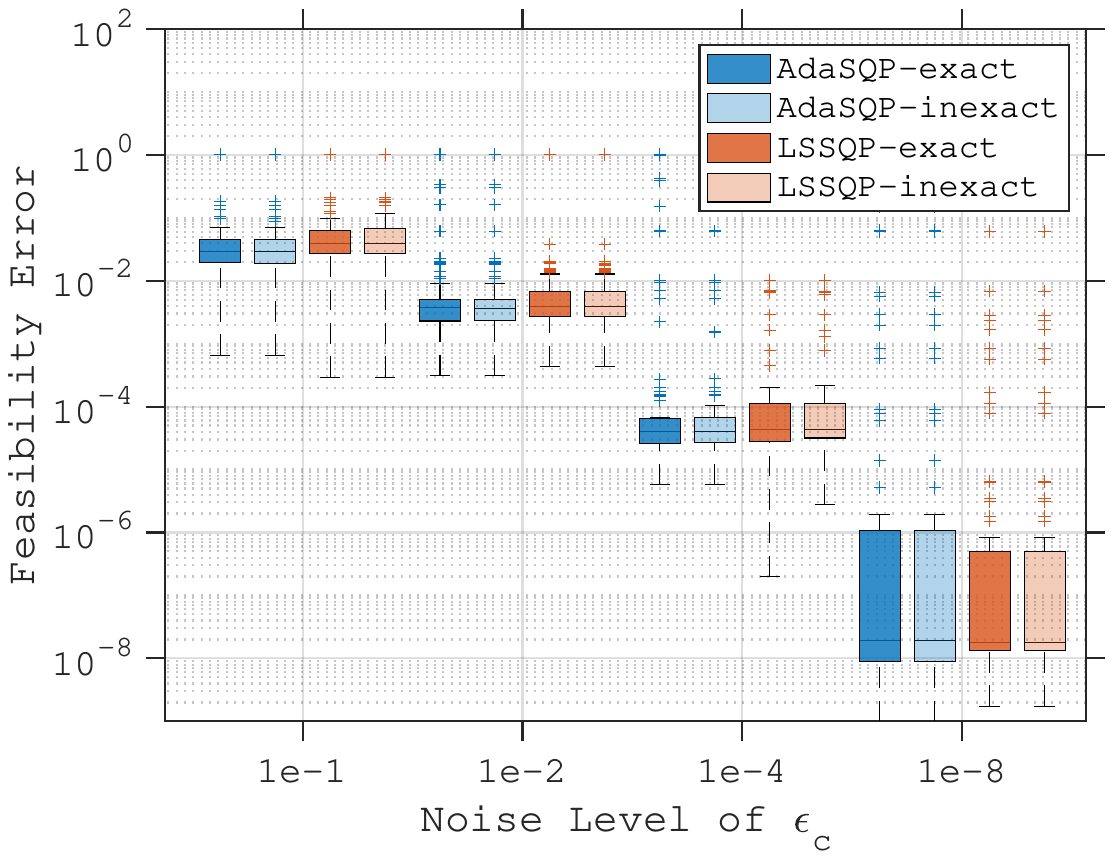}
\includegraphics[width=0.32\textwidth,clip=true,trim=10 180 50 150]{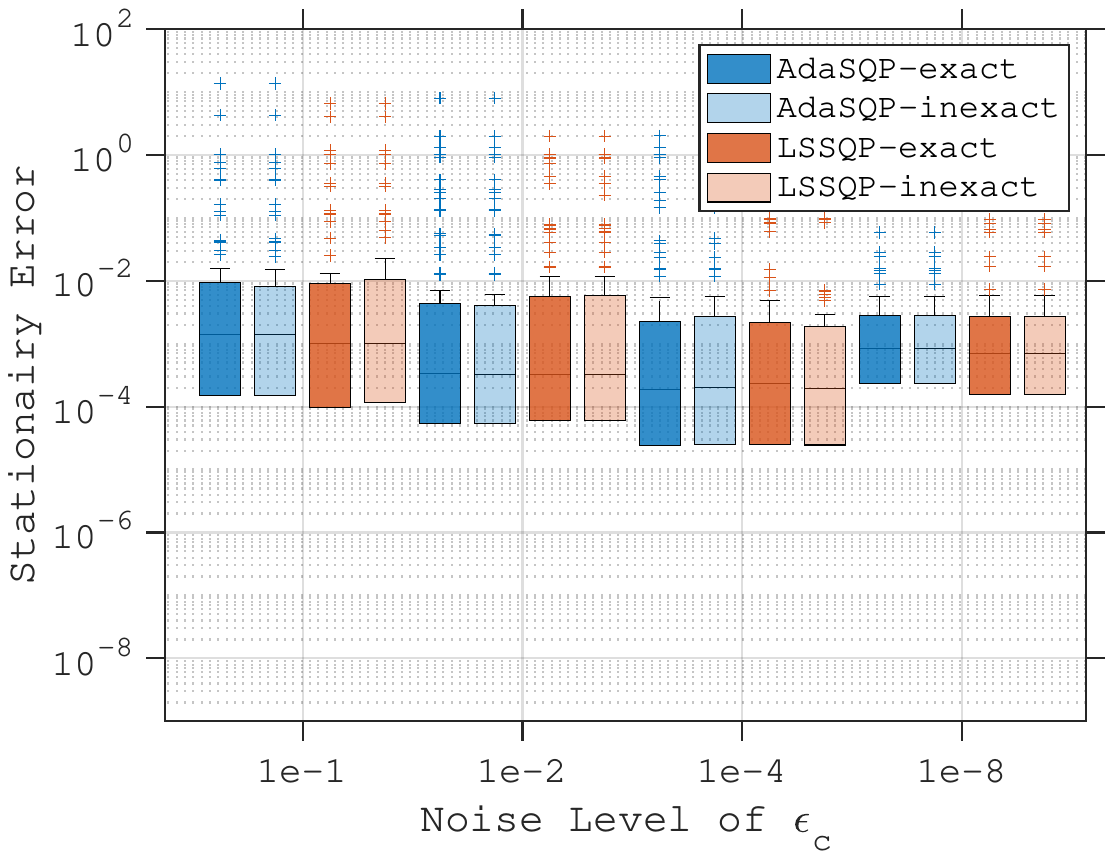}
\includegraphics[width=0.32\textwidth,clip=true,trim=10 180 50 150]{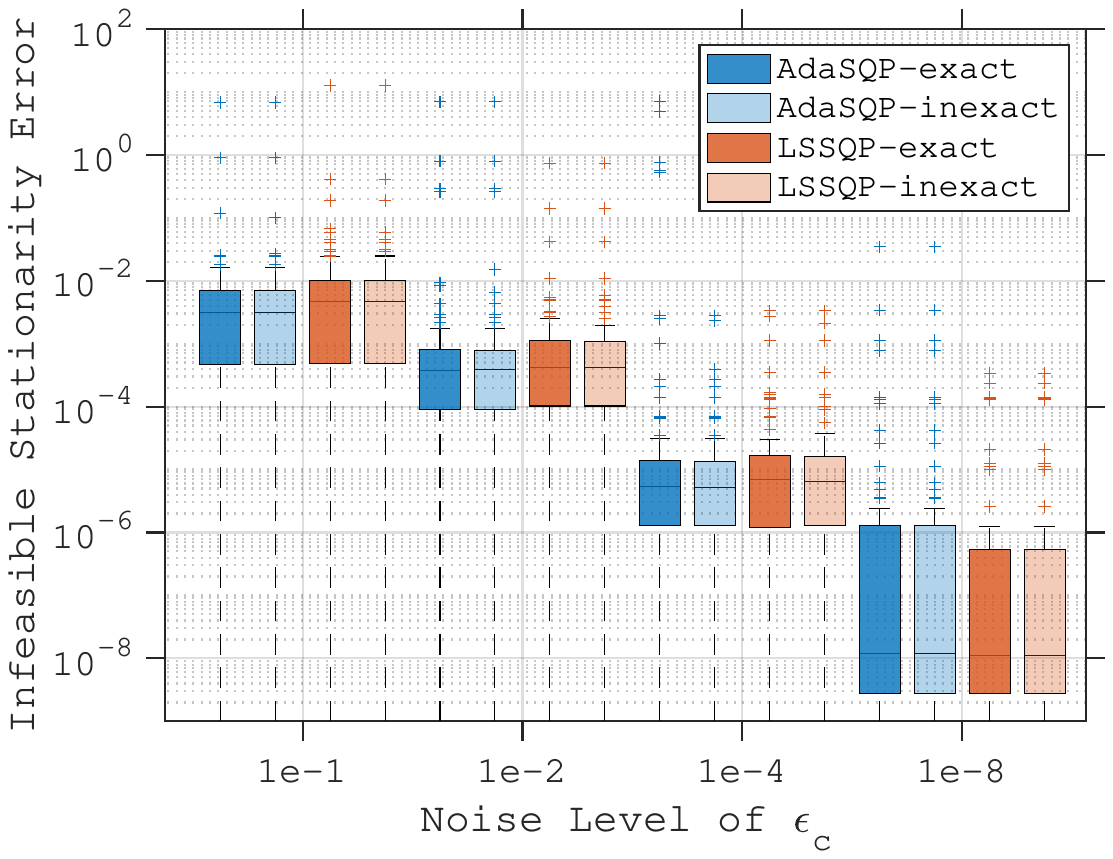}
\vspace{-0.15cm}

\caption{Box plots illustrating feasibility errors (left), stationarity errors (middle) and infeasible stationarity error (right) for exact and inexact \adasqpopt{} and \lssqpopt{} methods on CUTEst problems that violate the LICQ. First row: $\epsilon_c = 10^{-4}$ and $\epsilon_f \in \{10^{-1}, 10^{-2}, 10^{-4},  10^{-8}\}$. Second row: $\epsilon_f = 10^{-4}$ and $\epsilon_c \in \{10^{-1}, 10^{-2}, 10^{-4},  10^{-8}\}$.}
\label{fig.boxplot.NoLICQ}
\end{figure}

Next, we consider the same problems, noise levels and inexactness but in the setting for which the LICQ is violated. To do this in a controlled manner and without affecting the feasible region, we duplicated the last constraint (after adding noise) for each instance in the original set. For problems where the LICQ is violated, we adopt a similar approach to report the best iterate. In addition, we report the infeasible stationarity error $\|J_k^T c_k\|_{\infty}$ at that iterate. 
Figure \ref{fig.boxplot.NoLICQ} indicates that our method remains effective even when the LICQ is not satisfied. Furthermore, the best infeasible stationarity error decreases as $\epsilon_c$ is reduced.

\newpage

\begin{figure}[htbp]
        \centering    
\includegraphics[width=0.24\textwidth,clip=true,trim=10 180 50 150]{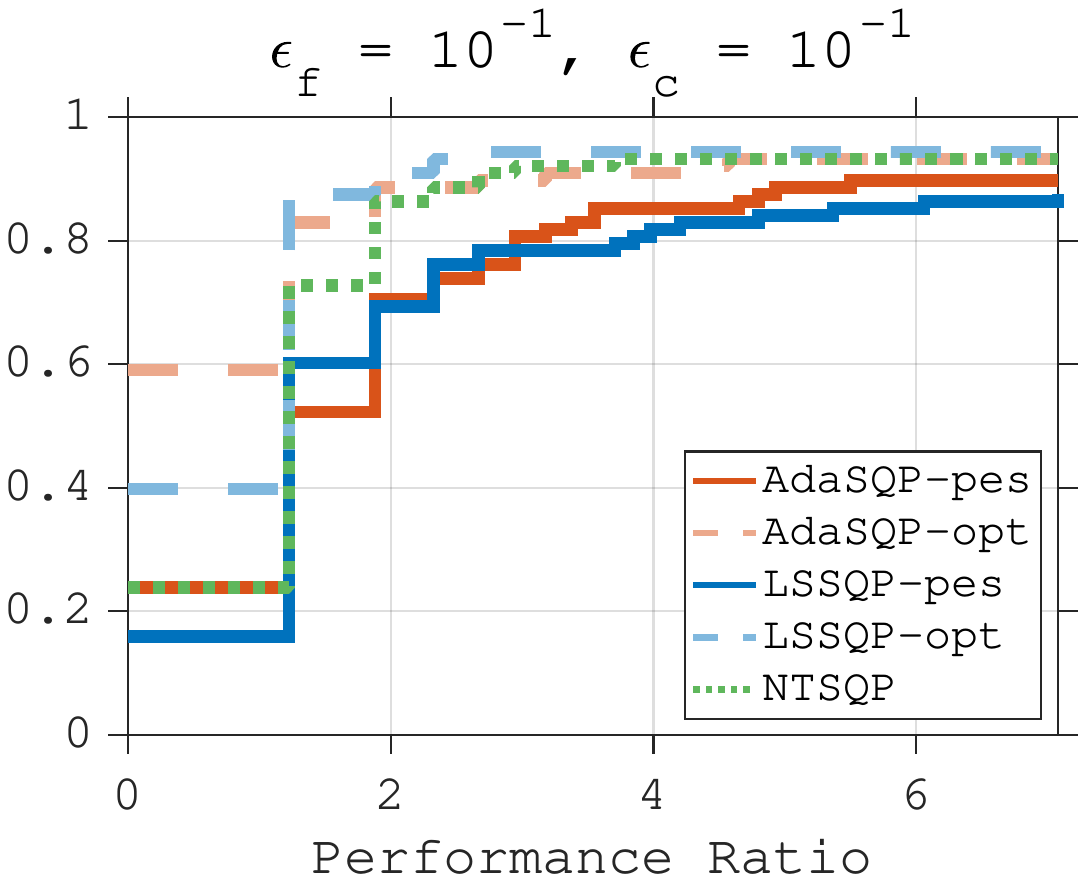}
\includegraphics[width=0.24\textwidth,clip=true,trim=10 180 50 150]{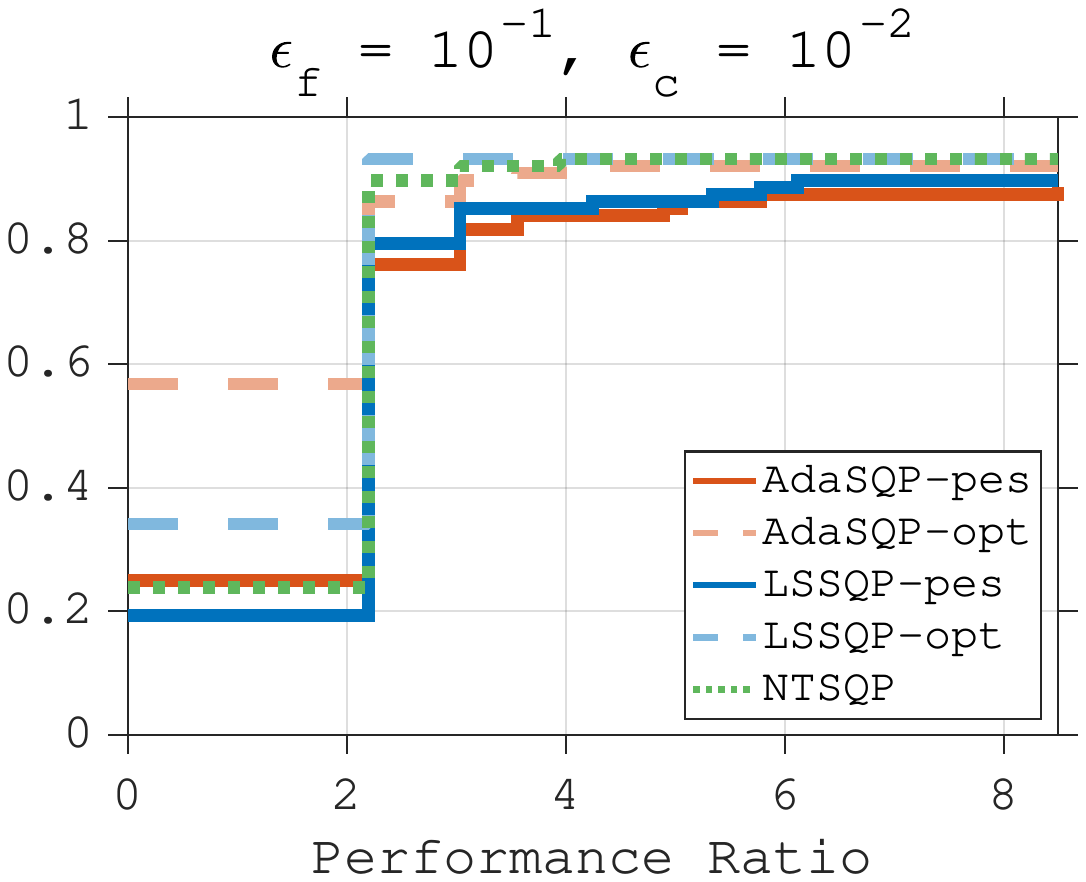}
\includegraphics[width=0.24\textwidth,clip=true,trim=10 180 50 150]{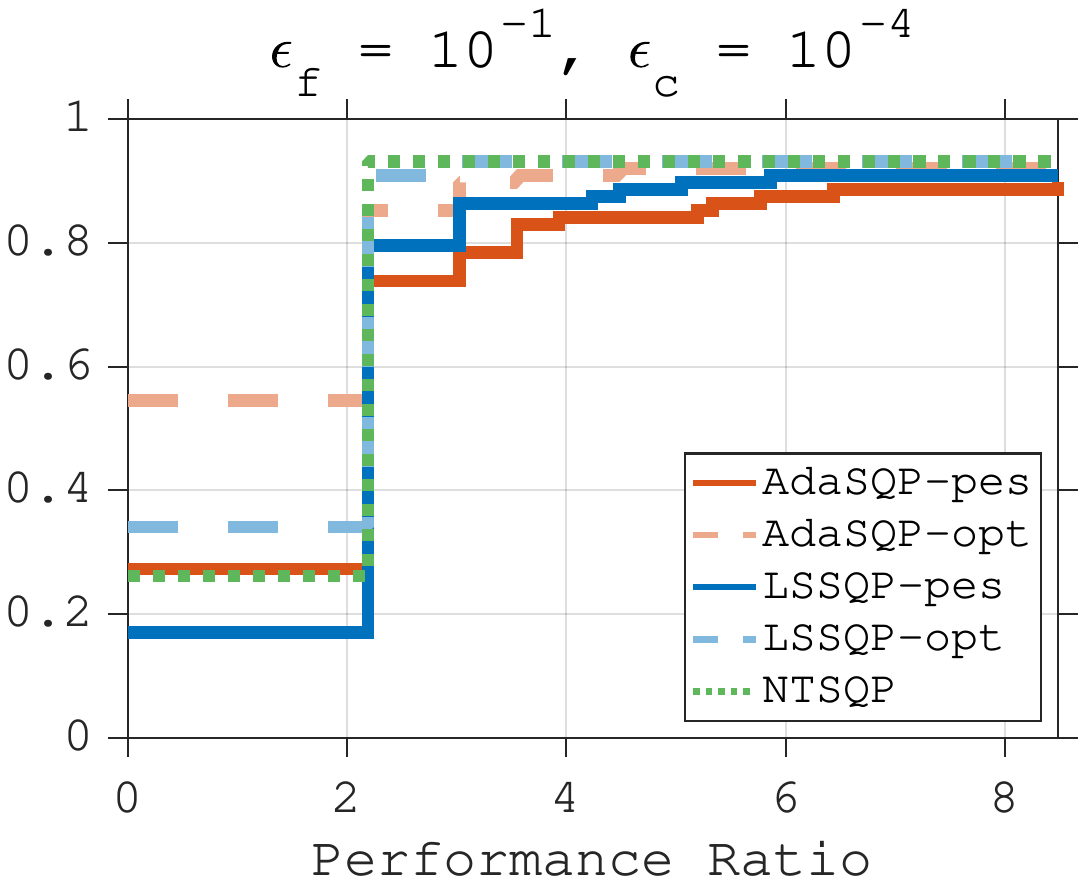}
\includegraphics[width=0.24\textwidth,clip=true,trim=10 180 50 150]{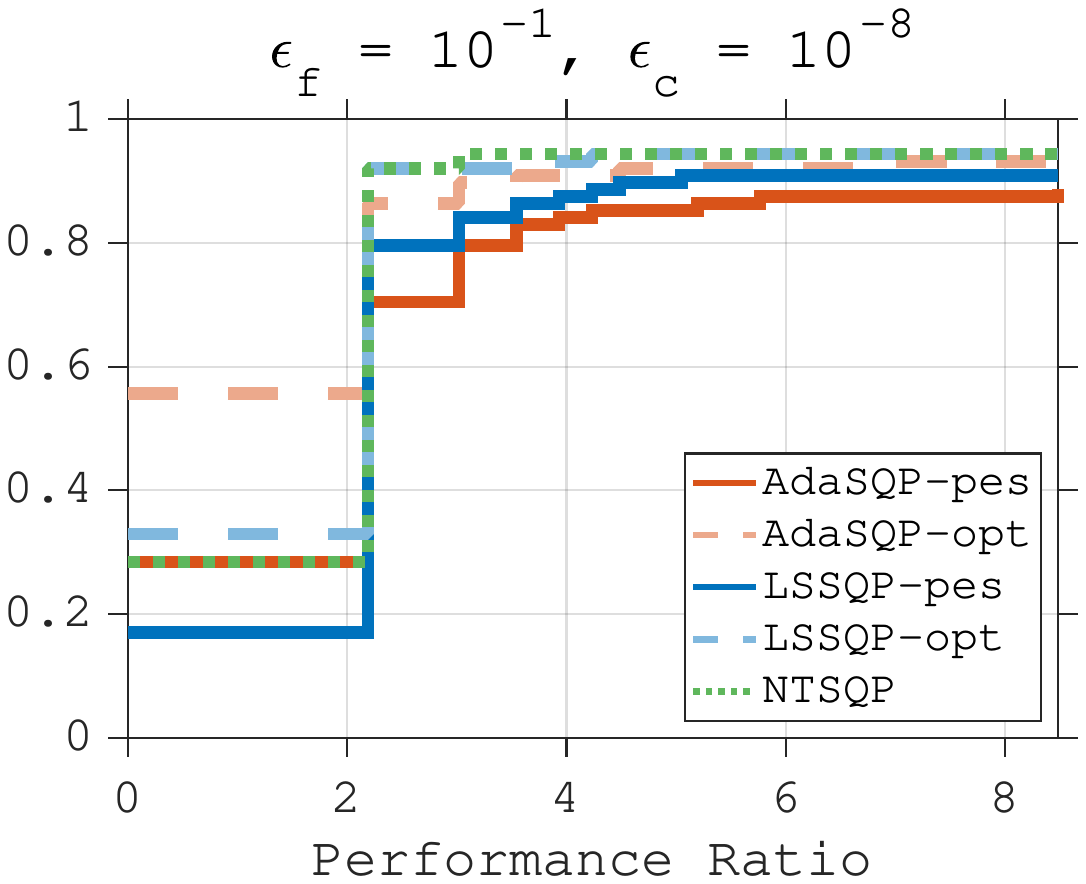}
\includegraphics[width=0.24\textwidth,clip=true,trim=10 180 50 150]{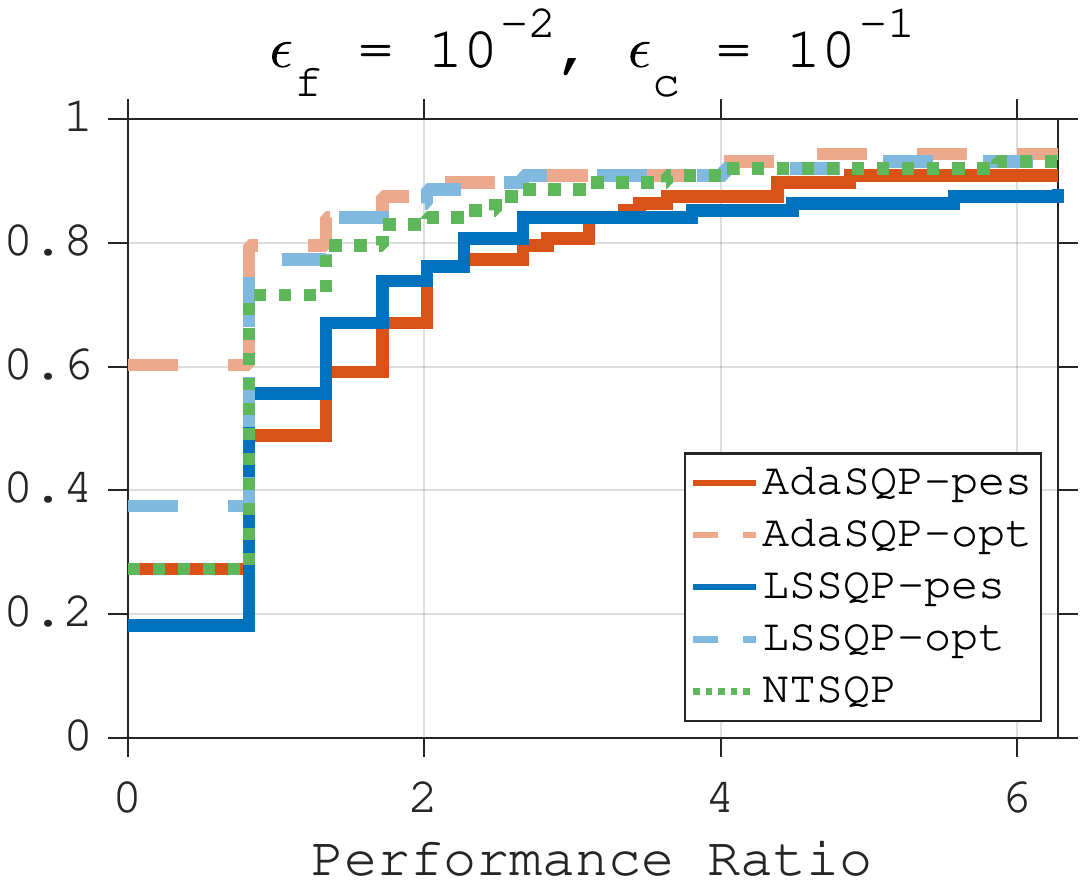}
\includegraphics[width=0.24\textwidth,clip=true,trim=10 180 50 150]{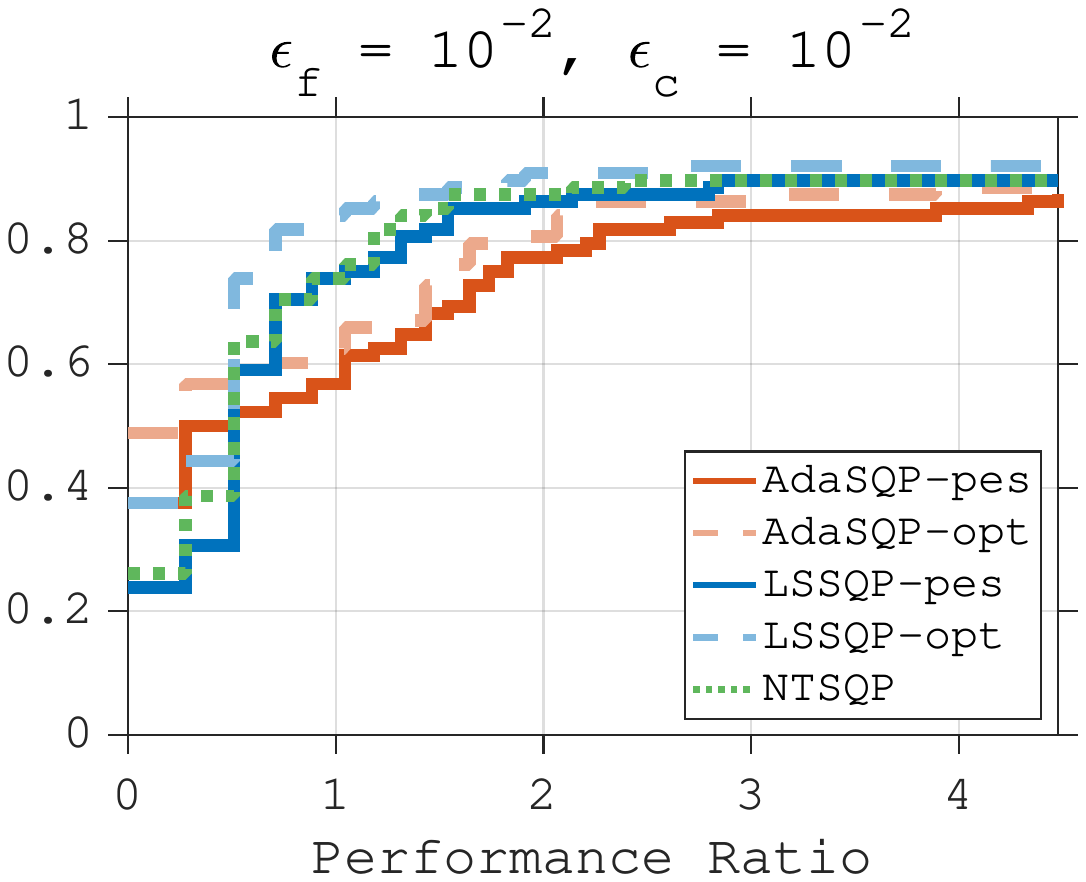}
\includegraphics[width=0.24\textwidth,clip=true,trim=10 180 50 150]{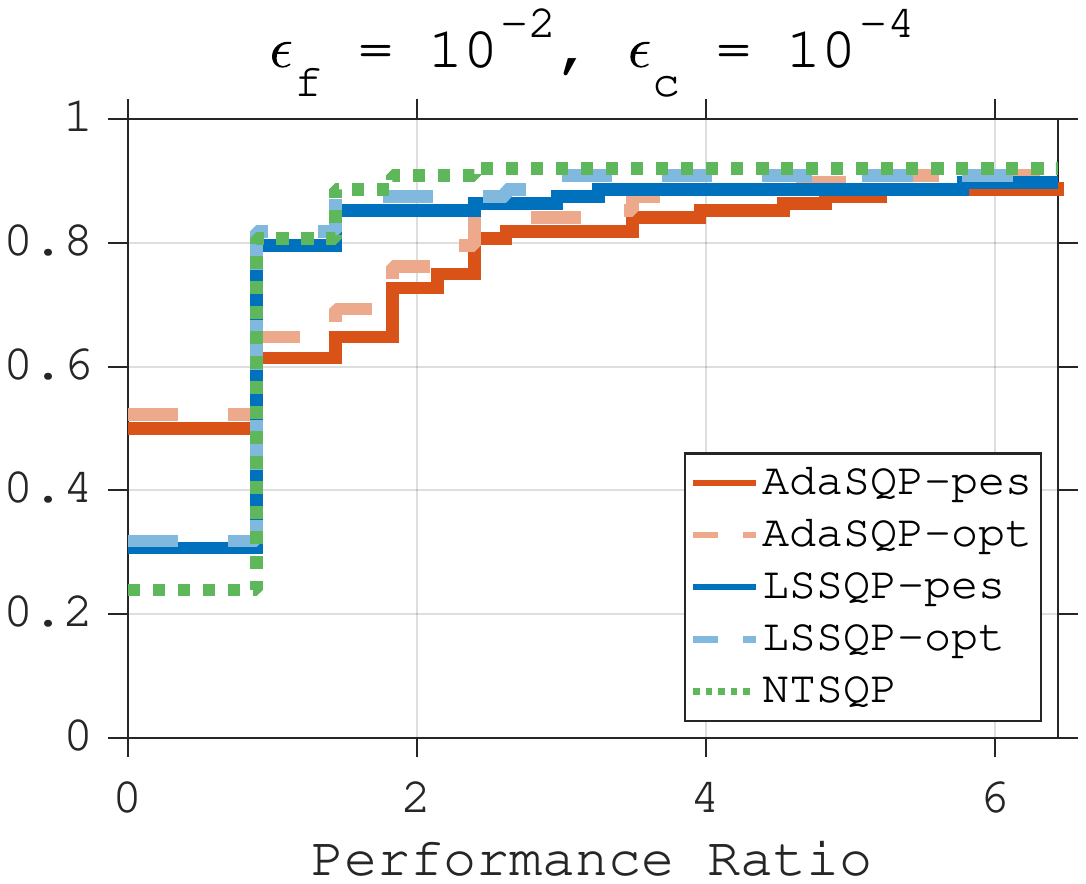}
\includegraphics[width=0.24\textwidth,clip=true,trim=10 180 50 150]{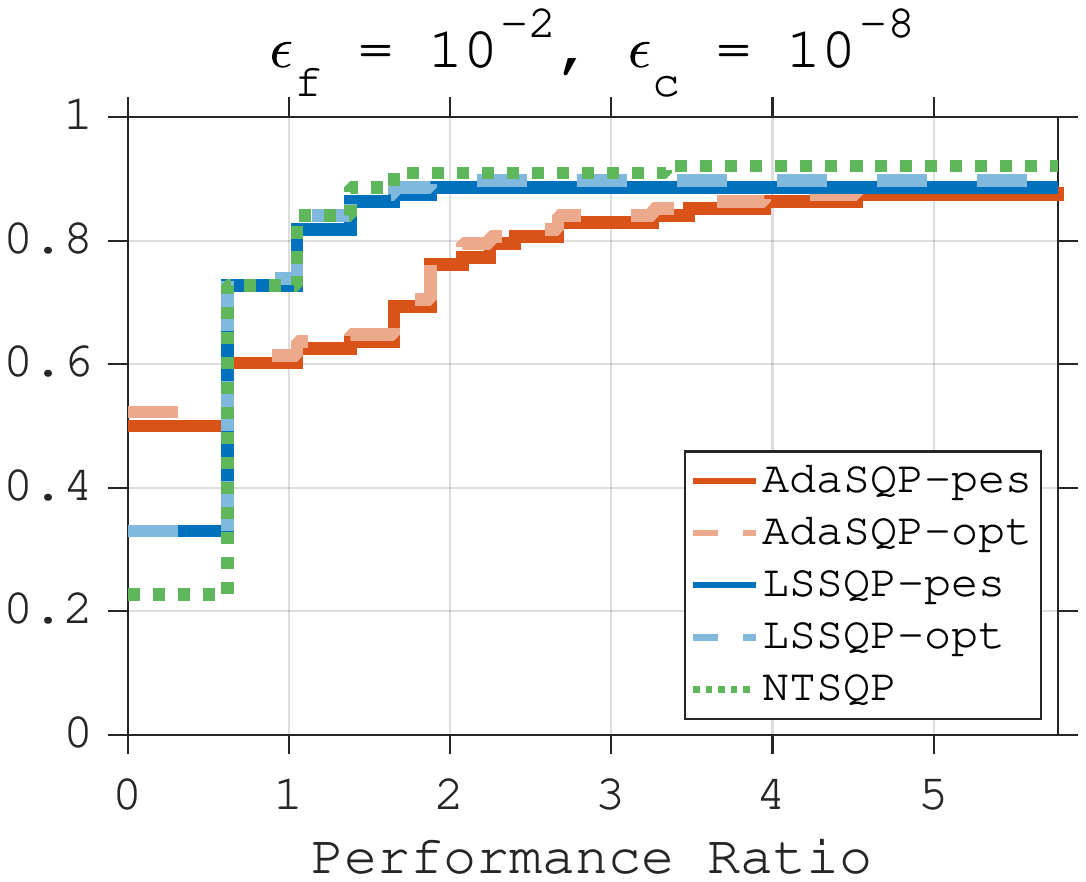}
\includegraphics[width=0.24\textwidth,clip=true,trim=10 180 50 150]{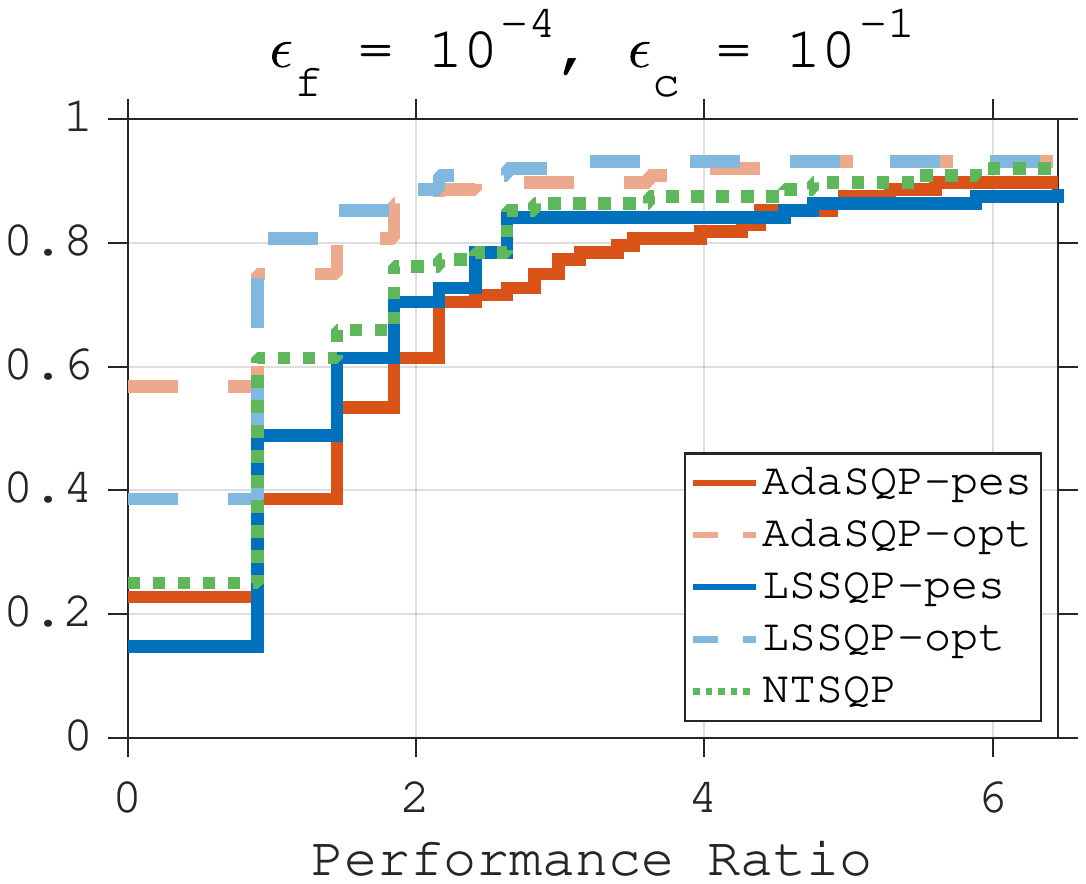}
\includegraphics[width=0.24\textwidth,clip=true,trim=10 180 50 150]{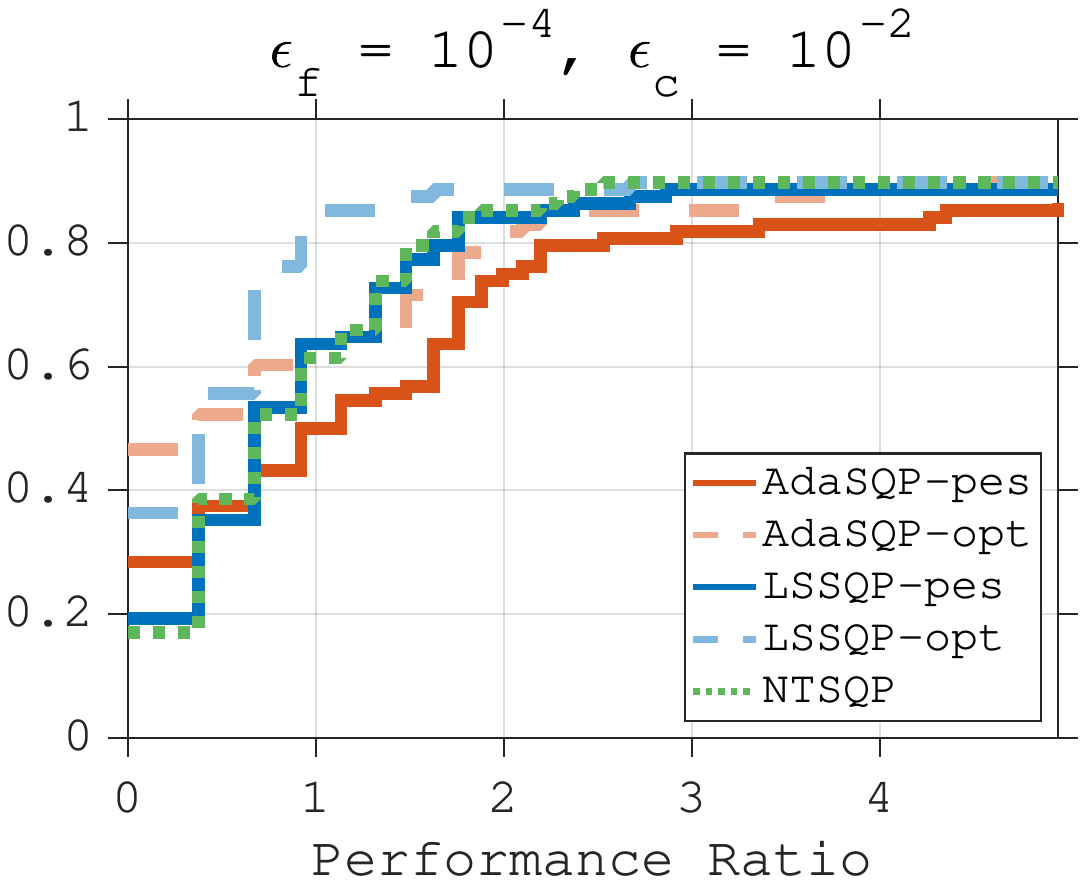}
\includegraphics[width=0.24\textwidth,clip=true,trim=10 180 50 150]{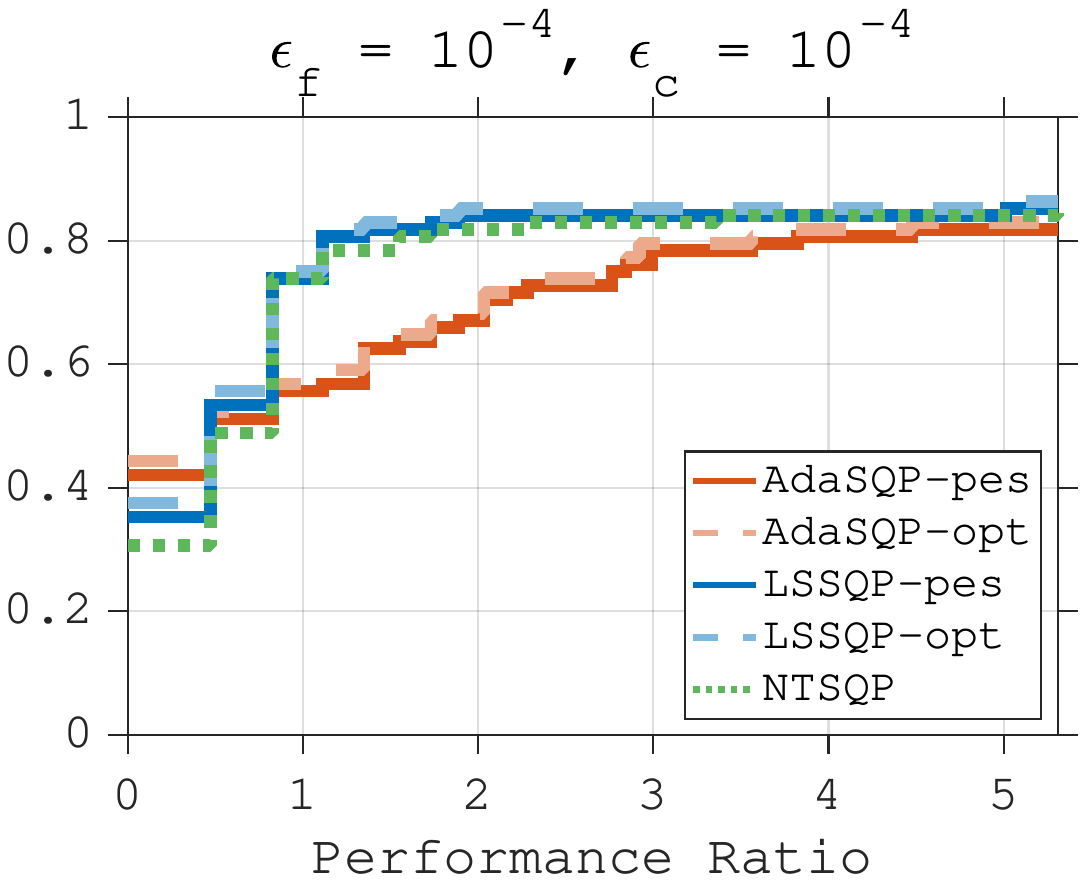}
\includegraphics[width=0.24\textwidth,clip=true,trim=10 180 50 150]{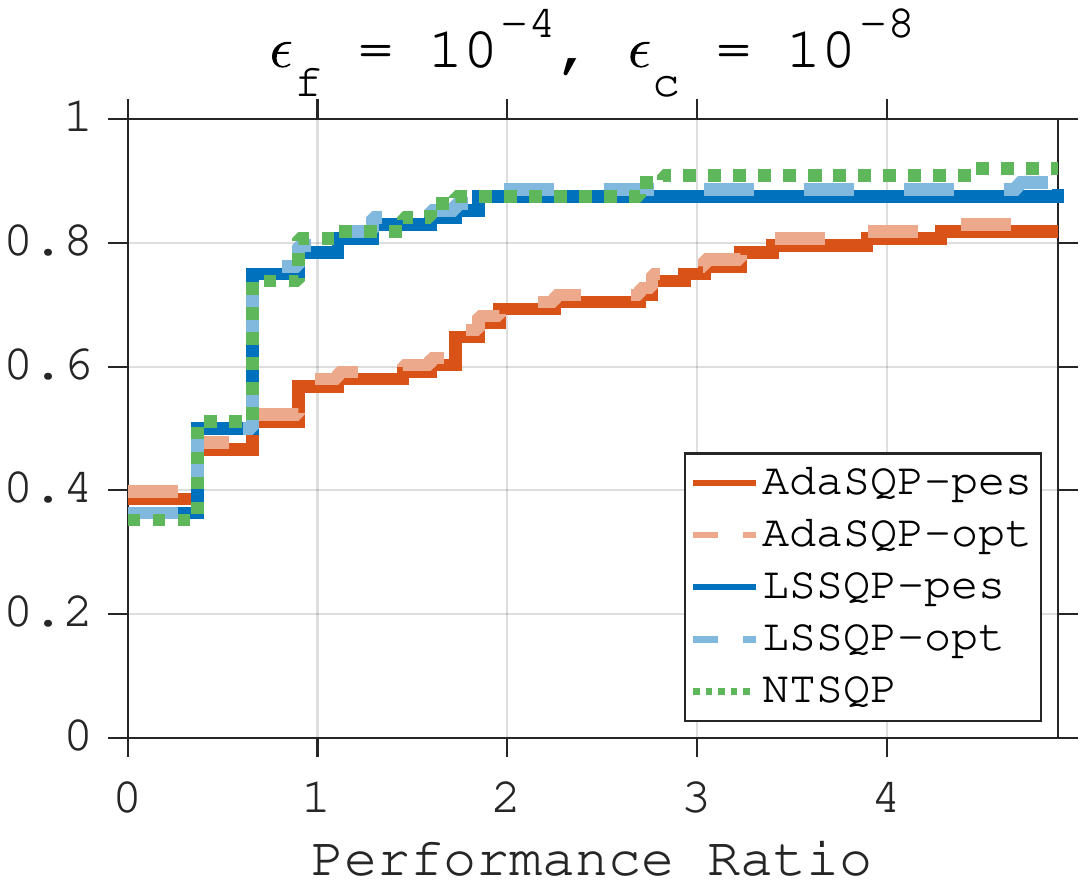}
\includegraphics[width=0.24\textwidth,clip=true,trim=10 180 50 150]{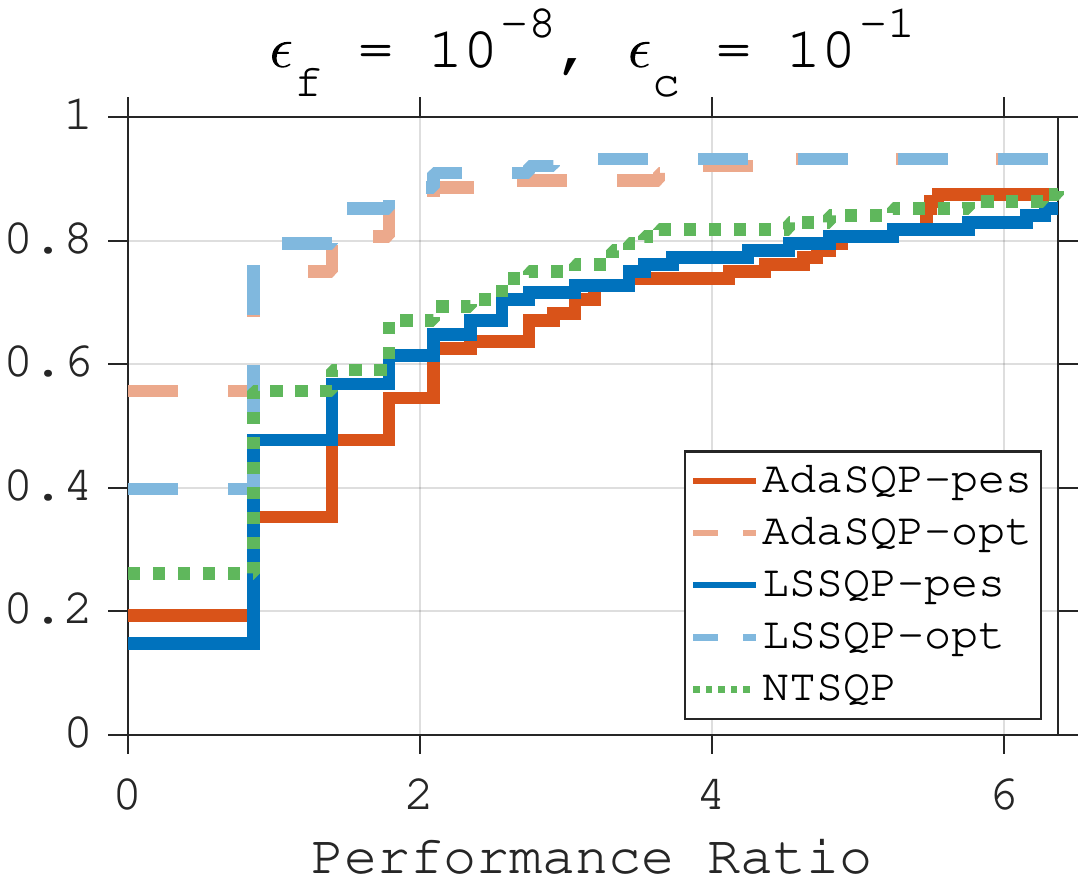}
\includegraphics[width=0.24\textwidth,clip=true,trim=10 180 50 150]{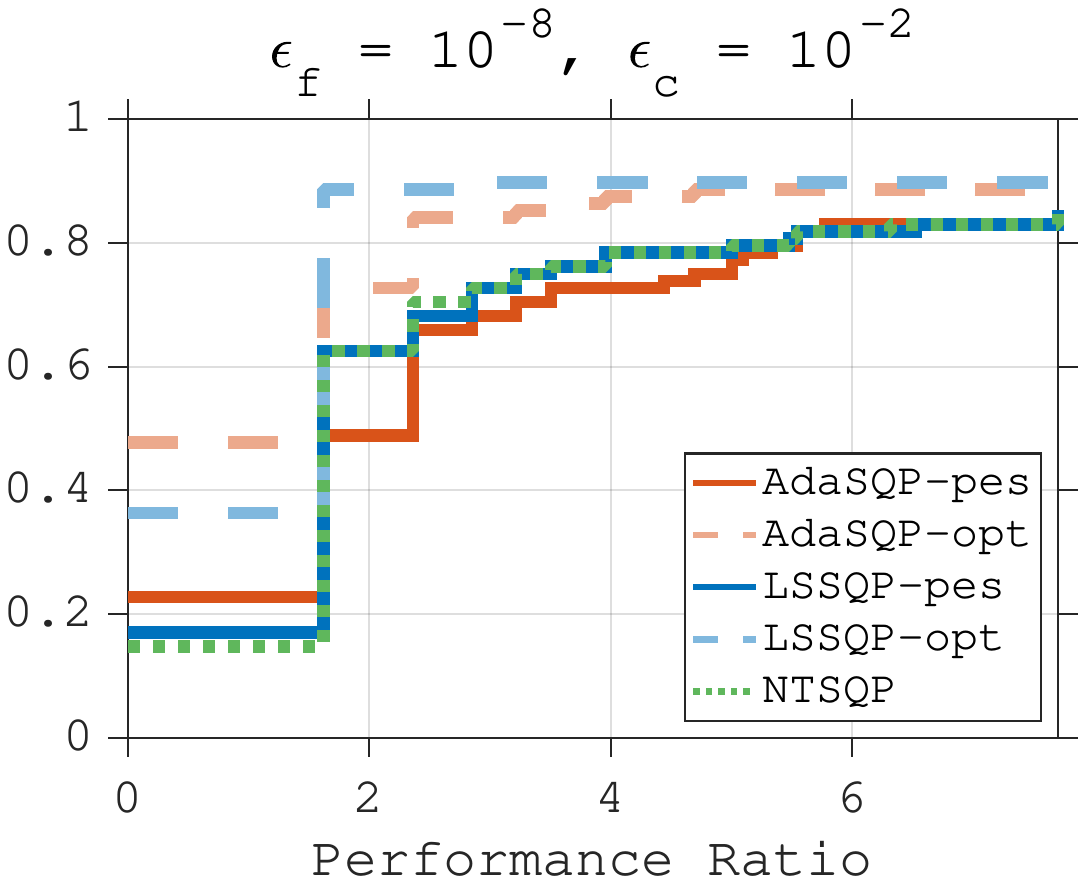}
\includegraphics[width=0.24\textwidth,clip=true,trim=10 180 50 150]{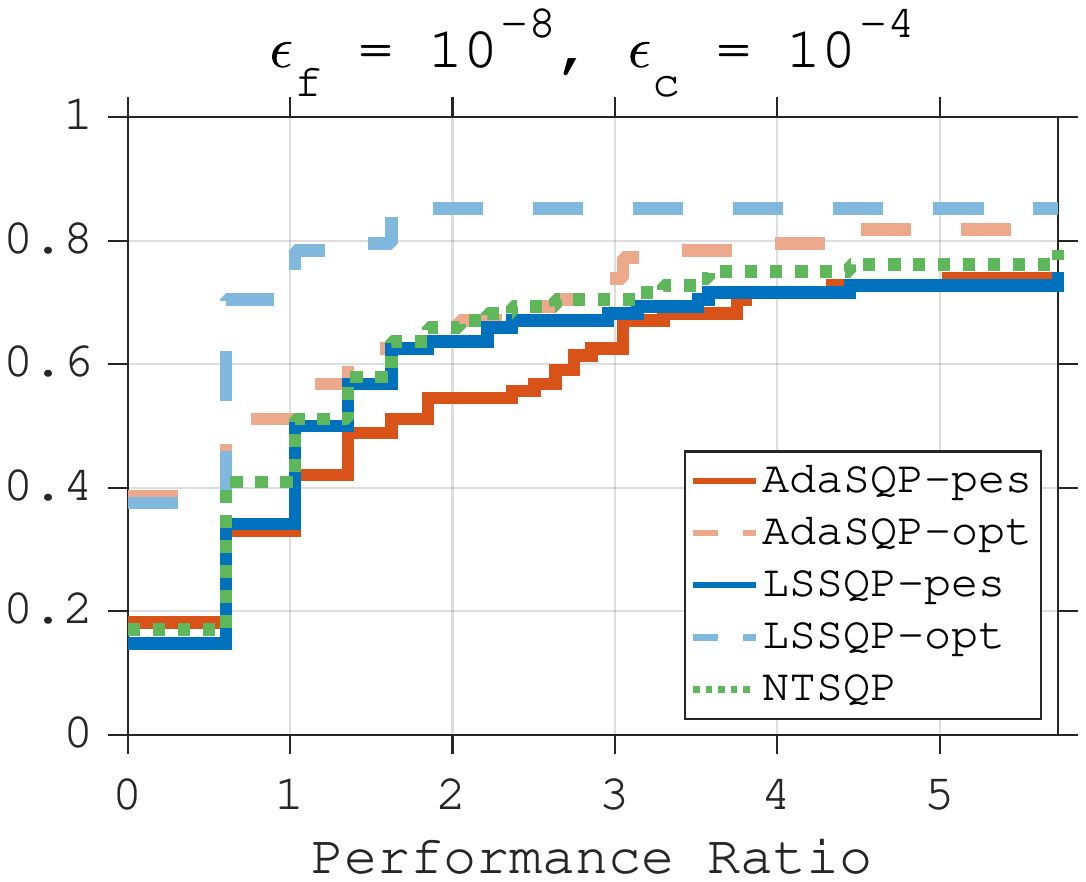}
\includegraphics[width=0.24\textwidth,clip=true,trim=10 180 50 150]{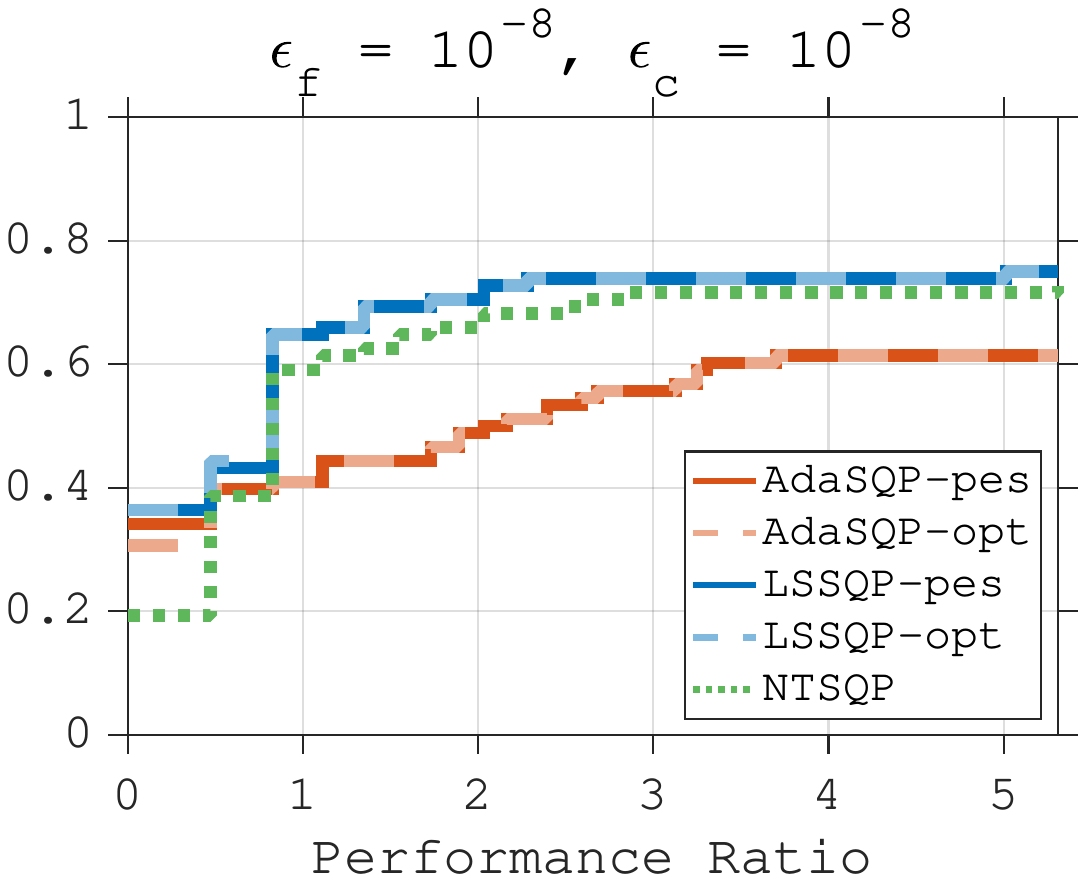}
\caption{Dolan-Mor\'e  performance profiles 
comparing \adasqppes{}, \adasqpopt{}, \lssqppes{}, \lssqpopt{} and \ntsqp{} on CUTEst collection of test problems that satisfy the LICQ 
in terms of number of \textbf{function evaluations} for $\epsilon_c \in \{ 10^{-1}, 10^{-2}, 10^{-4}, 10^{-8}\}$ (from \textbf{left} to \textbf{right}) and 
$\epsilon_f \in \{ 10^{-1}, 10^{-2}, 10^{-4}, 10^{-8}\}$ (from \textbf{top} to \textbf{bottom}). } 
\label{fig.DMplot.LICQ.fun}
\end{figure}

\begin{figure}[htbp]
        \centering    
\includegraphics[width=0.24\textwidth,clip=true,trim=10 180 50 150]{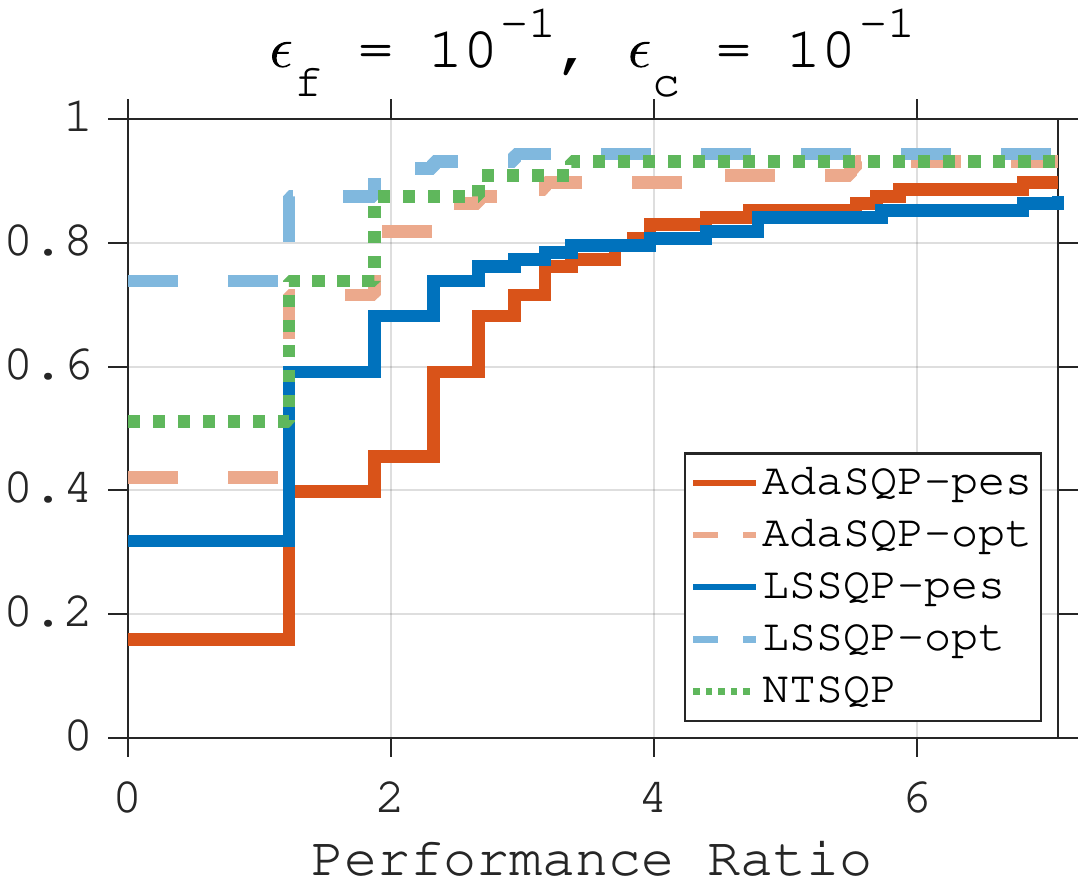}
\includegraphics[width=0.24\textwidth,clip=true,trim=10 180 50 150]{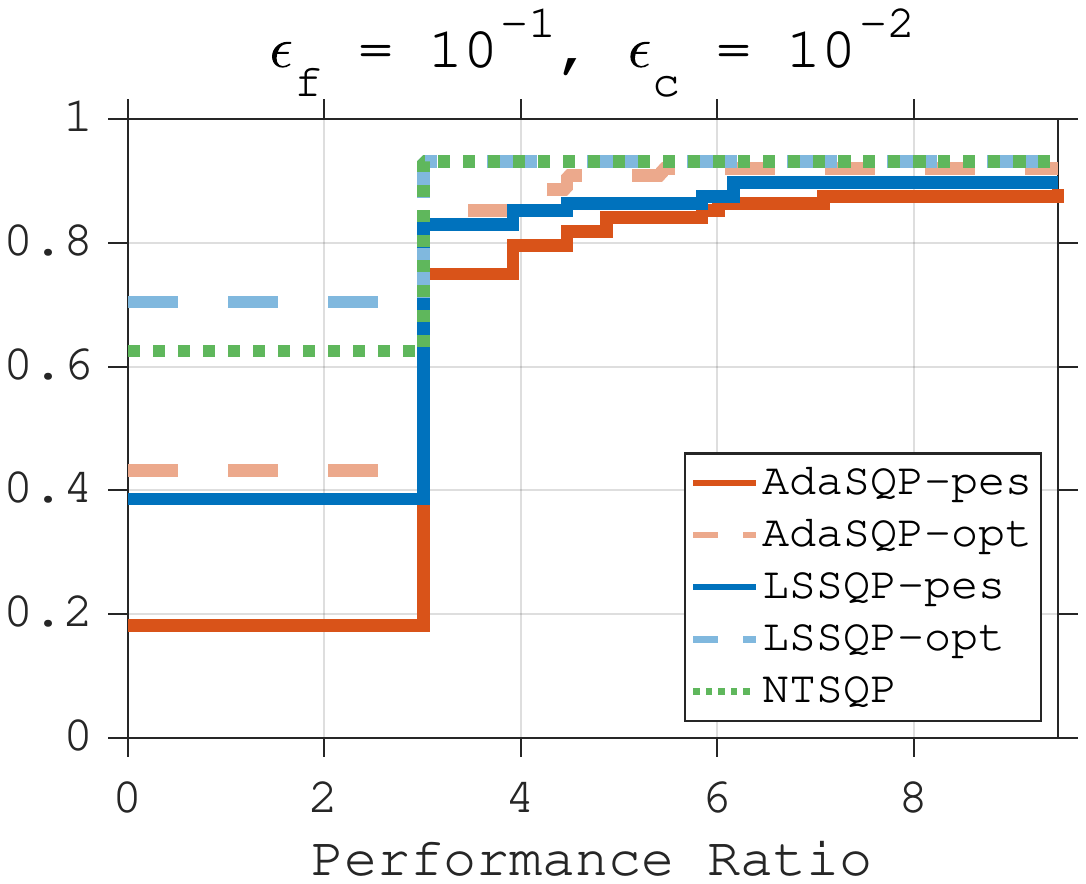}
\includegraphics[width=0.24\textwidth,clip=true,trim=10 180 50 150]{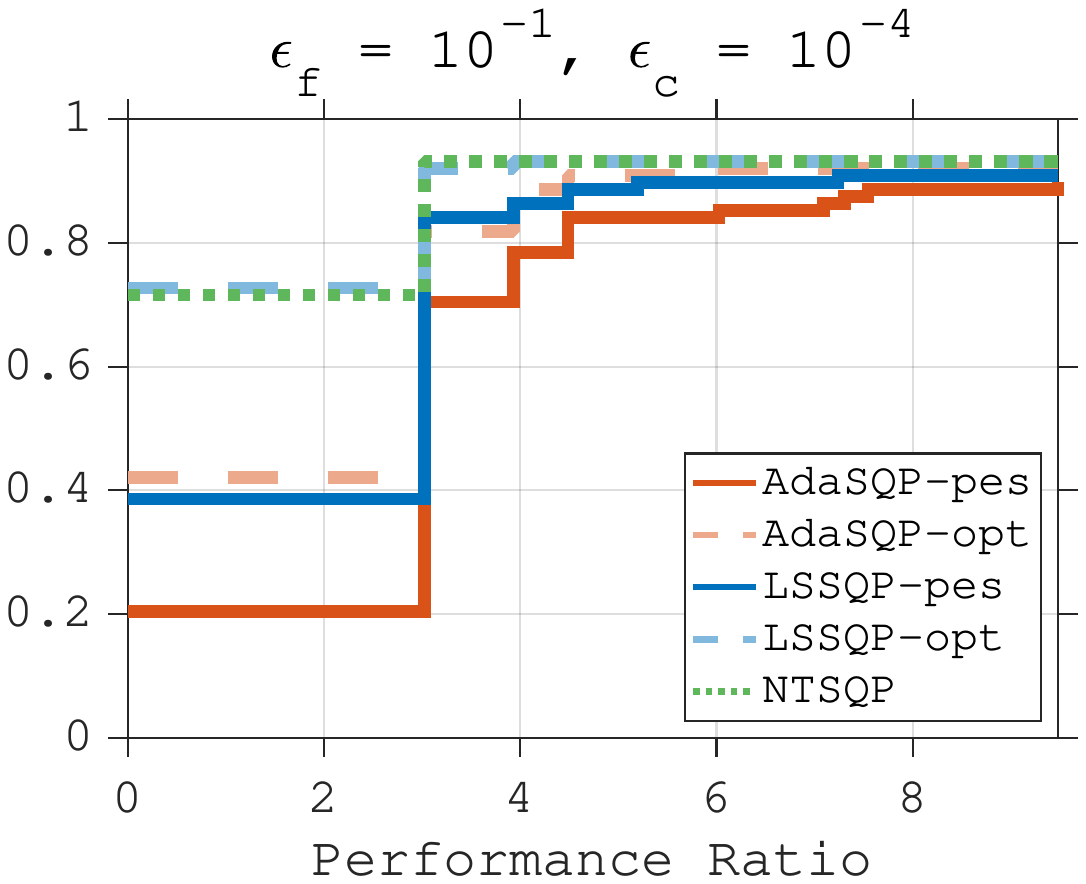}
\includegraphics[width=0.24\textwidth,clip=true,trim=10 180 50 150]{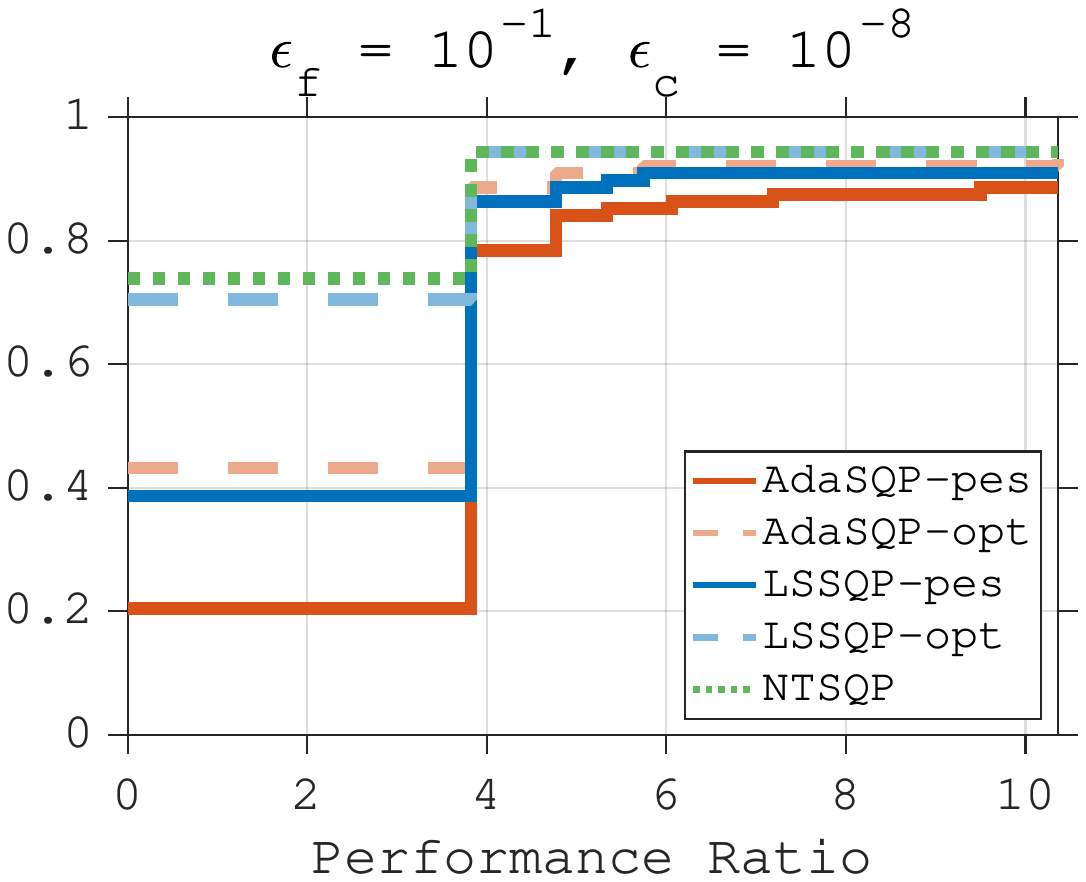}
\includegraphics[width=0.24\textwidth,clip=true,trim=10 180 50 150]{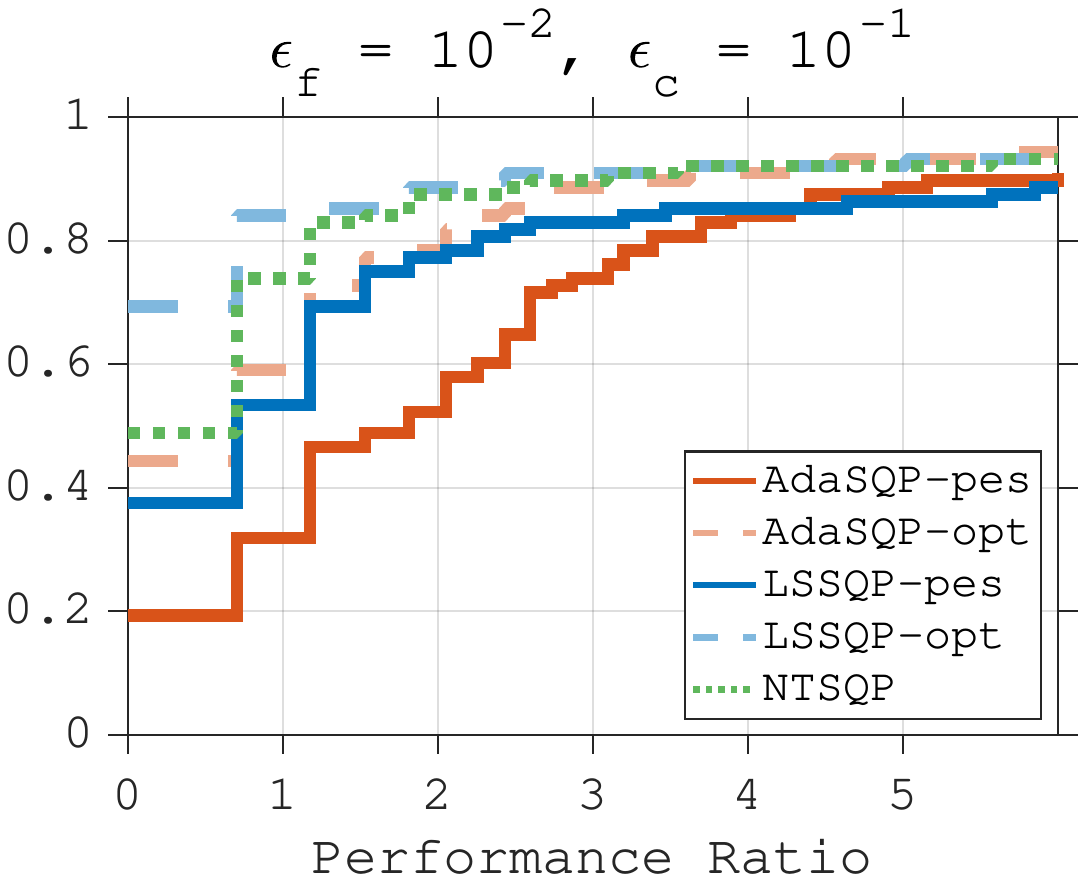}
\includegraphics[width=0.24\textwidth,clip=true,trim=10 180 50 150]{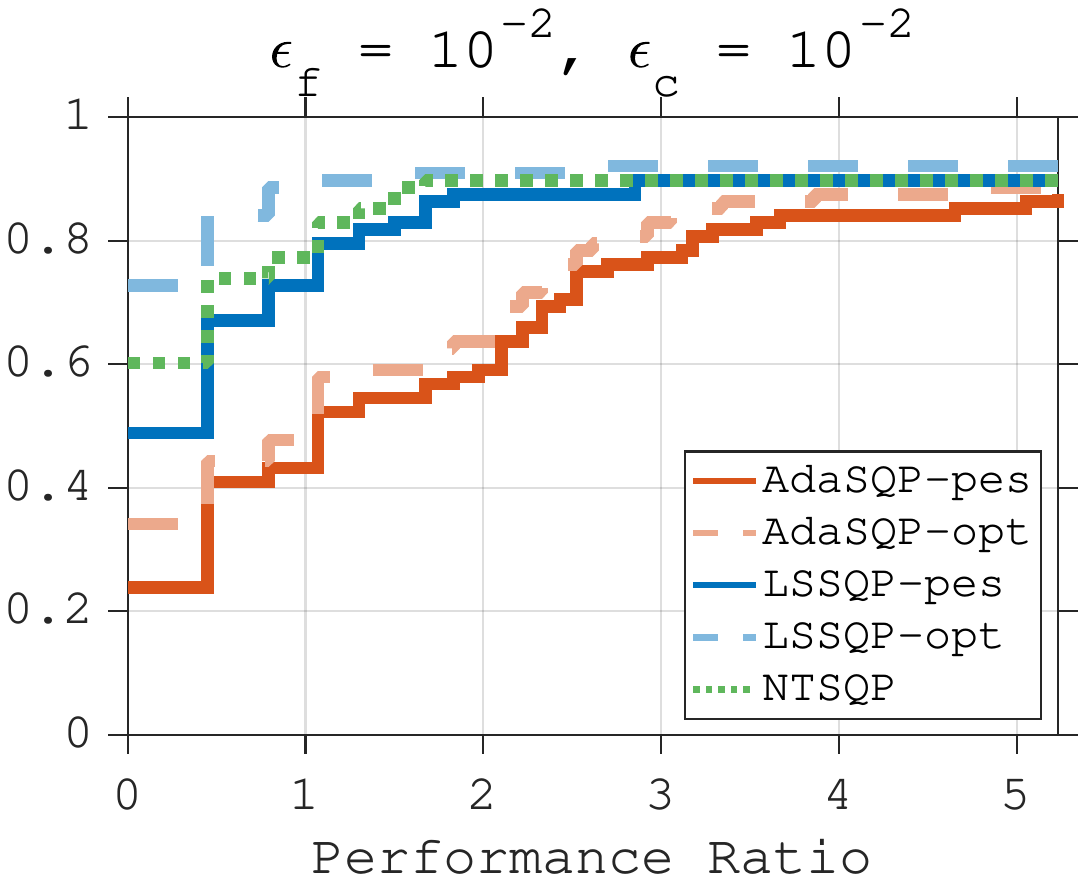}
\includegraphics[width=0.24\textwidth,clip=true,trim=10 180 50 150]{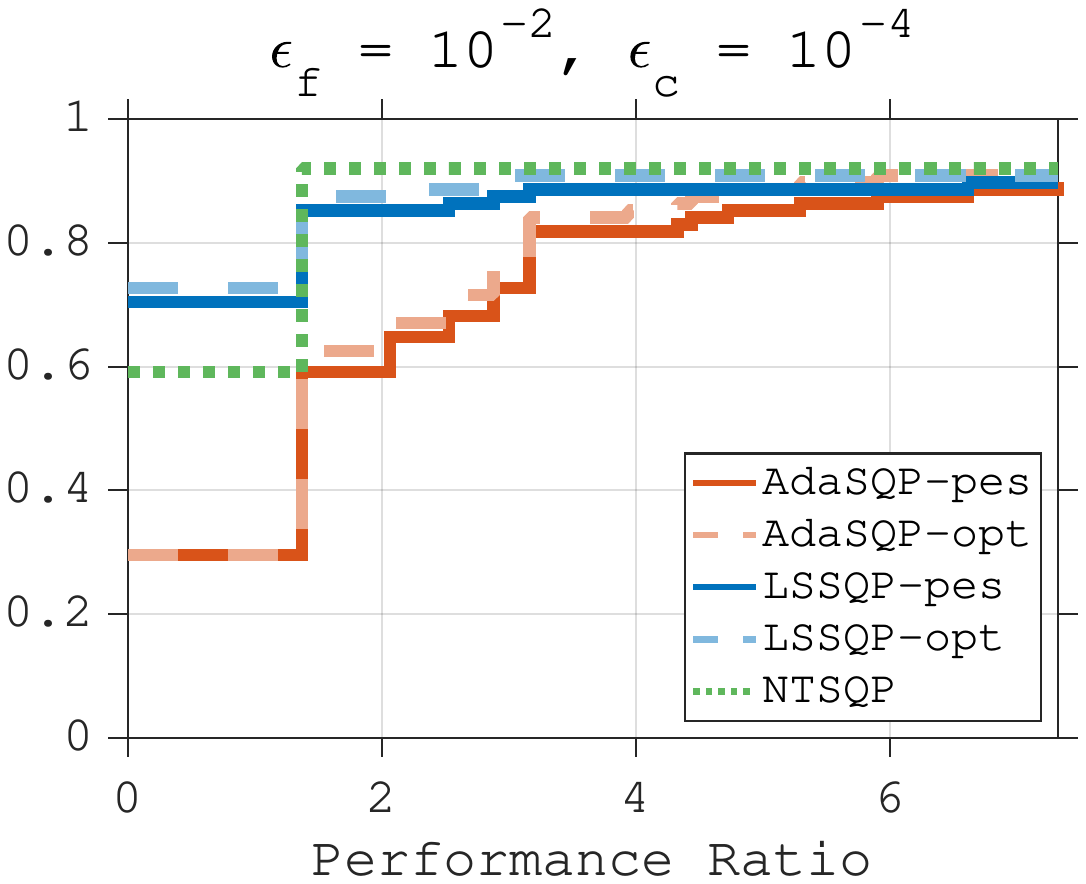}
\includegraphics[width=0.24\textwidth,clip=true,trim=10 180 50 150]{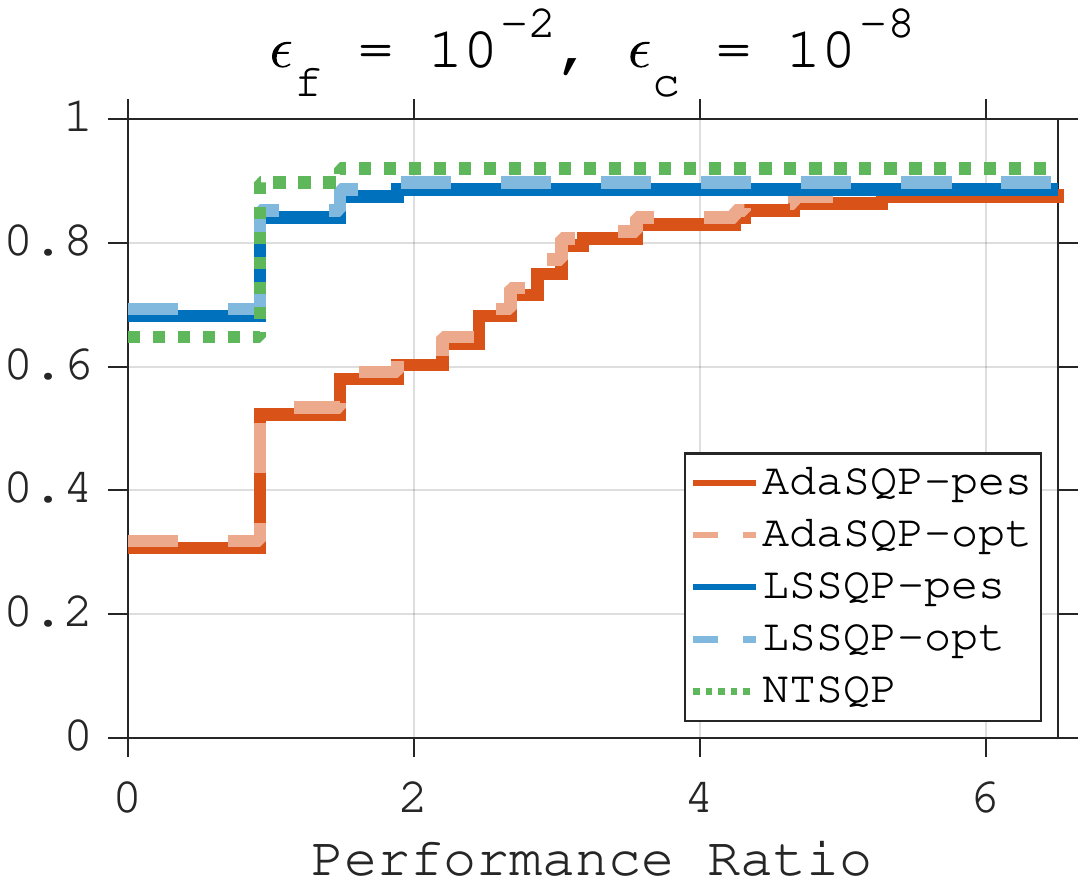}
\includegraphics[width=0.24\textwidth,clip=true,trim=10 180 50 150]{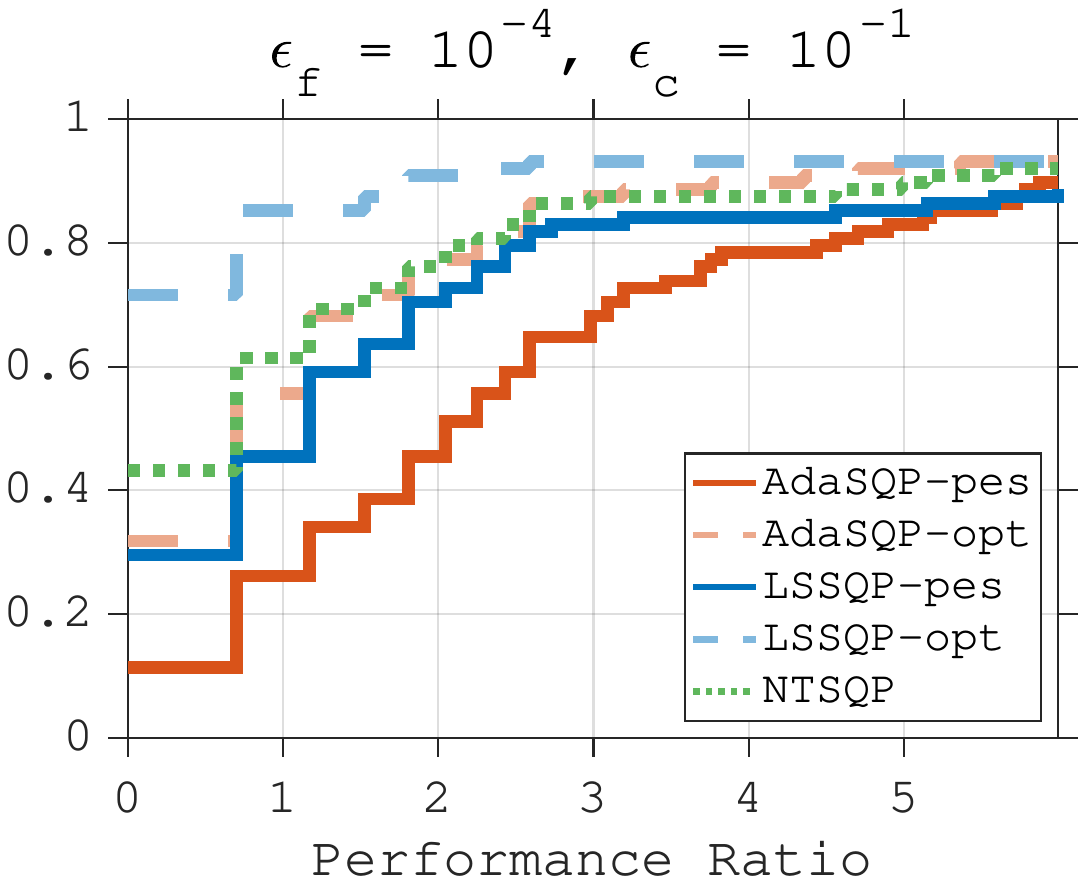}
\includegraphics[width=0.24\textwidth,clip=true,trim=10 180 50 150]{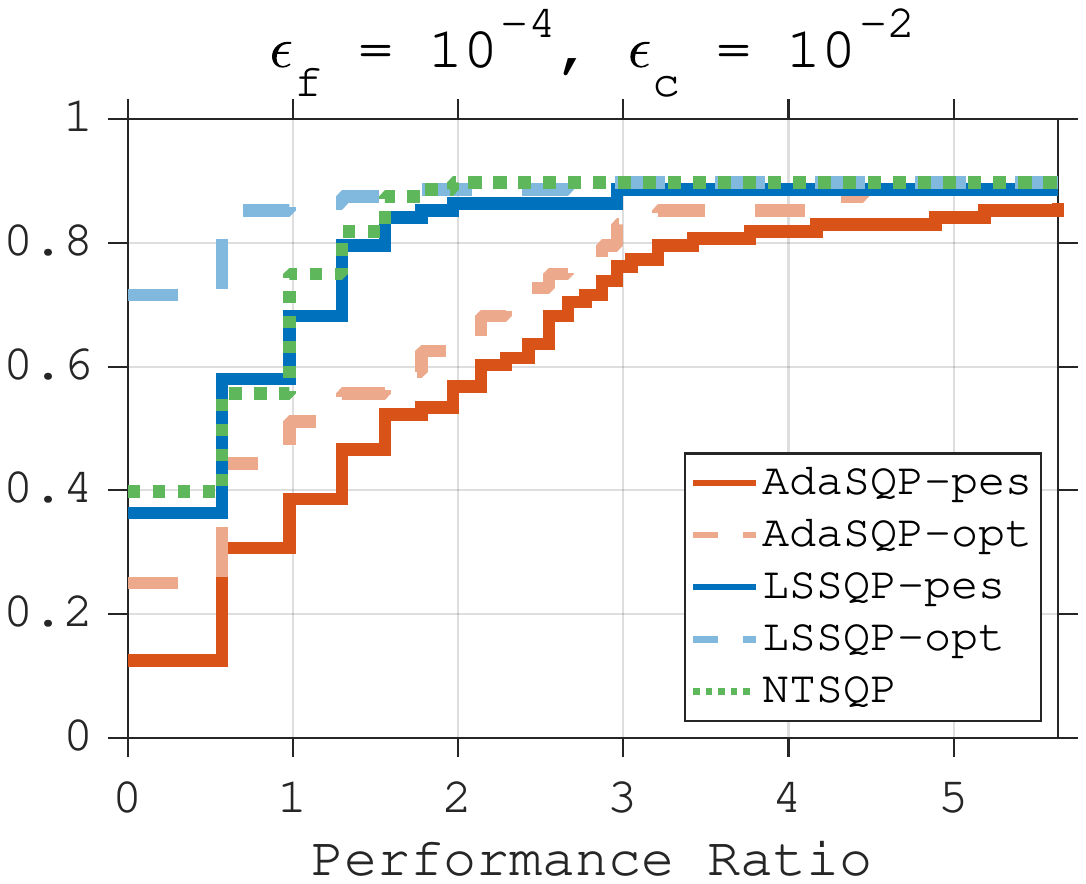}
\includegraphics[width=0.24\textwidth,clip=true,trim=10 180 50 150]{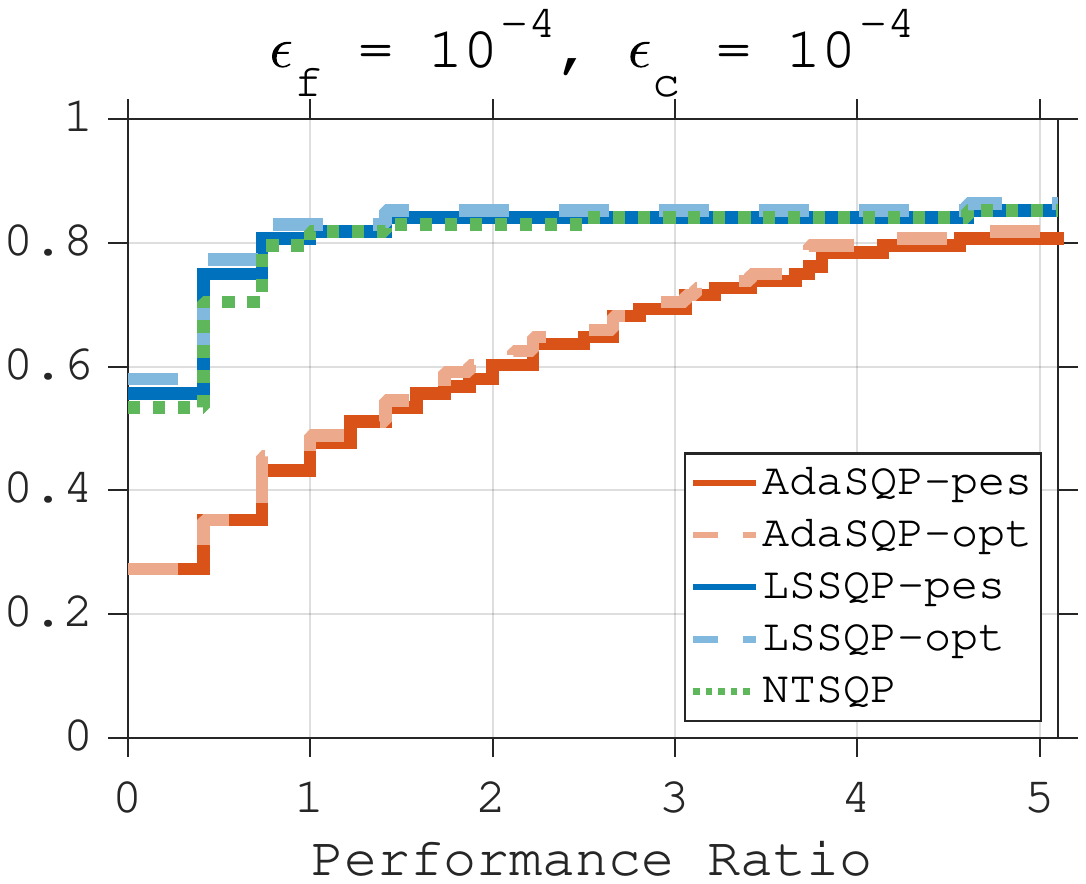}
\includegraphics[width=0.24\textwidth,clip=true,trim=10 180 50 150]{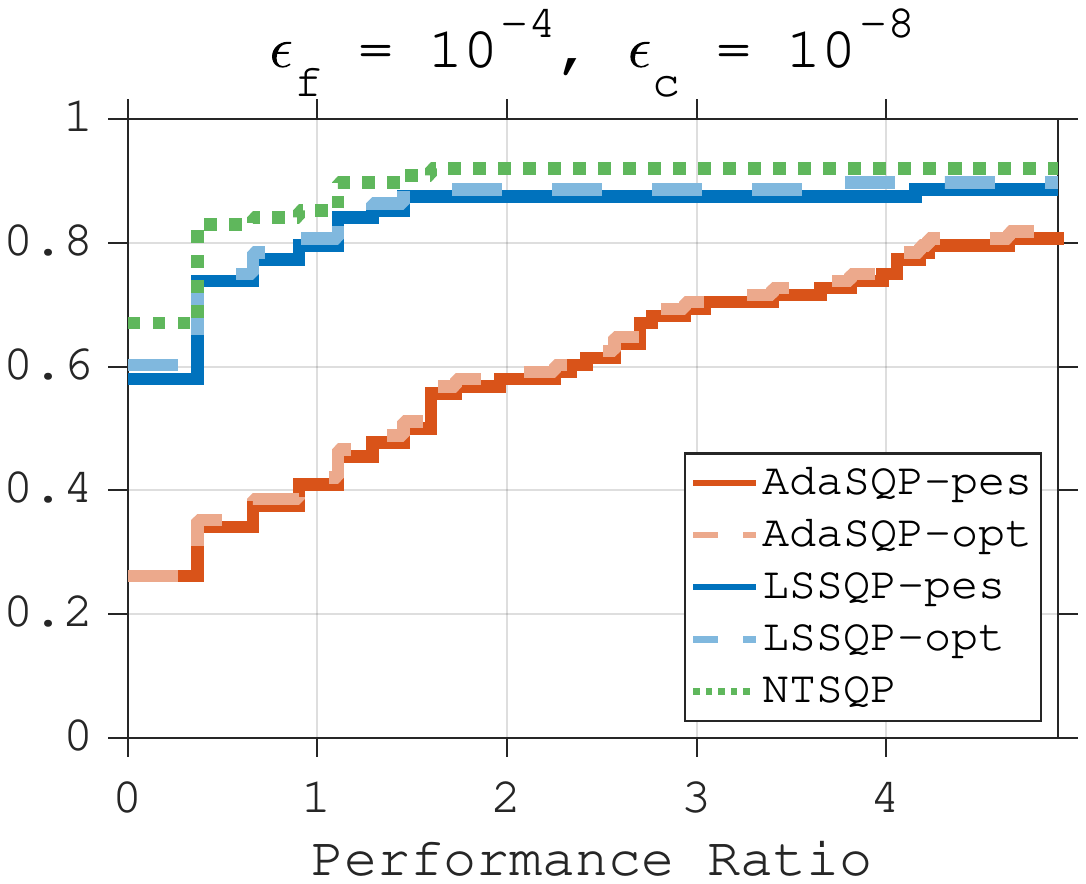}
\includegraphics[width=0.24\textwidth,clip=true,trim=10 180 50 150]{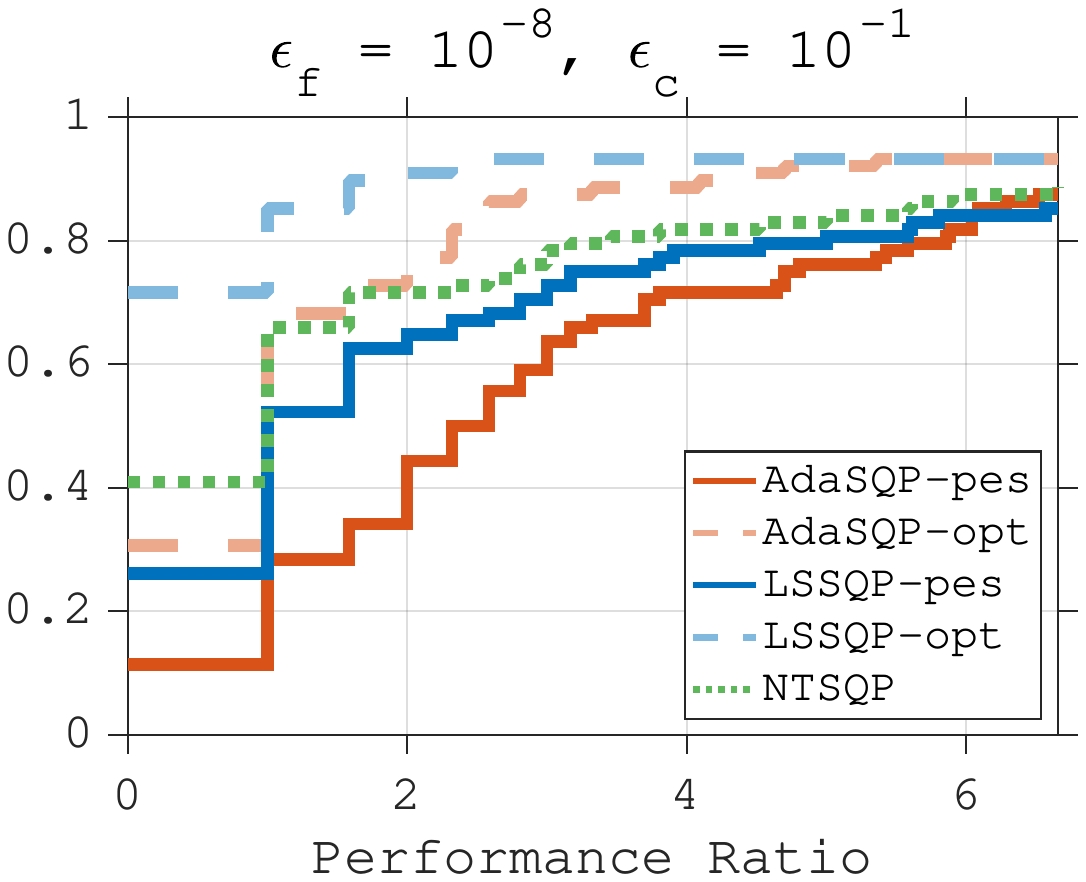}
\includegraphics[width=0.24\textwidth,clip=true,trim=10 180 50 150]{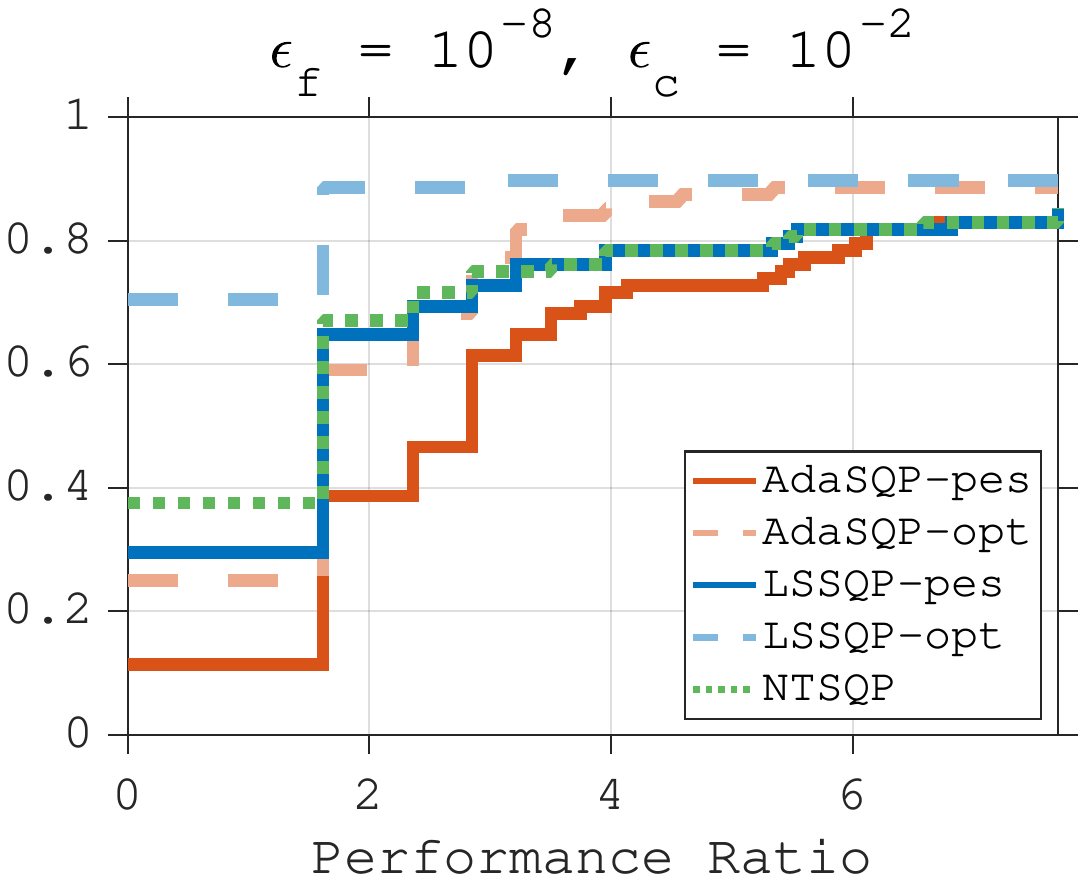}
\includegraphics[width=0.24\textwidth,clip=true,trim=10 180 50 150]{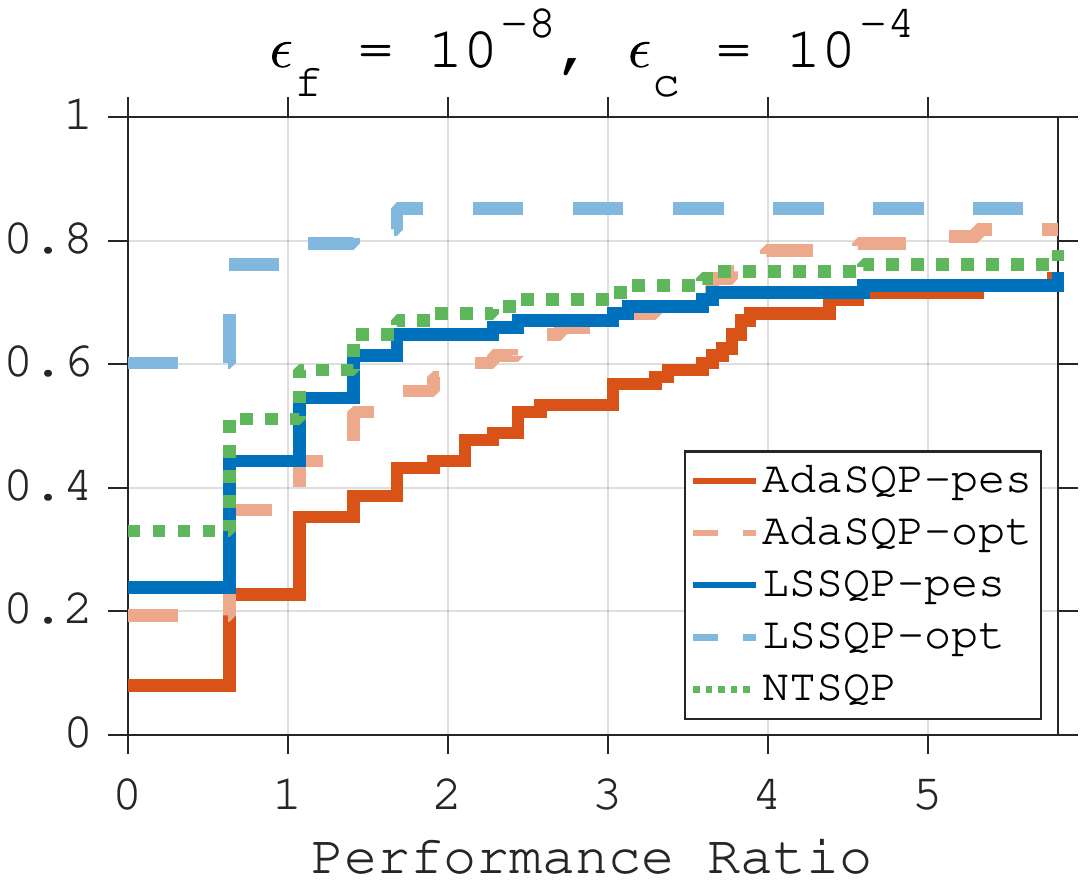}
\includegraphics[width=0.24\textwidth,clip=true,trim=10 180 50 150]{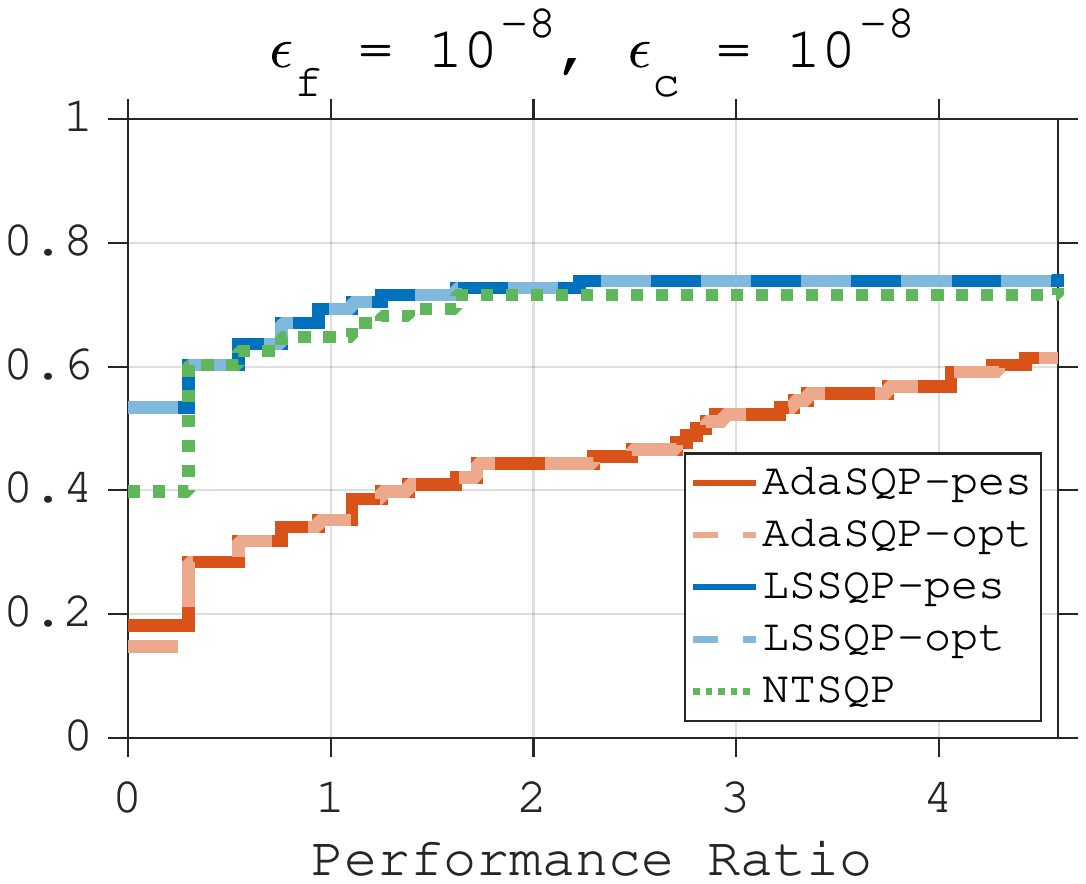} 
\caption{Dolan-Mor\'e  performance profiles 
comparing \adasqppes{}, \adasqpopt{}, \lssqppes{}, \lssqpopt{} and \ntsqp{} on CUTEst collection of test problems that satisfy the LICQ 
in terms of \textbf{MINRES iterations} for $\epsilon_c \in \{ 10^{-1}, 10^{-2}, 10^{-4}, 10^{-8}\}$ (from \textbf{left} to \textbf{right}) and 
$\epsilon_f \in \{ 10^{-1}, 10^{-2}, 10^{-4}, 10^{-8}\}$ (from \textbf{top} to \textbf{bottom}). }
\label{fig.DMplot.LICQ.minres}
    \end{figure}

\begin{figure}[htbp]
        \centering    
\includegraphics[width=0.24\textwidth,clip=true,trim=10 180 50 150]{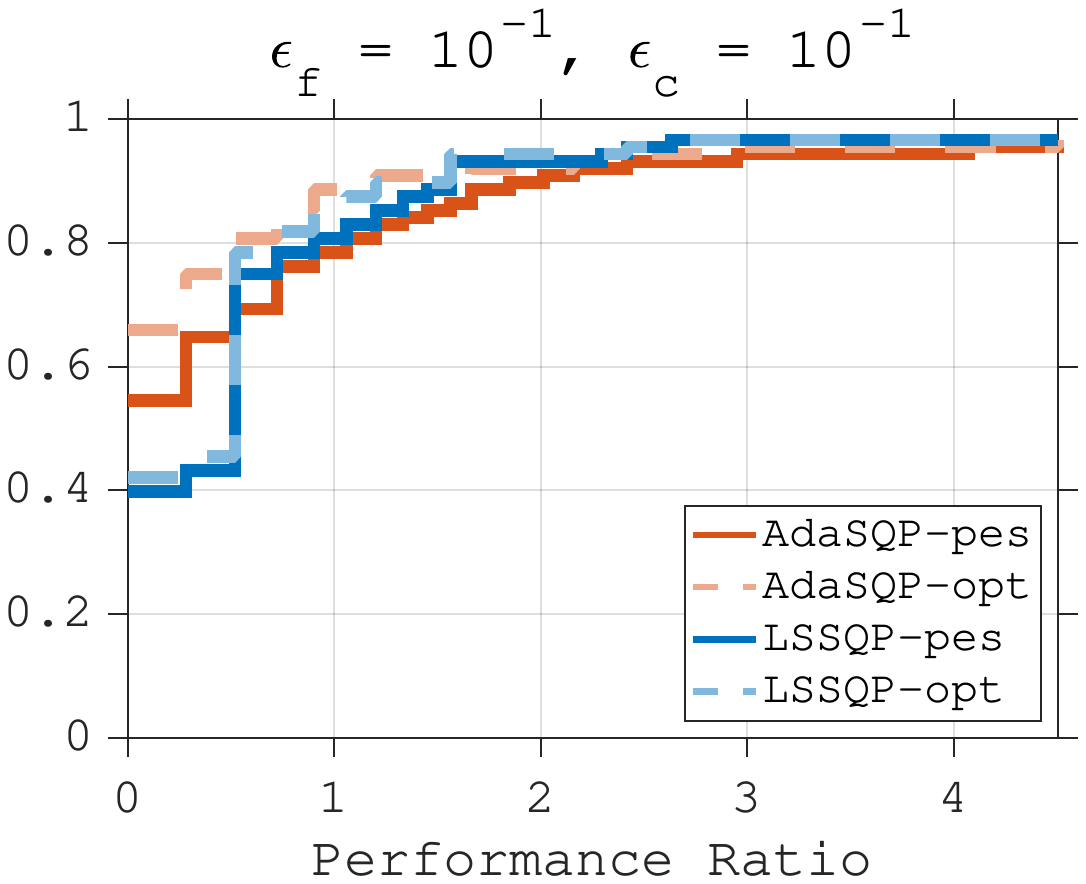}
\includegraphics[width=0.24\textwidth,clip=true,trim=10 180 50 150]{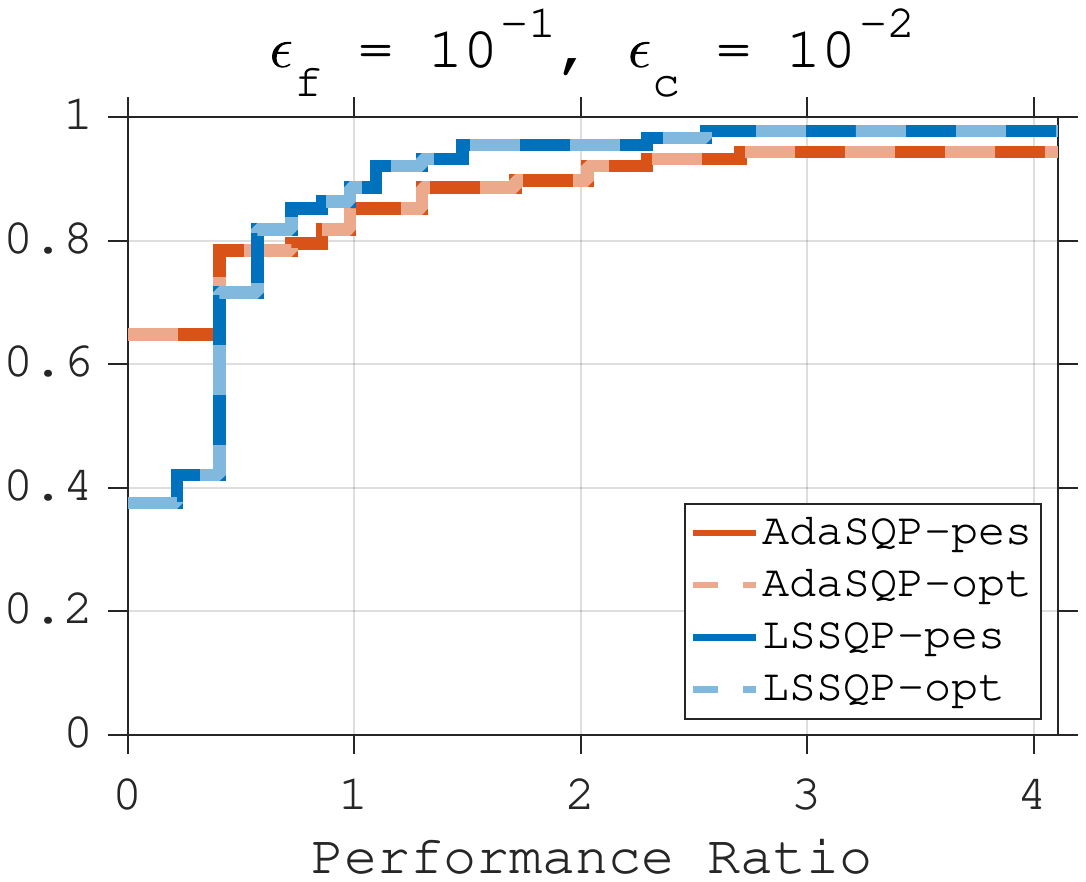}
\includegraphics[width=0.24\textwidth,clip=true,trim=10 180 50 150]{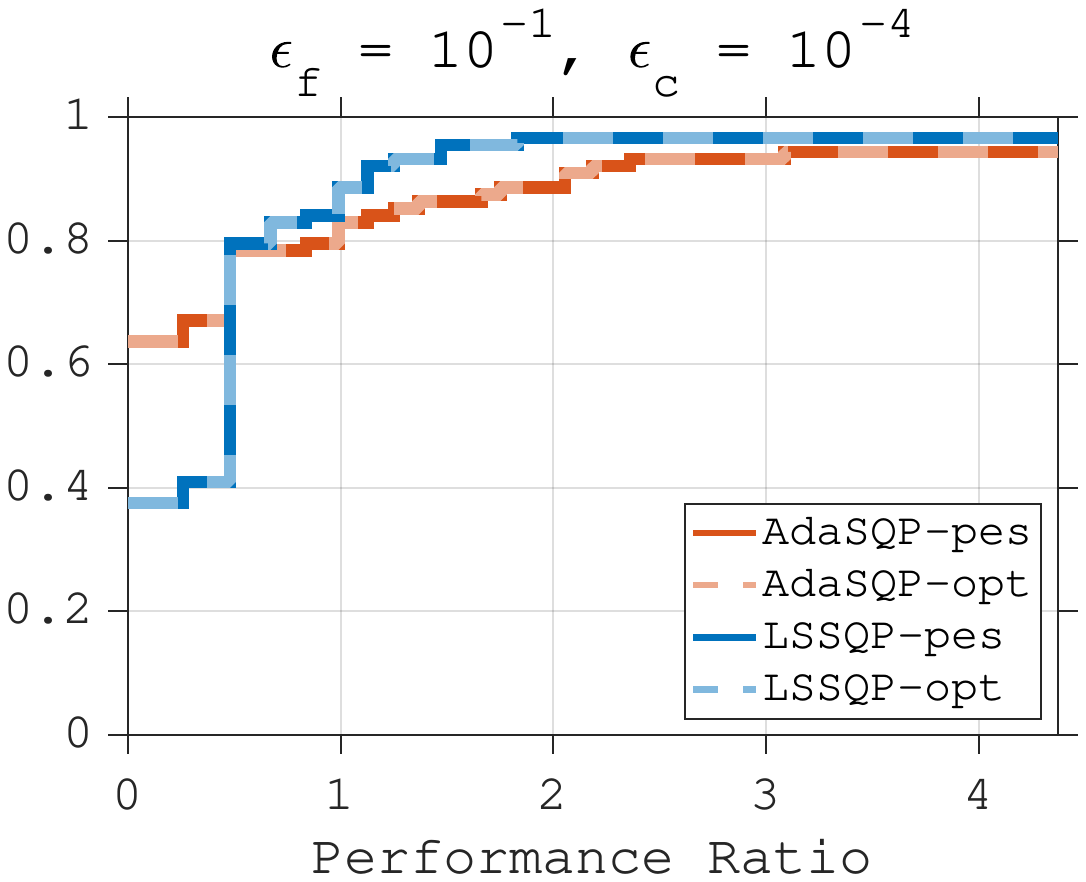}
\includegraphics[width=0.24\textwidth,clip=true,trim=10 180 50 150]{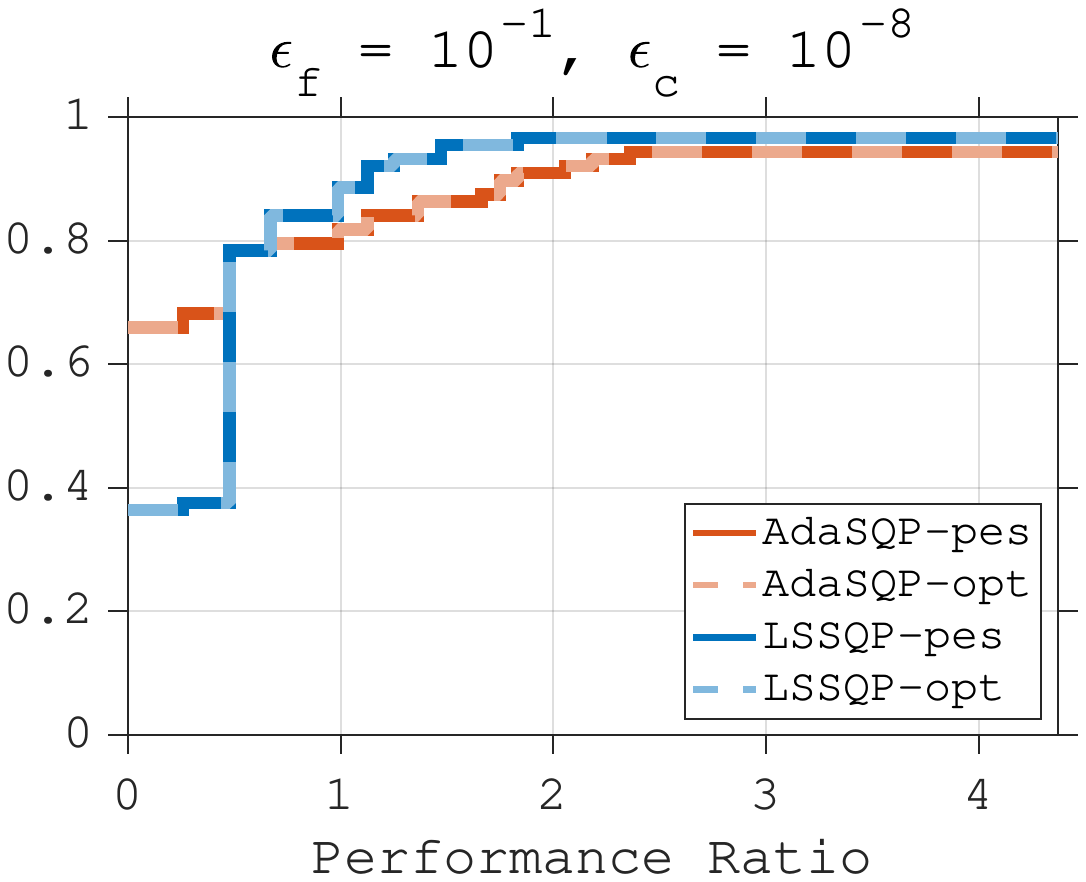}
\includegraphics[width=0.24\textwidth,clip=true,trim=10 180 50 150]{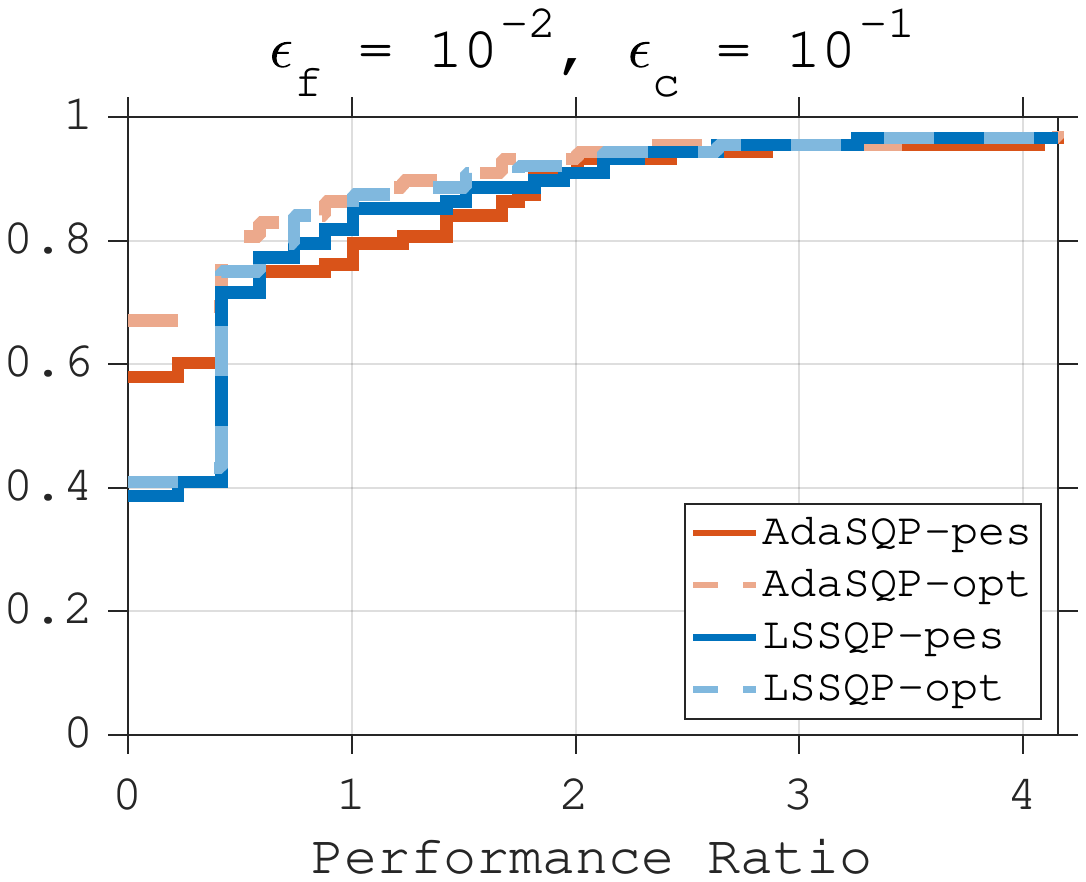}
\includegraphics[width=0.24\textwidth,clip=true,trim=10 180 50 150]{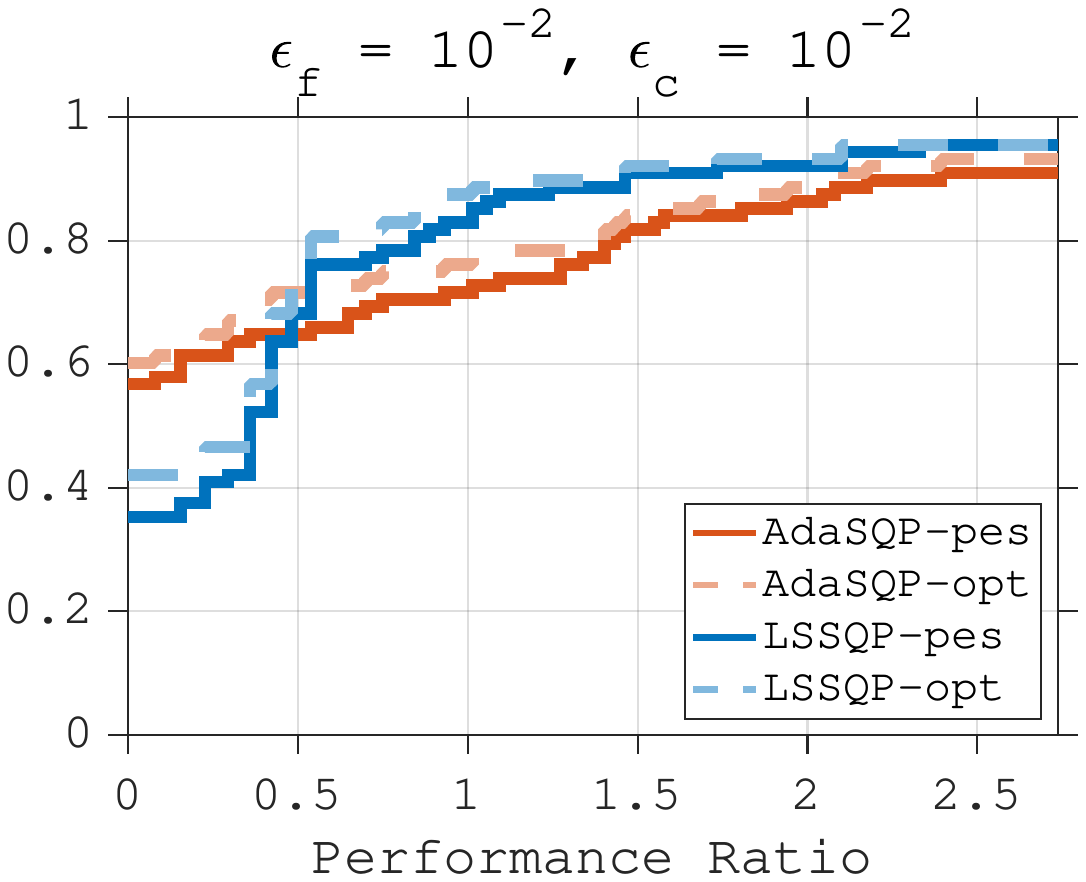}
\includegraphics[width=0.24\textwidth,clip=true,trim=10 180 50 150]{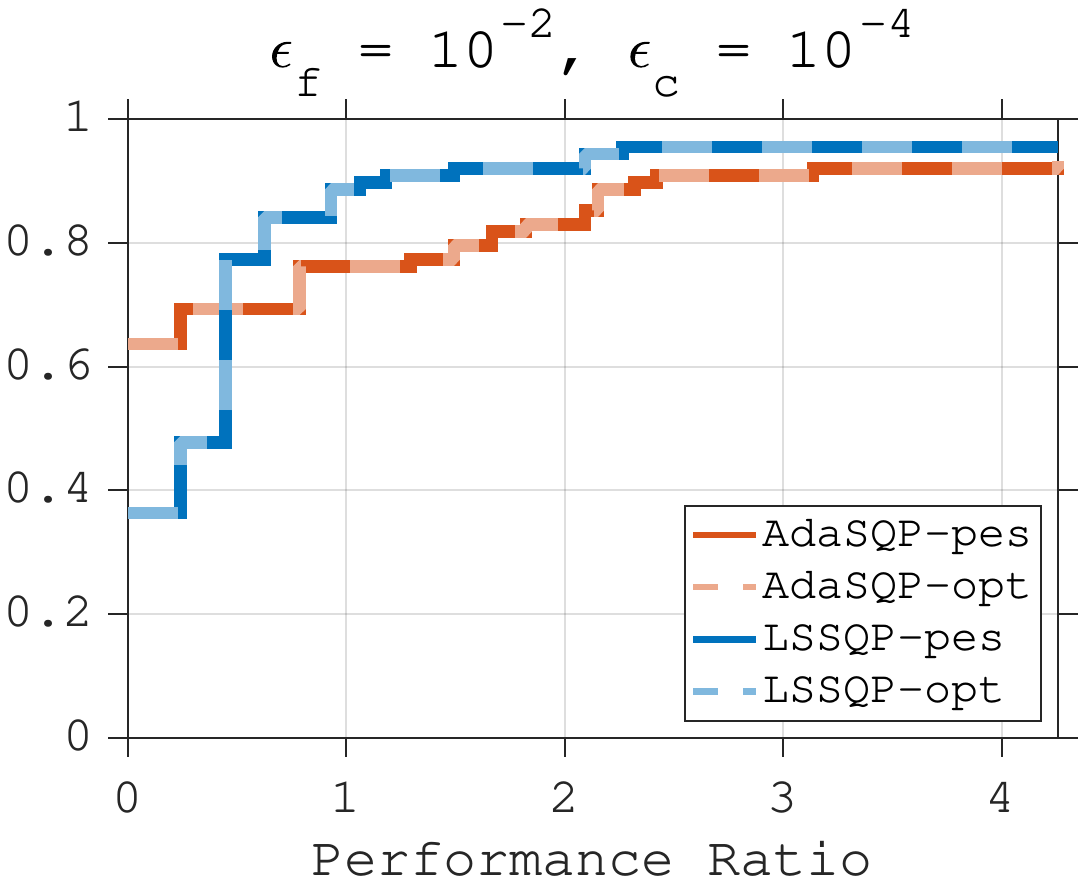}
\includegraphics[width=0.24\textwidth,clip=true,trim=10 180 50 150]{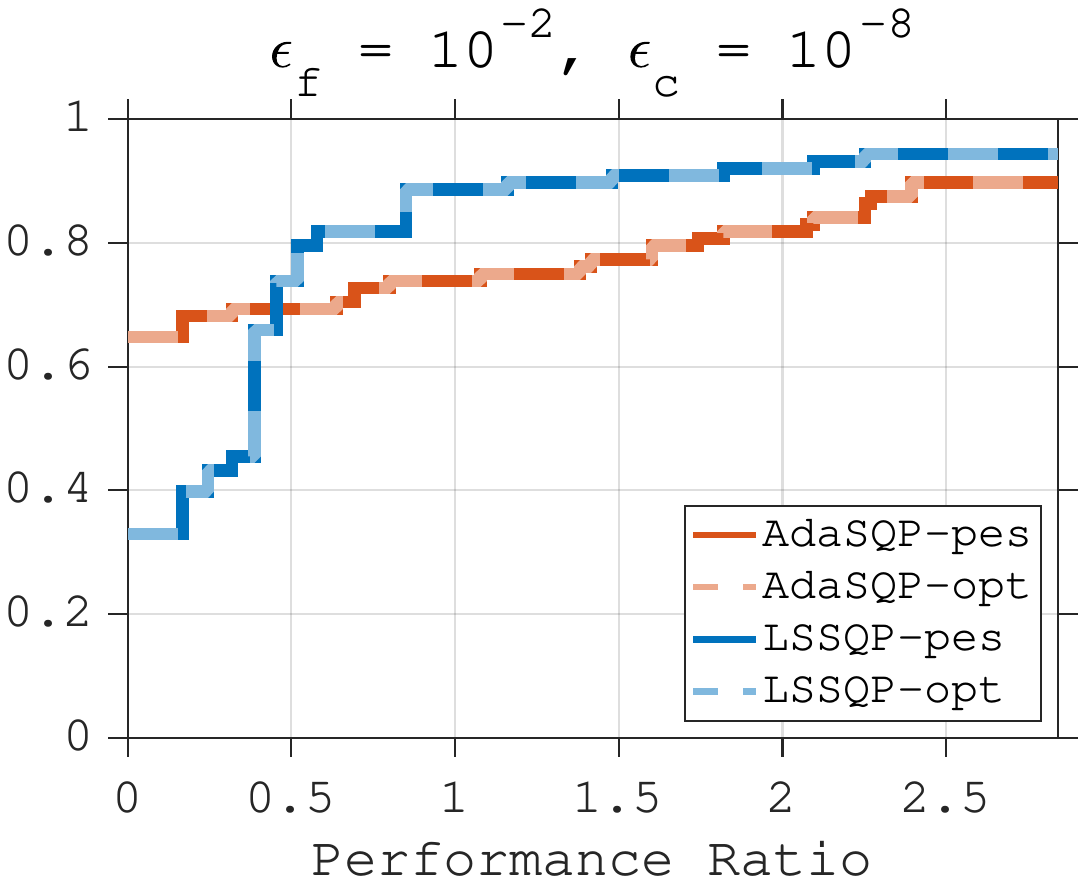}
\includegraphics[width=0.24\textwidth,clip=true,trim=10 180 50 150]{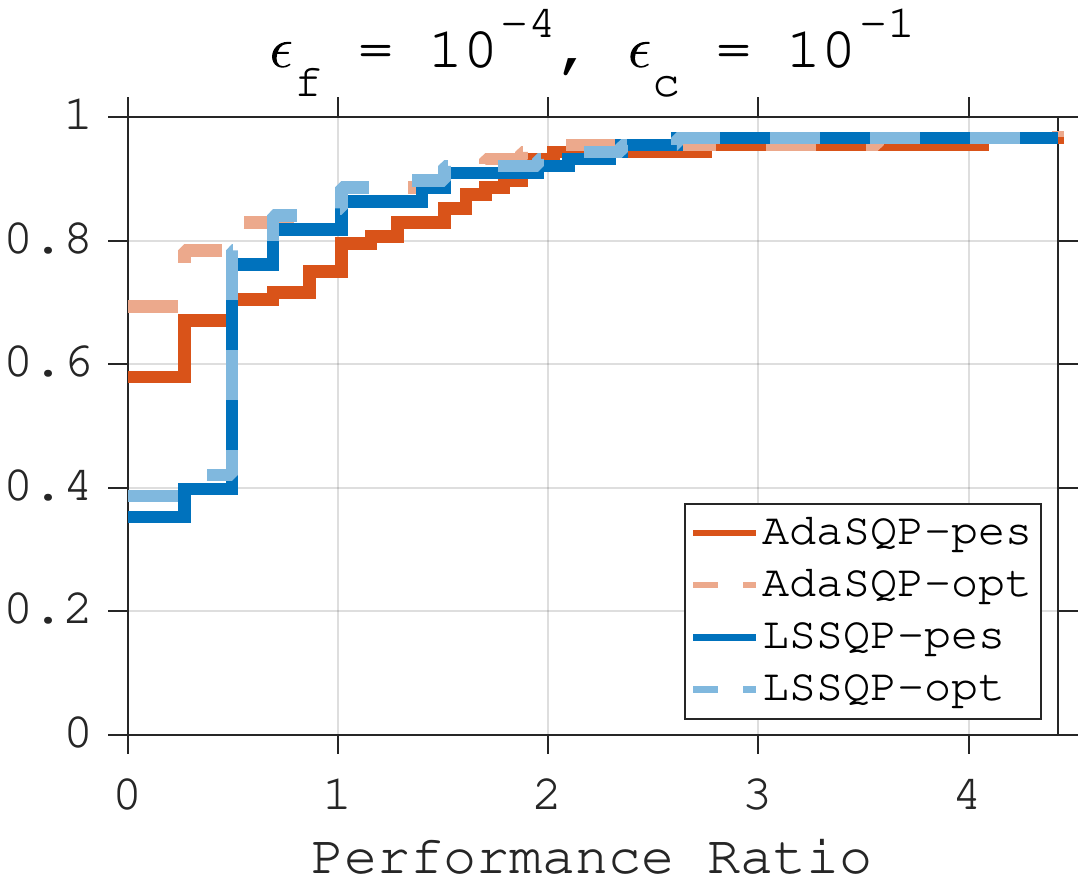}
\includegraphics[width=0.24\textwidth,clip=true,trim=10 180 50 150]{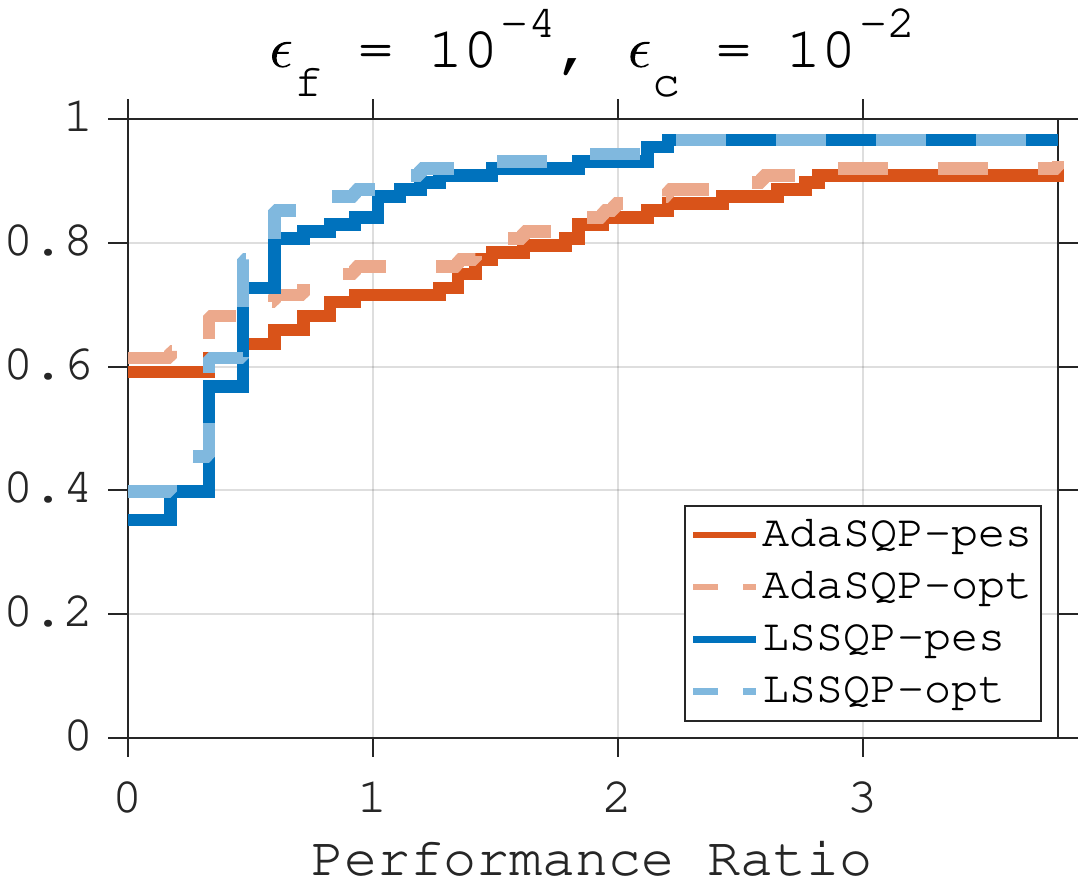}
\includegraphics[width=0.24\textwidth,clip=true,trim=10 180 50 150]{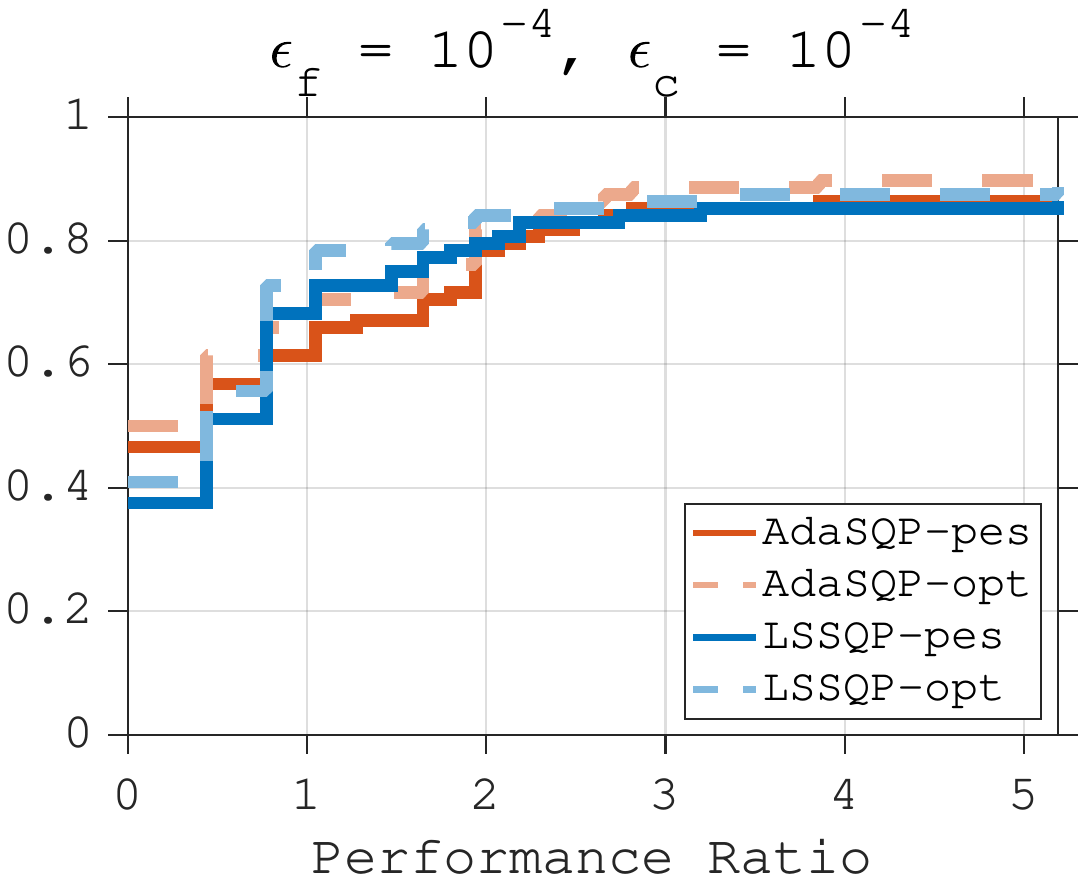}
\includegraphics[width=0.24\textwidth,clip=true,trim=10 180 50 150]{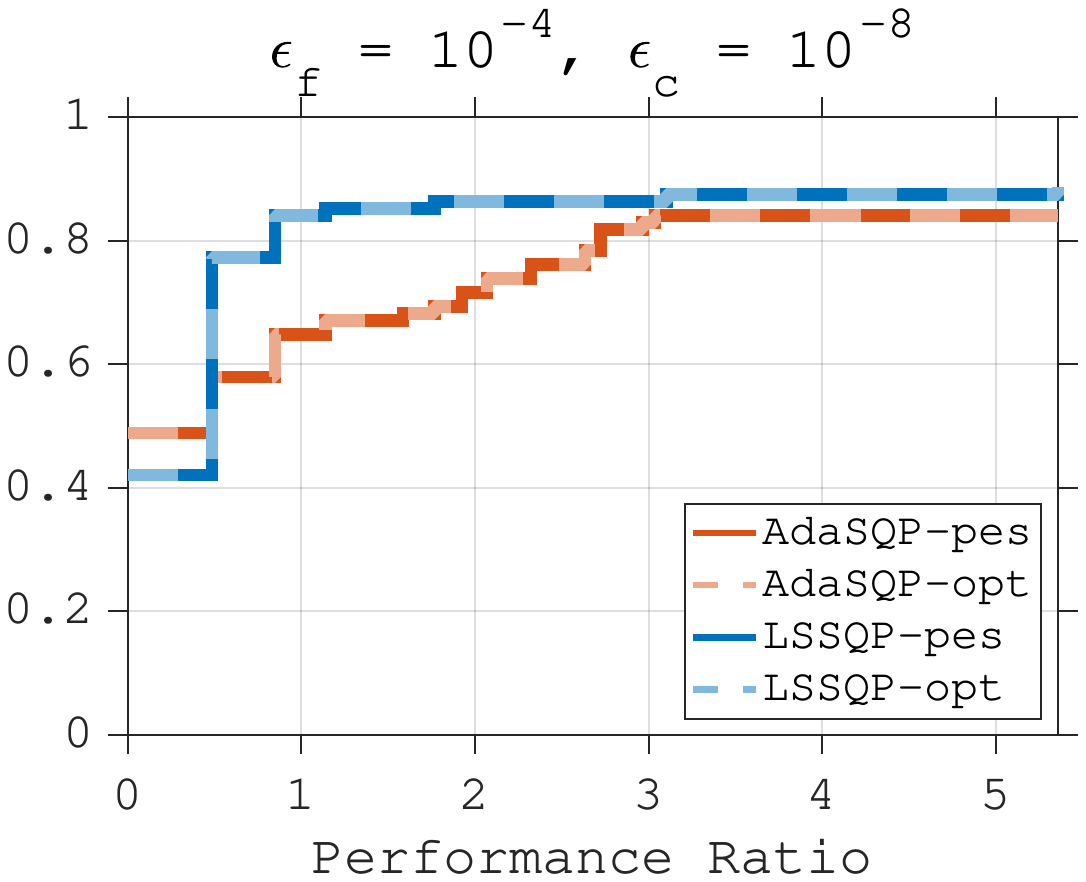}
\includegraphics[width=0.24\textwidth,clip=true,trim=10 180 50 150]{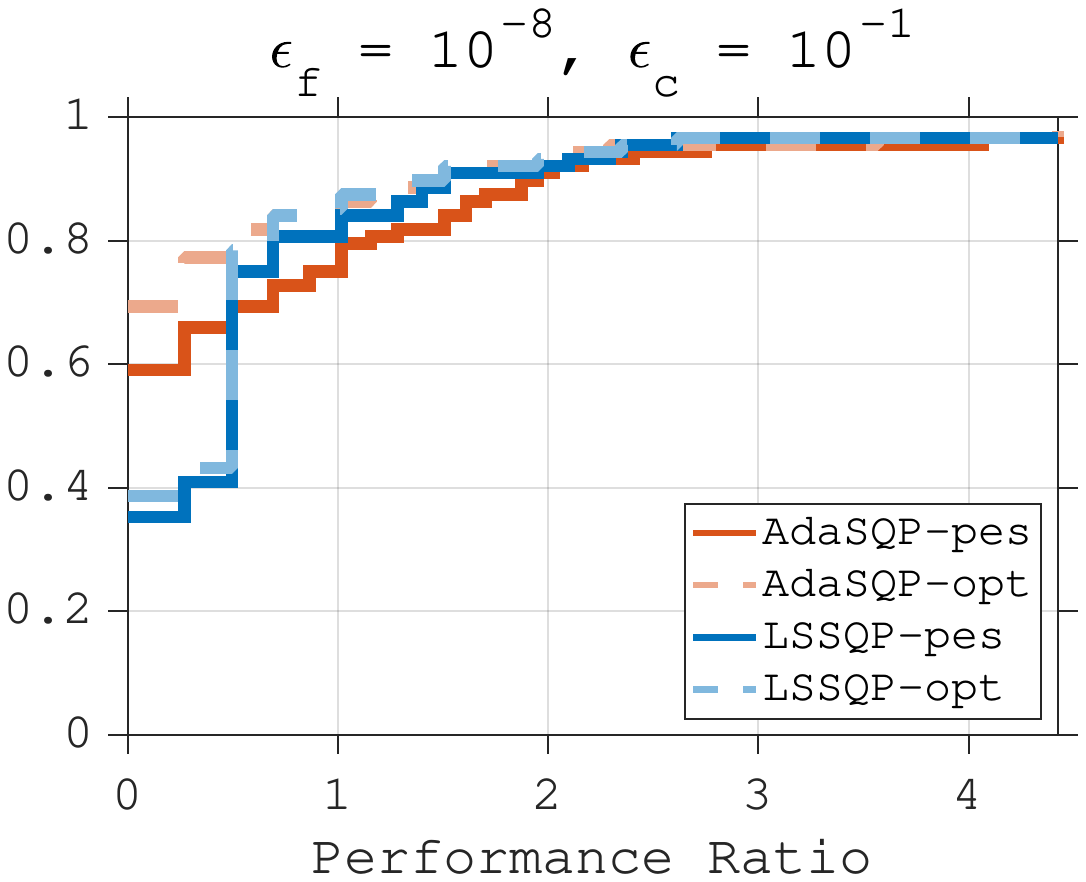}
\includegraphics[width=0.24\textwidth,clip=true,trim=10 180 50 150]{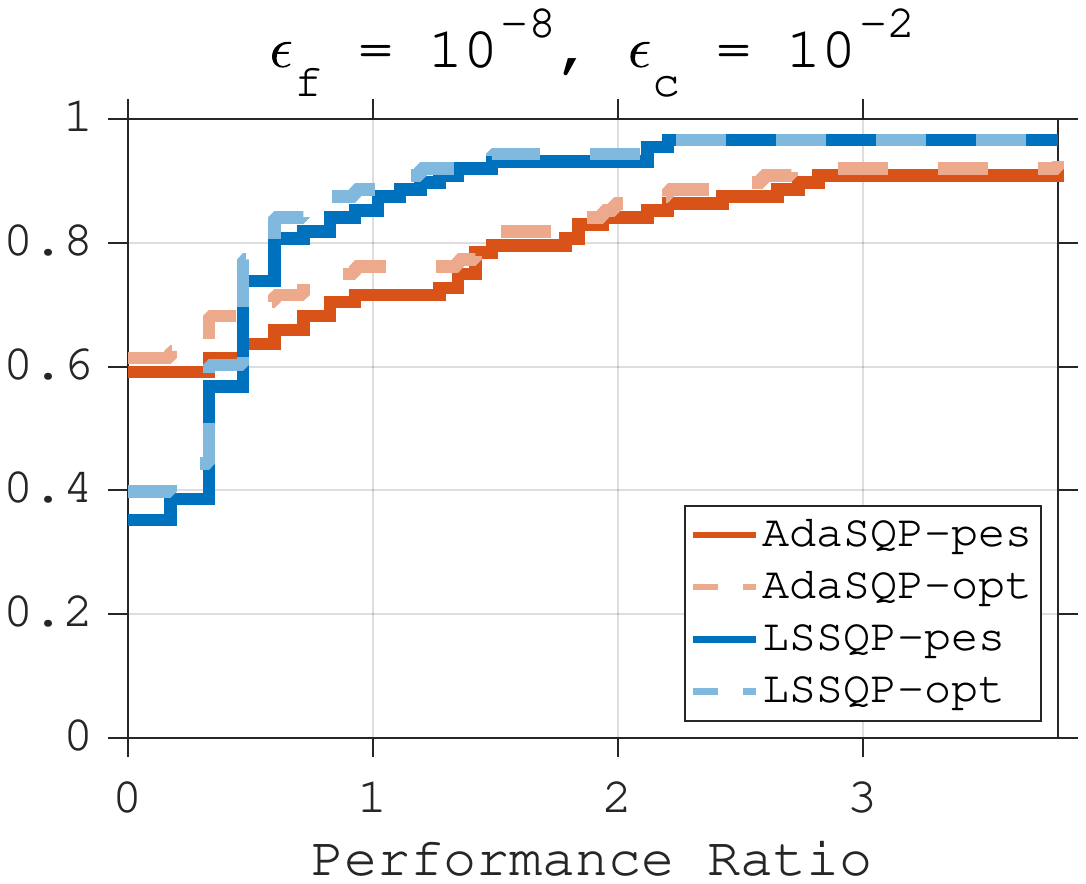}
\includegraphics[width=0.24\textwidth,clip=true,trim=10 180 50 150]{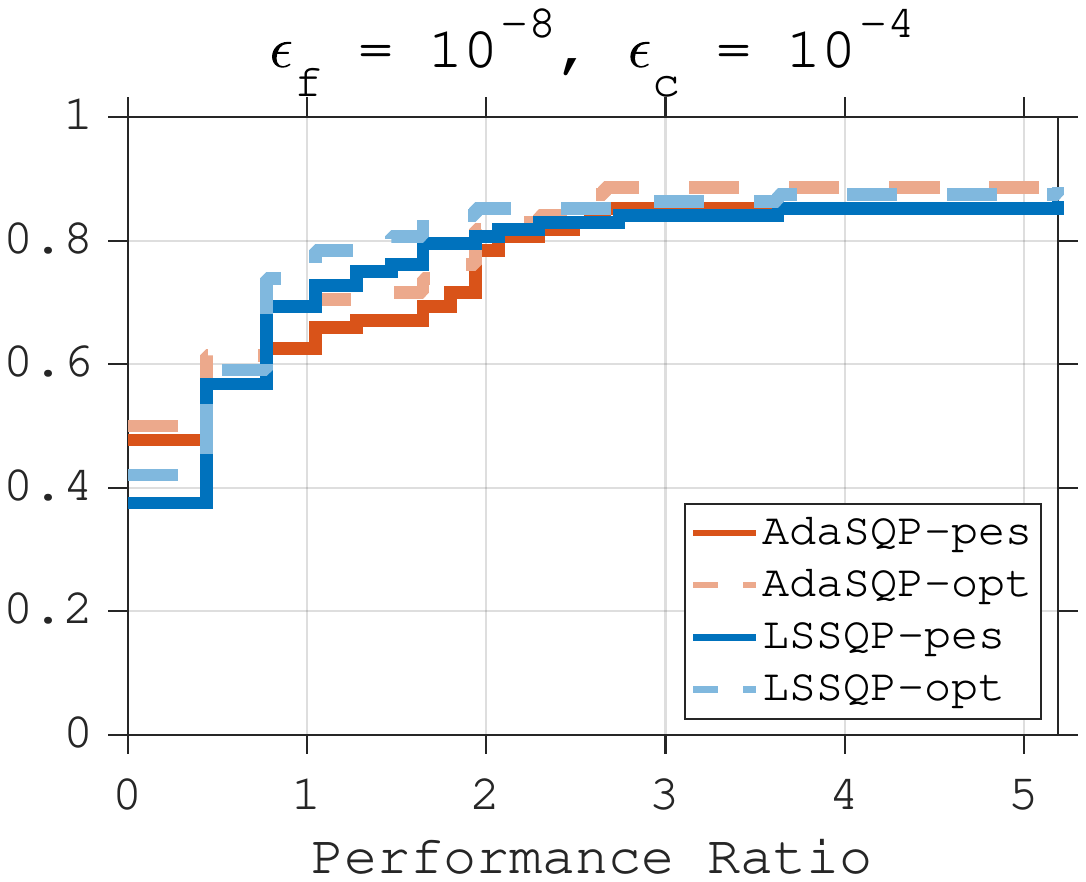}
\includegraphics[width=0.24\textwidth,clip=true,trim=10 180 50 150]{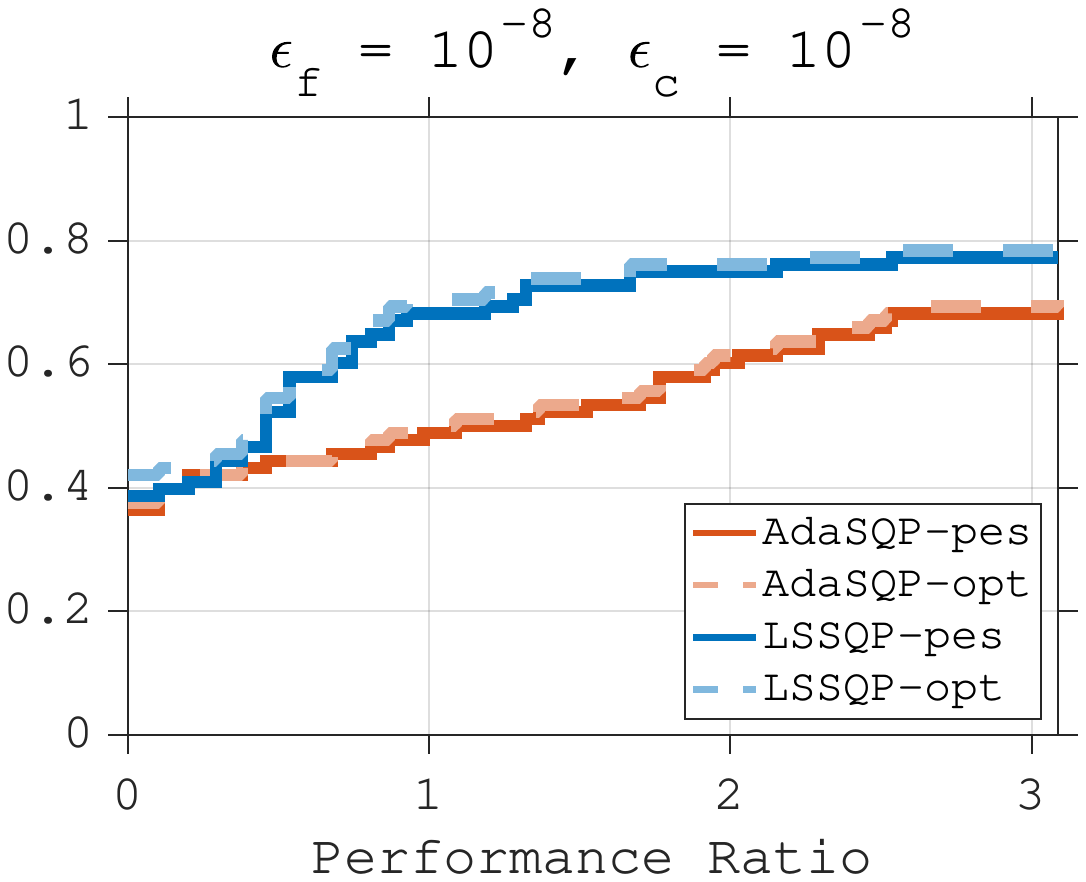}
\caption{Dolan-Mor\'e  performance profiles 
comparing \adasqppes{}, \adasqpopt{}, \lssqppes{}, and \lssqpopt{}  
on CUTEst collection of test problems in the absence of the LICQ in terms of \textbf{function evaluations} for $\epsilon_c \in \{ 10^{-1}, 10^{-2}, 10^{-4}, 10^{-8}\}$ (from \textbf{left} to \textbf{right}) and 
$\epsilon_f \in \{ 10^{-1}, 10^{-2}, 10^{-4}, 10^{-8}\}$ (from \textbf{top} to \textbf{bottom}). }
\label{fig.DMplot.NoLICQ.fun}
\end{figure}

\begin{figure}[htbp]
        \centering    
\includegraphics[width=0.24\textwidth,clip=true,trim=10 180 50 150]{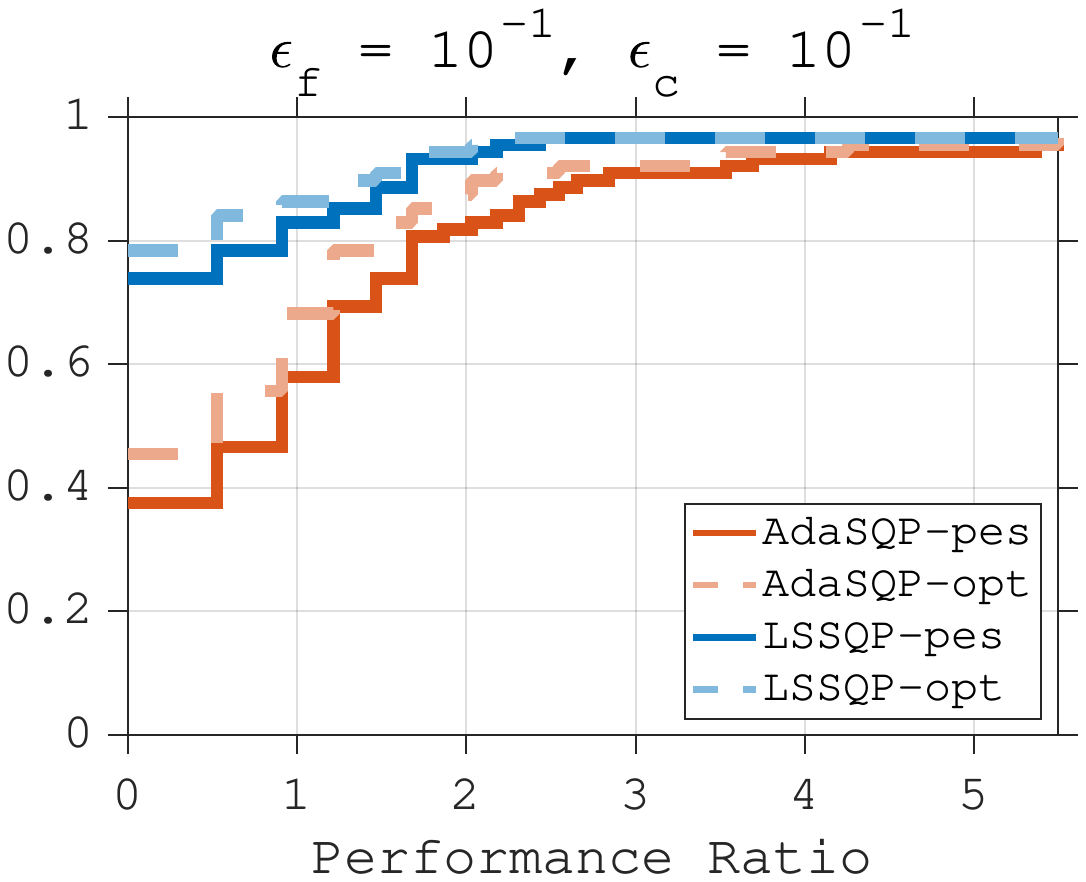}
\includegraphics[width=0.24\textwidth,clip=true,trim=10 180 50 150]{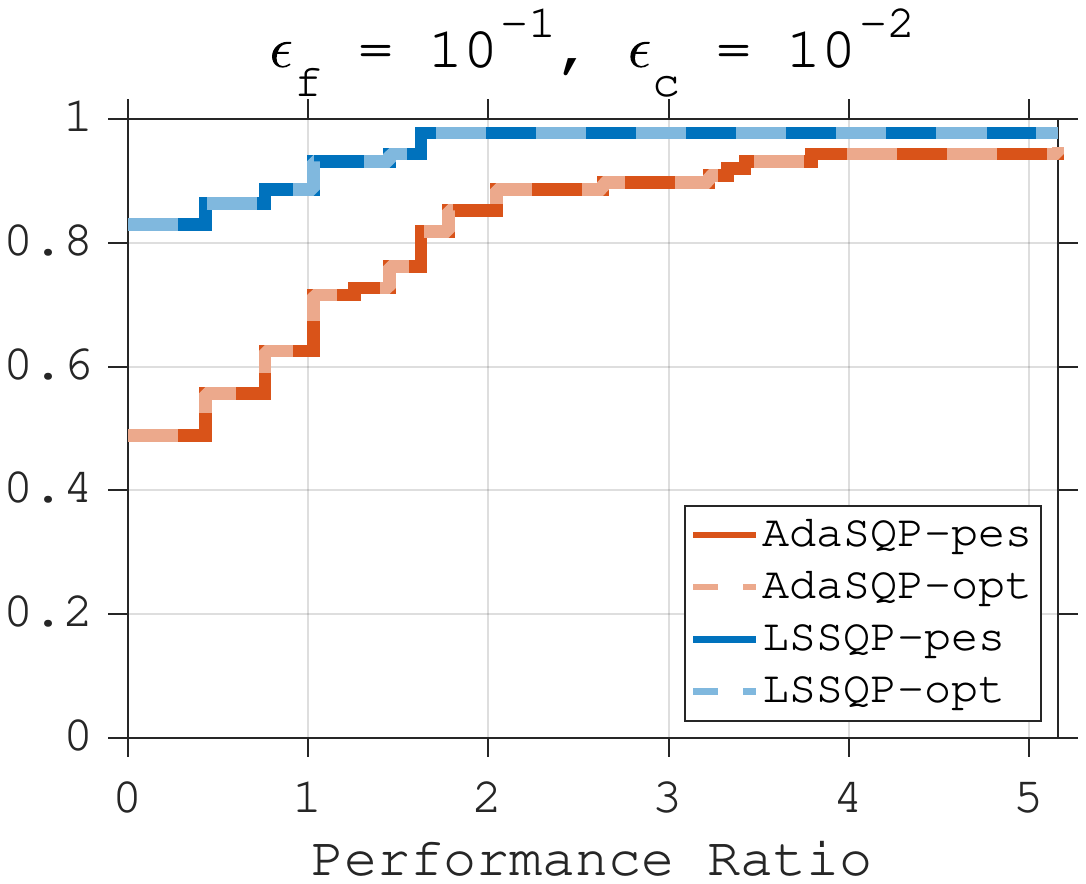}
\includegraphics[width=0.24\textwidth,clip=true,trim=10 180 50 150]{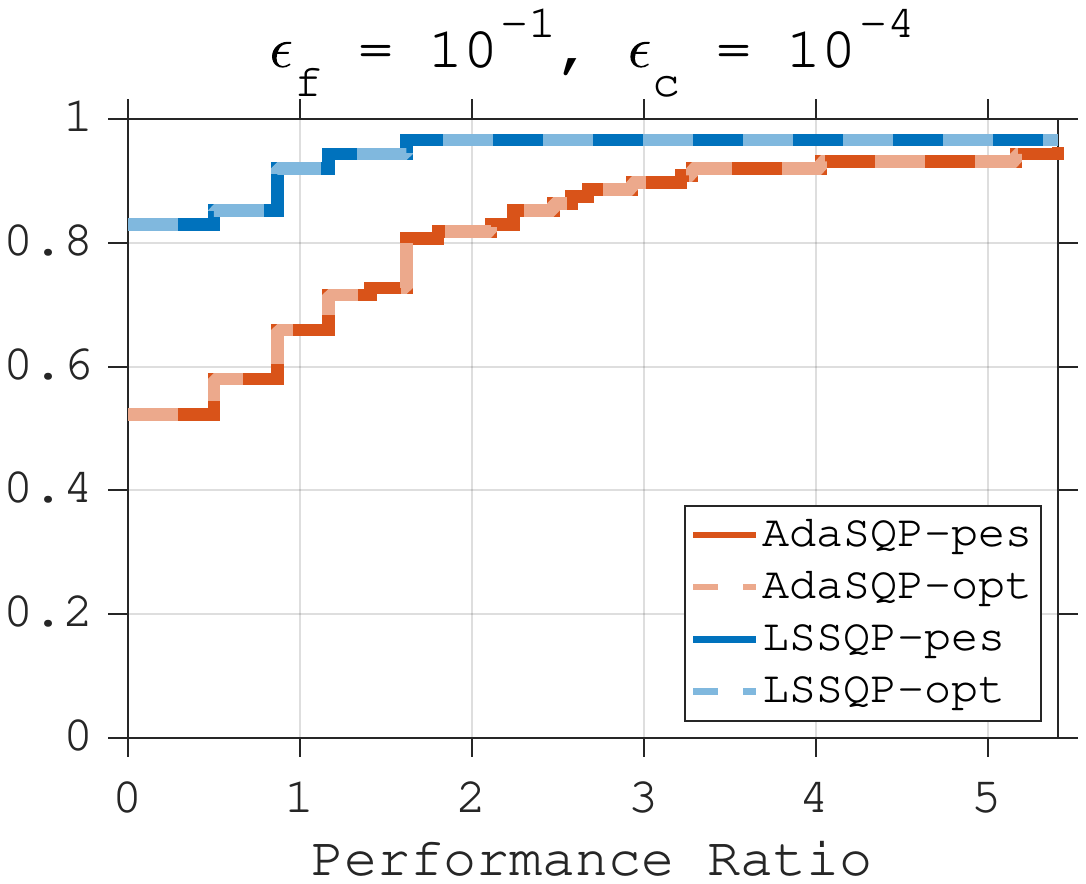}
\includegraphics[width=0.24\textwidth,clip=true,trim=10 180 50 150]{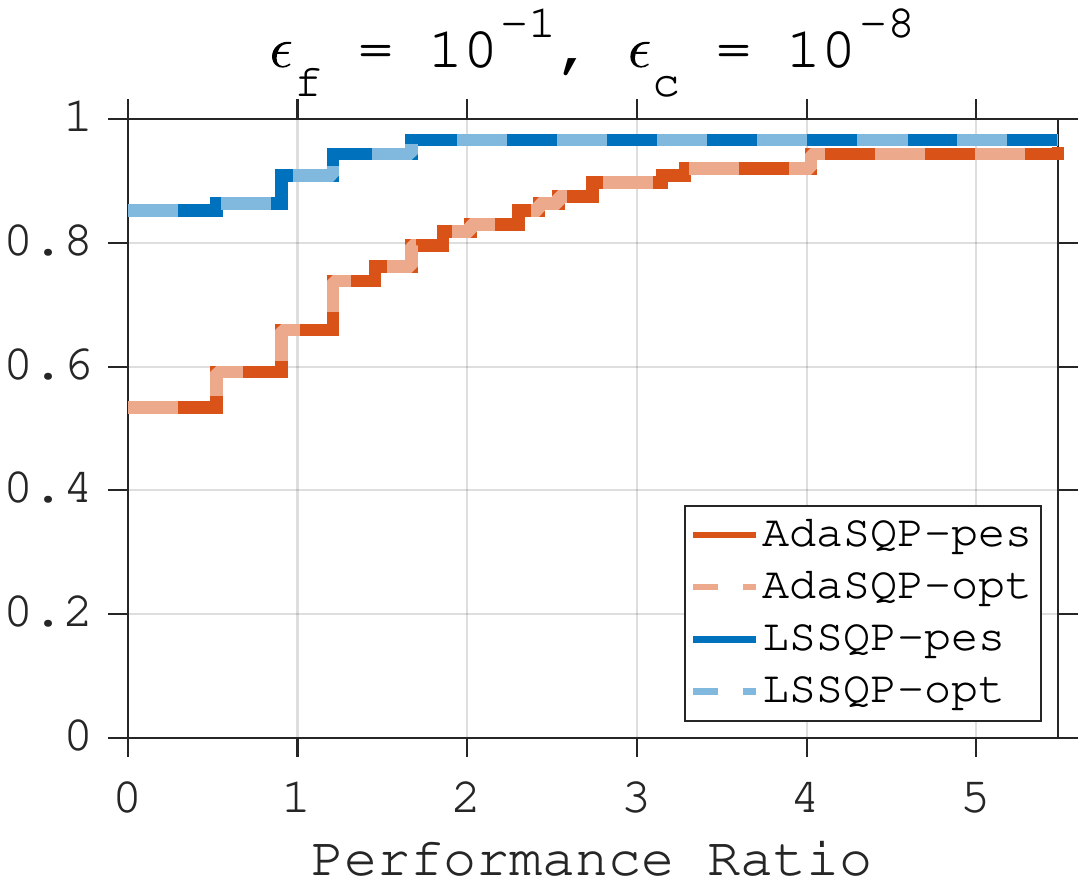}
\includegraphics[width=0.24\textwidth,clip=true,trim=10 180 50 150]{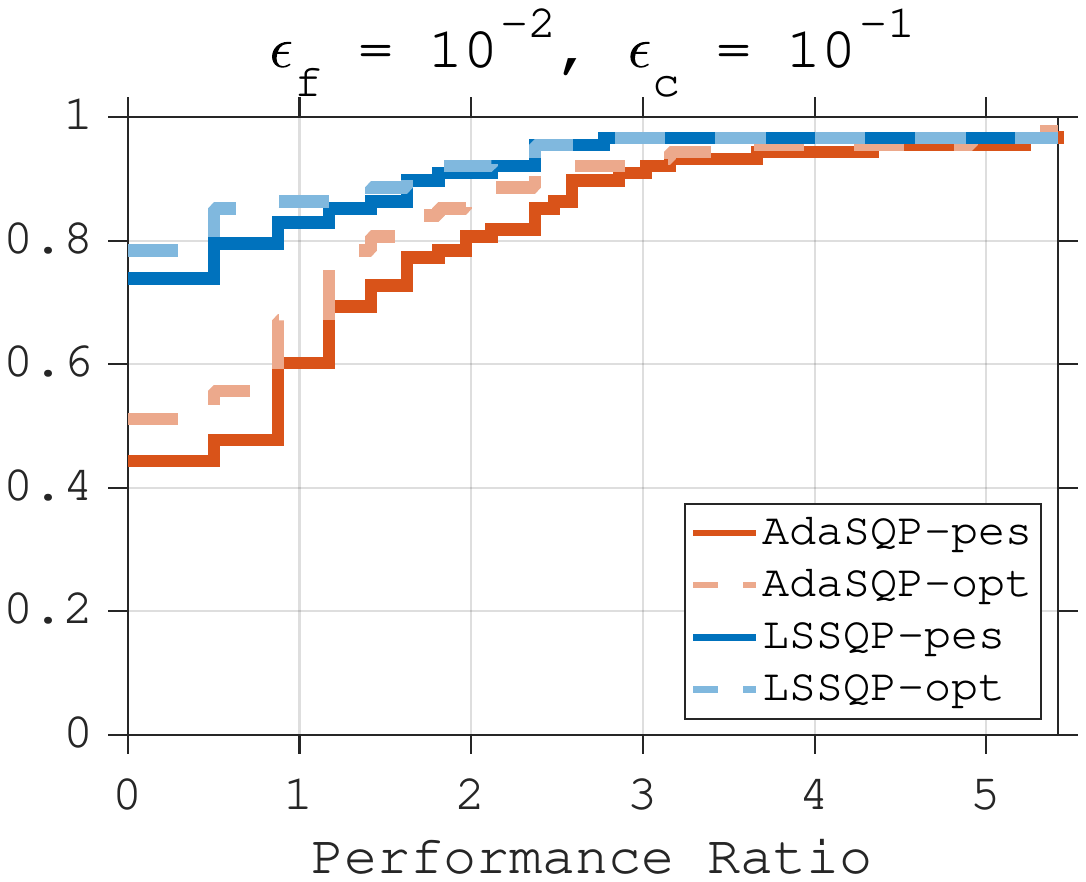}
\includegraphics[width=0.24\textwidth,clip=true,trim=10 180 50 150]{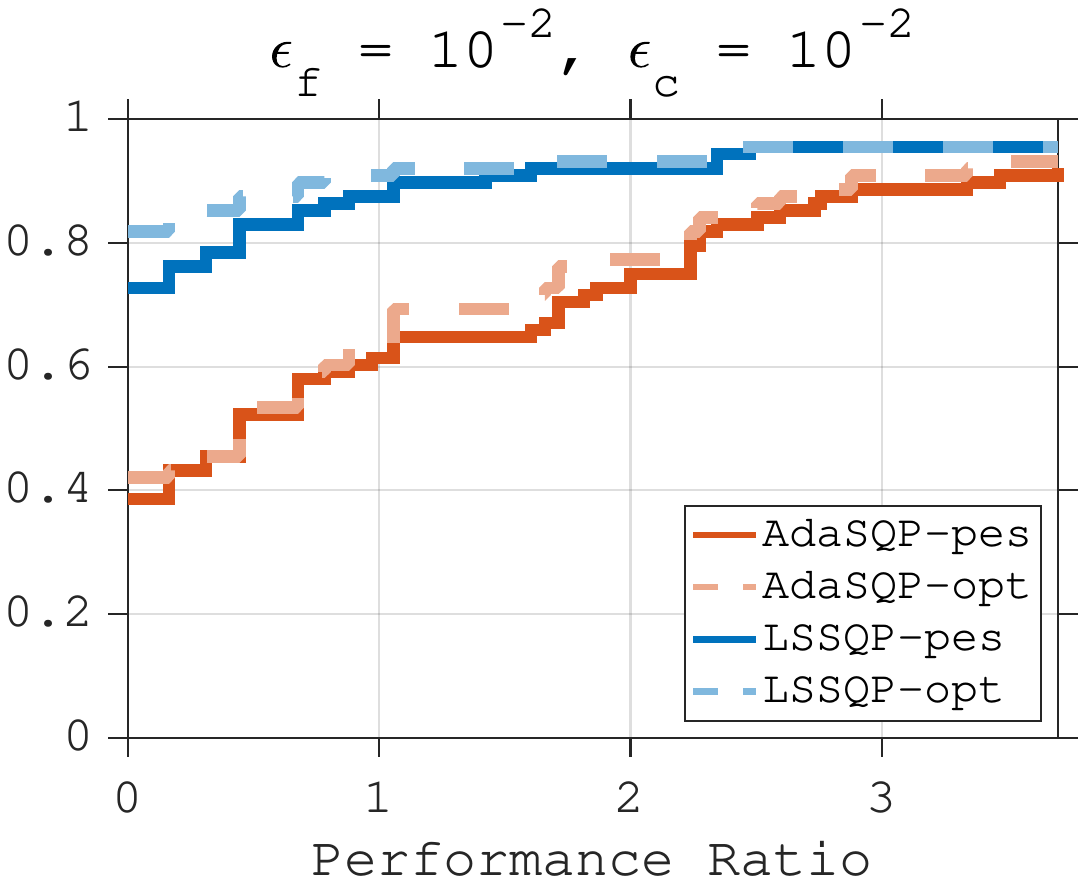}
\includegraphics[width=0.24\textwidth,clip=true,trim=10 180 50 150]{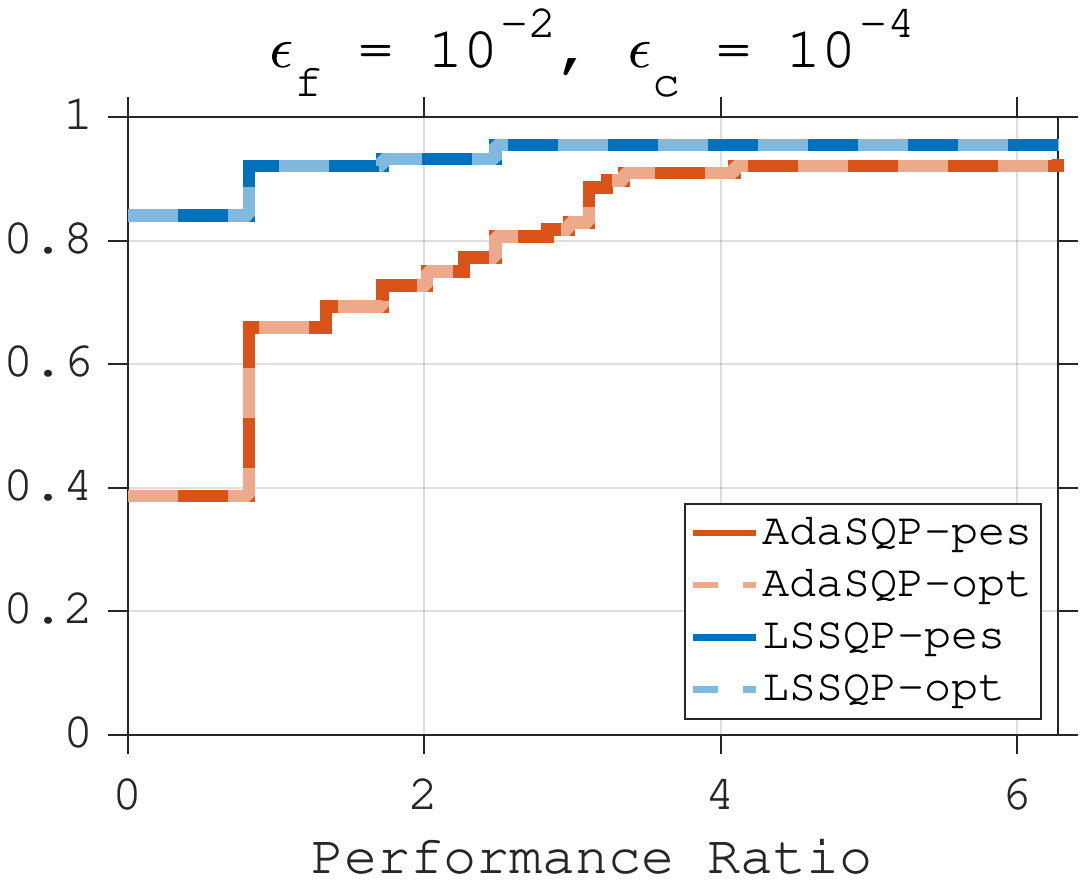}
\includegraphics[width=0.24\textwidth,clip=true,trim=10 180 50 150]{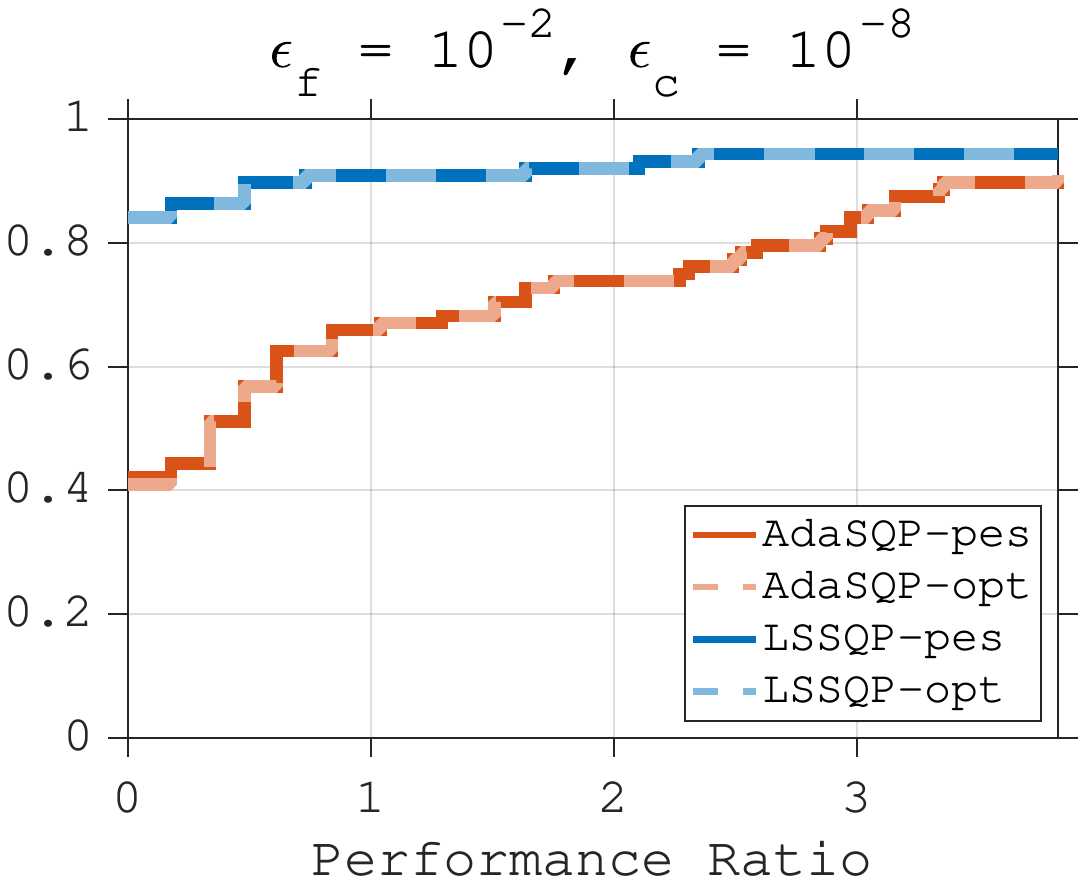}
\includegraphics[width=0.24\textwidth,clip=true,trim=10 180 50 150]{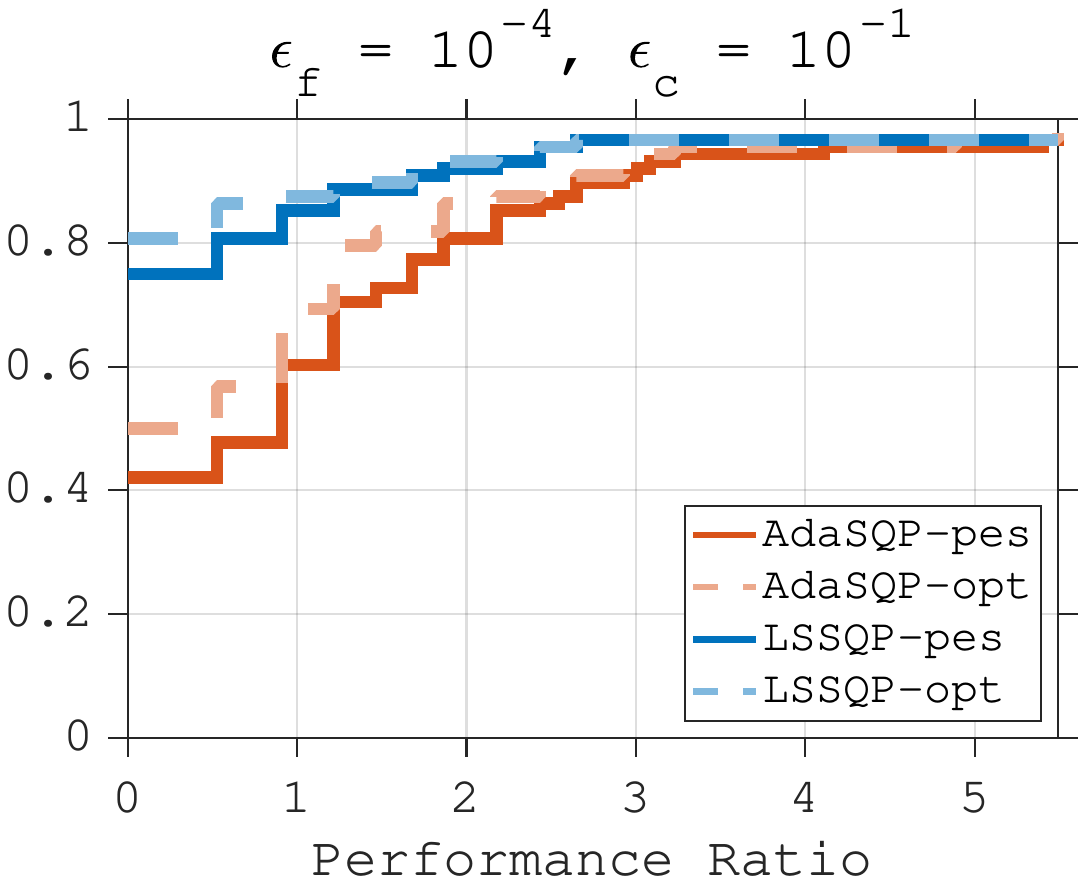}
\includegraphics[width=0.24\textwidth,clip=true,trim=10 180 50 150]{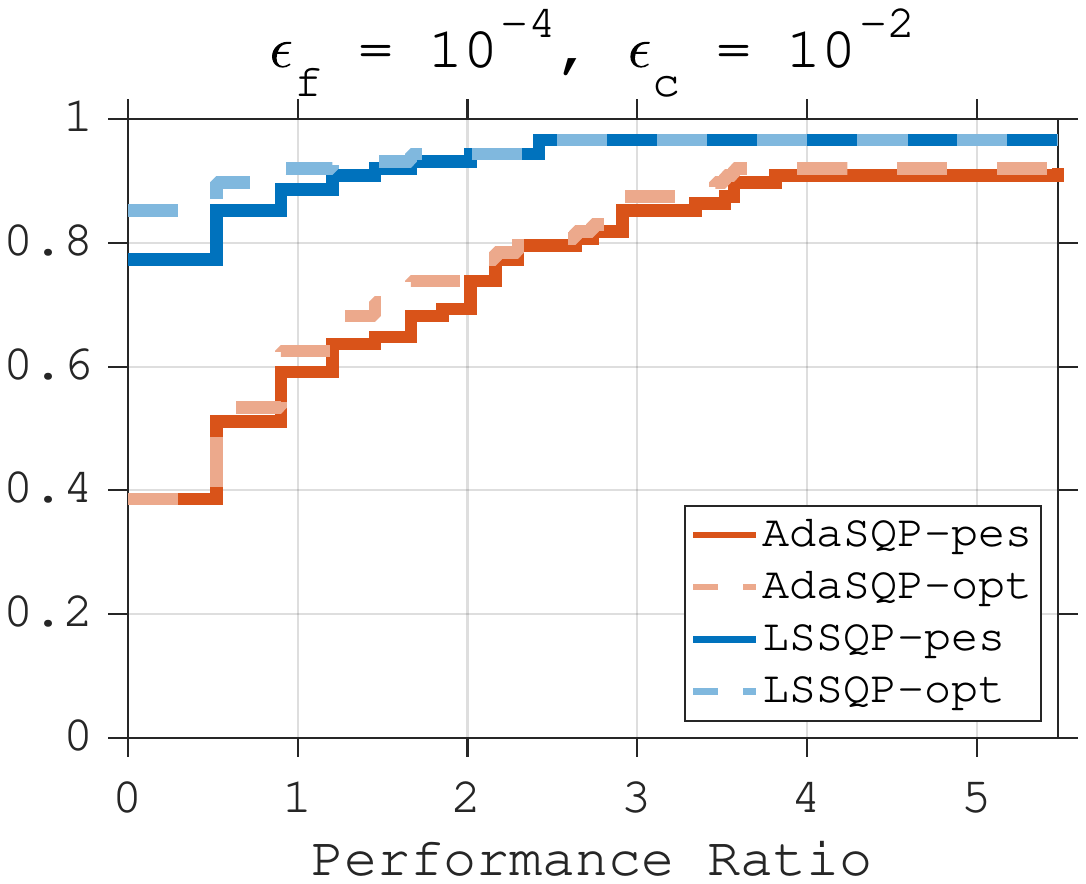}
\includegraphics[width=0.24\textwidth,clip=true,trim=10 180 50 150]{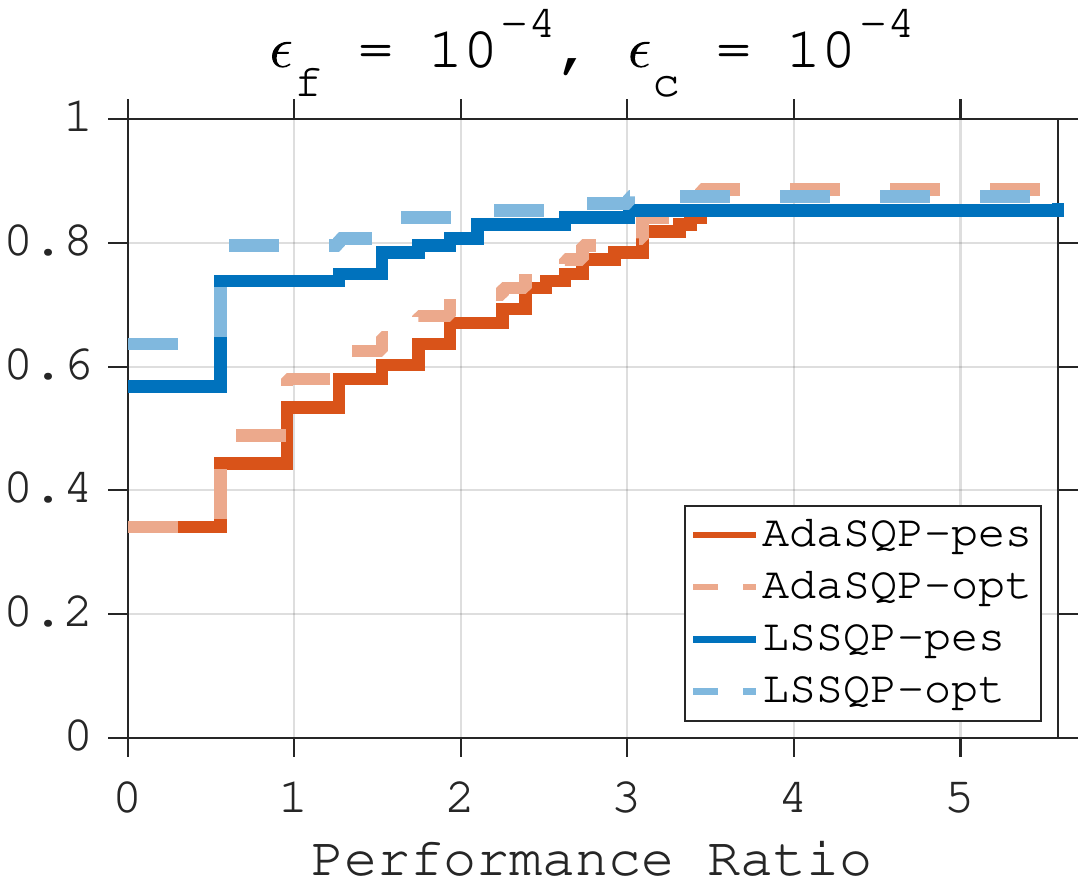}
\includegraphics[width=0.24\textwidth,clip=true,trim=10 180 50 150]{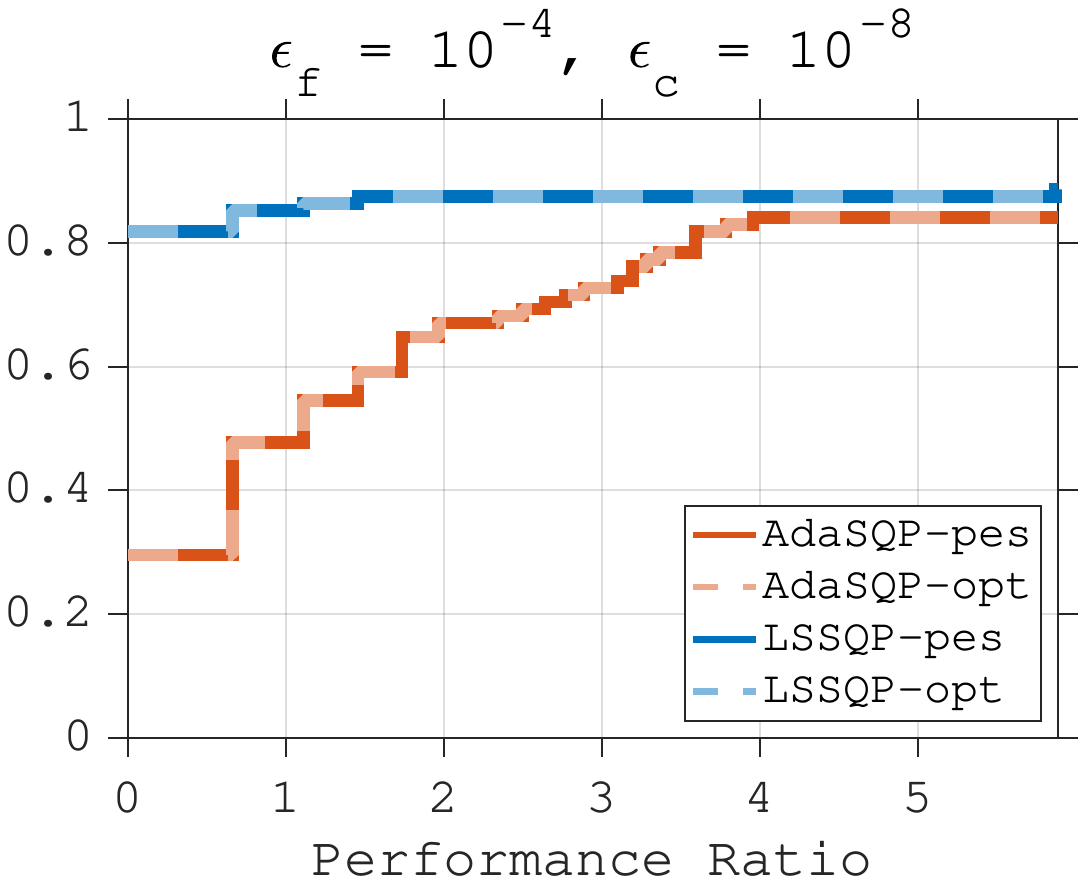}
\includegraphics[width=0.24\textwidth,clip=true,trim=10 180 50 150]{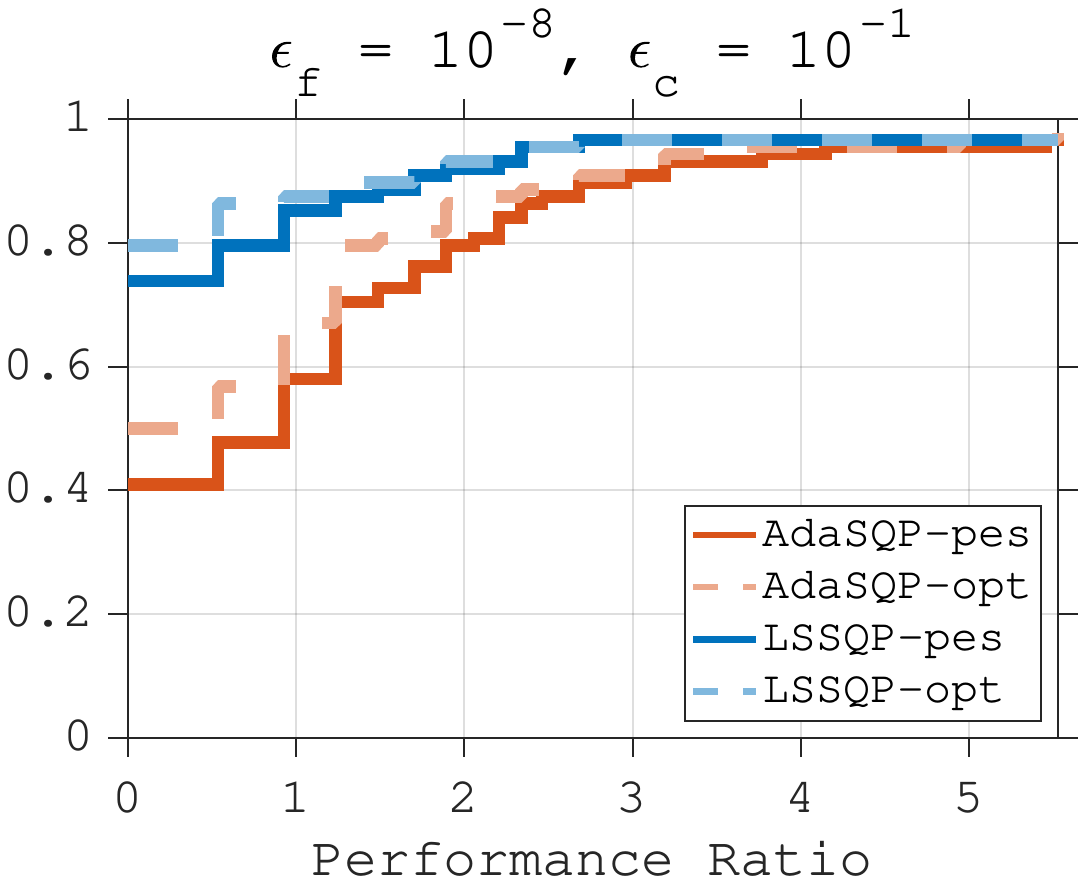}
\includegraphics[width=0.24\textwidth,clip=true,trim=10 180 50 150]{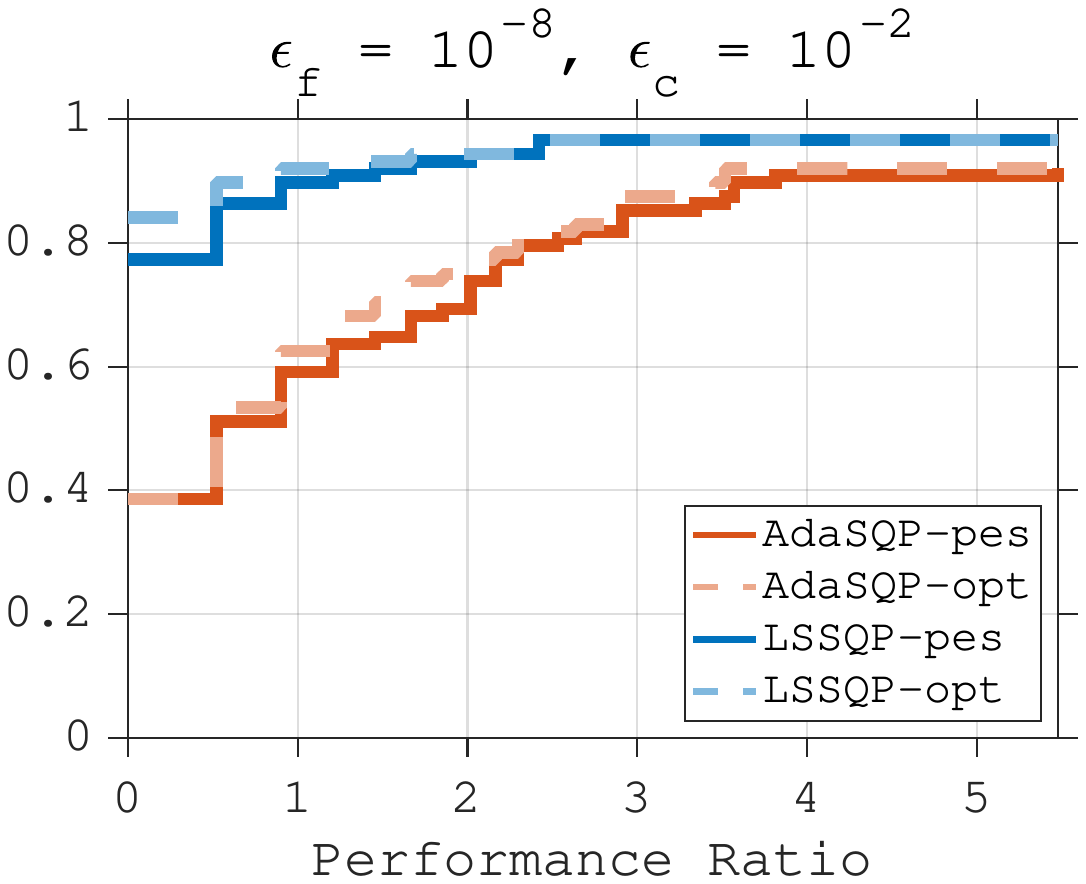}
\includegraphics[width=0.24\textwidth,clip=true,trim=10 180 50 150]{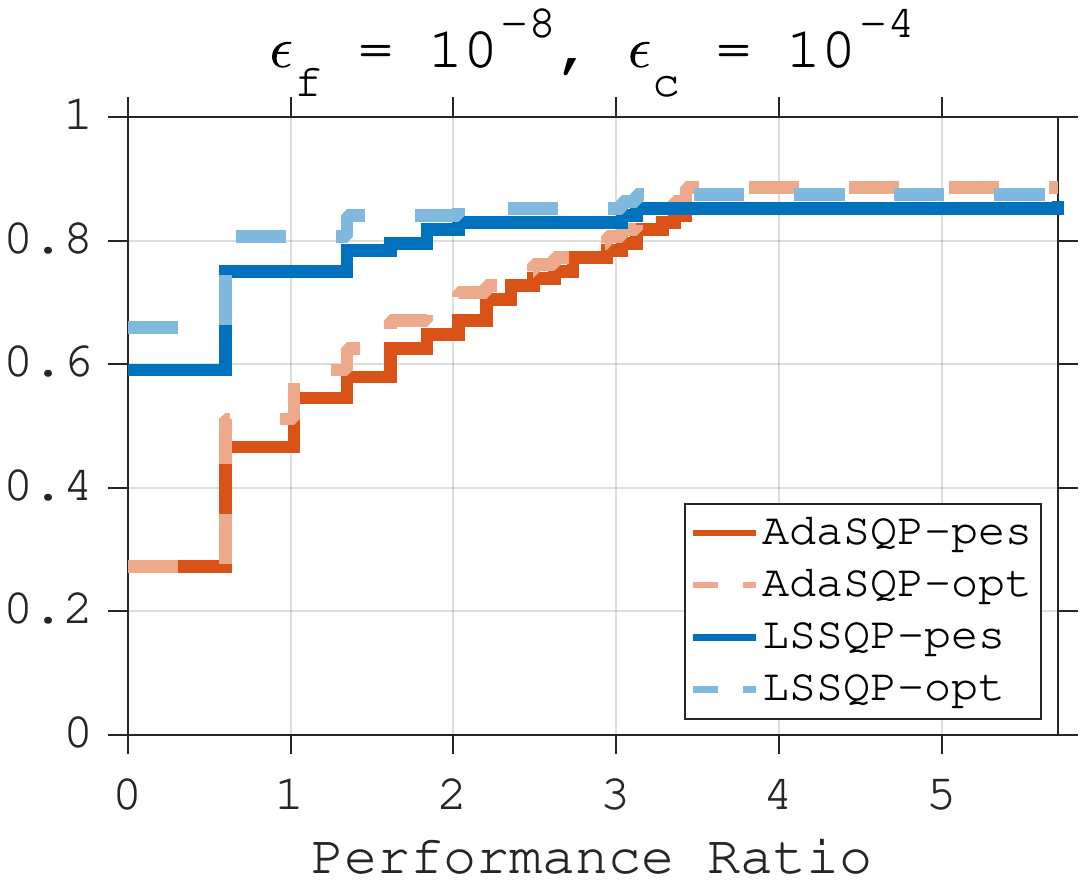}
\includegraphics[width=0.24\textwidth,clip=true,trim=10 180 50 150]{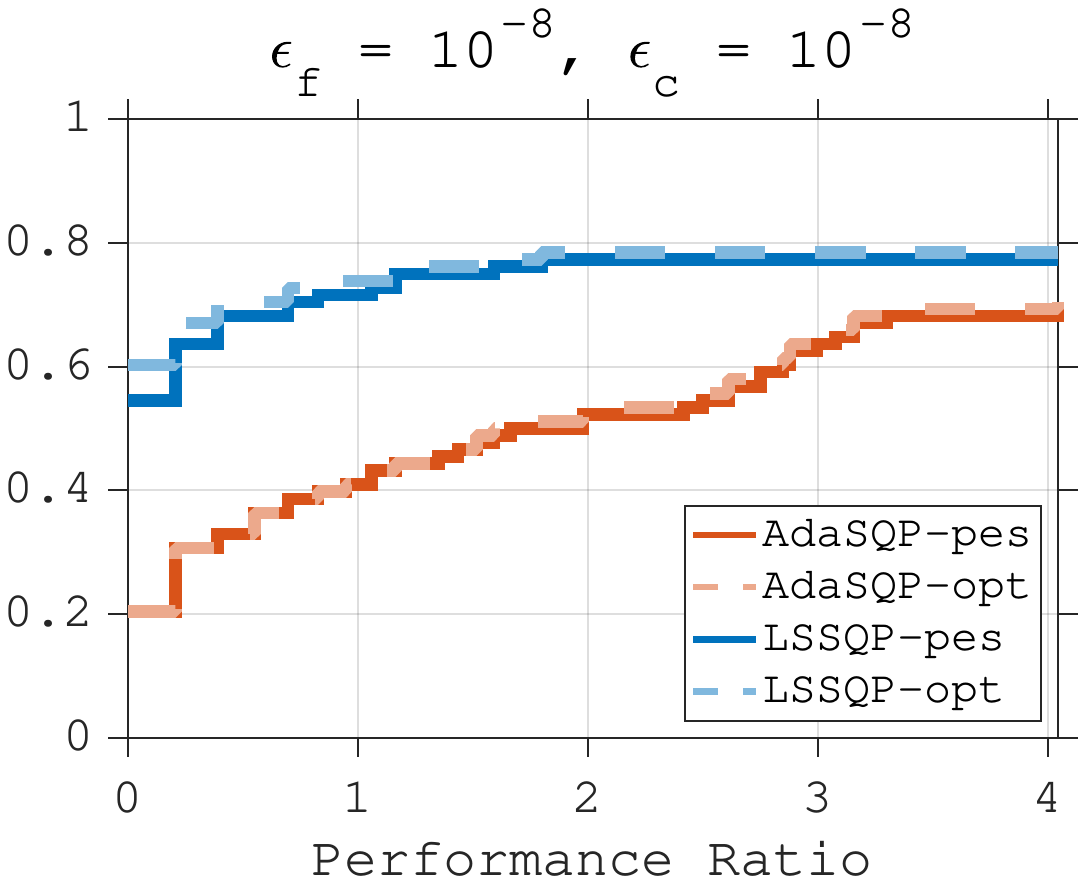}
\caption{Dolan-Mor\'e  performance profiles 
comparing \adasqppes{}, \adasqpopt{}, \lssqppes{}, and \lssqpopt{} on CUTEst collection of test problems in the absence of the LICQ in terms of \textbf{MINRES iterations} for $\epsilon_c \in \{ 10^{-1}, 10^{-2}, 10^{-4}, 10^{-8}\}$ (from \textbf{left} to \textbf{right}) and 
$\epsilon_f \in \{ 10^{-1}, 10^{-2}, 10^{-4}, 10^{-8}\}$ (from \textbf{top} to \textbf{bottom}). }
\label{fig.DMplot.NoLICQ.minres}
\end{figure}

\begin{figure}[htbp]
        \centering    
\includegraphics[width=0.24\textwidth,clip=true,trim=10 180 50 150]{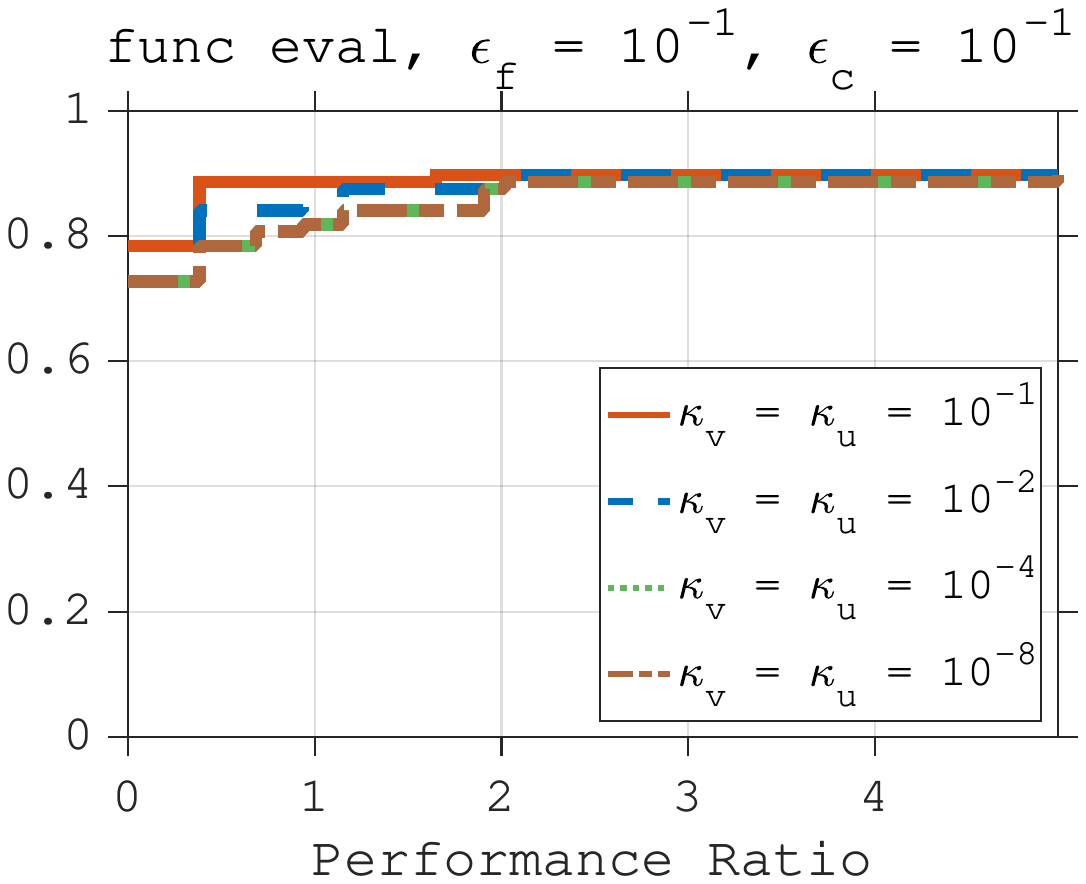}
\includegraphics[width=0.24\textwidth,clip=true,trim=10 180 50 150]{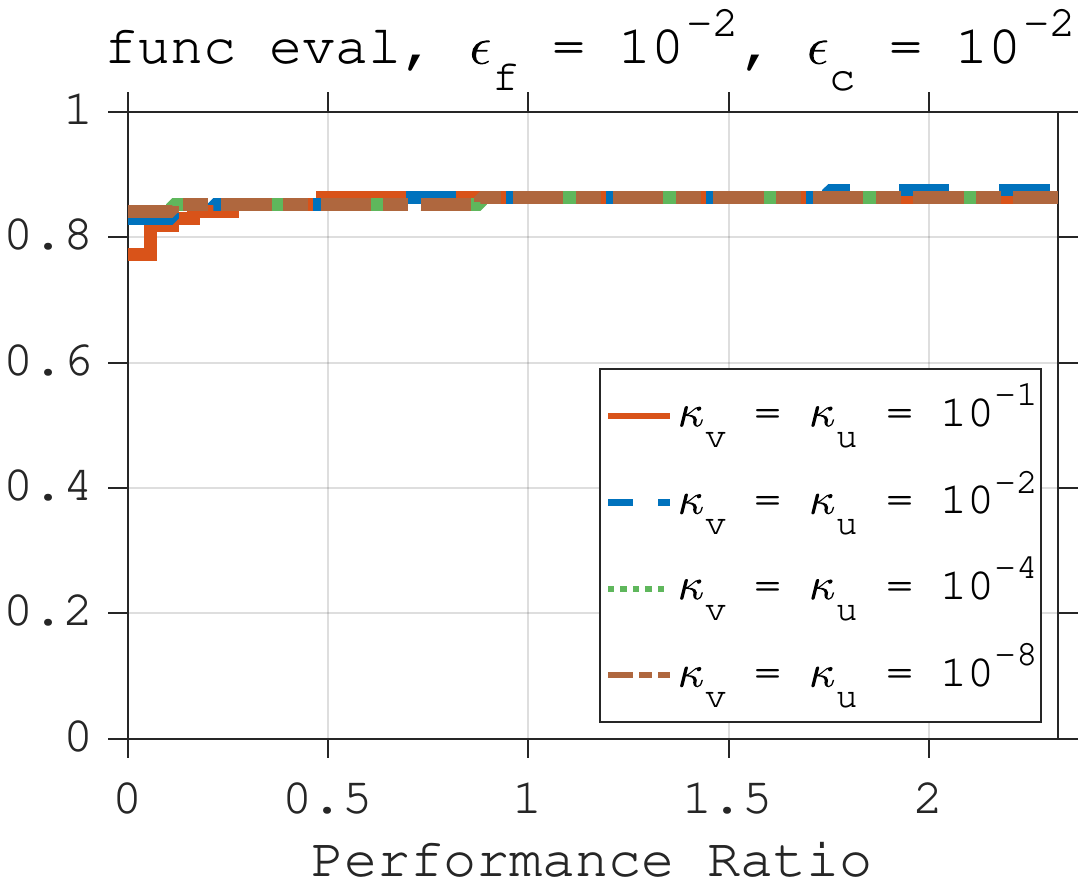}
\includegraphics[width=0.24\textwidth,clip=true,trim=10 180 50 150]{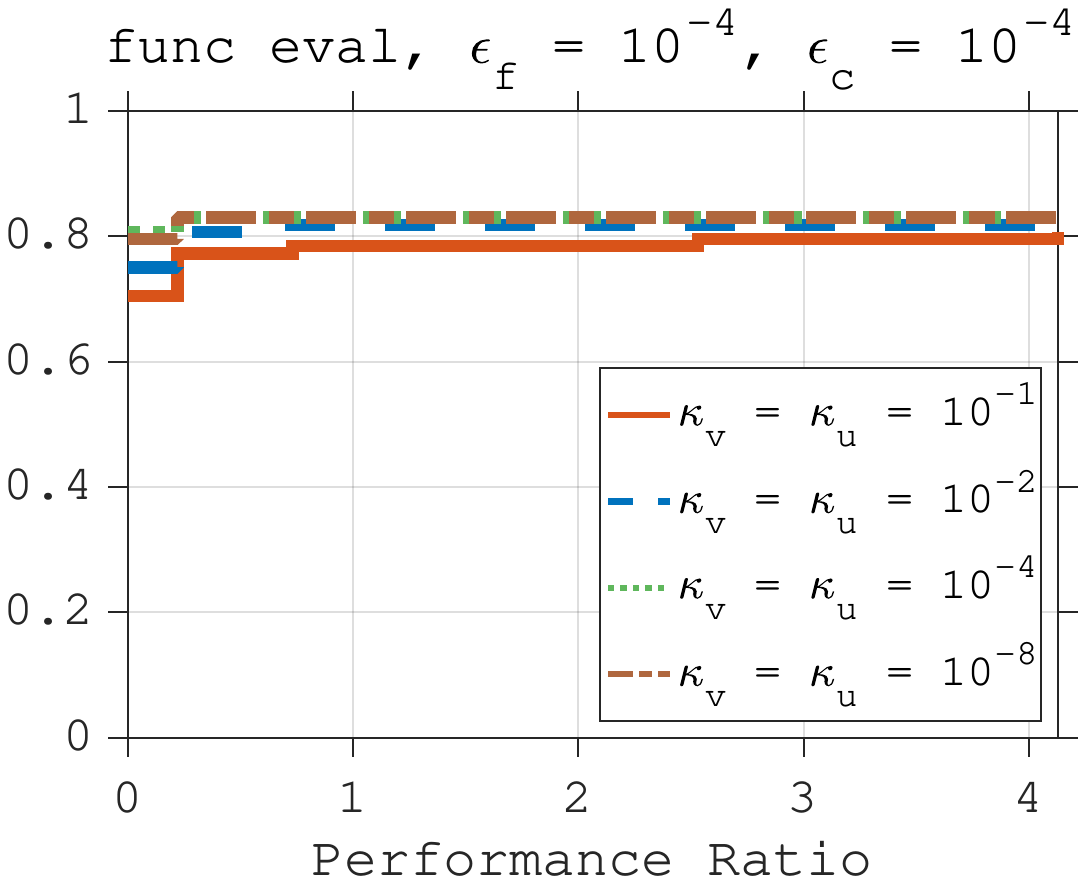}
\includegraphics[width=0.24\textwidth,clip=true,trim=10 180 50 150]{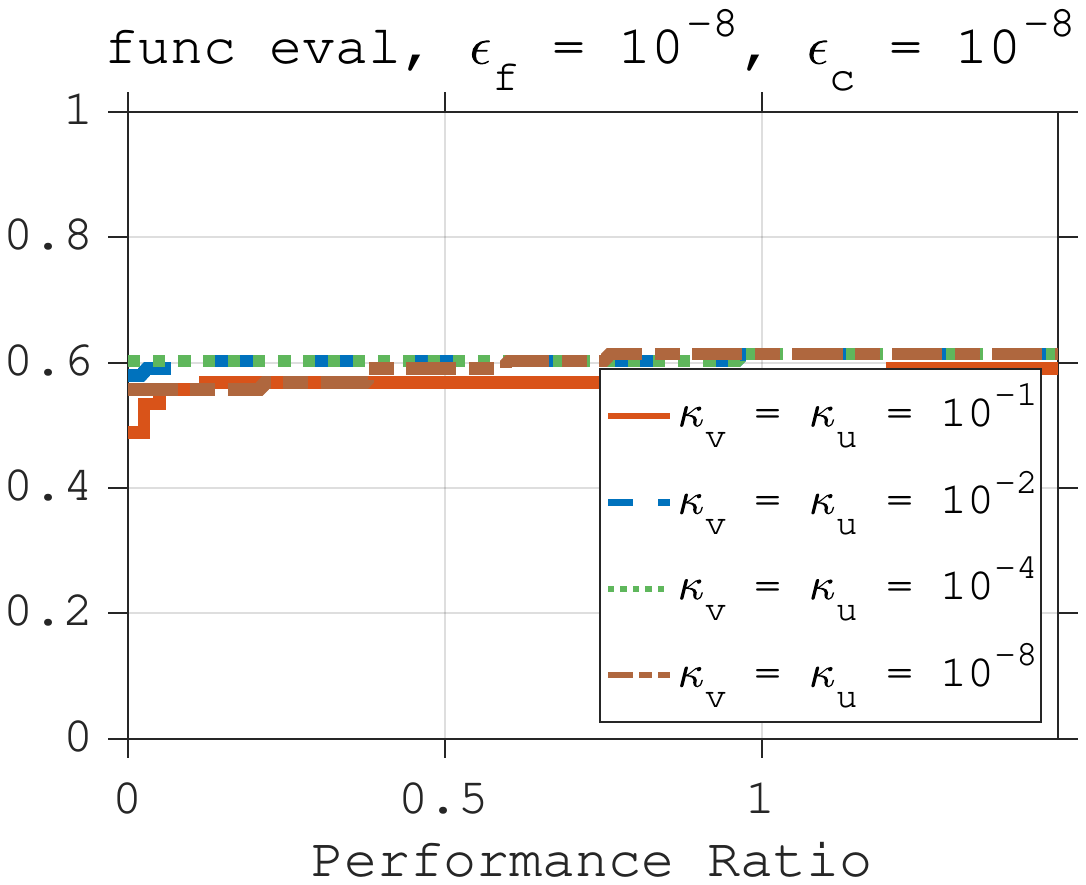}
\includegraphics[width=0.24\textwidth,clip=true,trim=10 180 50 150]{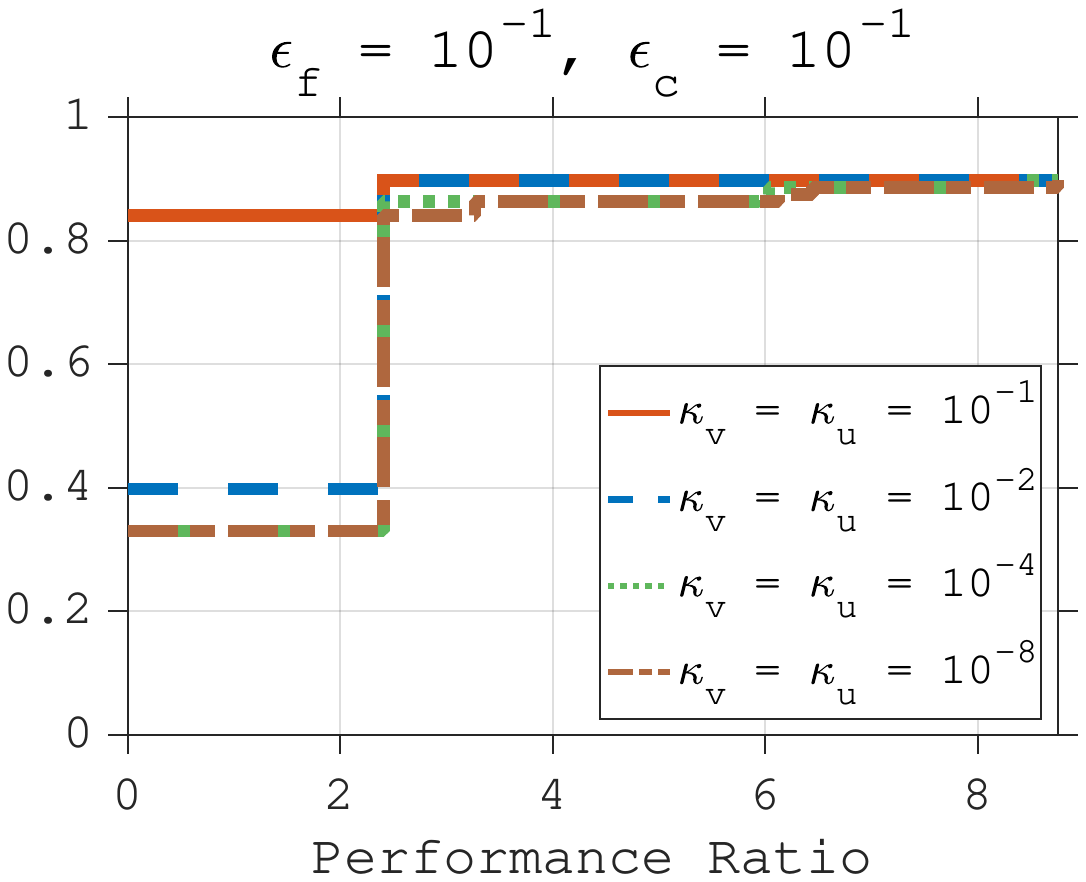}
\includegraphics[width=0.24\textwidth,clip=true,trim=10 180 50 150]{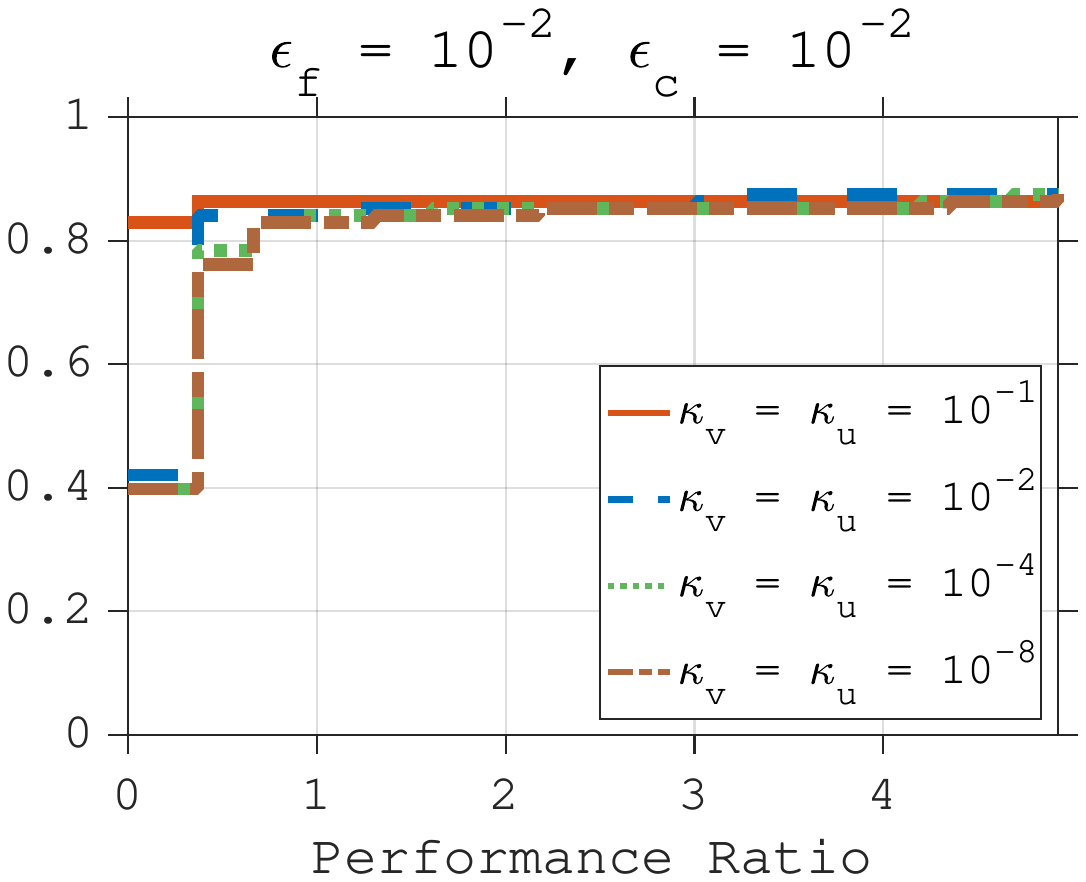}
\includegraphics[width=0.24\textwidth,clip=true,trim=10 180 50 150]{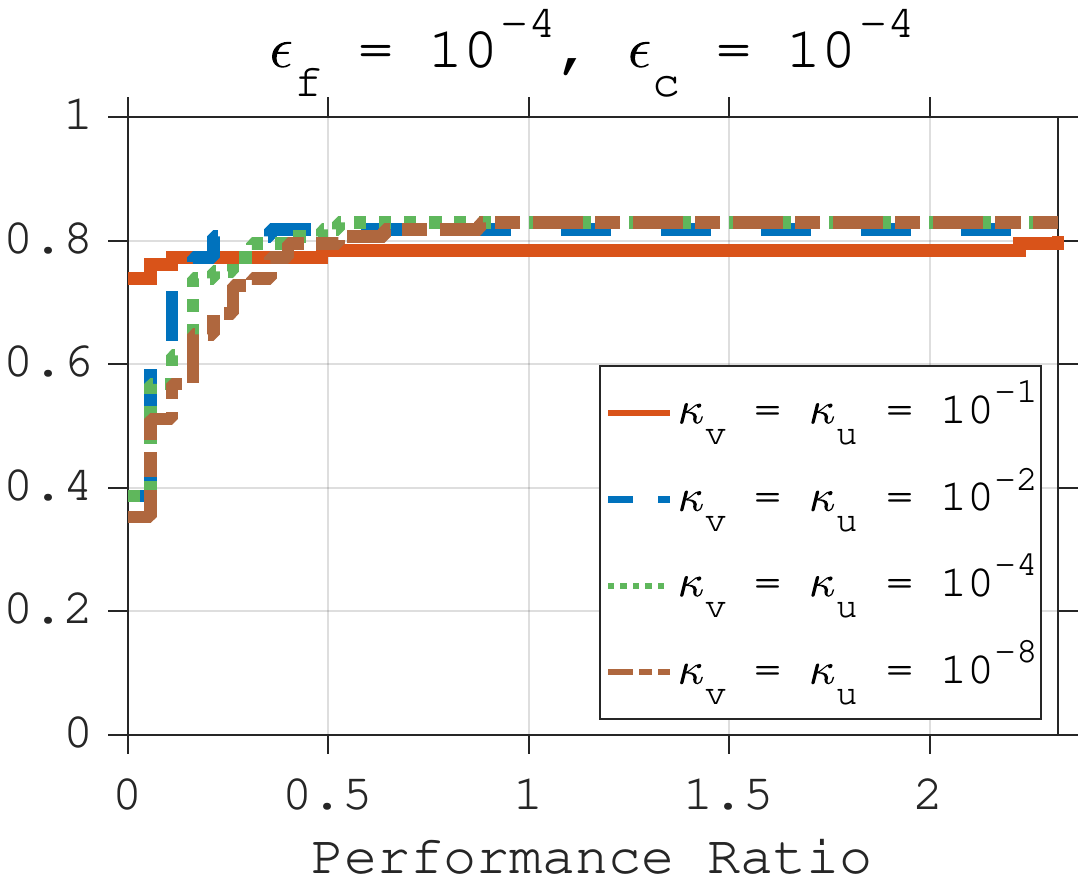}
\includegraphics[width=0.24\textwidth,clip=true,trim=10 180 50 150]{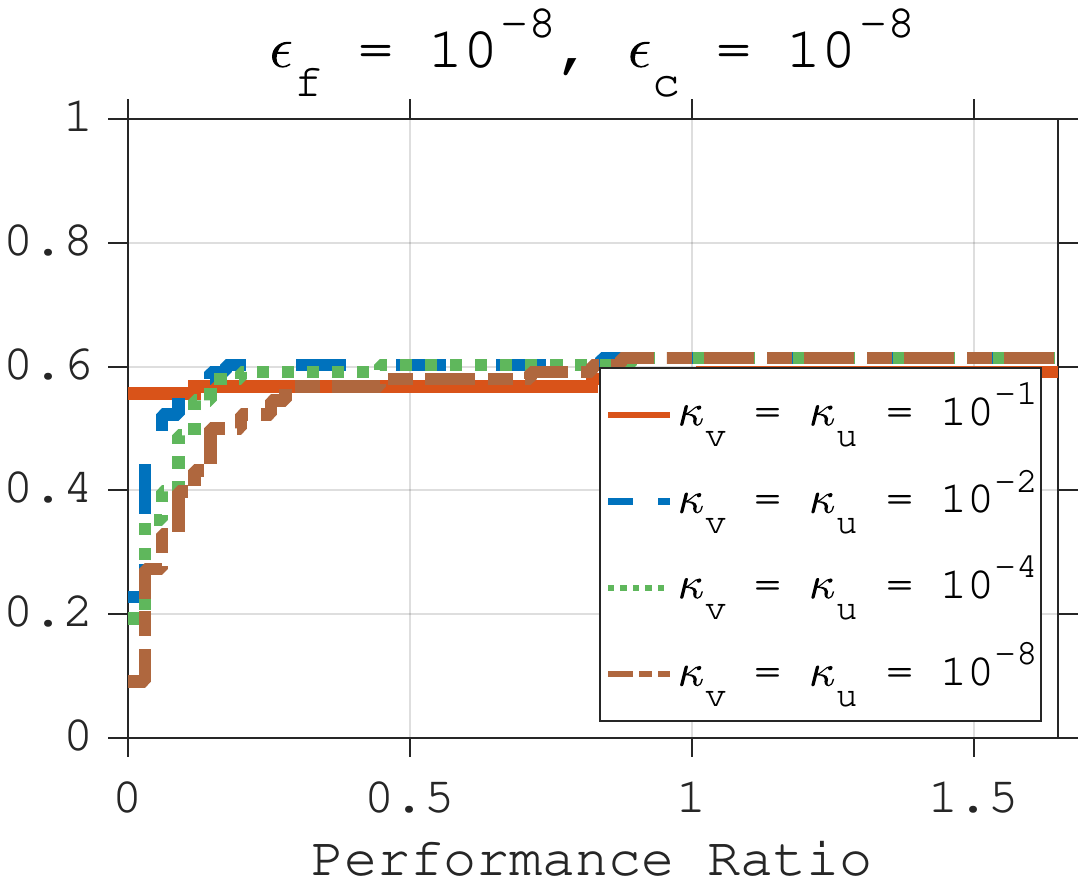}
\caption{Dolan-Mor\'e  performance profiles 
comparing different inexactness levels for \adasqpopt{} on CUTEst collection of test problems in terms of \textbf{function evaluations} (\textbf{top row}) and \textbf{CG and MINRES iterations} (\textbf{bottom row}) for $\epsilon_f = \epsilon_c \in \{ 10^{-1}, 10^{-2}, 10^{-4}, 10^{-8}\}$ (from \textbf{left} to \textbf{right}).}
\label{fig.DMplot.exactness.ada}
\end{figure}

\begin{figure}[htbp]
        \centering    
\includegraphics[width=0.24\textwidth,clip=true,trim=10 180 50 150]{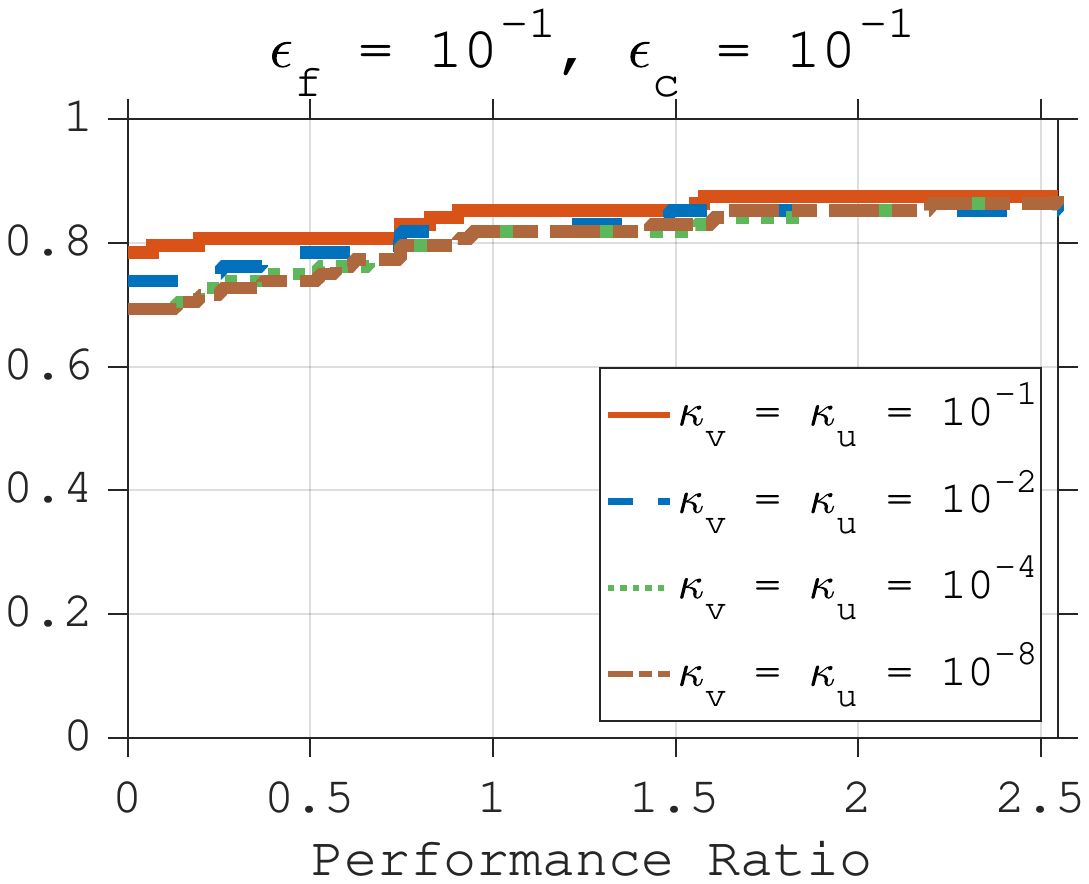}
\includegraphics[width=0.24\textwidth,clip=true,trim=10 180 50 150]{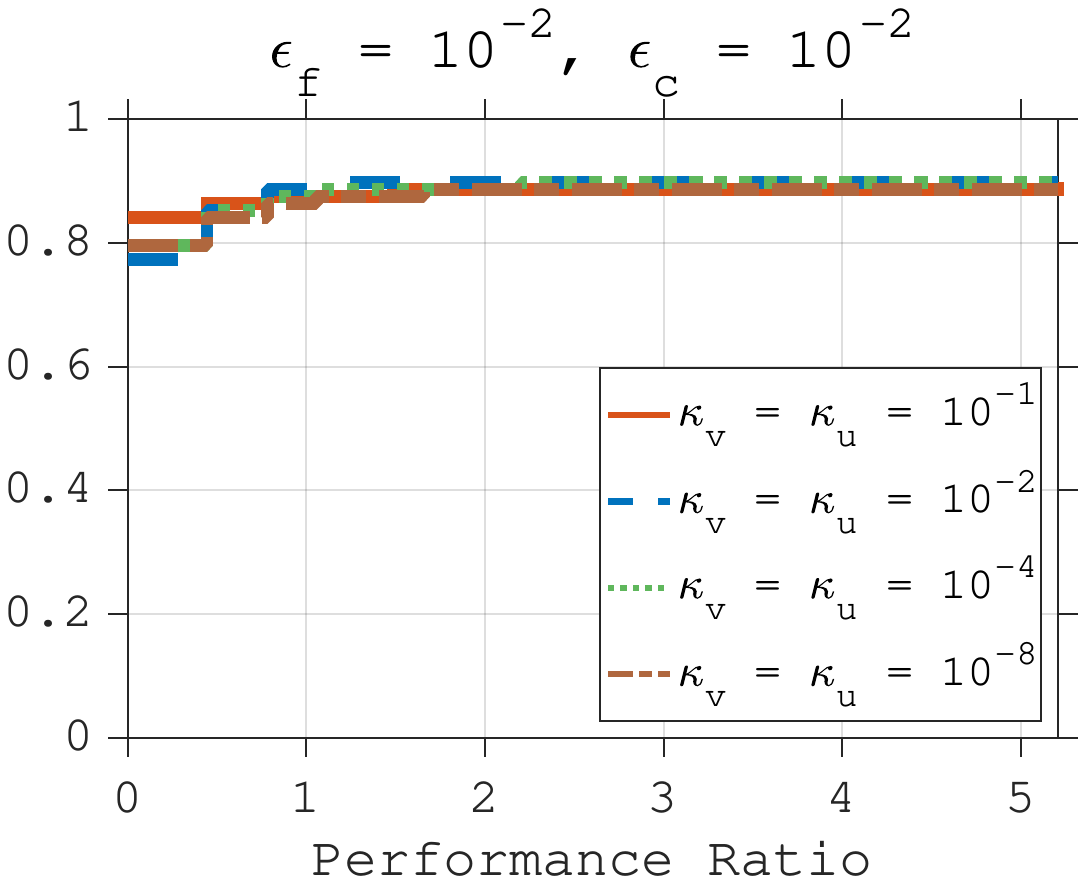}
\includegraphics[width=0.24\textwidth,clip=true,trim=10 180 50 150]{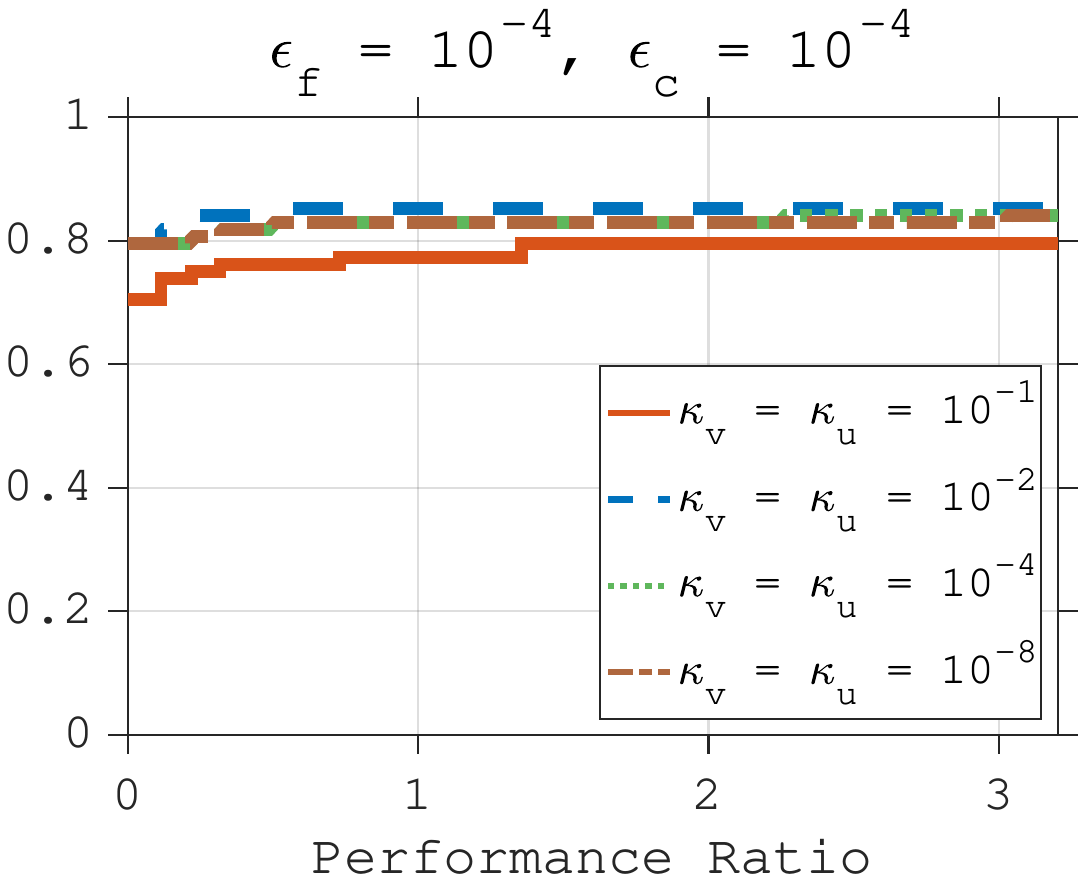}
\includegraphics[width=0.24\textwidth,clip=true,trim=10 180 50 150]{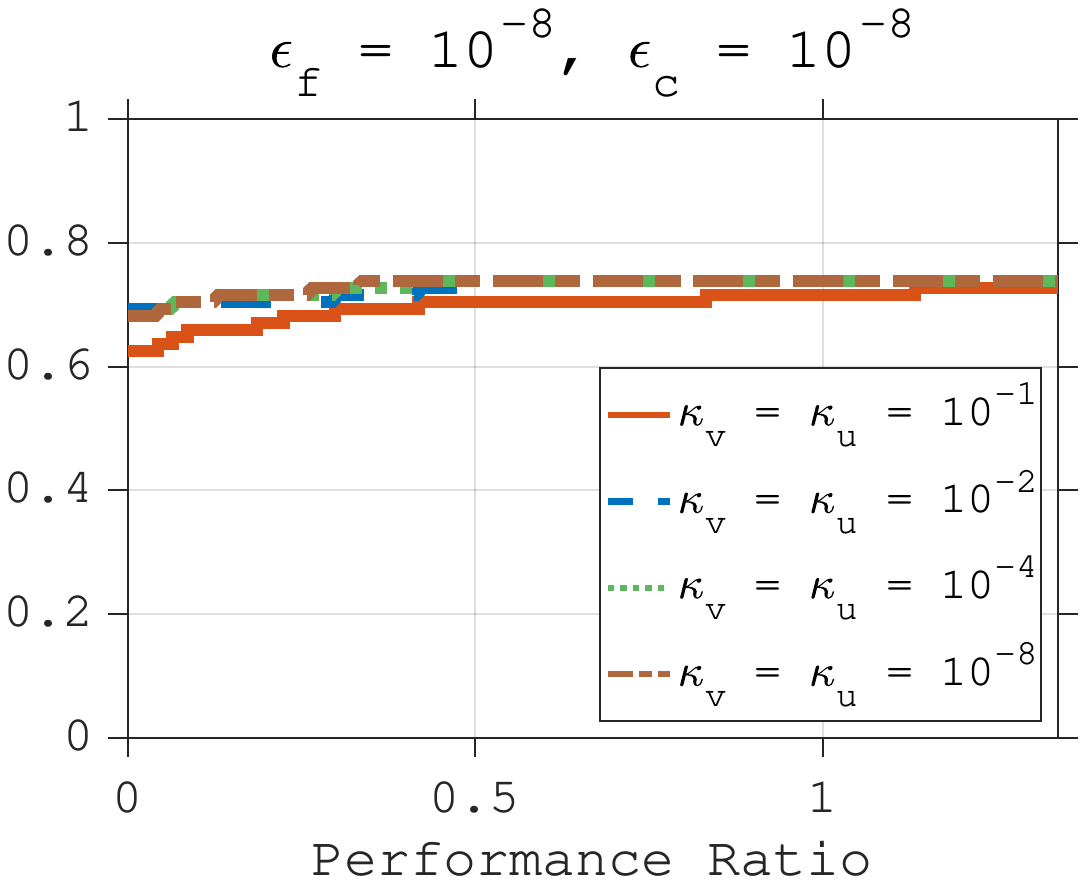}
\includegraphics[width=0.24\textwidth,clip=true,trim=10 180 50 150]{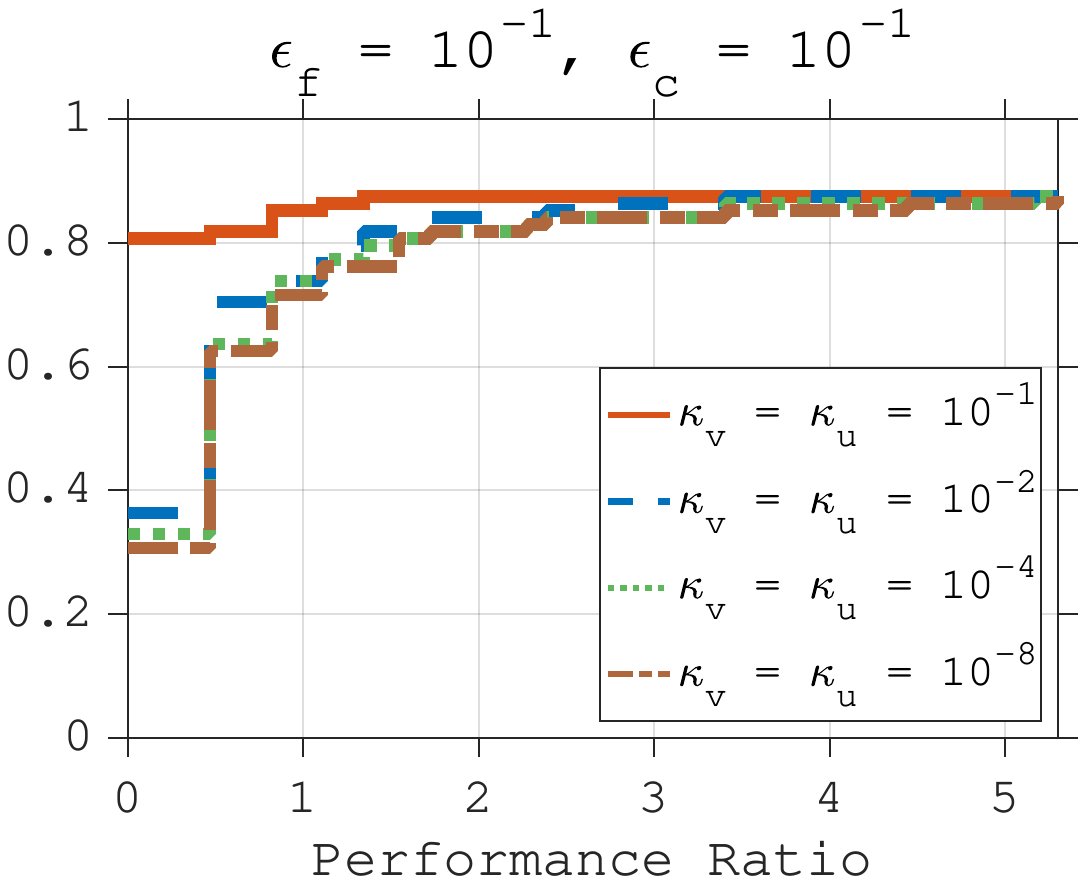}
\includegraphics[width=0.24\textwidth,clip=true,trim=10 180 50 150]{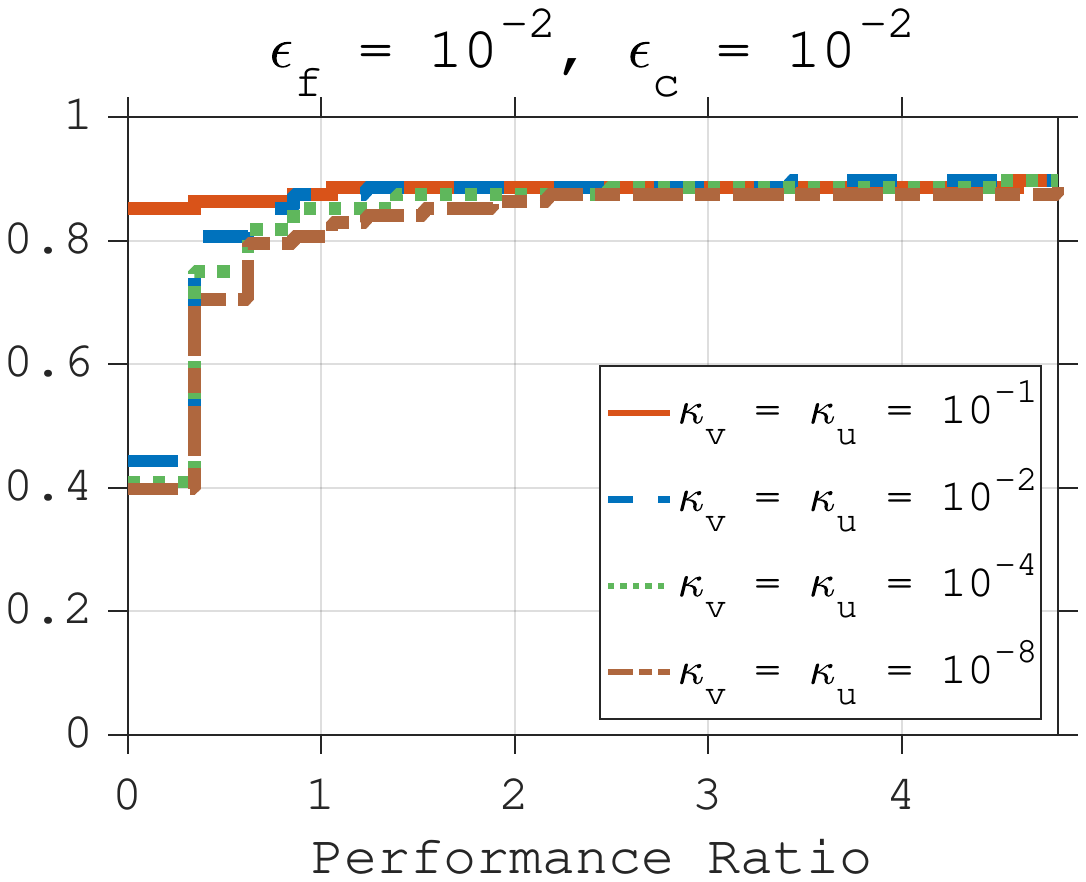}
\includegraphics[width=0.24\textwidth,clip=true,trim=10 180 50 150]{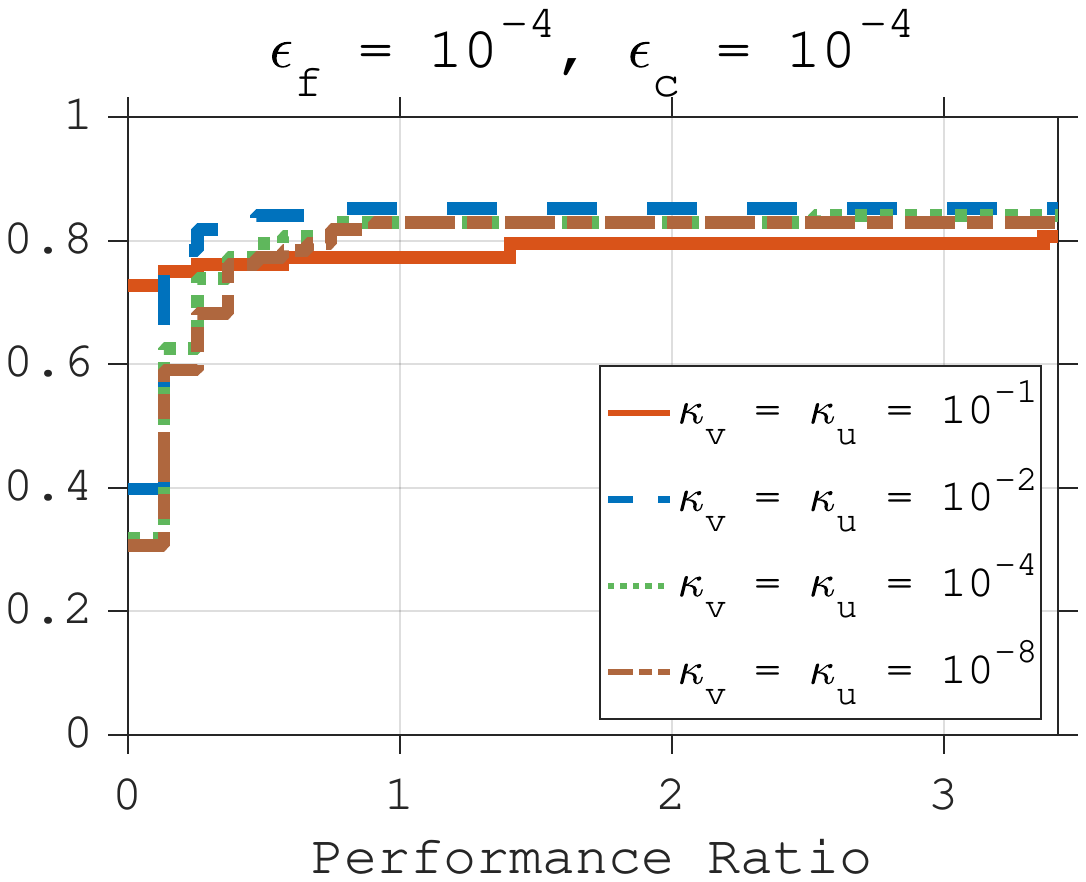}
\includegraphics[width=0.24\textwidth,clip=true,trim=10 180 50 150]{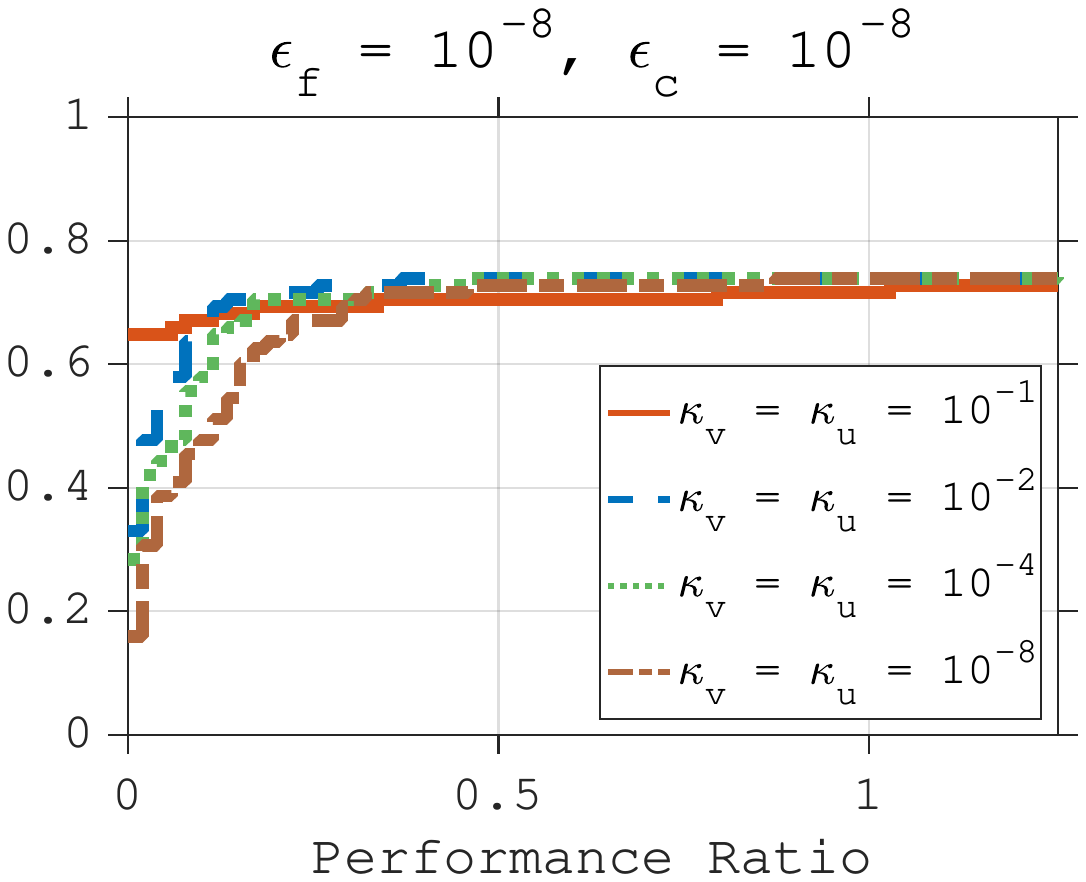}
\caption{Dolan-Mor\'e  performance profiles 
comparing different inexactness levels for \lssqpopt{} on CUTEst collection of test problems in terms of \textbf{function evaluations} (\textbf{top row}) and \textbf{CG and MINRES iterations} (\textbf{bottom row}) for $\epsilon_f = \epsilon_c \in \{ 10^{-1}, 10^{-2}, 10^{-4}, 10^{-8}\}$ (from \textbf{left} to \textbf{right}).}
\label{fig.DMplot.exactness.ls}
\end{figure}

In the remainder of our experiments, all algorithms are run with a budget of 1,000 iterations and a total of 10,000 function evaluations, where the total number of function evaluations is defined as the sum of the number of function evaluations $\bar{f}(\cdot)$ and two times the number of gradient evaluations $\bar{g}(\cdot)$. An algorithm is considered to have successfully solved a problem if an approximate stationary point is found, defined as
  \begin{align} 
 \label{eq.stationary.satisfied}
     \left\|c\left(x_k\right)\right\|_{\infty} \leq 2\max\{\epsilon_c , \epsilon_f \} \quad \text{and} \quad \left\|\nabla f\left(x_k\right)+\nabla c\left(x_k\right) y_k\right\|_{\infty} \leq 2(\epsilon_g + \|y_k\|_{\infty} \epsilon_J),
 \end{align}
where $y_k  = \arg \min_{y \in \mathbb{R}^m}\left\|\nabla f\left(x_k\right)+\nabla c\left(x_k\right) y\right\|^2$. If our proposed algorithm terminates early because of $\left\|c\left(x_k\right)\right\|_{\infty} \leq \epsilon_o$ and $\Delta \lbar(x_{k},\bar \tau_{k}, \bar d_{k}) \leq \epsilon_o$, it is deemed to have successfully solved the problem. For problems satisfying the LICQ, we consider both the pessimistic ($\epsilon_o = 0$) and optimistic ($\epsilon_o = \epsilon_c$) variants of \adasqp{} and \lssqp{}, and compare our inexact  methods with $\kappa_v = \kappa_u = 10^{-2}$ against the Noise-Tolerant SQP (\ntsqp{}) method proposed in \cite{oztoprak2023constrained}, where in each iteration the following linear system 
\begin{align*}
    \left[\begin{array}{cc}
H_k & \bar J_k^T \\
 \bar J_k & 0
\end{array}\right]\left[\begin{array}{l}
\bar  d_k \\
\bar y_k
\end{array}\right]=-\left[\begin{array}{l}
\bar g_k \\
\bar c_k
\end{array}\right] 
\end{align*}
is solved exactly using the MINRES procedure until 
 \begin{equation}
 \label{eq.inexact.terminate.u}
    \left\|\left[\begin{array}{c}
\bar\rho_{k,t} \\
\bar r_{k,t}
\end{array}\right]\right\|_{\infty} \leq 10^{-10}\max \left\{  \left\|\left[\begin{array}{c}
\bar {g}_k \\
\bar c_k
\end{array}\right]\right\|_{\infty}, 10^{-2}\right\}, 
\end{equation} 
to obtain the search direction. We present results in the form of Dolan-Mor\'e performance profiles \cite{dolan2002benchmarking} in terms of function evaluations and MINRES iterations. While our proposed methods incur additional computational cost to solve problem \eqref{eq.nolicq.tr.sub.p1}, this cost can be mitigated by combining \eqref{eq.nolicq.tr.sub.p1} and \eqref{eq.nolicq.sub.p2} into a single subproblem, similar to the \ntsqp{} algorithm when LICQ is satisfied. Therefore, we focus on MINRES iterations to enable  a fair comparison. In the first set of results, the LICQ is satisfied  (Figures~\ref{fig.DMplot.LICQ.fun} and \ref{fig.DMplot.LICQ.minres}) and in the second set of results the LICQ is violated (Figures~\ref{fig.DMplot.NoLICQ.fun} and \ref{fig.DMplot.NoLICQ.minres}). We also showcase the sensitivity of our method to the inexactness parameters $\kappa_u = \kappa_v \in \{10^{-1}, 10^{-2}, 10^{-4}, 10^{-8}\}$ (Figures~\ref{fig.DMplot.exactness.ada} and \ref{fig.DMplot.exactness.ls}).

Figures~\ref{fig.DMplot.LICQ.fun} and \ref{fig.DMplot.LICQ.minres} illustrate the performance of the optimistic and pessimistic inexact variants of \adasqp{} and \lssqp{}, and \ntsqp{}, in terms of function evaluations and MINRES iterations, respectively, on problems where the LICQ holds and under varying noise levels. In general, across all noise levels and the two metrics of comparison, the optimistic variants consistently perform as well as or better than their pessimistic counterparts and \ntsqp{}. In terms of function evaluations, \adasqp{} is more efficient than \lssqp{}. This is due to the fact that \lssqp{} requires at least one additional function evaluation at every iteration to find a suitable step size. For certain noise instances (e.g., $\epsilon_f = \epsilon_c = 10^{-8}$), \adasqp{} is able to solve fewer problem instances compared to \lssqp{}. This occurs because the step size for \adasqp{} may be too conservative, which prevents the algorithm from achieving the desired feasibility and stationarity within the budget. In terms of MINRES iterations, \lssqp{} is more efficient than \adasqp{}. This is due to the fact that \lssqp{} ensures that the step size yields sufficient reduction in the merit function, whereas \adasqp{} tends to adopt a more conservative step size. It is interesting to note that the optimistic variants solve a greater number of problem instances compared to the pessimistic methods in some noise settings (e.g., $\epsilon_f = 10^{-8}, \epsilon_c = 10^{-2}$). This can be attributed to the early termination condition and suggests that this condition can be utilized as a reliable termination condition for the algorithms. Although all methods solve fewer problem instances as the noise levels $\epsilon_f$ and $\epsilon_c$ decrease, they remain robust to changes in noise levels.

To gain further insights into the pairwise comparison between the optimistic and pessimistic methods, we report the number of instances, out of 88 CUTEst problems, where the optimistic method terminates earlier than the pessimistic method, as well as the number of instances where the pessimistic method terminates earlier than the optimistic method. As shown in Table~\ref{tab:opt_improve_table}, for both \adasqp{} and \lssqp{}, the optimistic variants terminated earlier for most of the problem instances and noise levels. The difference is more pronounced at high (objective function or constraint) noise levels. That said, there are a few problem instances where the optimistic method performs worse than the pessimistic method. This occurs because, even when $\|\bar{c}_k\| \leq \epsilon_c$, it may still be beneficial to continue to reduce the feasibility error in some cases since the simulated noise could be much smaller than the upper bound $\epsilon_c$.

\begin{table}[H]
\centering
\begin{tabular}{c ccccc}
\toprule
\multirow{2}{*}{$\epsilon_c$} & \multicolumn{4}{c}{$\epsilon_f$} \\ 
\cmidrule(l){2-5}
 & $10^{-1}$ & $10^{-2}$ & $10^{-4}$ & $10^{-8}$ \\  
\cmidrule(r){1-1} \cmidrule(l){2-5}
$10^{-1}$ & ($\pmb{54}/1$, $\pmb{45}/0$) & ($\pmb{47}/1$, $\pmb{43}/1$) & ($\pmb{52}/0$, $\pmb{47}/0$)  & ($\pmb{57}/0$, $\pmb{50}/0$)  \\
$10^{-2}$ & ($\pmb{44}/0$, $\pmb{37}/0$)    & ($\pmb{24}/1$, $\pmb{27}/0$)    & ($\pmb{36}/2$, $\pmb{38}/0$) & ($\pmb{43}/1$, $\pmb{44}/0$) \\
$10^{-4}$ & ($\pmb{43}/0$, $\pmb{40}/0$)        & ($11/0$, $7/0$)        & ($8/7$, $4/0$)   & ($\pmb{37}/3$, $\pmb{36}/0$) \\
$10^{-8}$ & ($\pmb{46}/0$, $\pmb{39}/0$)        & ($4/0$, $3/0$)        & ($1/0$, $2/0$) & ($2/9$, $2/0$)  \\
\bottomrule
\end{tabular}
\caption{Number of problem instances for which the optimistic/pessimistic variants is better among the 88 CUTEst problems for both \adasqp{} (first entry of the tuple) and \lssqp{} (second entry of the tuple).} 
\label{tab:opt_improve_table}
\end{table}

Figures~\ref{fig.DMplot.NoLICQ.fun} and \ref{fig.DMplot.NoLICQ.minres} illustrate the performance of the  optimistic and pessimistic inexact variants of \adasqp{} and \lssqp{} on the same 88 CUTEst problems (same as in Figures~\ref{fig.DMplot.LICQ.fun} and \ref{fig.DMplot.LICQ.minres}) but in the setting where the LICQ is violated in a controlled way (as described above). We do not compare to \ntsqp{} since this method was designed for the setting in which LICQ holds. As shown in Figures~\ref{fig.DMplot.NoLICQ.fun} and \ref{fig.DMplot.NoLICQ.minres},  \adasqp{} is more efficient than \lssqp{} when the evaluation metric is total function evaluations while 
\lssqp{} is more efficient when considering MINRES iterations. This observations are similar to those for the setting in which the LICQ is satisfied. Compared to the setting in which the LICQ is satisfied (Figures~\ref{fig.DMplot.LICQ.fun} and \ref{fig.DMplot.LICQ.minres}), 
the performance gap between the pessimistic and optimistic variants is smaller in this setting.

\newpage

Finally, we compare the performance of the \adasqpopt{} and \lssqpopt{} methods for different levels of inexactness in the subproblem solutions. Specifically, we set $\kappa_u = \kappa_v \in \{10^{-1}, 10^{-2}, 10^{-4}, 10^{-8}\}$. For conciseness, we present Dolan-Moré performance profiles in terms of both function evaluations and the total number of CG and MINRES iterations, focusing on the case where $\epsilon_c = \epsilon_f \in \{10^{-1}, 10^{-2}, 10^{-4}, 10^{-8}\}$. Figures \ref{fig.DMplot.exactness.ada} (\adasqpopt{}) and \ref{fig.DMplot.exactness.ls} (\lssqpopt{}) demonstrate that our proposed algorithms are highly robust to changes in $\kappa_u$ and $\kappa_v$. In particular, the most inexact method ($\kappa_u = \kappa_v = 10^{-1}$) fails on only one additional problem instance compared to the other exactness levels when $\epsilon_c = \epsilon_f = 10^{-4}$. That said, in general, this inexactness level is more efficient than the other levels in terms of the total number of CG and MINRES iterations.

\newpage

\section{Final Remarks and Future Work}\label{sec.remarks}

This paper introduces a sequential quadratic programming (SQP) algorithm to solve noisy nonlinear optimization problems with noisy nonlinear equality constraints, addressing the challenges posed by bounded noise in the objective and constraint functions as well as potential rank deficiency in the constraint Jacobians. Our approach incorporates a step decomposition strategy, inexact subproblem solutions,  two flexible step size selection schemes, and an optimistic early termination scheme to enhance computational efficiency and robustness. We have established convergence to a neighborhood of a first-order stationary point in the presence of full-rank constraint Jacobians and convergence to a neighborhood of an infeasible stationary point in the rank-deficient case, with both neighborhood sizes dictated by the noise levels in the problem.  Numerical experiments conducted on the CUTEst test set demonstrate the efficiency and robustness of the proposed method. 


{\small
\bibliographystyle{plain}
\bibliography{references}}

\newpage
\appendix

\end{document}